\documentclass[
    oneside, 
    fontsize=12pt,          % default font size 12pt
    paper=a4,               % DIN A4 page format
    numbers=noenddot,       % remove dots behind chapter numbers (e.g. 1.5 not 1.5.)
    listof=totoc,           % add list of figures, tables, etc. to ToC
    listof=entryprefix,     % add entry name to figures, tables, etc.
    listof=nochaptergap,    % no chapter gap for figures, tables, etc.
    bibliography=totoc,     % add bibliography to ToC but without a chapter number
    parskip=half            % half line spacing between paragraphs
]{scrbook}

\usepackage{cmap}               % PDF character encoding
\usepackage[T1]{fontenc}        % 8-bit font encoding
\usepackage[utf8]{inputenc}     % UTF-8 input encoding

\usepackage[draft]{todonotes}   % notes showed
\newcommand{\docAuthor}{Patrick Jung}
\newcommand{\docMatriculationNumber}{22382810}

\newcommand{\docCity}{Höchstadt}

\newcommand{\docTitle}{Optimizing initial Feature-Mapping Variables from given Designs via Tracking}
\newcommand{\docSubTitle}{}     % leave empty to remove
\newcommand{\docType}{Master Thesis}
\newcommand{\docSupervisor}{Prof. Dr. Michael Stingl} % Main supervisor
\newcommand{\docSecondSupervisor}{Dr. Fabian Wein}    % Co-supervisor
\newcommand{\docCoSupervisor}{Dr. Lukas Pflug}
\newcommand{\docStartDate}{01.04.2025}
\newcommand{\docDeadline}{30.09.2025}

\usepackage{scrhack}            % better KOMA adaptions
\usepackage[]{xcolor}      % color support
\usepackage{chngcntr}           % for renumbering stuff
\usepackage{amsmath}
\usepackage{listings}
\usepackage{subcaption}
\usepackage{booktabs}
\usepackage{sidecap}
\usepackage{graphicx}
\usepackage{svg}
\usepackage{pgf}
\usepackage{pgfplots}
\pgfplotsset{compat=1.18}
\usepackage{lmodern}
\usepackage{amsmath}
\usepackage{amssymb}
\usepackage{array}
\usepackage{placeins}
\usepackage{multirow}
\usepackage[ruled,vlined,linesnumbered]{algorithm2e}
\usepackage[
    colorlinks=false,
   	linkcolor=black,
   	citecolor=black,
  	filecolor=black,
	urlcolor=black,
    bookmarks=true,
    bookmarksopen=true,
    bookmarksopenlevel=3,
    bookmarksnumbered,
    plainpages=false,
    pdfpagelabels=true,
    hyperfootnotes,
    pdftitle ={\docTitle},
    pdfauthor={\docAuthor},
    pdfcreator={\docAuthor}
]{hyperref}

\usepackage[onehalfspacing]{setspace}       % default 1.5 spacing
\addtokomafont{chapter}{\normalfont\bfseries} %\sffamily\large\bfseries
\addtokomafont{section}{\normalfont\bfseries} %\sffamily\normalsize\bfseries
\addtokomafont{subsection}{\normalfont\bfseries} 
\addtokomafont{subsubsection}{\normalfont\bfseries} %\sffamily\normalsize\mdseries
\addtokomafont{caption}{\normalfont\mdseries} %\sffamily\normalsize\mdseries

\usepackage{relsize}

\renewcommand{\texttt}{\ttfamily\smaller\relax}

\usepackage[
	left=4cm,
	right=2cm,
	top=2cm,
	bottom=2cm,
        includehead,
        includefoot,
]{geometry}

\usepackage[
]{scrlayer-scrpage}
\clearpairofpagestyles                  % clear default settings
\addtokomafont{pagefoot}{\normalfont}
\ofoot*{\thepage}                       % page number
\ohead*{\leftmark}                      % chapter name

\graphicspath{{pictures/}}  % default path

\usepackage{multirow}   % multi row and multi column table functionality
\usepackage{pbox}

\counterwithout{table}{chapter}

\usepackage{beramono}   % use a typewriter font which supports bold characters
\usepackage{enumitem}
\usepackage{chngcntr}
\counterwithout{figure}{chapter} % keine Rücksetzung pro Kapitel

\definecolor{codegreen}{rgb}{0,0.6,0}
\definecolor{codegray}{rgb}{0.5,0.5,0.5}
\definecolor{codepurple}{rgb}{0.5,0,0.33}
\definecolor{codepurblue}{rgb}{0.16,0.0,1.0}
\definecolor{backcolour}{rgb}{0.95,0.95,0.92}

\lstdefinestyle{codestyle}{
    backgroundcolor=\color{backcolour},   
    commentstyle=\color{codegreen},
    keywordstyle=\bfseries\color{codepurple},
    numberstyle=\tiny\color{codegray},
    stringstyle=\color{codepurblue},
    basicstyle=\scriptsize\ttfamily,
    breakatwhitespace=false,
    breaklines=true,
    captionpos=b,
    keepspaces=true,
    numbers=left,
    numbersep=5pt,
    showspaces=false,
    showstringspaces=false,
    showtabs=false,
    tabsize=2,
    escapeinside={(*@}{@*)}
}
\KOMAoptions{toc=chapterentrydotfill}       % dotted lines for chapters
\addtokomafont{chapterentry}{\normalfont}   % use normal font for chapter entries
\setuptoc{toc}{totoc}                       % add ToC to ToC

\DeclareTOCStyleEntry[beforeskip=0cm]{chapter}{chapter}
\DeclareTOCStyleEntry[beforeskip=0cm]{section}{section}
\DeclareTOCStyleEntry[beforeskip=0cm]{default}{subsection}
\DeclareTOCStyleEntry[beforeskip=0cm]{default}{subsubsection}
\BeforeStartingTOC[lof]{}
\BeforeStartingTOC[lot]{}
\BeforeStartingTOC[lol]{}

\usepackage{csquotes}
\RequirePackage[backend=biber,
sorting=nyt,
style=authoryear,
sortcites=true,
maxcitenames=2,
mincitenames=1,
maxbibnames=100,
date=year,
url=false,
urldate=long,
isbn=false,
firstinits=true,
doi=false,
terseinits=false,
autocite=inline,
uniquename=false,
dashed=false,
uniquelist=false,
block=none, 
hyperref=true, 
]{biblatex}

\usepackage[printonlyused]{acronym}

\usepackage[title,titletoc]{appendix}

\usepackage[nameinlink, noabbrev]{cleveref} 
\usepackage{amsthm}
\usepackage{tcolorbox}
\usepackage{float}
\usepackage{caption}

\theoremstyle{remark}

\newcommand{\Ccal}{\mathcal{C}}

\newcommand{\R}{\mathbb{R}}
\newcommand{\N}{\mathbb{N}}

\newcommand{\vecsym}[1]{\boldsymbol{#1}}     % z. B. \vecsym{x} → fett kursiv
\newcommand{\pp}[2]{\frac{\partial #1}{\partial #2}}     % z. B. \pp{f}{x}
\newcommand{\ppp}[3]{\frac{\partial^2 #1}{\partial #2 \partial #3}} % gemischte
\newcommand{\ppsq}[2]{\frac{\partial^2 #1}{\partial #2^2}}          % 2. partielle
\DeclareMathOperator{\dist}{dist}

\newcommand{\abs}[1]{\left|#1\right|}
\newcommand{\norm}[1]{\left\|#1\right\|}

\newcommand{\transpose}{^\mathsf{T}}

\DeclareMathOperator{\sech}{sech}

\newcommand{\mycomment}[1]{}                      % load preamble
\begin{document}

\begin{titlepage}
% HFU Logo
\begin{center}
    \begin{figure}[ht]
        \centering
        \includegraphics[height=2cm]{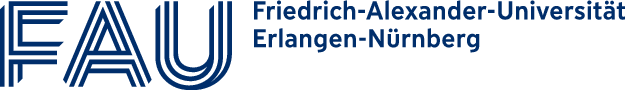}
    \end{figure}
\end{center}

\begin{center}

    {\fontsize{22}{26} \selectfont \textbf{\docTitle}}\\[5mm]
    {\fontsize{18}{22} \selectfont \docSubTitle}
    \vspace{1cm}
    {\fontsize{18}{22} \selectfont \docType}
    \vspace{4cm}
    
    \begin{tabular}{ll}

        Handed in by:       & \docAuthor   
            \\\\
        Matriculation number: & \docMatriculationNumber 
        \\\\
        First supervisor:      & \docSupervisor    \\\\
        Second supervisor: & \docSecondSupervisor \\\\
        Advisor:    & \docCoSupervisor  \\\\	
        Editing time:       & \docStartDate – \docDeadline	
    \end{tabular}
\end{center}

    \centering
    Chair of Continuous Optimization | Cauerstraße 11 | D-91058 Erlangen | www.math.fau.de/kontinuierliche-optimierung

\end{titlepage}           % title page

% Roman numbering
\frontmatter
\pagenumbering{Roman}

\chapter*{Abstract}
\addcontentsline{toc}{chapter}{Abstract}
\markboth{Abstract}{}

A feature--mapping framework for the inverse reconstruction of density--based topology optimization results is proposed. In contrast to SIMP, whose voxelized outputs are difficult to interpret or to reuse in parametric design, feature--mapping methods employ high--level geometric primitives that are mapped onto a fixed analysis grid. In the present formulation, capsule--shaped bars defined by endpoints and a radius are used as primitives. Closed--form signed distances and smooth transition functions with analytic derivatives up to second order are obtained, yielding differentiable pseudo--densities that can be consistently aggregated through smooth operators. This permits gradient--based optimization with exact Hessians. 

Several extensions are incorporated to increase robustness. Asymmetric transition functions propagate sensitivities into void regions, a reward--only objective guides features toward the target during initialization and geometric safeguards prevent degenerate configurations. The reconstruction is carried out in a staged optimization process progressing from exploration to bridging and convergence, with optional refinement schemes that adaptively add, remove or merge features based on residuals and geometric criteria.

Numerical studies on canonical SIMP benchmarks, including the five--bar and cantilever layouts, confirm that the framework reproduces target geometries with high fidelity while requiring only a moderate number of features. The staged optimization pipeline with $p$--norm and softmax aggregation consistently yields sharp reconstructions; heuristic pruning improves compactness by removing redundant features and additive refinement restores coverage in underrepresented regions. Comparisons further show that the availability of exact Hessians accelerates convergence and increases robustness relative to quasi--Newton updates. The resulting designs remain interpretable, efficient and closely aligned with the prescribed targets. These results demonstrate that analytic feature mappings can reliably reconstruct voxel--based SIMP outputs and thereby establish a bridge to explicit parametric models. The framework complements existing feature–mapping approaches such as geometry projection and moving morphable components by providing systematic strategies for initializing feature parameters. Through staged aggregation, refinement and heuristic mechanisms, it yields robust starting configurations that enable accurate and interpretable reconstructions of target density fields.            % Abstract

\tableofcontents                        % Contents
\listoffigures                          % List of Figures
\listoftables                           % List of Tables
%\listofalgorithms\markboth{List of Algorithms}{}
%\include{framework/abbreviations}       % Abbreviations

% Content
\chapter{Introduction\markboth{Introduction}{}}
\label{sec:introduction}

The aim of this thesis is to develop an inverse approximation method for results of density-based topology optimization. Rather than embedding geometric features directly into finite-element performance optimization, converged SIMP density fields are reconstructed from a compact set of explicit features with closed-form derivatives up to second order. In this way, voxelized outputs of classical topology optimization are mapped to sparse, interpretable and CAD-amenable models that are suitable for subsequent analysis or design refinement.

Topology optimization and particularly the SIMP paradigm, is widely used in structural mechanics. It handles topology changes without remeshing and supports gradient-based solvers. Yet the standard outcome of SIMP is a density field defined on a fixed grid. Boundaries are blurred by filtering and projection, geometric intent remains implicit and high-level constraints on shape or connectivity cannot easily be imposed. For industrial workflows, a substantial modeling gap remains between a SIMP layout and a parametric CAD model. The motivation of this work is therefore not to modify the SIMP formulation itself, but to define a separate reconstruction problem in which converged density fields are approximated by explicit parametric primitives.

Several methodological lines are closely related to this idea. Geometry-projection methods parameterize structures by bar-like primitives with semicircular ends and project them onto fixed grids using smooth window functions \parencite{Norato2015GeometryProjection}. This enables gradient-based optimization of compliance or other objectives with explicit control of feature dimensions. Wein and Stingl \parencite{WeinStingl2018ParamErsatz} combined such parametric representations with ersatz-material formulations, where tanh-based transitions, saturation functions for overlaps and slope and curvature regularizations were derived. Periodicity, symmetry and overhang constraints were further expressed smoothly on the level of geometric parameters. These approaches are closely related precursors, as they highlight the benefits of smooth geometry-to-density mappings and analytic sensitivities; however, they remain embedded in finite-element performance optimization, with a focus on compliance rather than on reconstructing a given density field.

The broader family of feature-mapping methods has been synthesized in the review of Wein, Dunning and Norato \parencite{wein2020review}. Two defining aspects are emphasized: the use of high-level geometric parameterizations and their mapping to non-body-fitted (fixed) grids for analysis. Within this framework, signed distance fields, smooth boundary transitions, differentiable aggregation rules and appropriate numerical quadrature constitute the key ingredients and their interplay governs regularity and stability of sensitivities. The review contrasts polynomial and sigmoid transitions and argues for analytic sensitivities in order to avoid failures of semi-analytical schemes near topology or active-set changes. This methodological basis underlies the present thesis. The pill-shaped features used here fall into the class of offset geometries with closed-form distances and first-and second-order derivatives are obtained in closed form in accordance with these recommendations. At the same time, the framework is extended by introducing asymmetric transitions, a reward-only exploration objective and staged optimization pipelines beyond what has been discussed previously.

Recent advances have also explored spline- and arc-based features. Greifenstein et al.\ \parencite{Greifenstein2023Spaghetti} introduced “spaghetti” features for continuous fiber-reinforced layouts, constructed from piecewise linear splines with smoothly rounded arcs, with an emphasis on differentiable signed distances and efficient analytic evaluation. A combination step for anisotropic material models, determining effective angles of overlapping features, was additionally proposed. The present work shares the use of explicit geometric primitives with closed-form distances, but differs in two respects: the objective is the tracking of prescribed SIMP densities rather than compliance of fiber composites and the chosen primitive is deliberately simple—a pill-shaped bar defined by two endpoints and a radius. This simplicity avoids iterative distance evaluations required for general splines and renders explicit second derivatives tractable, which in turn supports the use of Newton-type solvers.

At the interface to CAD, Shannon et al.\ \parencite{Shannon2023PostProcessing} investigated how feature-mapping outputs can be converted into fully parametric CAD models via systematic feature removal, enlargement, collinearity-based hierarchies and snapping operations to establish clean connections. While the focus lies on post-processing of moving morphable component layouts, the underlying motivation is identical: to obtain parametric, fabrication-oriented models from grid-based representations. Direct CAD output is not produced here, but the sparse and interpretable feature sets obtained in this work can be exported to such pipelines. The heuristic pruning and greedy extension strategies developed below serve as optimization-driven counterparts to these cleaning procedures.

Other branches of literature have proposed related feature-based formulations. Moving morphable bars and components (MMB/MMC) \parencite{Hoang2017MMB,Zhang2016MMC} represent structures by parametric bars or rectangles that can move, rotate and vary in thickness, with an emphasis on explicit length-scale control and manufacturability. Zhang et al.\ \parencite{Zhang2017CBS} introduced closed B-spline primitives combined by Boolean operations, offering an alternative library-based representation. 

Within this landscape, the contributions of the thesis are summarized as follows. First, the reconstruction of SIMP density fields by pill-shaped bars is formulated as a differentiable optimization problem in which signed distances, smooth transition mappings and aggregation operators form a closed analytic chain; all first-and second-order sensitivities are derived in closed form, enabling Newton-type solvers with exact Hessians. Second, an asymmetric transition extension is introduced that increases sensitivities in undercovered regions and improves convergence in staged runs. Third, the standard tracking objective is complemented by a reward-only objective used during exploration to promote initial target overlap and stable feature orientation. Fourth, a staged optimization process is designed that progresses from exploration to bridging and convergence, with optional heuristic pruning/merging and a greedy residual-based addition scheme that incrementally augments the feature set until no further improvement is achieved. Together, these elements provide a systematic mapping from voxel-based densities to explicit parametric models.

The scope of this work is deliberately focused. The study is restricted to two-dimensional density fields on Cartesian grids, with density tracking as the central objective rather than compliance or other performance criteria. Manufacturing constraints are not modeled explicitly, apart from a minimum segment length condition to prevent degeneracy. Nevertheless, the chosen primitives and the analytic framework are general and extensions to three dimensions or richer constraint sets are possible. The work should therefore be understood as a methodological contribution that demonstrates feasibility and explores the benefits of analytic second-order sensitivities, rather than as a final application-ready design tool. The remainder of the thesis is structured as follows: \Cref{sec:feature_geometry} introduces pill-shaped bars and defines their geometric representationas well as their the signed distance formulation and \Cref{sec:feature_derivatives} derives analytical first-and second-order sensitivities; \Cref{sec:transition_function} discusses smooth transition mappings, their properties and extensions to asymmetric variants; \Cref{sec:aggregation} presents aggregation operators for multiple features; \Cref{sec:grid_mapping} describes the projection of continuous pseudo-densities onto a fixed finite element grid; \Cref{sec:optimization_problem} formulates reconstruction objectives, bounds and constraints and summarizes the overall optimization model; \Cref{sec:optimization_strategy} details the staged optimization strategy, including heuristics and iterative refinement; \Cref{sec:numerical_results} reports numerical experiments on canonical benchmarks such as the five-bar and cantilever layouts; and finally, \Cref{sec:outlook} summarizes the contributions, discusses limitations and outlines directions for future research.

\mainmatter
\chapter{Geometric feature definition and distance representation\markboth{Geometric feature definition and distance representation}{}}
\label{sec:feature_geometry}

In geometry-based topology optimization, compact parametric features are employed as basic modeling entities embedded in the design domain. In this thesis, a class of capsule-shaped primitives, hereafter referred to as \emph{pill features}, is used. Each pill feature consists of a straight core segment with semicircular ends and therefore maintains a constant thickness along its medial axis. This geometry is well suited to represent the bar-like members that frequently appear in SIMP layouts.

\subsubsection*{Definition of the pill feature}
Let \(P,Q \in \R^{2}\) denote two distinct endpoints of a line segment and let \(r \in \R_{>0}\) denote a radius parameter. The pill feature \(\mathcal{F}(P,Q,r)\) is defined as the set of all points whose Euclidean distance to the segment \(\overline{PQ}\) does not exceed \(r\):
\begin{equation}
\mathcal{F}(P,Q,r) := \left\{\, \vecsym{x} \in \R^{2} \;\middle|\; \dist\bigl(\vecsym{x},\overline{PQ}\bigr) \le r \,\right\}.
\label{eq:feature_definition}
\end{equation}

Through this formulation the pill shape arises naturally. In the central region the distance is determined by the perpendicular offset to the segment, while near the ends the relevant contribution is provided by the distances to \(P\) and \(Q\). The resulting set is compact, has continuous boundary curvature and maintains a constant cross-sectional width of \(2r\) in directions orthogonal to the segment axis.

\begin{figure}[tbp]
 \centering
 \includegraphics[width=0.6\linewidth]{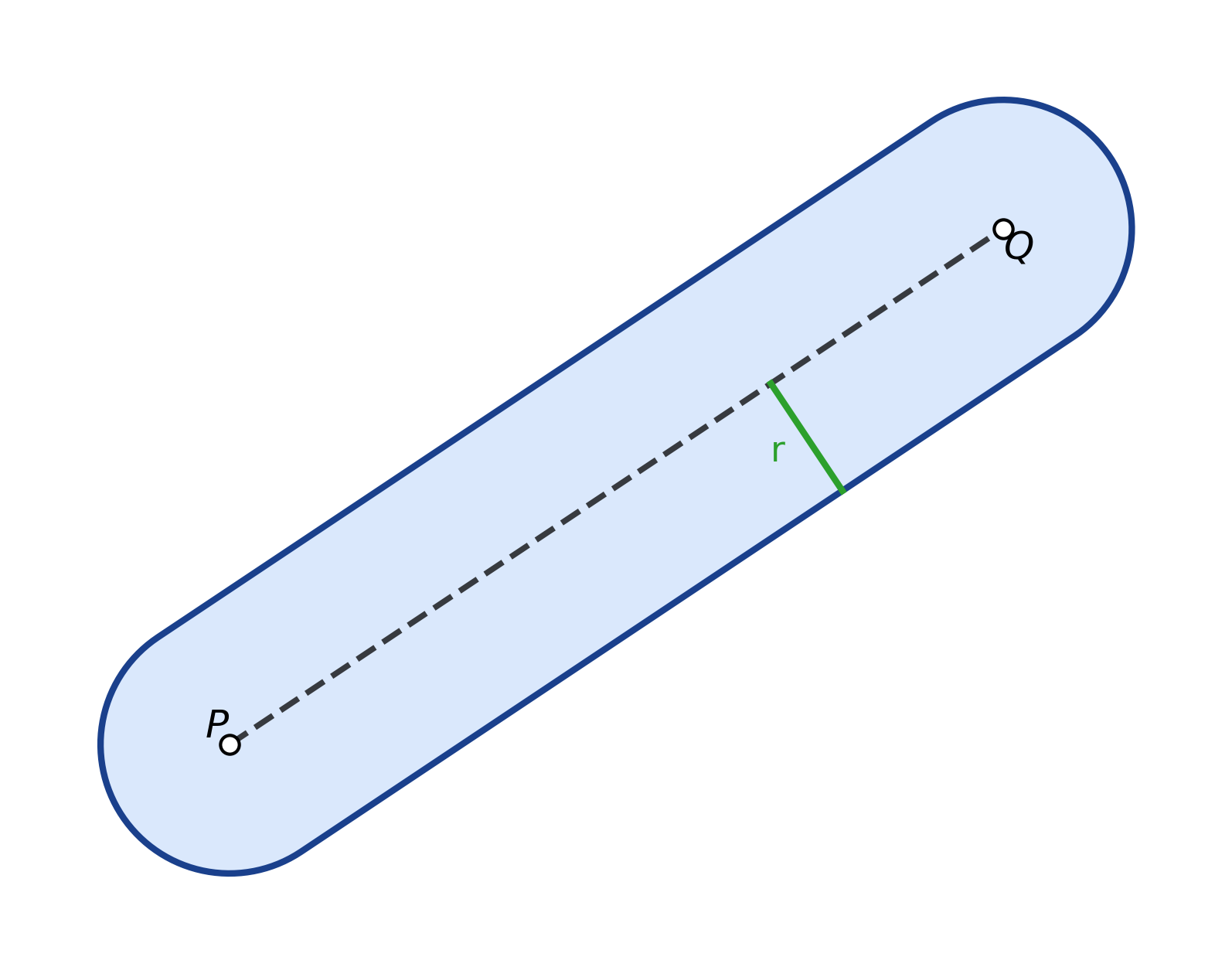}
 \caption{Illustration of pill feature \(\mathcal{F}(P,Q,r)\).}
 \label{fig:feature_geometry}
\end{figure}

\subsubsection*{Distance representation}
\begin{figure}[tbp]
 \centering
 \includegraphics[width=0.48\linewidth]{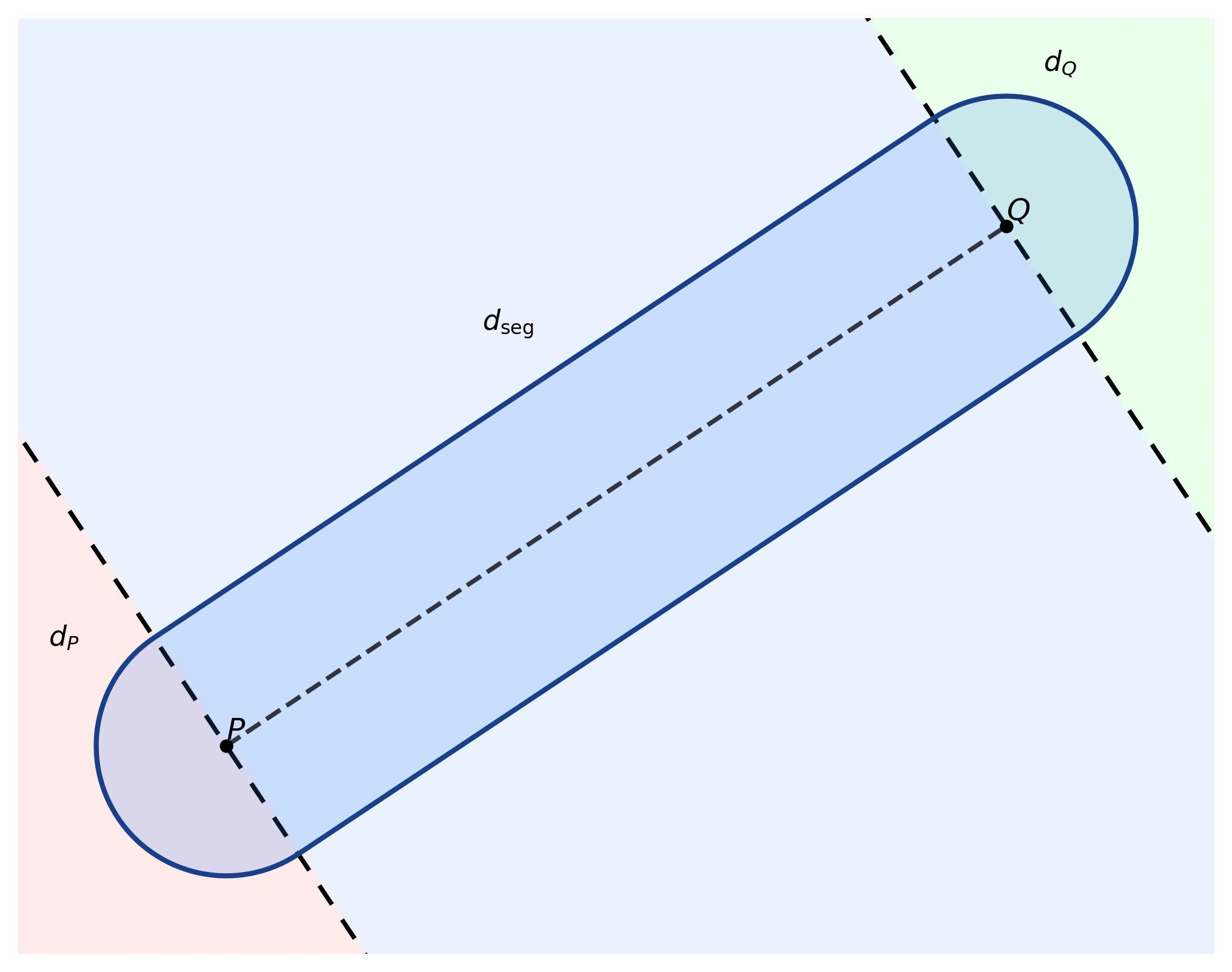}
 \includegraphics[width=0.48\linewidth]{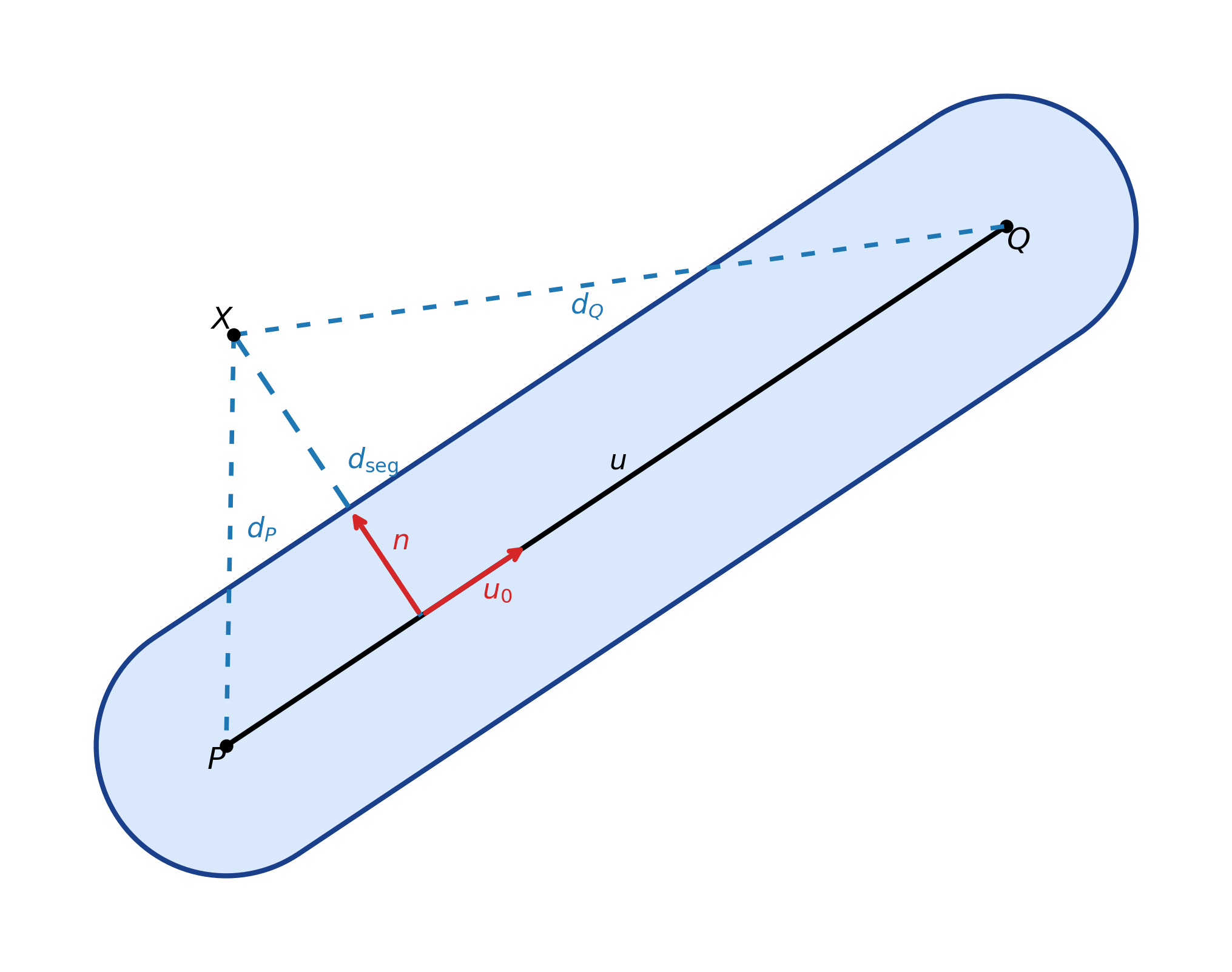}
 \caption{Left: partition of the plane into regions where $d_{\mathrm{seg}}$ (blue) or the
 endpoint distances $d_{\mathrm{pt}}(\cdot,P)$, $d_{\mathrm{pt}}(\cdot,Q)$ (red/green) are minimal.
 Right: distance decomposition for a sample point $\vecsym{x}$, showing
 $d_{\mathrm{seg}}$, $d_{\mathrm{pt}}(\vecsym{x},P)$, $d_{\mathrm{pt}}(\vecsym{x},Q)$ together with
 the segment direction $u_0$ and normal $n$.} 
 \label{fig:distance_components}
\end{figure}

To evaluate the pill definition in \cref{eq:feature_definition}, the distance \(\dist(\vecsym{x},\overline{PQ})\) must be computed explicitly. The closest point on the segment is determined by three possible cases: the orthogonal projection of \(\vecsym{x}\) falls (i) inside the segment, (ii) to the left of \(P\) or (iii) to the right of \(Q\). In the first case the perpendicular offset to the supporting line is relevant, while in the latter two cases the distances to the endpoints dominate. This decomposition is illustrated in \cref{fig:distance_components}.

Formally, let 
\[
u = Q-P \in \R^2, \qquad
u_0 = \tfrac{u}{\|u\|} \in \R^2, \qquad
n = (-u_{0,2},\,u_{0,1})\transpose \in \R^2,
\]
where $u_{0,1}$ and $u_{0,2}$ denote the first and second components of $u_0$.
By construction, $u_0$ and $n$ are unit vectors and mutually orthogonal. Throughout this thesis, the notation $\|\cdot\|$ denotes the Euclidean norm in $\R^2$ and $|\cdot|$ the absolute value in $\R$, unless stated otherwise. The relevant distance candidates are given by the perpendicular offset to the supporting line and the distances to the two endpoints:

\begin{align}
d_{\mathrm{seg}}(\vecsym{x}) &= |(\vecsym{x}-Q)\cdot n|, 
\label{eq:dseg} \\
d_{\mathrm{pt}}(\vecsym{x},P) &= \|\vecsym{x}-P\|, \qquad 
d_{\mathrm{pt}}(\vecsym{x},Q) = \|\vecsym{x}-Q\|.
\label{eq:dpt}
\end{align}
Here $d_{\mathrm{seg}}$ coincides with the true distance to the finite segment only if the orthogonal projection of $\vecsym{x}$ onto the supporting line falls within the bounds
\[
0 \le (\vecsym{x}-P)\cdot u_0 \le \|Q-P\|.
\]
Otherwise, one of the endpoint distances $d_{\mathrm{pt}}$ becomes dominant. This case distinction is naturally resolved by combining all candidates into a single expression,
\begin{equation}
d(\vecsym{x}) = \min \{\, d_{\mathrm{pt}}(\vecsym{x},P),\ d_{\mathrm{pt}}(\vecsym{x},Q),\ d_{\mathrm{seg}}(\vecsym{x}) \,\},
\label{eq:dmin}
\end{equation}
which provides the unsigned distance to the pill.

From this unsigned field the signed distance function follows directly as
\begin{equation}
d_{\mathrm{signed}}(\vecsym{x}) = d(\vecsym{x}) - r,
\label{eq:dsigned}
\end{equation}
which takes negative values inside the pill, positive values outside and vanishes exactly on its boundary, here denoted by $\Gamma(P,Q,r)$.

\begin{figure}[tbp]
 \centering
 \includegraphics[width=0.48\linewidth]{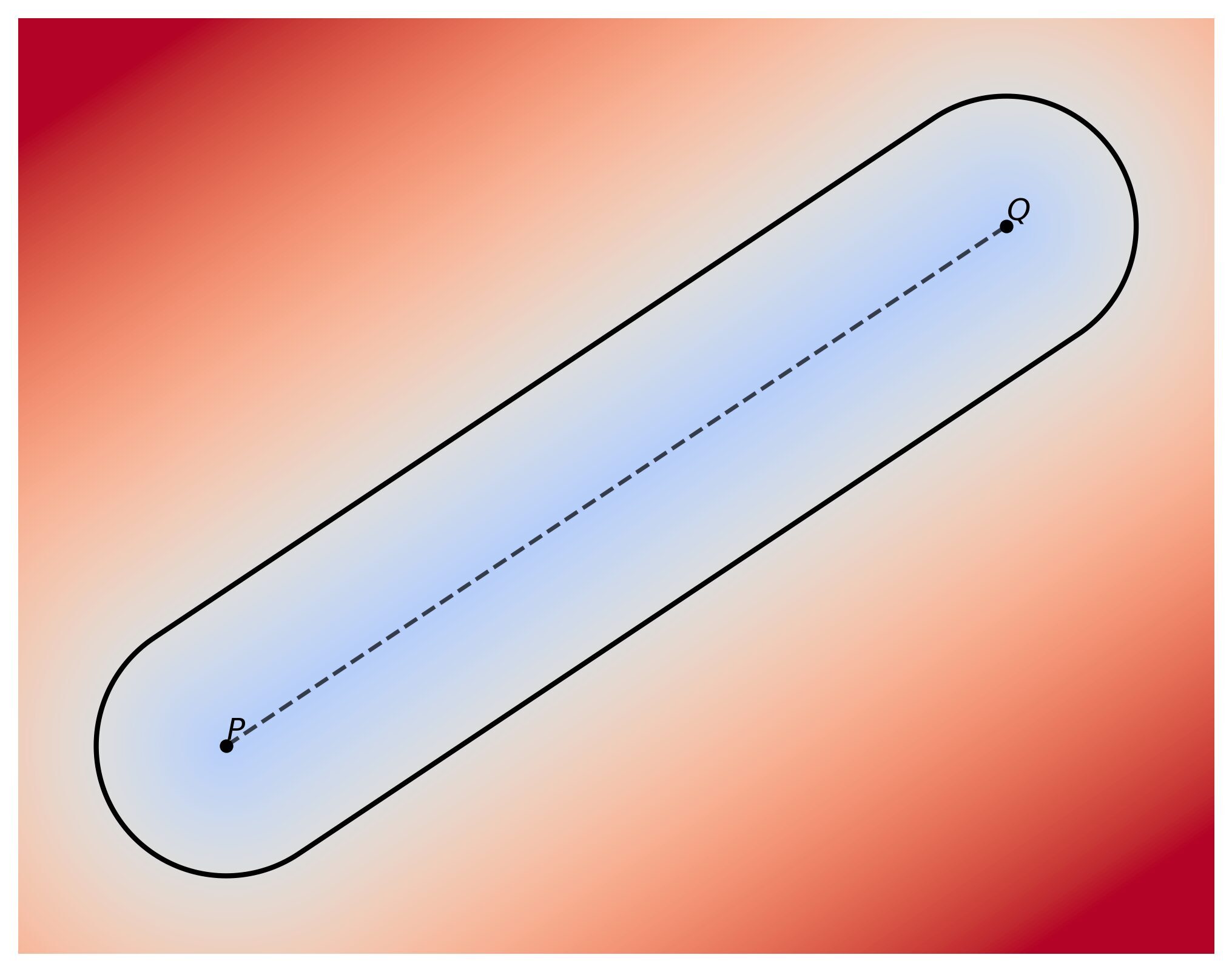}
 \caption{Signed distance field of a pill, with negative values(blue) inside, zero on the boundary and positive values(red) outside.}
 \label{fig:signed_distance_field}
\end{figure}

For subsequent use, each pill is represented by the parameter vector
\[
\vecsym{z} = (p_x,p_y,q_x,q_y,r)\transpose \in \R^{5},
\]
which collects the coordinates of the two endpoints together with the radius. This compact form serves as the design variable vector in the reconstruction problem. With the geometry and its distance representation established, the following chapter addresses the analytical computation of derivatives of the distance function up to second order.

\chapter{Analytical Derivatives of the Pill Distance\markboth{Analytical Derivatives of the Pill Distance}{}}
\label{sec:feature_derivatives}

In the preceding chapter the geometry of a pill feature and its distance representation were established. On this basis, the present chapter is devoted to the analytical derivation of sensitivities. Explicit first- and second-order derivatives of the distance function with respect to the design parameters are developed, providing quantities that are required in gradient- and Newton-type optimization methods. Since all features share the same structure, the analysis can be restricted to a single feature without loss of generality and the resulting expressions apply independently at the parametric level.
\section{Derivatives of Distance Functions}
\label{sec:derivatives_distance}

The unsigned distance to a pill has been introduced in \cref{eq:dmin} as the minimum of two point distances and a segment distance. Since the two point contributions share the same analytical structure, it suffices to derive the derivatives for a single endpoint, say \(P\). The corresponding expressions for \(Q\) follow by symmetry. Together with the derivatives of the segment distance, this provides the full set of building blocks required for the analytical treatment of the pill distance. Consider a fixed query point $\vecsym{x} = (x_1, x_2)\transpose \in \R^2$ and the endpoint $P = (p_x, p_y)\transpose$. The corresponding point distance is
\begin{equation}
 d_{\mathrm{pt}}(\vecsym{x},P)=\sqrt{(x_1-p_x)^2+(x_2-p_y)^2}.
\end{equation}
The point distance depends only on the coordinates of $P$, so derivatives with respect to $q_x,q_y$ and $r$ vanish. The nonzero first derivatives are
\begin{align}
\frac{\partial d_{\mathrm{pt}}}{\partial p_x} &=
-\frac{x_1-p_x}{\sqrt{(x_1-p_x)^2+(x_2-p_y)^2}}, &
\frac{\partial d_{\mathrm{pt}}}{\partial p_y} &=
-\frac{x_2-p_y}{\sqrt{(x_1-p_x)^2+(x_2-p_y)^2}}.
\end{align}

By differentiating these expressions once more with respect to $p_x$ and $p_y$, the nonzero second derivatives (for $\vecsym{x}\neq P$) are obtained as
\begin{align}
\frac{\partial^2 d_{\mathrm{pt}}}{\partial p_x^2} &=
\frac{(x_2-p_y)^2}{\bigl[(x_1-p_x)^2+(x_2-p_y)^2\bigr]^{3/2}}, &
\frac{\partial^2 d_{\mathrm{pt}}}{\partial p_y^2} &=
\frac{(x_1-p_x)^2}{\bigl[(x_1-p_x)^2+(x_2-p_y)^2\bigr]^{3/2}}, \\
\frac{\partial^2 d_{\mathrm{pt}}}{\partial p_x\,\partial p_y} &=
-\,\frac{(x_1-p_x)(x_2-p_y)}{\bigl[(x_1-p_x)^2+(x_2-p_y)^2\bigr]^{3/2}}.
\end{align}

All remaining mixed derivatives with respect to other parameters vanish. For completeness, these entries can be written in matrix form as
\[
\nabla^2_{(p_x,p_y)}\, d_{\mathrm{pt}}(\vecsym{x},P)
= \frac{1}{\norm{\vecsym{x}-P}}\,I_2 \;-\; \frac{(\vecsym{x}-P)(\vecsym{x}-P)^\top}{\norm{\vecsym{x}-P}^{3}},
\]
which is symmetric and positive semidefinite (and becomes singular at $\vecsym{x}=P$).
The same structure applies to the second endpoint $Q=(q_x,q_y)\transpose$, with $p_x,p_y$ replaced by $q_x,q_y$, so that the gradients and Hessians of both point distances are available as building blocks for the analytic derivatives of the full pill distance.

\section{Derivatives of the Segment Distance}

The segment distance introduced in \cref{eq:dseg} is now expressed in parametrized form with respect to the design vector $\vecsym{z}=(p_x,p_y,q_x,q_y,r)\transpose \in \R^5$. For a query point $\vecsym{x}\in\R^2$ and endpoints $P=(p_x,p_y)\transpose$, $Q=(q_x,q_y)\transpose$, the perpendicular distance to the supporting line becomes
\begin{equation}
 d_{\mathrm{seg}}(\vecsym{x},\vecsym{z}) =
 \abs{\frac{(x_1-q_x)(p_y-q_y) + (x_2-q_y)(q_x-p_x)}
       {\norm{(p_x-q_x,\;p_y-q_y)}}}.
 \label{eq:segment_distance_param}
\end{equation}
Compared to the compact form in \cref{eq:dseg}, this representation makes the explicit dependence on the optimization variables visible and is therefore better suited for derivative calculations. In order to prepare the differentiation, the expression in \cref{eq:segment_distance_param} is rewritten in quotient form with a numerator and a denominator,
\[
\begin{aligned}
N(\vecsym{x},\vecsym{z}) &= (x_1 - q_x)(p_y - q_y) + (x_2 - q_y)(q_x - p_x), \\[0.3em]
D(\vecsym{z}) &= \sqrt{(p_x - q_x)^2 + (p_y - q_y)^2}
               = \bigl\|(p_x - q_x,\; p_y - q_y)\bigr\|.
\end{aligned}
\]

so that
\[
d_{\mathrm{seg}}(\vecsym{x},\vecsym{z}) 
= \frac{\abs{N(\vecsym{x},\vecsym{z})}}{D(\vecsym{z})}.
\]
Since the denominator $D(\vecsym{z})$ represents the Euclidean distance between $P$ and $Q$, it is strictly positive for $P\neq Q$. The absolute value therefore applies only to the numerator $N(\vecsym{x},\vecsym{z})$. Writing 
\[
d_{\mathrm{seg}}(\vecsym{x},\vecsym{z}) = \frac{\abs{N(\vecsym{x},\vecsym{z})}}{D(\vecsym{z})}
\]
separates the linear and quadratic contributions of the endpoint coordinates and prepares the expression for differentiation. Gradients and Hessians can then be derived by considering $N$ and $D$ individually and combining their contributions through the quotient rule. For a generic design parameter $z_i$ this takes the form
\[
\frac{\partial d_{\mathrm{seg}}}{\partial z_i}(\vecsym{x},\vecsym{z}) 
= \frac{\operatorname{sgn}(N(\vecsym{x},\vecsym{z}))\,\frac{\partial N}{\partial z_i}(\vecsym{x},\vecsym{z}) \cdot D(\vecsym{z}) 
   - \lvert N(\vecsym{x},\vecsym{z}) \rvert \,\frac{\partial D}{\partial z_i}(\vecsym{z})}
    {D(\vecsym{z})^2}.
\]

The required partial derivatives of $N$ and $D$ with respect to the endpoint coordinates are:

\small
\[
\begin{array}{c|ccccc}
 & p_x & p_y & q_x & q_y & r \\ \hline
\partial N/\partial(\cdot) 
 & -(x_2-q_y) 
 & (x_1-q_x) 
 & -(p_y-q_y)+(x_2-q_y) 
 & (p_x-q_x)-(x_1-q_x) 
 & 0 \\
\partial D/\partial(\cdot) 
 & \tfrac{p_x-q_x}{D} 
 & \tfrac{p_y-q_y}{D} 
 & \tfrac{q_x-p_x}{D} 
 & \tfrac{q_y-p_y}{D} 
 & 0
\end{array}
\]
\normalsize

All derivatives with respect to $r$ vanish identically. Substitution of these entries into the general derivative of generic design parameter $z_i$ yields the explicit expressions
\[
\begin{aligned}
\frac{\partial d_{\mathrm{seg}}}{\partial p_x} 
&= \frac{\operatorname{sgn}(N)\,\bigl(-(x_2-q_y)\bigr) D^2 - (p_x-q_x)\,N}{D^3}, \\[0.5em]
\frac{\partial d_{\mathrm{seg}}}{\partial p_y} 
&= \frac{\operatorname{sgn}(N)\,(x_1-q_x) D^2 - (p_y-q_y)\,N}{D^3}, \\[0.5em]
\frac{\partial d_{\mathrm{seg}}}{\partial q_x} 
&= \frac{\operatorname{sgn}(N)\,\bigl(-(p_y-q_y)+(x_2-q_y)\bigr) D^2 - (q_x-p_x)\,N}{D^3}, \\[0.5em]
\frac{\partial d_{\mathrm{seg}}}{\partial q_y} 
&= \frac{\operatorname{sgn}(N)\,\bigl((p_x-q_x)-(x_1-q_x)\bigr) D^2 - (q_y-p_y)\,N}{D^3}.
\end{aligned}
\]
Second-order derivatives of the segment distance follow by differentiating the gradient
expression once more. Recall
\[
\frac{\partial d_{\mathrm{seg}}}{\partial z_i} 
= \frac{\operatorname{sgn}(N)\,\frac{\partial N}{\partial z_i}\,D 
   - \abs{N}\,\frac{\partial D}{\partial z_i}}{D^{2}},
\]
where $N=N(\vecsym{x},\vecsym{z})$ and $D=D(\vecsym{z})$. 

Differentiating this expression with respect to a second parameter $z_k$ gives
\[
\frac{\partial^2 d_{\mathrm{seg}}}{\partial z_i\,\partial z_k} 
= \frac{\partial}{\partial z_k}\!\left(
  \frac{\operatorname{sgn}(N)\,\tfrac{\partial N}{\partial z_i}\,D}{D^2}\right)
 - \frac{\partial}{\partial z_k}\!\left(
  \frac{\abs{N}\,\tfrac{\partial D}{\partial z_i}}{D^2}\right).
\]

Each term is treated separately. For the first,
\[
\frac{\partial}{\partial z_k}\!\left(
  \frac{\operatorname{sgn}(N)\,\tfrac{\partial N}{\partial z_i}\,D}{D^2}\right)
= \frac{\operatorname{sgn}(N)\,\tfrac{\partial^2 N}{\partial z_i\partial z_k}\,D}{D^2}
+ \frac{\operatorname{sgn}(N)\,\tfrac{\partial N}{\partial z_i}\,\tfrac{\partial D}{\partial z_k}}{D^2}
- 2\,\frac{\operatorname{sgn}(N)\,\tfrac{\partial N}{\partial z_i}\,D\,\tfrac{\partial D}{\partial z_k}}{D^3}.
\]

For the second,
\[
\frac{\partial}{\partial z_k}\!\left(
  \frac{\abs{N}\,\tfrac{\partial D}{\partial z_i}}{D^2}\right)
= \frac{\operatorname{sgn}(N)\,\tfrac{\partial N}{\partial z_k}\,\tfrac{\partial D}{\partial z_i}}{D^2}
+ \frac{\abs{N}\,\tfrac{\partial^2 D}{\partial z_i \partial z_k}}{D^2}
- 2\,\frac{\abs{N}\,\tfrac{\partial D}{\partial z_i}\,\tfrac{\partial D}{\partial z_k}}{D^3}.
\]

Combining both contributions yields
\begin{align}
\frac{\partial^2 d_{\mathrm{seg}}}{\partial z_i \partial z_k}
&= \frac{\operatorname{sgn}(N)\,\tfrac{\partial^2 N}{\partial z_i\partial z_k}}{D}
- \frac{\operatorname{sgn}(N)\,\tfrac{\partial N}{\partial z_i}\,\tfrac{\partial D}{\partial z_k}
   + \operatorname{sgn}(N)\,\tfrac{\partial N}{\partial z_k}\,\tfrac{\partial D}{\partial z_i}}{D^2} \notag \\
&\quad + \frac{\abs{N}}{D^2}\left(\frac{2\,\tfrac{\partial D}{\partial z_i}\,\tfrac{\partial D}{\partial z_k}}{D}
- \tfrac{\partial^2 D}{\partial z_i \partial z_k}\right).
\label{eq:hessian_dseg_general}
\end{align}

This representation expresses the Hessian in terms of first and second derivatives of $N$ and $D$. The first derivatives were listed above; for completeness the second derivatives are summarized below, noting that most entries vanish due to the bilinear structure of $N$ and the quadratic form of $D$.

\[
\begin{array}{c|ccccc}
 & p_x & p_y & q_x & q_y & r \\ \hline
\partial^2 N/\partial(\cdot)\partial p_x & 0 & 0 & 0 & \;\;1 & 0 \\
\partial^2 N/\partial(\cdot)\partial p_y & 0 & 0 & -1 & 0 & 0 \\
\partial^2 N/\partial(\cdot)\partial q_x & 0 & -1 & 0 & 0 & 0 \\
\partial^2 N/\partial(\cdot)\partial q_y & 1 & 0 & 0 & 0 & 0 \\
\partial^2 N/\partial(\cdot)\partial r & 0 & 0 & 0 & 0 & 0
\end{array}
\]

The bilinear structure of $N$ implies that all pure second derivatives vanish and only the mixed terms $\partial^2 N/\partial p_x\partial q_y=+1$ and $\partial^2 N/\partial p_y\partial q_x=-1$ remain nonzero. 
\small
\[
\begin{array}{c|ccccc}
 & p_x & p_y & q_x & q_y & r \\ \hline
\partial^2 D/\partial(\cdot)\partial p_x 
 & \tfrac{(p_y-q_y)^2}{D^3} 
  & -\tfrac{(p_x-q_x)(p_y-q_y)}{D^3} 
  & -\tfrac{(p_y-q_y)^2}{D^3} 
  & \tfrac{(p_x-q_x)(p_y-q_y)}{D^3} 
  & 0 \\
  \partial^2 D/\partial(\cdot)\partial p_y 
  & -\tfrac{(p_x-q_x)(p_y-q_y)}{D^3} 
  & \tfrac{(p_x-q_x)^2}{D^3} 
  & \tfrac{(p_x-q_x)(p_y-q_y)}{D^3} 
  & -\tfrac{(p_x-q_x)^2}{D^3} 
  & 0 \\
  \partial^2 D/\partial(\cdot)\partial q_x 
  & -\tfrac{(p_y-q_y)^2}{D^3} 
  & \tfrac{(p_x-q_x)(p_y-q_y)}{D^3} 
  & \tfrac{(p_y-q_y)^2}{D^3} 
  & -\tfrac{(p_x-q_x)(p_y-q_y)}{D^3} 
  & 0 \\
  \partial^2 D/\partial(\cdot)\partial q_y 
  & \tfrac{(p_x-q_x)(p_y-q_y)}{D^3} 
  & -\tfrac{(p_x-q_x)^2}{D^3} 
  & -\tfrac{(p_x-q_x)(p_y-q_y)}{D^3} 
  & \tfrac{(p_x-q_x)^2}{D^3} 
  & 0 \\
  \partial^2 D/\partial(\cdot)\partial r 
   & 0 & 0 & 0 & 0 & 0
\end{array}
\]
\normalsize
With the first- and second-order derivatives of $N$ and $D$ tabulated above, all entries of the Hessian can now be assembled systematically by substitution into \cref{eq:hessian_dseg_general}. To illustrate the structure of the resulting expressions, three representative cases are given explicitly. The first case concerns the pure second derivative with respect to the same variable, here $p_x$, the second a mixed derivative with respect to both coordinates of the same endpoint $(p_x,p_y)$ and the third a mixed derivative across the two endpoints $(p_x,q_y)$. These examples cover the typical patterns, while all other entries follow analogously.

\begin{align}
 \frac{\partial^2 d_{\mathrm{seg}}}{\partial p_x^2}
 &= \frac{2\,\operatorname{sgn}(N)\,(x_2-q_y)(p_x-q_x)}{D^3}
   + \abs{N}\,\frac{2(p_x-q_x)^2-(p_y-q_y)^2}{D^5}, 
 \label{eq:dseg_pp}\\[6pt]
 \frac{\partial^2 d_{\mathrm{seg}}}{\partial p_x\,\partial p_y}
 &= \frac{\operatorname{sgn}(N)}{D^3}
   \Bigl[(x_2-q_y)(p_y-q_y)-(x_1-q_x)(p_x-q_x)\Bigr] \notag\\
 &\quad + \abs{N}\,\frac{3(p_x-q_x)(p_y-q_y)}{D^5}, 
 \label{eq:dseg_ppy}\\[6pt]
 \frac{\partial^2 d_{\mathrm{seg}}}{\partial p_x\,\partial q_y}
 &= \frac{\operatorname{sgn}(N)}{D}
   - \frac{\operatorname{sgn}(N)}{D^3}
    \Bigl[(x_2-q_y)(p_y-q_y)+(p_x-x_1)(p_x-q_x)\Bigr] \notag\\
 &\quad - \abs{N}\,\frac{3(p_x-q_x)(p_y-q_y)}{D^5}.
 \label{eq:dseg_pq}
\end{align}

These formulas exemplify the algebraic structure of the Hessian entries: they consist of sign-weighted bilinear terms originating from $N$ with denominators of order $D^3$, together with rational contributions proportional to $\abs{N}$ that involve quadratic factors of the endpoint coordinates and denominators of order $D^5$. By symmetry, the corresponding permuted indices follow directly and in combination with the tabulated derivatives this provides the complete Hessian of the segment distance function.

\section{Smoothness of the Combined Distance Function}
\label{sec:distance_smoothness}

The minimum of several smooth functions is in general continuous but not
differentiable at points where two or more arguments coincide. In the present
setting, the unsigned pill distance
\[
d(\vecsym{x},\vecsym{z}) := \min\{\, d_{\mathrm{pt}}(\vecsym{x},P),\ d_{\mathrm{seg}}(\vecsym{x},\vecsym{z}),\ d_{\mathrm{pt}}(\vecsym{x},Q)\,\},
\]
involves two point-to-point distances and the perpendicular segment distance.
The goal is to show that, despite the minimum operator, the combined field is
twice continuously differentiable away from a singular set. In particular,
\[
d \in C^2\bigl((\R^2 \times \R^5)\setminus \mathcal S\bigr),
\qquad 
\mathcal S := \{\vecsym{x}=P\}\cup\{\vecsym{x}=Q\}\cup\{N(\vecsym{x},\vecsym{z})=0\},
\]
where the last set consists of points lying on the supporting line, so that $d_{\mathrm{seg}}=\abs{N}/D$ inherits the kink of the absolute value.
Note that $\operatorname{sgn}(N)$ is undefined at $N=0$; all derivative expressions involving $\operatorname{sgn}(N)$ are therefore understood on $(\R^2\times\R^5)\setminus\mathcal S$.

Potential nonsmoothness can occur only along transition sets where two of the
contributing functions attain the same value. At these locations it must be
verified that not only the function values but also the first and second
derivatives of the competing terms coincide, ensuring continuous differentiability
of the overall minimum. Since $d$ is defined as the minimum of the two endpoint
distances and the segment distance, possible transitions involve one of the point
distances meeting the segment contribution. The analysis may therefore be
restricted to the endpoint $P$; the case of $Q$ follows analogously.

To characterize the transition set between the point-based distance $d_{\mathrm{pt}}(\vecsym{x},P)$ and the segment-based distance $d_{\mathrm{seg}}(\vecsym{x},\vecsym{z})$, the condition
\[
d_{\mathrm{pt}}(\vecsym{x},P) = d_{\mathrm{seg}}(\vecsym{x},\vecsym{z})
\]
is considered. In coordinates, the transition condition takes the form
\[
\sqrt{(x_1-p_x)^2+(x_2-p_y)^2}
= \left|\frac{(x_1-q_x)(p_y-q_y) + (x_2-q_y)(q_x-p_x)}
       {\sqrt{(p_x-q_x)^2+(p_y-q_y)^2}}\right|,
\]
where the denominator denotes the segment length $D$.
By squaring both sides, the absolute values cancel and one obtains
\[
(x_1-p_x)^2+(x_2-p_y)^2
= \frac{\bigl((x_1-q_x)(p_y-q_y) + (x_2-q_y)(q_x-p_x)\bigr)^2}{D^2}.
\]
The right-hand side can be recognized as the squared area formula of the
parallelogram spanned by the vectors $(x_1-q_x,\;x_2-q_y)$ and $(p_x-q_x,\;p_y-q_y)$,
divided by $D^2$. Hence the equality requires that the vector
\[
(x_1-p_x,\;x_2-p_y)
\]
is orthogonal to $(p_x-q_x,\;p_y-q_y)$ and has the same length as its projection in that orthogonal direction. It follows that the solution set is precisely the line through $P$ in the direction of the unit normal $(n_x,n_y)$ of the segment, i.e.
\[
(x_1,x_2) = (p_x, p_y) + s\,(n_x,n_y), \qquad s\in\R.
\]

Substitution of this parametrization into the defining equation then yields
\[
d_{\mathrm{pt}}(\vecsym{x},P)=\abs{s}, 
\qquad
d_{\mathrm{seg}}(\vecsym{x},\vecsym{z})=\abs{s},
\]
so that both distances coincide on all points of
\[
\mathcal{T}_P := \{\, (p_x+s\,n_x,\;p_y+s\,n_y) \mid s\in\R \,\}.
\]
 
For points of this form, the relations
\[
d_{\mathrm{pt}}(\vecsym{x},P) = \abs{s},
\qquad
d_{\mathrm{seg}}(\vecsym{x},\vecsym{z})
= \frac{\abs{s}\sqrt{n_x^2+n_y^2}\,\sqrt{(p_x-q_x)^2+(p_y-q_y)^2}}{D}
\]
are obtained. Since $(n_x,n_y)$ is normalized, the numerator simplifies to $\abs{s}D$, which cancels with the denominator and thus
\[
d_{\mathrm{seg}}(\vecsym{x},\vecsym{z})=\abs{s}
\]
follows.

The agreement of gradients and Hessians on the transition set is verified by
substituting $\mathcal{T}_P$ into the explicit formulas for $d_{\mathrm{pt}}$ with
respect to $p_x$ and $p_y$ and into the general expressions for $d_{\mathrm{seg}}$.
Along $\mathcal{T}_P$ one has
\[
(x_1,x_2)=(p_x,p_y)+s(n_x,n_y), \qquad s\in\R,
\]
together with $x_1-p_x=s\,n_x$, $x_2-p_y=s\,n_y$, 
$\lvert N(\vecsym{x},\vecsym{z})\rvert = \abs{s}D$ and
$\operatorname{sgn}(N(\vecsym{x},\vecsym{z}))=\operatorname{sgn}(s)$.
Since $D$ does not depend on the query point, these relations are sufficient to
evaluate all derivatives on $\mathcal{T}_P$.

Insertion into the point-based gradient yields
\[
\frac{\partial d_{\mathrm{pt}}}{\partial p_x}
=-\,\operatorname{sgn}(s)\,n_x,
\qquad
\frac{\partial d_{\mathrm{pt}}}{\partial p_y}
=-\,\operatorname{sgn}(s)\,n_y,
\]
while evaluation of the segment-based gradient with the tabulated expressions for
$\partial N/\partial(\cdot)$ and $\partial D/\partial(\cdot)$ gives
\[
\frac{\partial d_{\mathrm{seg}}}{\partial p_x}
=-\,\operatorname{sgn}(s)\,n_x,
\qquad
\frac{\partial d_{\mathrm{seg}}}{\partial p_y}
=-\,\operatorname{sgn}(s)\,n_y.
\]
Hence the gradients of both distance functions coincide everywhere on
$\mathcal{T}_P$ for $s\neq 0$.

A similar argument applies to the second derivatives. Substitution of
$x_1-p_x=s\,n_x$ and $x_2-p_y=s\,n_y$ into the explicit Hessian of the point
distance yields
\[
\nabla_{(p_x,p_y)}^2 d_{\mathrm{pt}}
= \frac{1}{\abs{s}}
\begin{pmatrix}
n_y^2 & -\,n_x n_y \\
-\,n_x n_y & n_x^2
\end{pmatrix}.
\]
Evaluation of the general Hessian expression for the segment distance with the
same substitutions and the tabulated derivatives of $N$ and $D$ leads to the
identical result,
\[
\nabla_{(p_x,p_y)}^2 d_{\mathrm{seg}}
= \frac{1}{\abs{s}}
\begin{pmatrix}
n_y^2 & -\,n_x n_y \\
-\,n_x n_y & n_x^2
\end{pmatrix}.
\]

It follows that both gradient and Hessian agree exactly on the transition set
$\mathcal{T}_P$ for all $s\neq 0$. At $s=0$, corresponding to $\vecsym{x}=P$, the
well-known singularity of order $1/\abs{s}$ appears in both expressions, so that
matching limits are obtained although a local $C^2$ extension cannot be defined
at that point. The second derivatives of $d_{\mathrm{pt}}$ and $d_{\mathrm{seg}}$ thus
coincide exactly along the interface and both Hessians are identical on
$\mathcal{T}_P\setminus\{P\}$. Combined with the agreement of first derivatives,
this confirms that the composed base distance is $C^2$ in the spatial and parametric variables
away from $\mathcal S$. This smoothness is sufficient for the construction of the transition
mappings in the following chapter. As argued at the beginning of the analysis, the discussion may be restricted to the endpoint $P$; the case of $Q$ follows analogously by symmetry of the construction. With this, the section is concluded and the development of transition functions is addressed next.
\chapter{Transition functions\markboth{Transition functions}{}}
\label{sec:transition_function}

With the smoothness of the signed distance field established in \cref{sec:distance_smoothness}, the focus turns to material representation. To transform distance values into differentiable pseudo-densities $\rho(\vecsym{x})\in[0,1]$, transition functions $\phi:\R\to[0,1]$ are employed. These mappings assign material fractions to signed distances, thereby regularizing the binary indicator of $\mathcal{F}(P,Q,r)$ and enabling gradient-based optimization on a fixed grid. The following sections introduce their definition, discuss representative families and analyze their derivatives for use in optimization.

\section{Definition and properties}
\label{sec:transition_def}

Transition functions convert the signed distance field into a smooth pseudo-density field. Let \(d_{\mathrm{signed}}(\vecsym{x})=d(\vecsym{x})-r\), abbreviated by \(d\) when unambiguous.
A transition function with half-width \(\delta>0\) is defined by
\begin{equation}
\phi_\delta(d)=
\begin{cases}
1 & d \le -\delta,\\
\phi(d) & -\delta<d<\delta,\\
0 & d \ge \delta,
\end{cases}
\label{eq:clipped_phi}
\end{equation}
where \(\phi:(-\delta,\delta)\to(0,1)\) is a smooth core mapping. Thus \(\phi_\delta\) coincides with the binary indicator outside the transition zone \([-\delta,\delta]\).

The core mapping is normalized at the interface by
\begin{equation}
\phi(0)=\tfrac12,
\label{eq:midvalue}
\end{equation}
and is assumed monotone (\(\phi'(d)\le 0\) on \((-\delta,\delta)\)).
No symmetry condition is imposed; many choices used here (e.g.\ tanh, smoothstep) are symmetric by construction, while asymmetric variants are introduced in \cref{sec:asym_scaled}.

For higher smoothness, \(\phi\in C^k(-\delta,\delta)\) is required for some \(k\ge 0\). If, in addition,
\begin{equation}
\phi^{(j)}(\pm\delta)=0,\qquad j=1,\dots,k,
\label{eq:boundary_conditions}
\end{equation}
then the composite field 
\[
\rho(\vecsym{x})=\phi_\delta\!\bigl(d_{\mathrm{signed}}(\vecsym{x})\bigr)
\]
satisfies \(\rho\in C^k(\R^2)\) and exhibits no artificial boundary gradients. 
The differentiability of $\phi_\delta$ is determined by the smoothness of $\phi$ on $(-\delta,\delta)$ together with the endpoint conditions~\eqref{eq:boundary_conditions}, ensuring the desired regularity of $\rho$ while retaining a strictly bounded transition zone.

\section{Tanh-based and polynomial transition functions}
\label{sec:transition_examples}

\begin{figure}[H]
\centering
\includegraphics[width=0.7\textwidth]{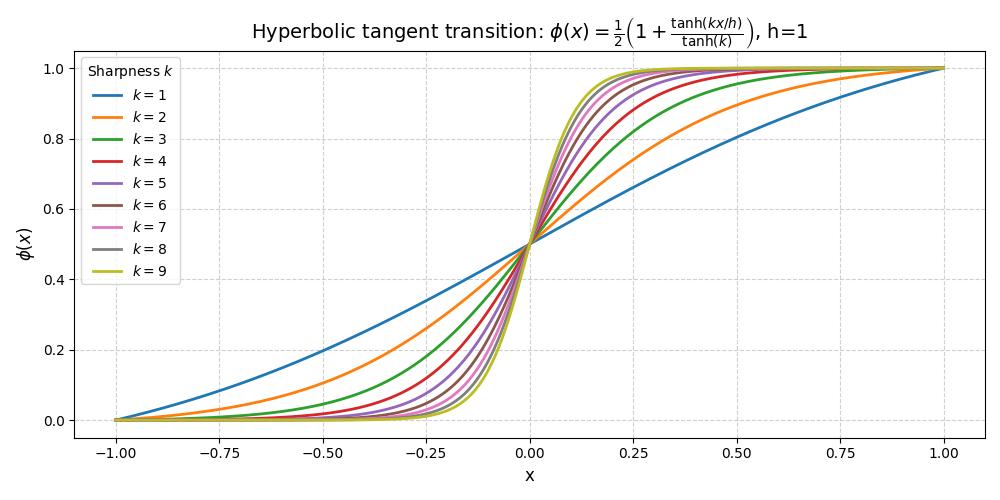}
\caption{Tanh-based transition functions $\phi_\delta(d)$ for different steepness values $\beta$.}
\label{fig:tanh_variants}
\end{figure}

Two representative choices of transition functions illustrate the design space
and highlight the trade-off between practical smoothness and exact boundary behavior. The first class is based on the hyperbolic tangent and, for a steepness
parameter $\beta>0$, is defined as
\begin{equation}
\phi(d)=\tfrac{1}{2}\,\bigl(1-\tanh(\tfrac{\beta d}{\delta})\bigr).
\label{eq:tanh_transition}
\end{equation}
This mapping satisfies $\phi(0)=\tfrac{1}{2}$, decreases monotonically and tends
asymptotically to unity and zero for $d\to-\infty$ and $d\to+\infty$. Although not
compactly supported, it becomes effectively flat outside $[-\delta,\delta]$ for
moderate $\beta$ and can therefore be combined with the clipping
rule~\eqref{eq:clipped_phi}. The parameter $\beta$ controls the steepness, with
larger values producing sharper transitions (cf.\ \cref{fig:tanh_variants}). The drawback of this construction is that the clipping introduces residual slope
mismatches at the boundaries, so that exact vanishing of derivatives is not
guaranteed.

\begin{figure}[H]
\centering
\includegraphics[width=0.7\textwidth]{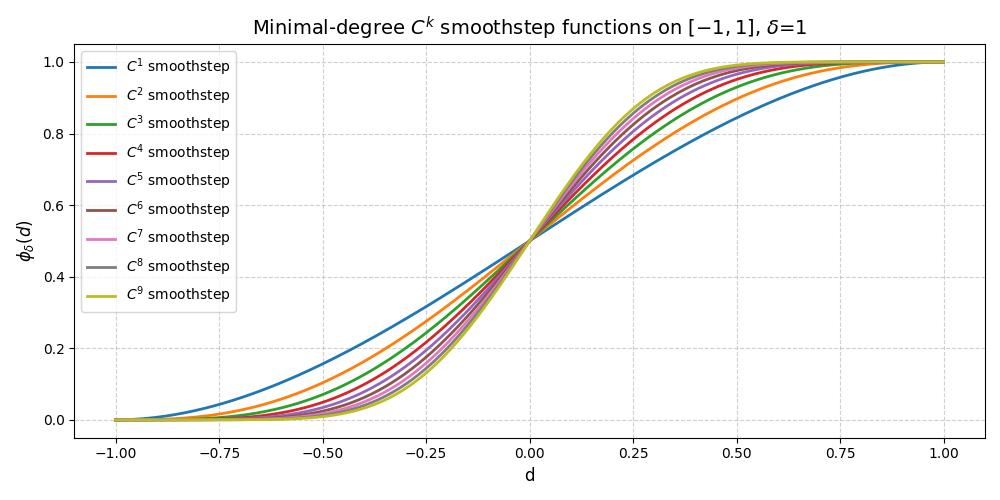}
\caption{Polynomial smoothstep functions $\phi_\delta(d)$ for different regularity classes $C^k$.}
\label{fig:smoothstep_variants}
\end{figure}

To achieve exact boundary behavior, a second family based on smoothstep polynomials
is introduced. With the normalized coordinate $t=(d+\delta)/(2\delta)\in[0,1]$, a
polynomial $p:[0,1]\to[0,1]$ defines $\phi(d):=p(t(d))$ provided the endpoint
conditions
\[
p(0)=1,\qquad p(1)=0,\qquad p^{(j)}(0)=p^{(j)}(1)=0 \quad\text{for } j=1,\dots,k
\]
are satisfied for some integer $k\ge 0$. These $2(k+1)$ conditions enforce a minimal
degree of $2k+1$, yielding smoothstep functions of class $C^k$. Unlike the tanh variant, their derivatives vanish exactly at $\pm\delta$, ensuring a smooth transition into the plateaus. Increasing $k$ improves regularity at the
boundary and flattens the tails, as shown in \cref{fig:smoothstep_variants}.

\subsubsection*{Example ($k=2$): construction of the quintic smoothstep}

For the case $k=2$, a quintic ansatz
\[
p(t)=a_5 t^5+a_4 t^4+a_3 t^3+a_2 t^2+a_1 t+a_0, \qquad t\in[0,1],
\]
is posed and the six endpoint conditions
\[
p(0)=1,\quad p(1)=0,\qquad
p'(0)=0,\quad p'(1)=0,\qquad
p''(0)=0,\quad p''(1)=0
\]
are enforced. Evaluation at $t=0$ fixes
\[
a_0=1,\qquad a_1=0,\qquad a_2=0.
\]
With these values inserted, the remaining conditions at $t=1$ become the linear system
\begin{equation}
\begin{aligned}
p(1)=0&:&\quad a_5+a_4+a_3&=-1,\\
p'(1)=0&:&\quad 5a_5+4a_4+3a_3&=0,\\
p''(1)=0&:&\quad 20a_5+12a_4+6a_3&=0.
\end{aligned}
\label{eq:smoothstep_quintic_system}
\end{equation}
Solving \eqref{eq:smoothstep_quintic_system} yields
\[
a_5=-6,\qquad a_4=15,\qquad a_3=-10,
\]
and therefore
\[
p(t)=-6t^5+15t^4-10t^3+1.
\]
Using $t=(d+\delta)/(2\delta)$ gives the explicit transition function
\[
\phi(d)=p\!\left(\frac{d+\delta}{2\delta}\right)
=-6\Bigl(\tfrac{d+\delta}{2\delta}\Bigr)^5
+15\Bigl(\tfrac{d+\delta}{2\delta}\Bigr)^4
-10\Bigl(\tfrac{d+\delta}{2\delta}\Bigr)^3+1,
\]
which is of class $C^2$, has compact support under the clipping rule \eqref{eq:clipped_phi} and exhibits
$\phi^{(j)}(\pm\delta)=0$ for $j=1,2$ as required.

\subsubsection*{General form ($k\in\N$)}

For a prescribed smoothness order $k\ge 0$, the unique minimal-degree polynomial of degree $2k+1$ that satisfies the endpoint conditions 
\[
p(0)=1,\ \ p(1)=0,\qquad
p^{(j)}(0)=0,\ \ p^{(j)}(1)=0\ \ \text{for } j=1,\dots,k
\]
can be written in closed form as
\begin{equation}
p(t)=\sum_{m=0}^{k}(-1)^m \binom{2k+1}{m}\binom{k+m}{m}\,t^{k+1+m}, 
\qquad t\in[0,1].
\label{eq:smoothstep_closedform}
\end{equation}
The corresponding transition function on the distance interval is defined by
\begin{equation}
\phi(d)=p\!\left(\frac{d+\delta}{2\delta}\right),\qquad d\in[-\delta,\delta],
\label{eq:smoothstep_general}
\end{equation}
and extended by clipping according to \eqref{eq:clipped_phi}. By construction, \eqref{eq:smoothstep_general} satisfies the boundary conditions \eqref{eq:boundary_conditions}; increasing $k$ raises the regularity and flattens the tails at the expense of a higher polynomial degree.

\subsection{Derivatives of transition functions}
\label{sec:derivatives_transition}

In the previous chapter, gradients and Hessians of the distance representation
were derived and the $C^2$ smoothness of the base field was established. In order
to preserve this level of regularity after the application of transition mappings,
the derivatives of $\phi$ must be analyzed and combined with those of the distance.
The composite field
\[
\rho(\vecsym{x};\vecsym{z})=\phi\!\bigl(d_{\mathrm{signed}}(\vecsym{x};\vecsym{z})\bigr)
\]
inherits its differentiability from both $d_{\mathrm{signed}}$ and $\phi$. Ensuring
that $\phi$ admits analytic derivatives up to second order is therefore crucial for
maintaining global $C^2$-smoothness in spatial and parametric variables. Let $d=d_{\mathrm{signed}}(\vecsym{x};\vecsym{z})$ and $z_i,z_j\in\{p_x,p_y,q_x,q_y,r\}$. Then
\begin{equation}
\pp{\rho}{z_i}=\phi'(d)\,\pp{d}{z_i},\qquad
\ppp{\rho}{z_i}{z_j}
=\phi''(d)\,\pp{d}{z_i}\,\pp{d}{z_j}
+\phi'(d)\,\ppp{d}{z_i}{z_j}.
\label{eq:rho_chain_rules_components}
\end{equation}
For diagonal entries ($i=j$) this specializes to
\begin{equation}
\ppsq{\rho}{z_i}
=\phi''(d)\,\bigl(\pp{d}{z_i}\bigr)^{2}
+\phi'(d)\,\ppsq{d}{z_i}.
\label{eq:rho_chain_rules_diagonal}
\end{equation}
For the radius one has $\pp{d}{r}=-1$, $\ppsq{d}{r}=0$ and $\ppp{d}{r}{z_j}=0$, hence
\begin{equation}
\pp{\rho}{r}=-\phi'(d),\qquad
\ppp{\rho}{r}{z_j}=-\,\phi''(d)\,\pp{d}{z_j},\qquad
\ppsq{\rho}{r}=\phi''(d).
\label{eq:rho_chain_rules_r}
\end{equation}

Both the tanh-based mapping \eqref{eq:tanh_transition} and the polynomial smoothstep \eqref{eq:smoothstep_general} admit analytic derivatives up to second order. Their explicit forms are collected below in order to assess smoothness and support properties.

For the tanh mapping, introduction of the normalized coordinate $\xi=\beta d/\delta$ with steepness parameter $\beta>0$ yields
\[
\tanh(\xi)=\frac{e^\xi-e^{-\xi}}{e^\xi+e^{-\xi}}, 
\qquad 
\sech^2(\xi)=\frac{4}{(e^\xi+e^{-\xi})^2}.
\]
Using these identities, the derivatives of $\phi$ become
\begin{equation}
\pp{\phi}{d}
= -\frac{\beta}{2\delta}\,\frac{4}{\bigl(e^{\xi}+e^{-\xi}\bigr)^{2}},
\qquad
\ppsq{\phi}{d}
= \Bigl(\frac{\beta}{\delta}\Bigr)^{\!2}\,
\frac{4\,(e^{\xi}-e^{-\xi})}{\bigl(e^{\xi}+e^{-\xi}\bigr)^{3}}.
\label{eq:tanh_derivatives}
\end{equation}
Both expressions decay exponentially as $|d|\to\infty$, so sensitivities outside the clipping band $[-\delta,\delta]$ are practically negligible.
The parameter $\beta$ directly controls the steepness: larger values sharpen the central transition, but small residual slopes remain in the clipped regions.
For the smoothstep construction of order $k$, the closed form \eqref{eq:smoothstep_general} can be differentiated directly. The first and second derivatives of $p(t)$ are
\begin{equation}
p'(t)=\sum_{m=0}^{k}(-1)^m\binom{2k+1}{m}\binom{k+m}{m}\,(k+1+m)\,t^{k+m},
\label{eq:smoothstep_polynomial_general_d1}
\end{equation}
\begin{equation}
p''(t)=\sum_{m=0}^{k}(-1)^m\binom{2k+1}{m}\binom{k+m}{m}\,(k+1+m)(k+m)\,t^{k+m-1}.
\label{eq:smoothstep_polynomial_general_d2}
\end{equation}
Inserting $\phi'(d)=p'(t)/(2\delta)$ and $\phi''(d)=p''(t)/(2\delta)^2$ into \eqref{eq:rho_chain_rules_r}
yields explicit sensitivities of $\rho$ with respect to all design parameters for arbitrary smoothness order~$k$.
\begin{figure}[H]
 \centering
 \begin{minipage}{.18\textwidth}\centering
  \includegraphics[width=\linewidth]{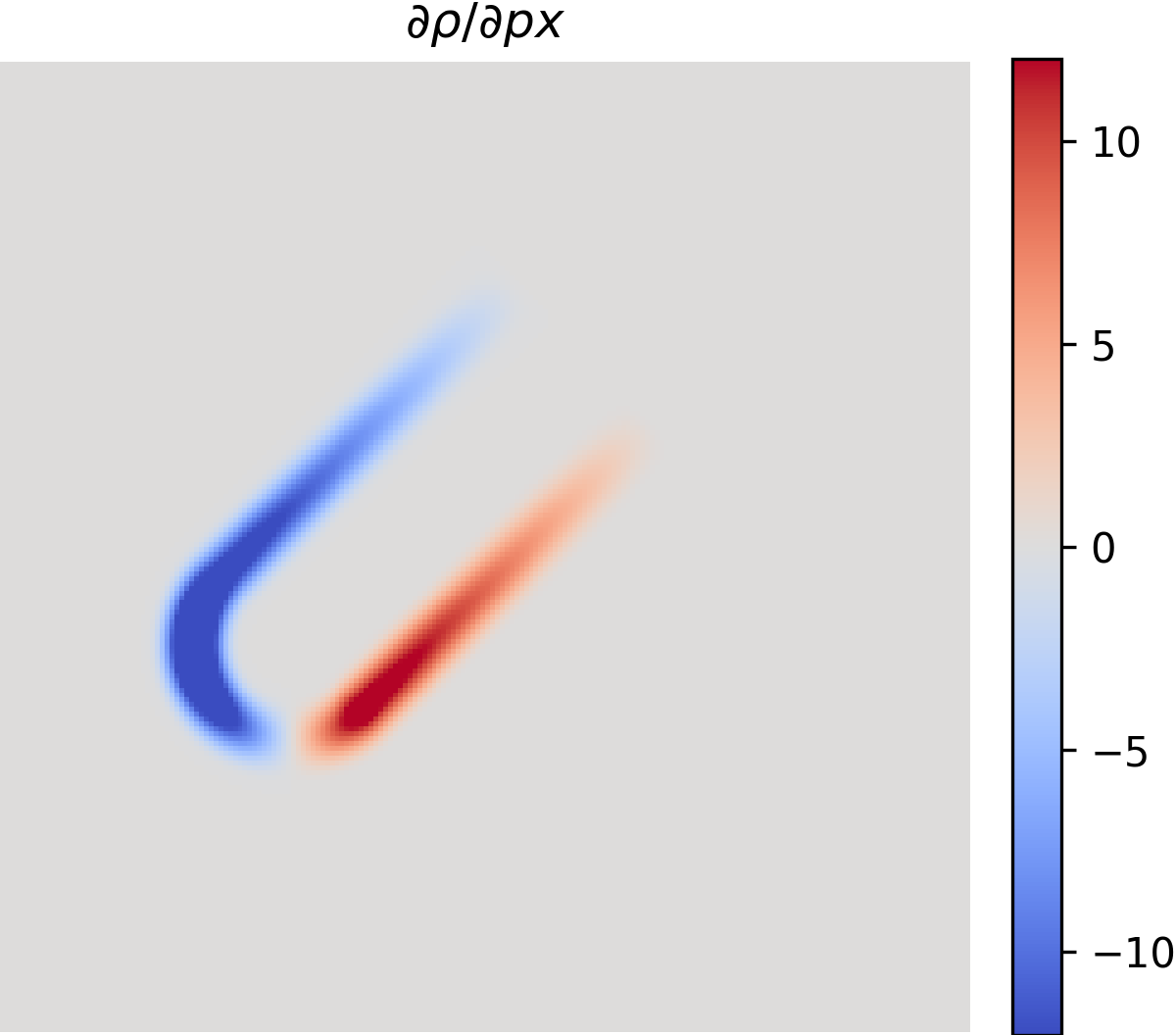}
  \caption*{$\partial\rho/\partial p_x$}
 \end{minipage}\hfill
 \begin{minipage}{.18\textwidth}\centering
  \includegraphics[width=\linewidth]{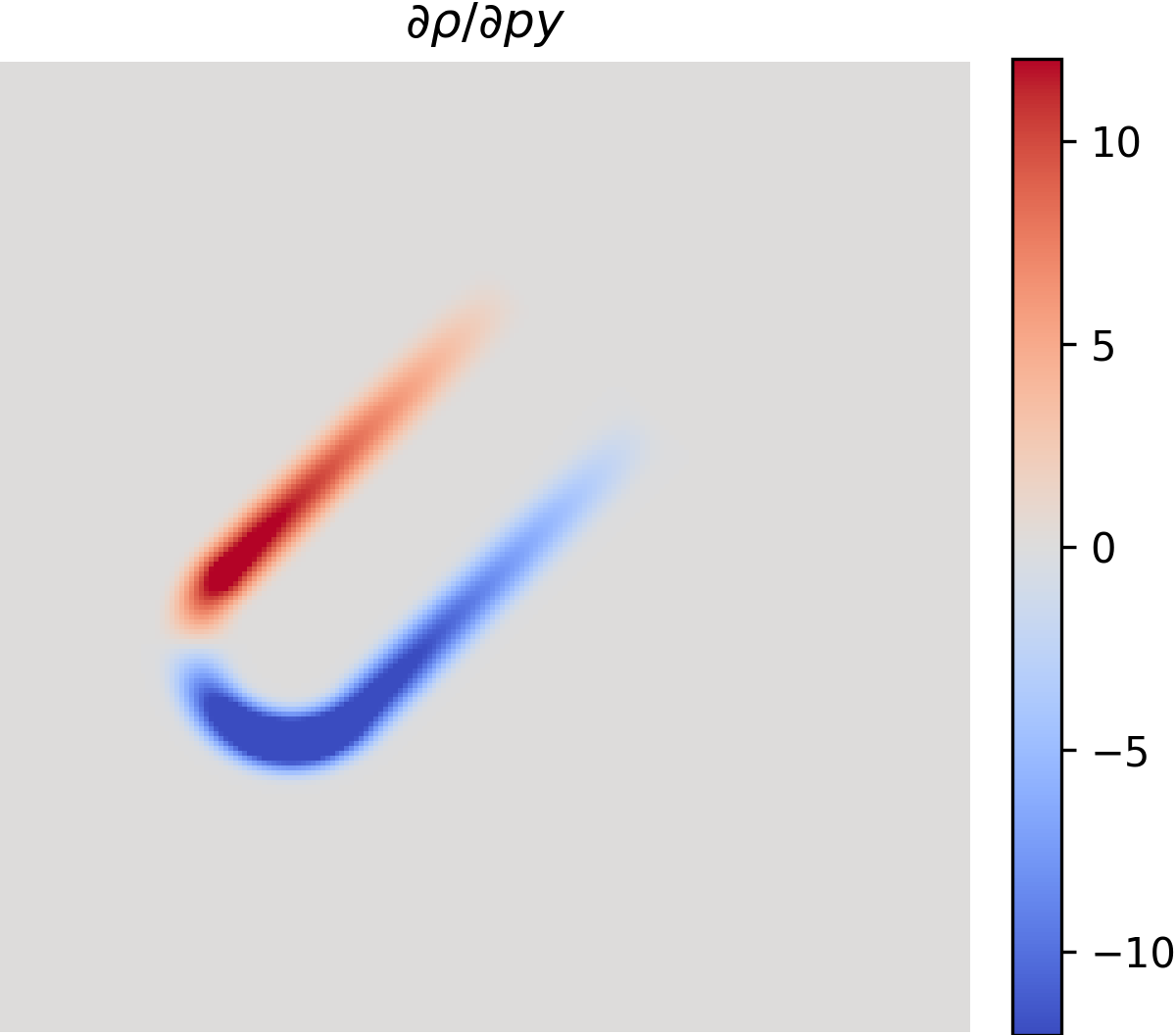}
  \caption*{$\partial\rho/\partial p_y$}
 \end{minipage}\hfill
 \begin{minipage}{.18\textwidth}\centering
  \includegraphics[width=\linewidth]{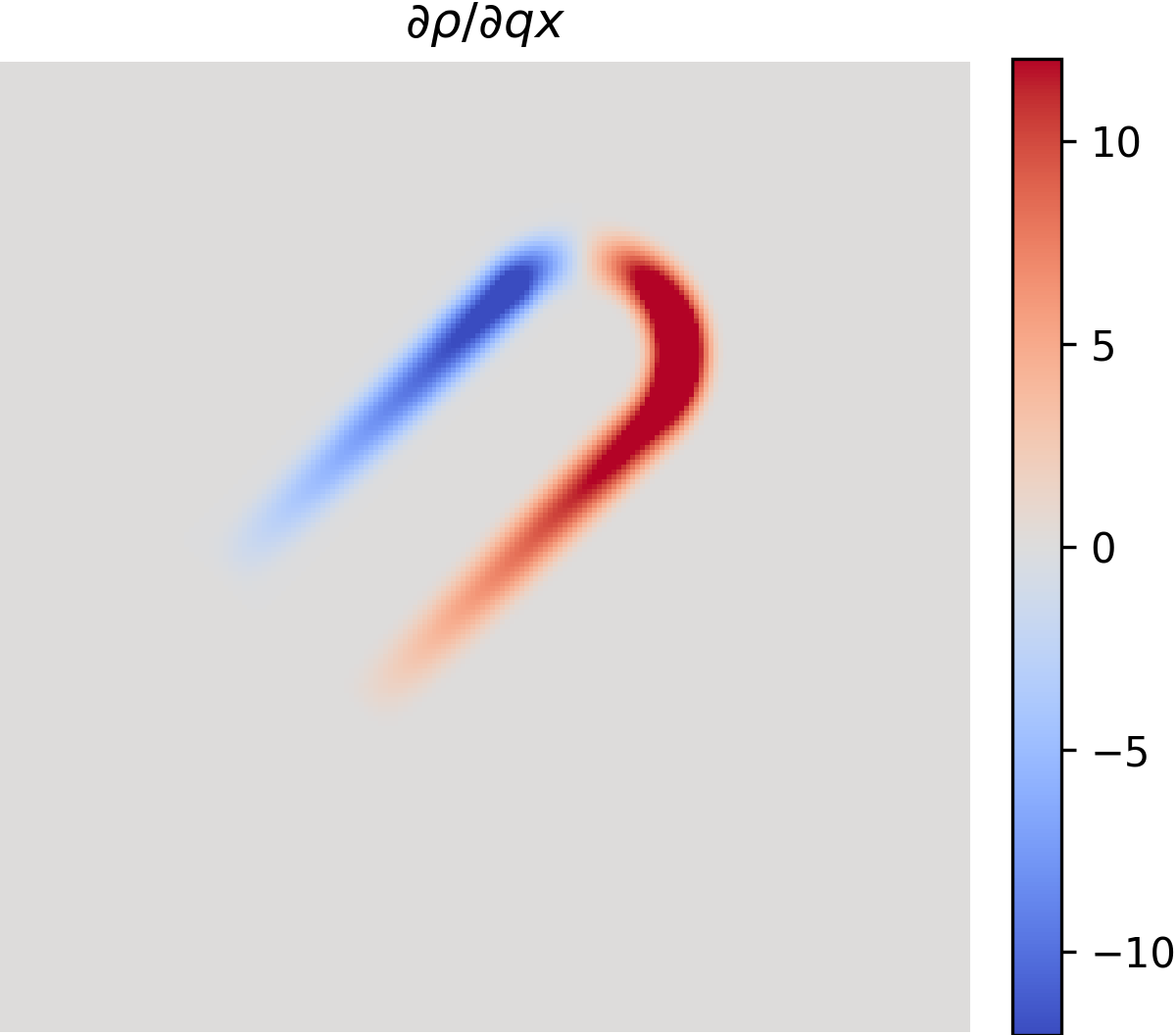}
  \caption*{$\partial\rho/\partial q_x$}
 \end{minipage}\hfill
 \begin{minipage}{.18\textwidth}\centering
  \includegraphics[width=\linewidth]{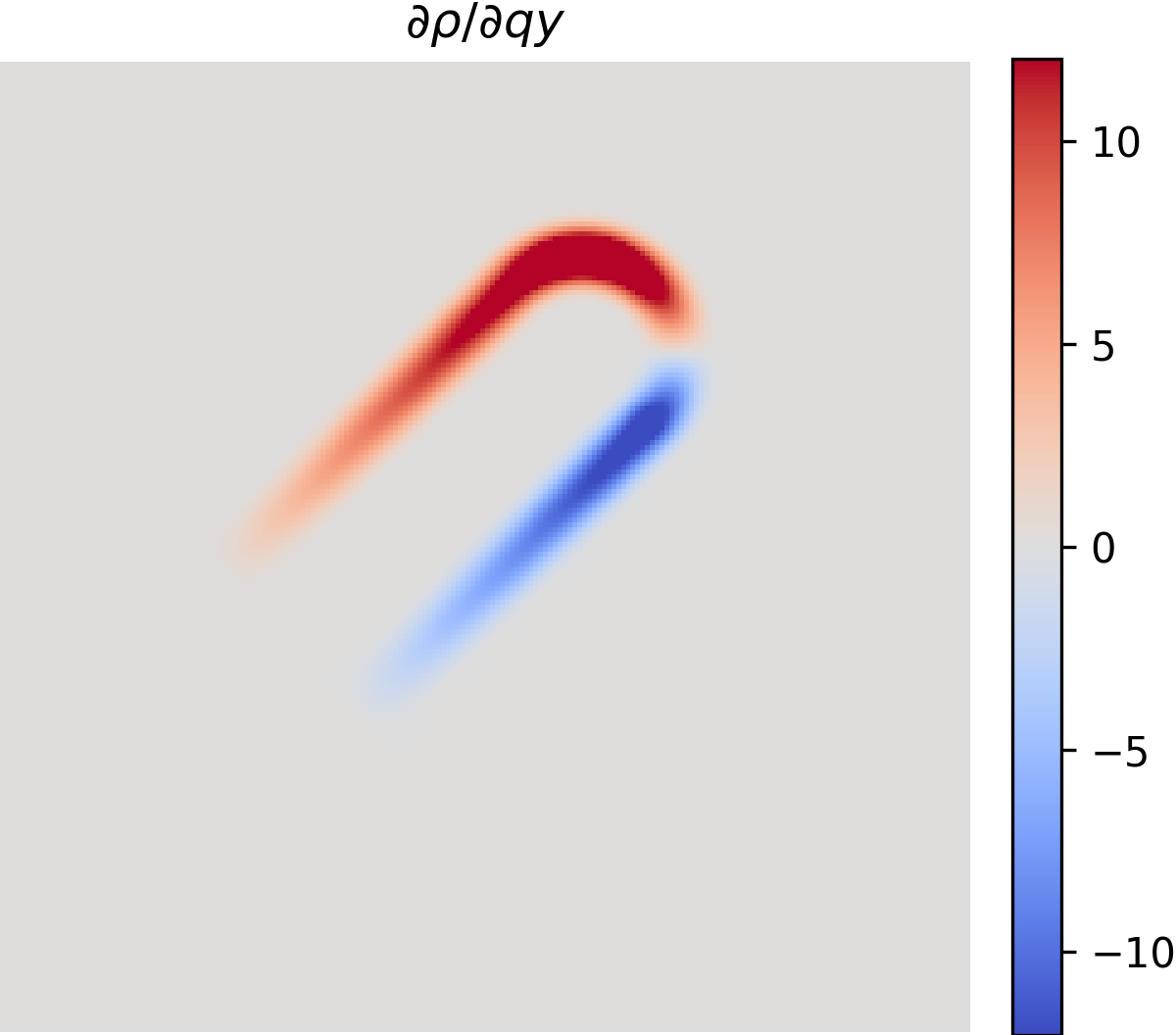}
  \caption*{$\partial\rho/\partial q_y$}
 \end{minipage}\hfill
 \begin{minipage}{.18\textwidth}\centering
  \includegraphics[width=\linewidth]{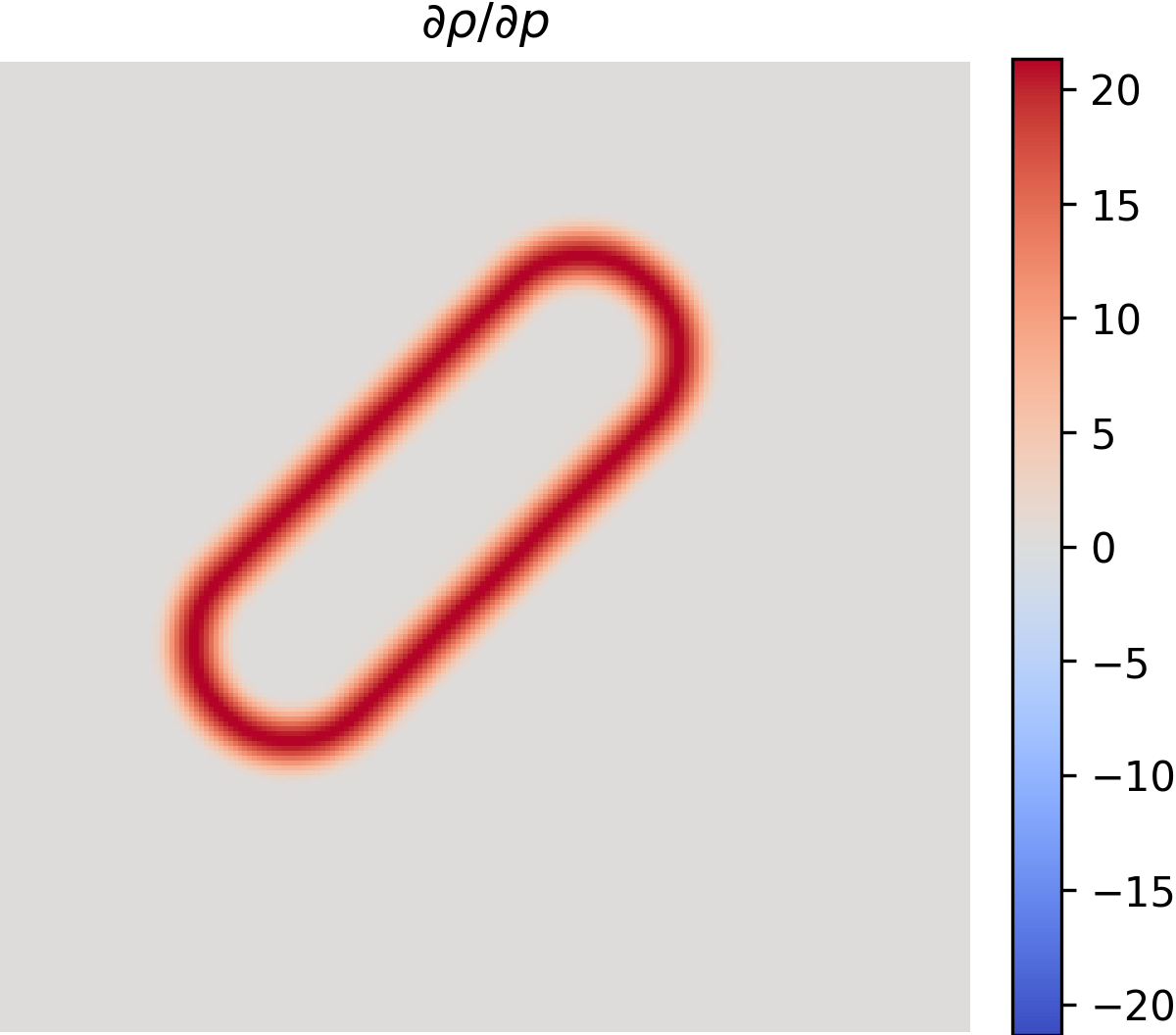}
  \caption*{$\partial\rho/\partial r$}
 \end{minipage}
 \caption{First-order sensitivities of $\rho$ for the smoothstep transition with $k=2$.}
 \label{fig:heatmaps_grad}
\end{figure}
\begin{figure}[H]
 \centering
 \begin{minipage}{.18\textwidth}\centering
  \includegraphics[width=\linewidth]{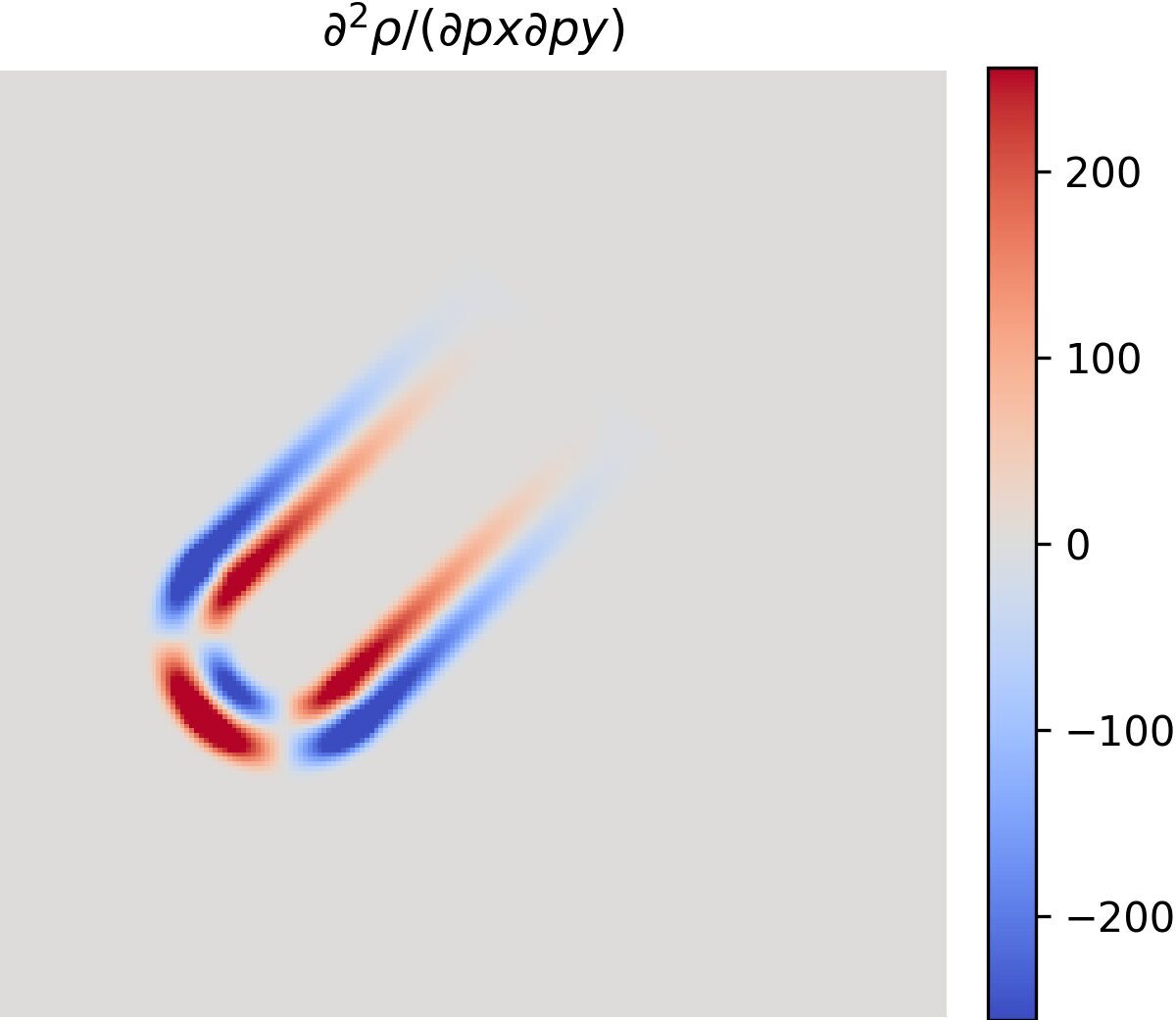}
  \caption*{$\partial^2\rho/(\partial p_x\partial p_y)$}
 \end{minipage}\hfill
 \begin{minipage}{.18\textwidth}\centering
  \includegraphics[width=\linewidth]{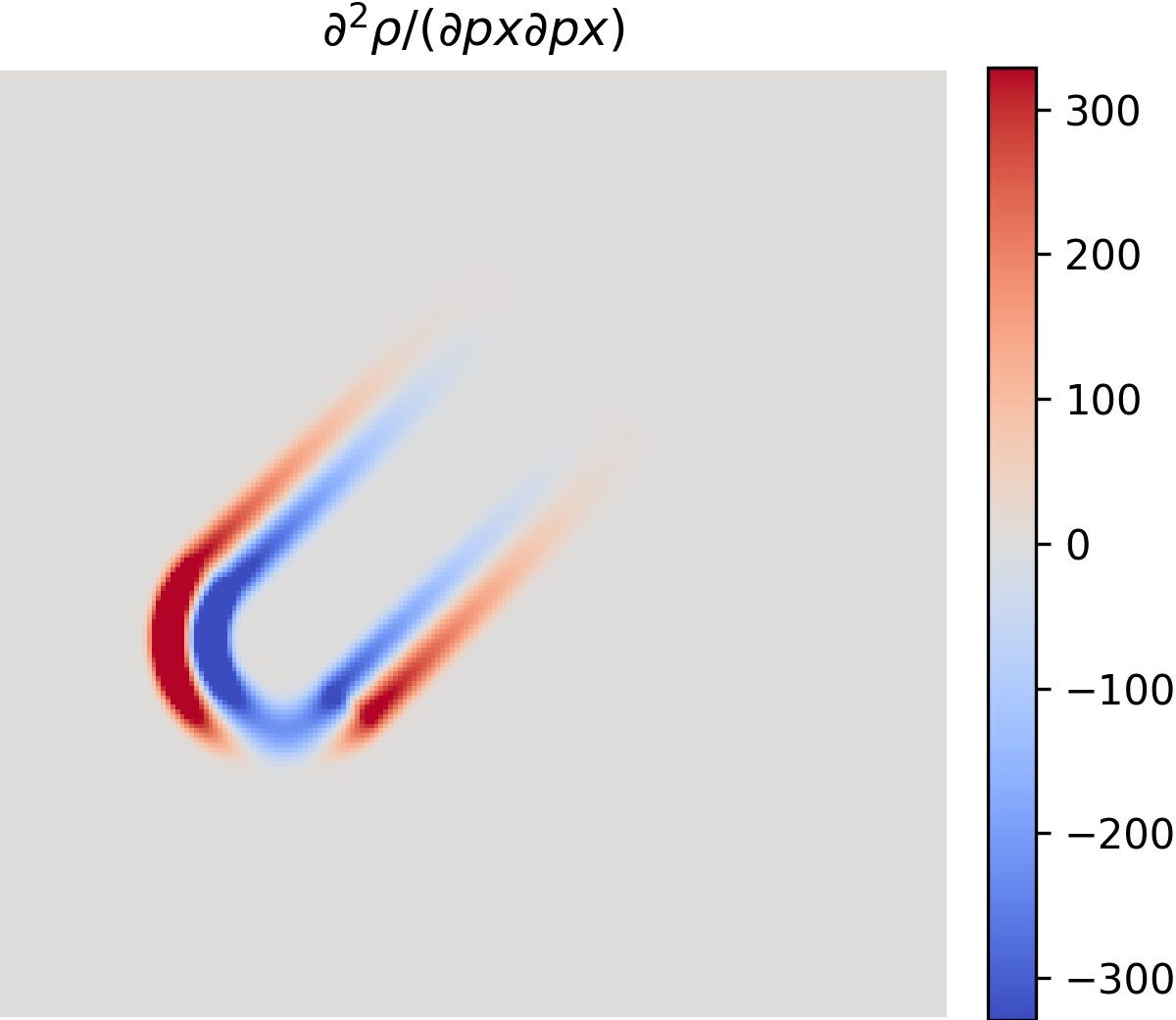}
  \caption*{$\partial^2\rho/(\partial p_x^2)$}
 \end{minipage}\hfill
 \begin{minipage}{.18\textwidth}\centering
  \includegraphics[width=\linewidth]{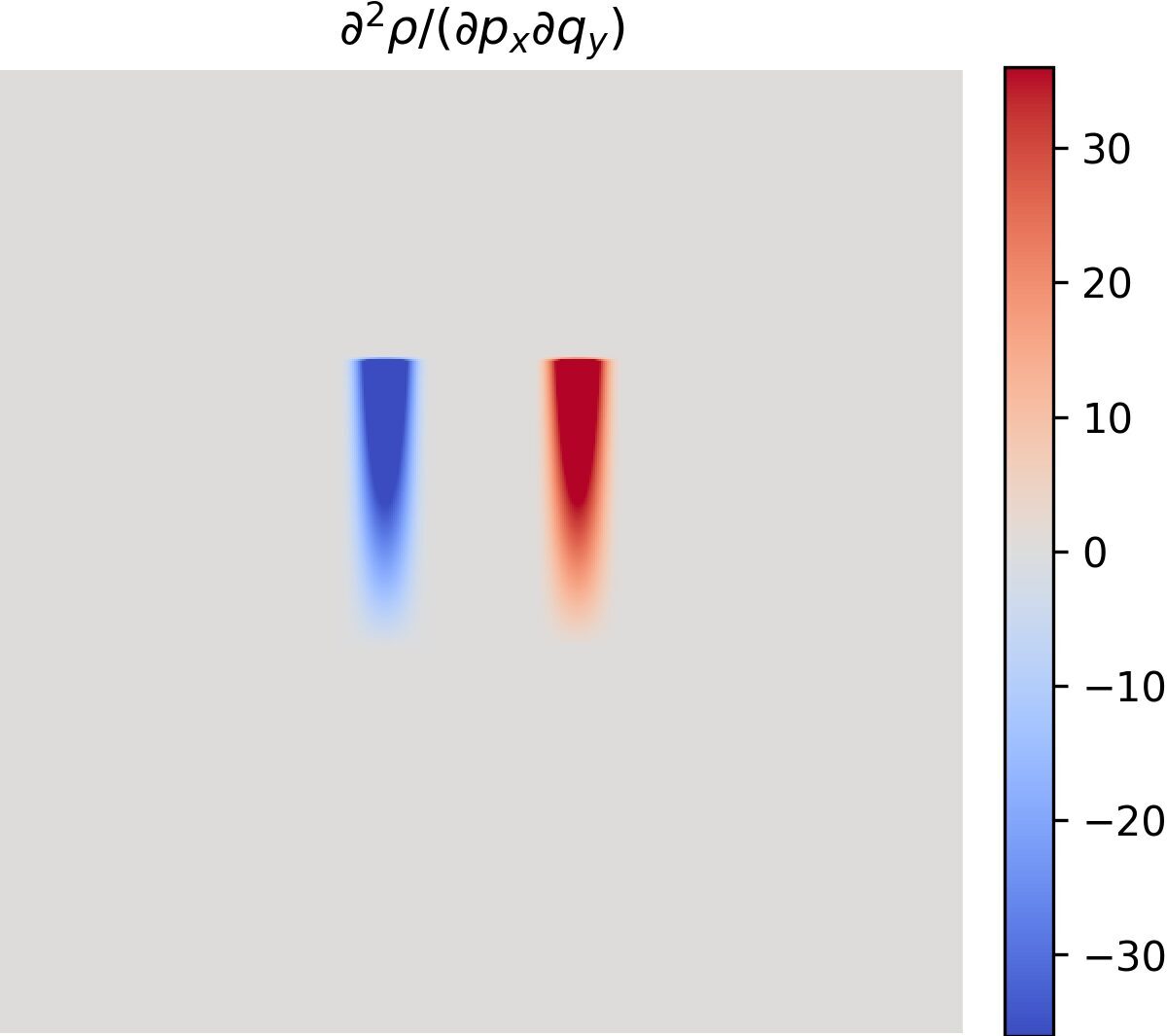}
  \caption*{$\partial^2\rho/(\partial p_x\partial q_y)$}
 \end{minipage}\hfill
 \begin{minipage}{.18\textwidth}\centering
  \includegraphics[width=\linewidth]{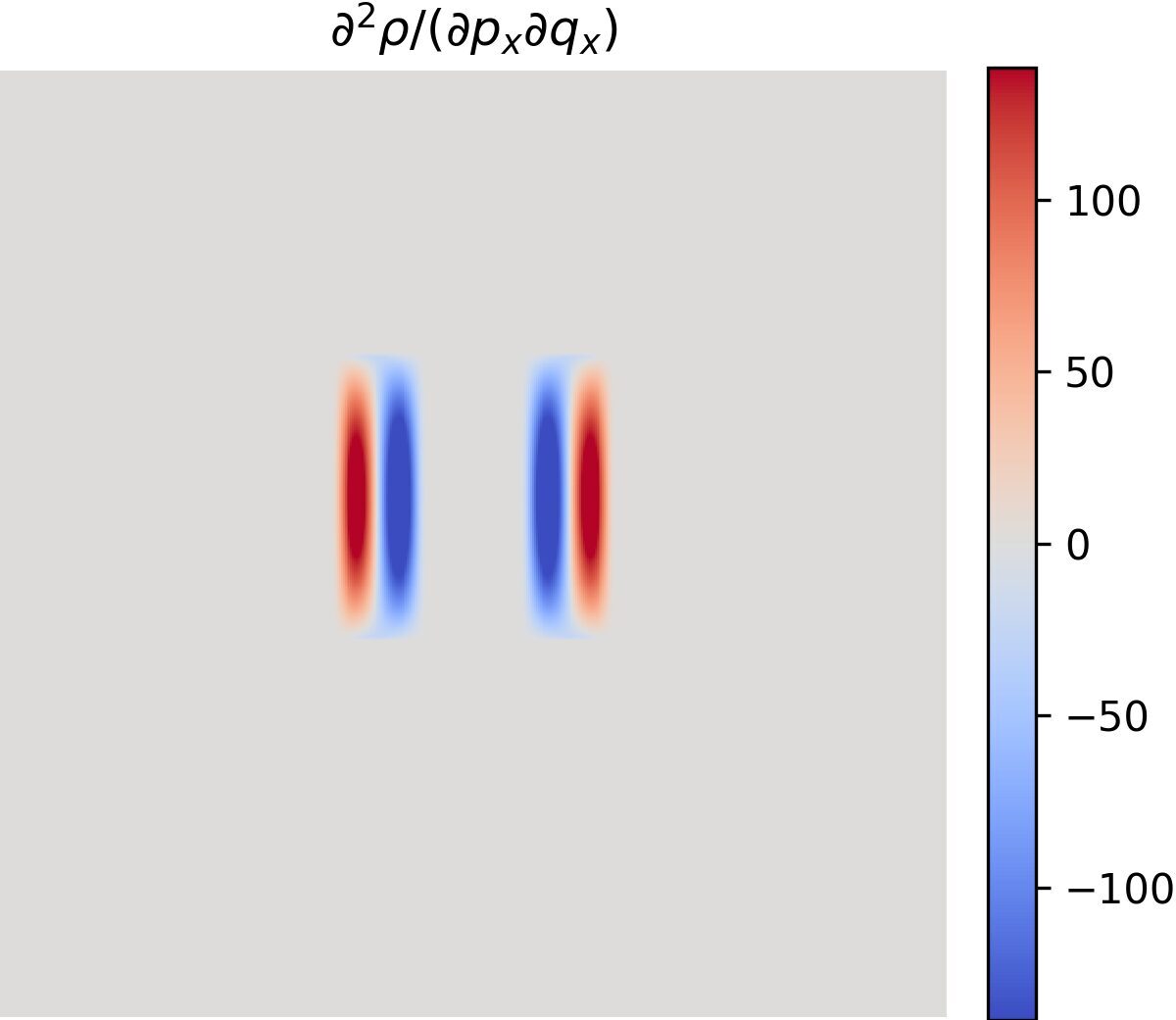}
  \caption*{$\partial^2\rho/(\partial p_x\partial q_x)$}
 \end{minipage}\hfill
 \begin{minipage}{.18\textwidth}\centering
  \includegraphics[width=\linewidth]{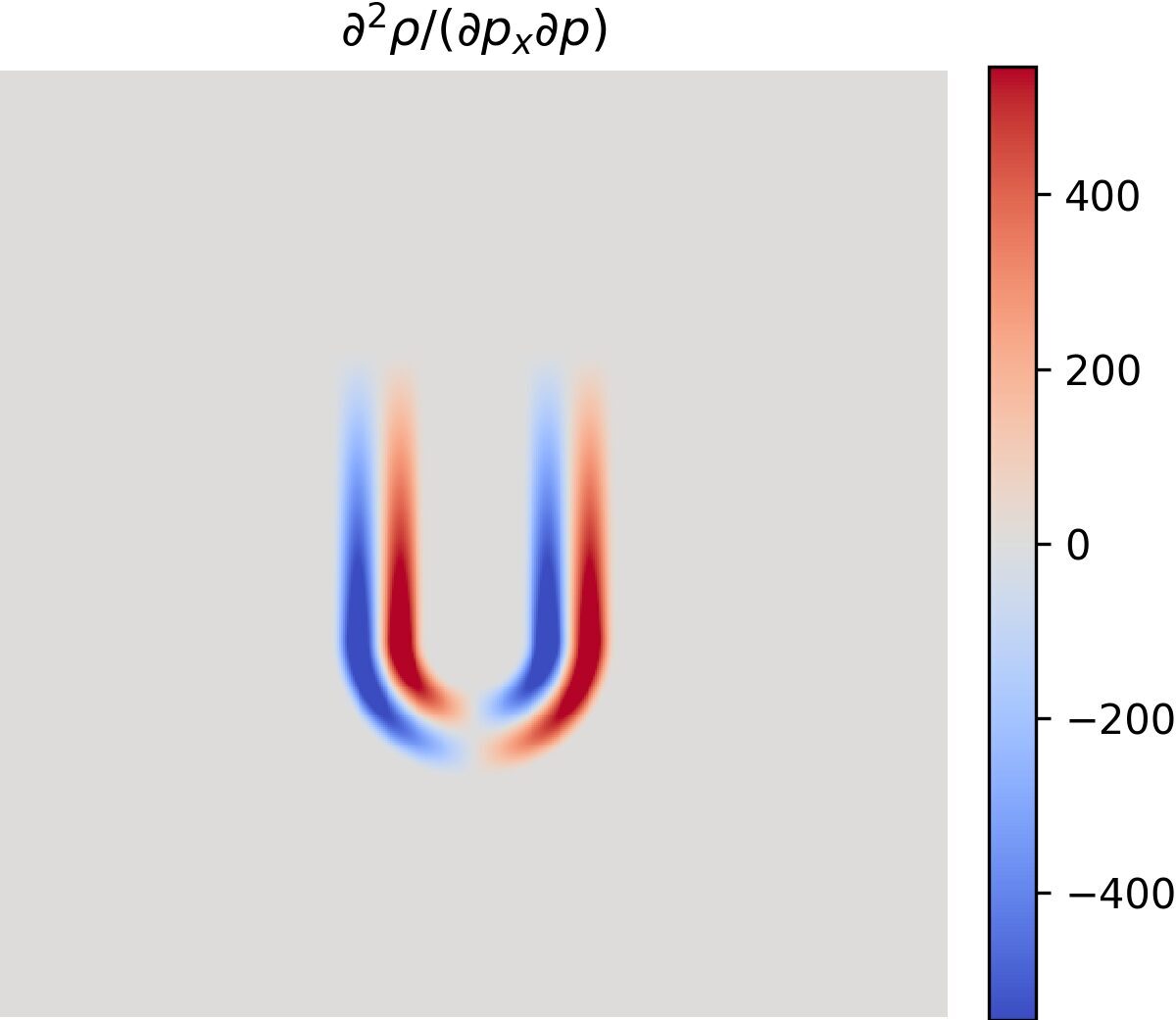}
  \caption*{$\partial^2\rho/(\partial p_x\partial r)$}
 \end{minipage}
 \caption{Representative Hessian entries of $\rho$ for the smoothstep transition with $k=2$.}
 \label{fig:heatmaps_hess2}
\end{figure}

The sensitivity patterns observed in the plots confirm the analytic structure of the derivatives. Nonzero contributions occur exclusively within the transition zone $\{\,\abs{d(\vecsym{x};\vecsym{z})}<\delta\,\}$, while both the interior plateau ($d\le -\delta$) and the exterior plateau ($d\ge \delta$) remain unaffected, since $\phi'(d)=0$ in these regions. Endpoint translations give rise to antisymmetric lobes across the local normal direction: a displacement in $+x$ or $+y$ increases $\rho$ on one side of the band and decreases it on the opposite side, with the strongest magnitudes occurring near the curved endcaps of the segment. Variations of the radius act uniformly and produce a strictly positive annulus aligned with the transition band, reflecting the isotropic expansion of the pill. The second-order derivatives inherit these structures in a more complex form: mixed entries alternate in sign and generate four-lobe patterns around the endcaps, while terms involving the radius combine the annular contribution with the antisymmetric structure of positional derivatives. In all cases, sensitivities outside the transition zone are negligible and sign symmetries cause substantial cancellation, ensuring localization and smoothness of the gradient and Hessian fields. This confirms that the transition functions preserve the $C^2$ trend established for the distance field and therefore provide a consistent basis for subsequent optimization.

A closer inspection of the second-order derivatives reveals how these structures arise from the analytic form. For the diagonal entries associated with endpoint translations, the dominant contribution is 
\[
\phi''(d)\,\Bigl(\tfrac{\partial d}{\partial p_i}\Bigr)^2,
\]
which changes sign between the inner and outer flank of the transition band. This explains the observed separation into positive and negative curvature zones across the band. Since the factor $(\partial d/\partial p_i)^2$ amplifies directions where the local normal has a large component in the corresponding coordinate, the maxima are aligned with the horizontal axis for $\partial^2/\partial p_x^2$ and with the vertical axis for $\partial^2/\partial p_y^2$. In contrast, along the straight segment where the normal is constant, the values remain symmetric on both sides. The secondary term \[
\phi'(d)\,\frac{\partial^2 d}{\partial p_i^2}
\]
becomes relevant near the endcaps, where the distance Hessian is nonzero and induces an inward shift of the negative contribution, consistent with the curvature of the circular arc.

The mixed second-order derivatives, e.g.\ $\partial^2/(\partial p_x\partial p_y)$, follow directly from the structure 
\[
\phi''(d)\,\frac{\partial d}{\partial p_x}\,\frac{\partial d}{\partial p_y}.
\]
Their sign depends on the product of the normal components $n_x n_y$ and therefore alternates between quadrants around the endcap, giving rise to the characteristic four-lobe pattern with alternating signs. The additional $\phi'(d)$-term again modulates these values close to the cap but does not alter the sign structure.

Cross-derivatives between the two endpoints reduce to the segment contribution. When the same coordinate appears at both ends (e.g.\ $\partial^2/(\partial p_x\partial q_x)$), the product $(\partial d/\partial p_x)(\partial d/\partial q_x)$ is positive on the outer flank and negative on the inner flank of the band. For mixed coordinates (e.g.\ $\partial^2/(\partial p_x\partial q_y)$) the sign is inverted, since $n_x$ and $n_y$ enter with opposite roles. These analytic relations explain the systematic alternation of sign in the Hessian blocks and confirm that all nontrivial contributions remain confined to the transition zone.

\subsection{Asymmetric transition}
\label{sec:asym_scaled}

\begin{figure}[H]
\centering
\includegraphics[width=0.8\textwidth]{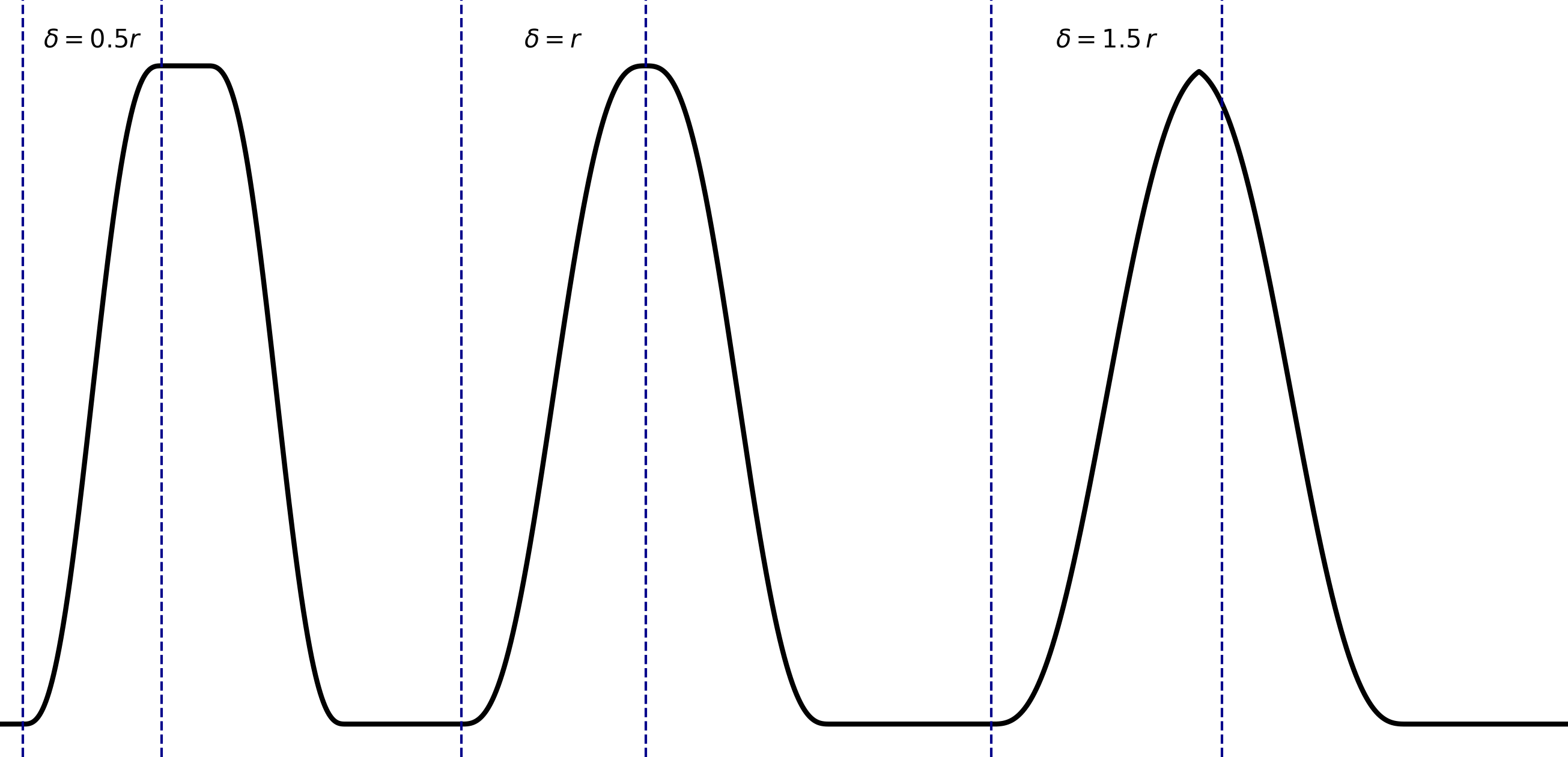}
\caption{Oversized symmetric transition with $\delta>r$. The interior plateau does not reach $\rho=1$, producing slope discontinuities at the segment.}
\label{fig:oversized_transition}
\end{figure}

In some cases, it is desirable to increase the exterior sensitivity while preserving a bounded pill radius. A naive enlargement of the symmetric transition zone beyond the pill radius $r$ leads to artifacts: the interior plateau fails to saturate at $\rho=1$ and mismatched one-sided derivatives violate $C^2$-continuity at the medial axis. This degeneracy is illustrated in \cref{fig:oversized_transition}, where oversized symmetric zones with $\delta>r$ produce visible kinks along the segment.

To avoid this, the interior half-width is restricted to $h=\delta/2$ while only the exterior side is extended. Denoting the additional reach of the right flank by $\mathrm{ext}\ge0$ and using the symmetric $C^2$ smoothstep basis $\phi$ on $[-\delta,\delta]$, the piecewise parametrization becomes
\[
t(d)=
\begin{cases}
\dfrac{-d+h}{2h}, & -h \le d \le 0,\\[0.4em]
\dfrac{-d+h+\mathrm{ext}}{2(h+\mathrm{ext})}, & 0<d<h+\mathrm{ext},
\end{cases}
\qquad
\rho(\vecsym{x})=\phi(t(d)),
\]
with clipping $\rho=1$ for $d\le-h$ and $\rho=0$ for $d\ge h+\mathrm{ext}$. Continuity holds at $d=0$ since $t(0^-)=t(0^+)=\tfrac12$, but the derivatives differ:
\[
\rho'(0^-)= -\tfrac{\phi'(1/2)}{2h},\qquad 
\rho'(0^+)= -\tfrac{\phi'(1/2)}{2(h+\mathrm{ext})}.
\]
For $\mathrm{ext}>0$ the slopes on both sides of the medial axis differ, so that $\rho'(d)$ exhibits a discontinuity at $d=0$. The profile is therefore continuous but no longer differentiable at the interface, i.e.\, only $C^0$. This provides the desired asymmetric reach while at the same time eliminating the plateau defect of oversized symmetric transitions, albeit at the expense of reduced smoothness.

\begin{figure}[H]
 \centering
 \includegraphics[width=.30\textwidth]{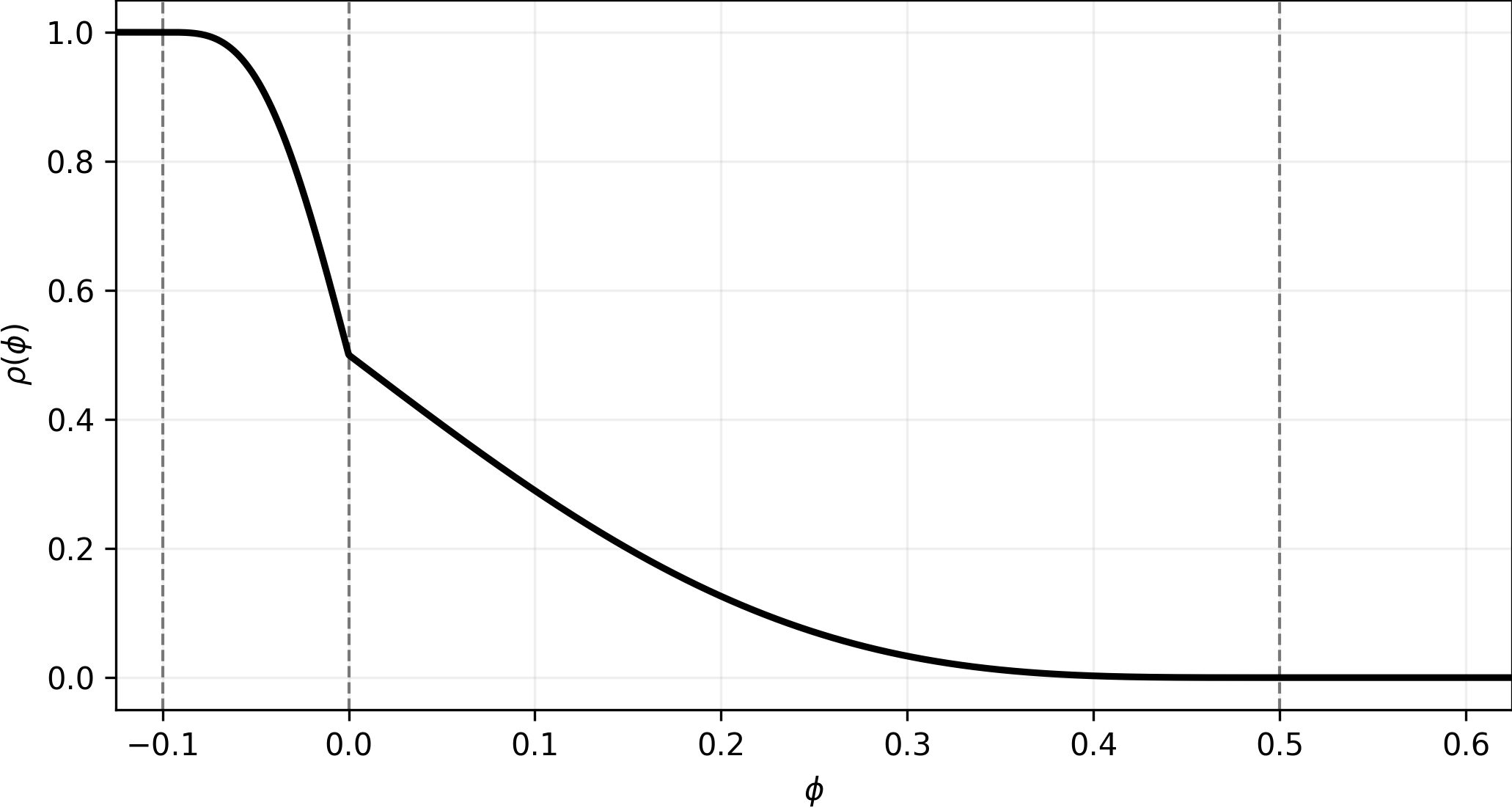}
 \includegraphics[width=.30\textwidth]{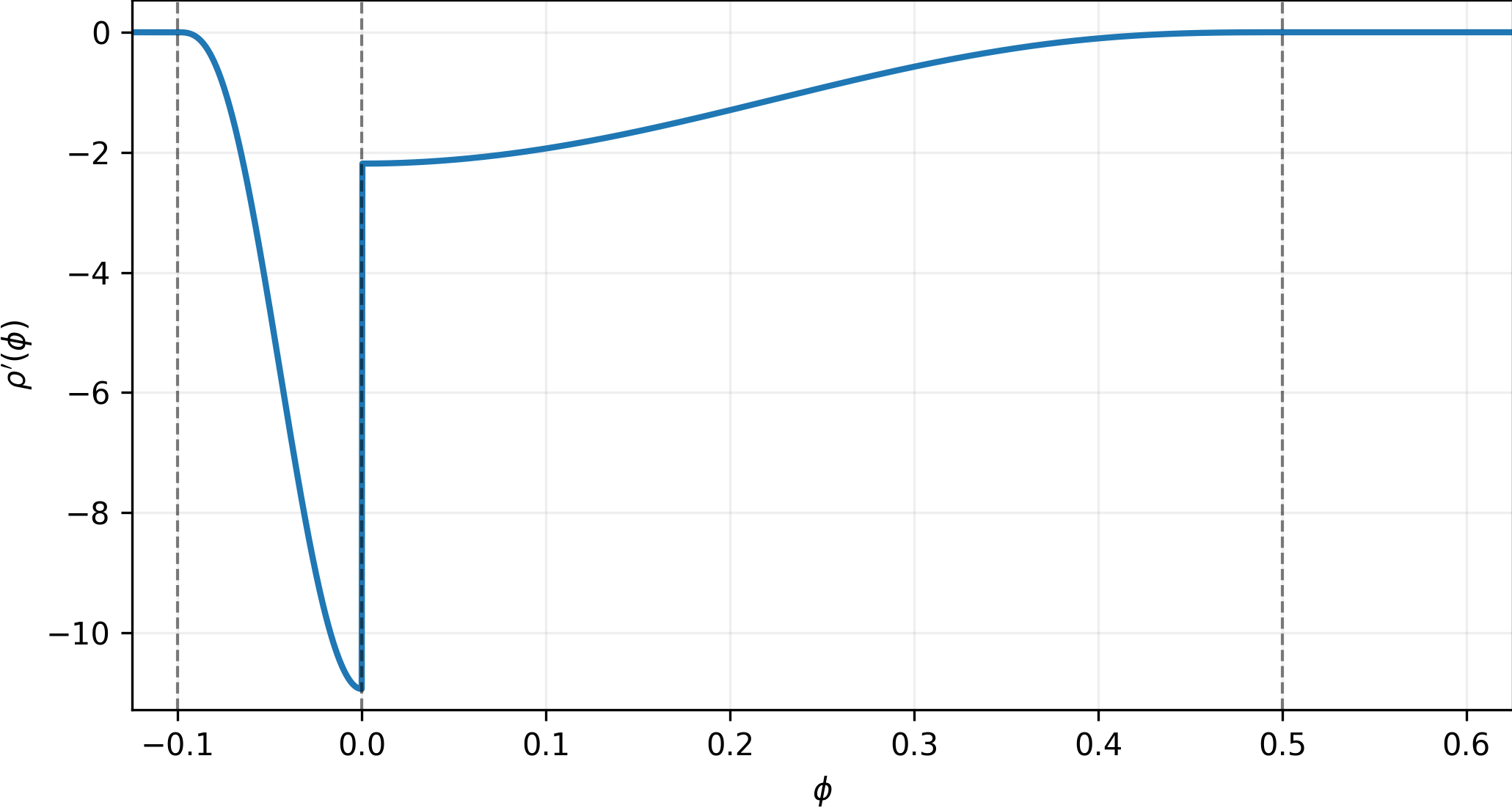}
 \includegraphics[width=.30\textwidth]{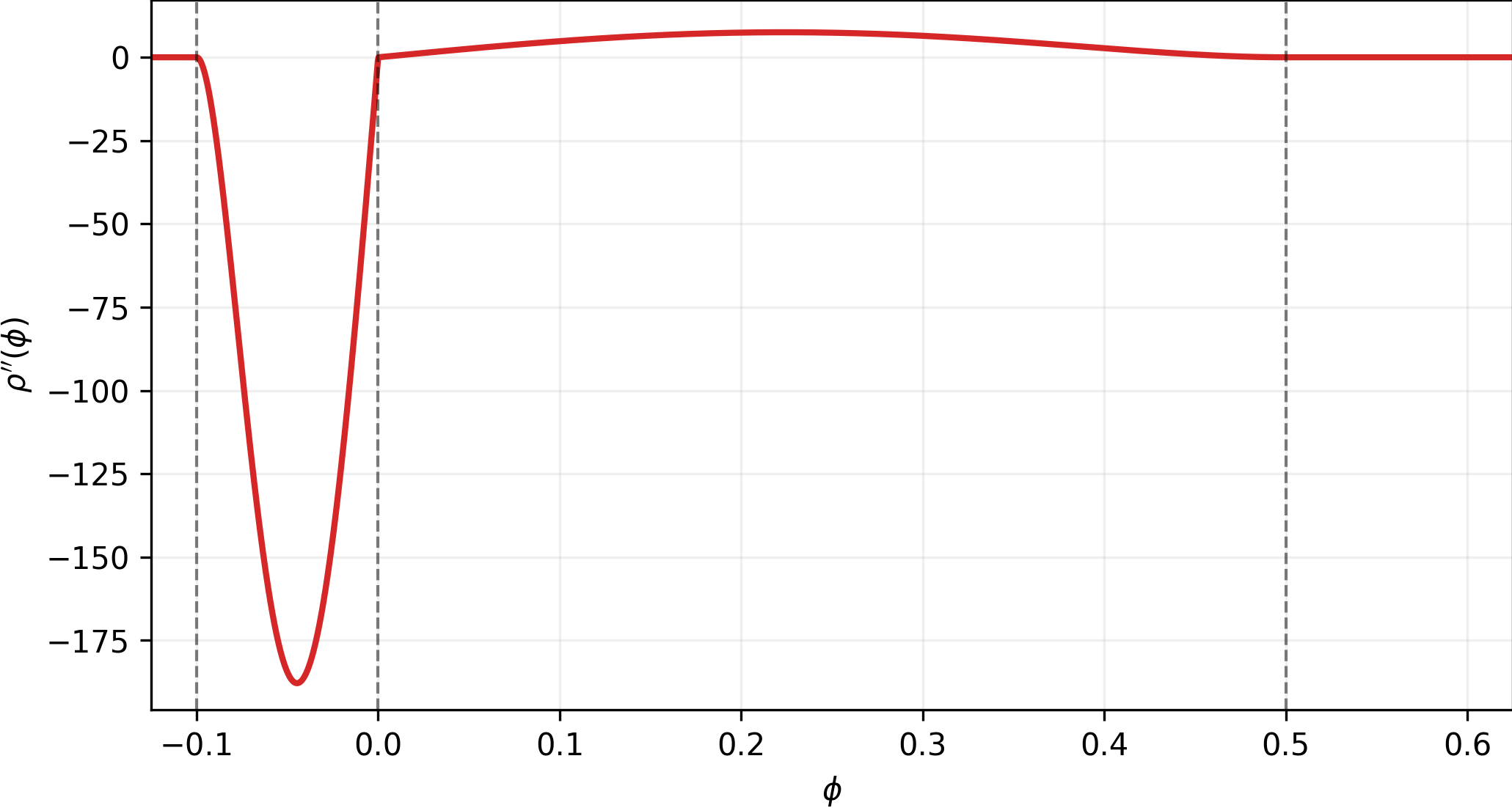}
 \caption{Asymmetric transition with scaled right flank: density $\rho$ (left), first derivative $\rho'$ (center) and second derivative $\rho''$ (right). The discontinuity of $\rho'$ and the spike in $\rho''$ at $d=0$ reflect the reduction to $C^0$ regularity.}
 \label{fig:asymm_scaled_derivs}
\end{figure}

Figure~\ref{fig:asymm_scaled_derivs} illustrates the resulting behavior of the function and its derivatives. The density $\rho$ remains monotone and extends farther into the exterior, while the first derivative shows a jump at $d=0$: the interior slope is steeper than the exterior slope because $h<h+\mathrm{ext}$. The second derivative does not exhibit a discontinuity but forms a cusp at $d=0$, indicating the loss of $C^1$ regularity. Away from the medial axis, both $\rho'$ and $\rho''$ decay smoothly and vanish at the clipping points $d=-h$ and $d=h+\mathrm{ext}$, so that compact support is retained.

\begin{figure}[H]
\centering
\includegraphics[width=0.6\textwidth]{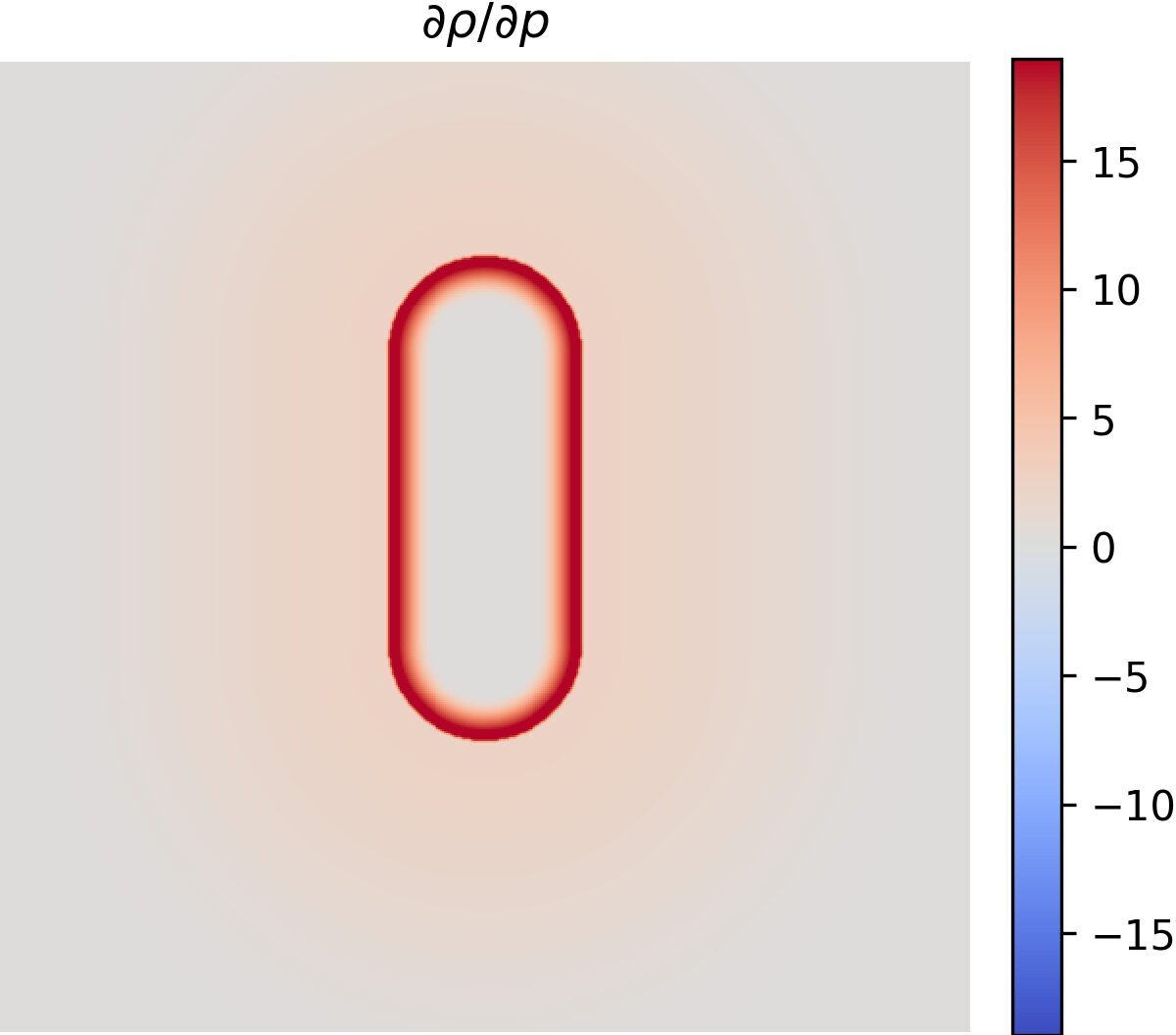}
\caption{Radius sensitivity $\partial \rho/\partial r$ for the scaled-flank transition. Nonzero values extend into a wider exterior band compared to the symmetric case.}
\label{fig:asymm_scaled_r}
\end{figure}

The impact of this construction on sensitivities is seen in the derivative with respect to the pill radius $r$. Compared to the symmetric case, $\partial\rho/\partial r$ now extends farther into the exterior region, confirming that the enlarged flank provides gradient information over a wider band. Figure~\ref{fig:asymm_scaled_r} shows this effect: the sensitivity forms an annulus that is strictly positive but reaches further outward than in the symmetric construction, thereby enhancing responsiveness of the geometry to radius updates.

Overall, the scaled-flank approach avoids the plateau defect of oversized symmetric zones but reduces smoothness at the medial axis. In optimization, the discontinuity of the derivative may lead to artificial gradient concentration along the pill axis, yet the extended support of $\partial\rho/\partial r$ improves the effectiveness of radius updates. The method thus represents a trade-off: it sacrifices $C^2$ regularity in favor of increased exterior sensitivity, making it suitable when strict locality inside the pill is less critical than extending influence to the surroundings.
The transition mappings developed in this chapter complete the projection of individual features into smooth pseudo-density fields. In practice, however, design domains contain many such features whose contributions must be combined. The next chapter therefore addresses the aggregation of multiple pseudo-densities into a global field while preserving locality, smoothness and analytic differentiability.

\chapter{Aggregation of Multiple Features\markboth{Aggregation of Multiple Features}{}}
\label{sec:aggregation}

In realistic designs, many geometric features coexist and must be represented jointly. While each pill can be projected into a smooth pseudo–density field, the resulting fields need to be combined into a single aggregation density that remains differentiable and suitable for gradient–based optimization. Formally, let $\Omega\subset\R^2$ be the design domain and let $\rho_m:\Omega\to[0,1]$, $m=1,\dots,n$, denote the pseudo–densities of the individual features. An \emph{aggregation operator}
\[
A:[0,1]^n \to \R_{\ge 0}, \qquad 
A(\vecsym{x}) = A\bigl(\rho_1(\vecsym{x}),\dots,\rho_n(\vecsym{x})\bigr), \;\; \vecsym{x}\in\Omega,
\]
defines how these contributions interact. The operator $A$ must ensure continuity and boundedness of the global field, interpolate smoothly in regions of partial overlap, preserve locality without introducing spurious long–range interactions and remain interpretable in terms of feature accumulation, competition or suppression. Moreover, differentiability with respect to the design parameters is essential in order to provide analytic gradients and Hessians for optimization.

\section{Aggregation methods}
\label{sec:aggregation_methods}

Various constructions of the aggregation operator $A$ are possible, each encoding different interaction patterns between overlapping features. The following sections present representative examples, together with their analytic derivatives.

\subsection*{Linear sum}

The most direct baseline is the linear superposition
\[
A(\vecsym{x}) := \sum_{m=1}^n \rho_m(\vecsym{x};\vecsym{z}^{(m)}),
\]
where each pseudo–density $\rho_m$ is associated with its own parameter block 
\[
\vecsym{z}^{(m)}=(p_x^{(m)},p_y^{(m)},q_x^{(m)},q_y^{(m)},r^{(m)}),
\]
and $z_i^{(m)}$ denotes the $i$–th component of this block. The contributions of the features are accumulated additively, resulting in values in $[0,n]$ without global rescaling. A division by $n$ would artificially suppress isolated features, while an adaptive rescaling by the number of local overlaps would introduce discontinuities. The linear sum is therefore left unscaled.

Differentiation with respect to the design parameters proceeds feature–wise. For any pill indices $a,b\in\{1,\dots,n\}$ and component indices $i,j\in\{1,\dots,5\}$, the sensitivities read
\[
\pp{A}{z_i^{(a)}} = \pp{\rho_a}{z_i^{(a)}}, \qquad
\ppp{A}{z_i^{(a)}}{z_j^{(b)}} = 
\begin{cases}
\ppp{\rho_a}{z_i^{(a)}}{z_j^{(a)}}, & a=b,\\
0, & a\neq b.
\end{cases}
\]
Hence the gradient of the aggregated field is simply the collection of the individual pill gradients and the Hessian exhibits a block–diagonal structure with one $5\times 5$ block per pill. Cross–feature couplings do not occur, so the operator is fully transparent and computationally inexpensive. The principal limitation is the absence of saturation: in regions of full overlap the density may exceed unity, whereby the natural binary interpretation $\rho\in[0,1]$ is no longer preserved.

\subsection*{$p$–norm aggregation}

A more selective alternative is provided by the $p$–norm union,
\[
A(\vecsym{x}) := \Biggl(\sum_{m=1}^n \rho_m(\vecsym{x};\vecsym{z}^{(m)})^p\Biggr)^{1/p}, 
\qquad p>1,
\]
which interpolates smoothly between linear accumulation and maximum selection. In the limit $p\to 1$ the result reduces to the linear sum, whereas $p\to\infty$ yields the pointwise maximum $\max_m \rho_m$. Selectivity therefore increases with $p$: dominant features prevail, while weaker overlaps are attenuated.

Differentiation with respect to the design parameters follows from the chain rule. Defining
\[
S(\vecsym{x})=\sum_{m=1}^n \rho_m(\vecsym{x};\vecsym{z}^{(m)})^p,
\]
one obtains for any component $z_i^{(a)}$ of the parameter block of pill $a$
\[
\pp{A}{z_i^{(a)}} 
= S^{\frac{1}{p}-1}\,\rho_a^{\,p-1}\,\pp{\rho_a}{z_i^{(a)}}.
\]
For second derivatives with respect to components $z_i^{(a)}$ and $z_j^{(b)}$ it follows that
\begin{align}
 \ppp{A}{z_i^{(a)}}{z_j^{(b)}} &=
 S^{\frac{1}{p}-1}\,
 \delta_{ab}\Bigl[
 \rho_a^{\,p-1}\ppp{\rho_a}{z_i^{(a)}}{z_j^{(a)}}
 +(p-1)\rho_a^{\,p-2}\pp{\rho_a}{z_i^{(a)}}\pp{\rho_a}{z_j^{(a)}}
 \Bigr] \notag \\[0.4em]
 &\quad
 + (1-p)\,S^{\tfrac{1}{p}-2}
 \Bigl(\rho_a^{\,p-1}\pp{\rho_a}{z_i^{(a)}}\Bigr)
 \Bigl(\rho_b^{\,p-1}\pp{\rho_b}{z_j^{(b)}}\Bigr), 
 \label{eq:pnorm_hessian}
\end{align}
 
where $\delta_{ab}$ is the Kronecker delta. Hence, if $a=b$, the Hessian collects both the intrinsic second derivatives of $\rho_a$ and the quadratic contribution from its gradient. If $a\neq b$, only the last term contributes, reflecting smooth cross–feature couplings, which vanish in the linear case $p=1$ but appear for $p>1$.
The $p$–norm can therefore be regarded as a differentiable, maximum-biased operator that generalizes the linear sum, preserves locality and introduces a tunable degree of saturation in regions of overlap.

\subsection*{Softmax aggregation}

Another smooth maximum-biased construction is given by the log–sum–exp or softmax, operator
\[
A(\vecsym{x}) := \tfrac{1}{\beta}\log\Biggl(\sum_{m=1}^n 
e^{\beta\,\rho_m(\vecsym{x};\vecsym{z}^{(m)})}\Biggr), 
\qquad \beta>0.
\]
As $\beta\to 0$, the operator approaches the arithmetic mean $\tfrac{1}{n}\sum_{m=1}^n \rho_m$, whereas in the limit $\beta\to\infty$ it converges to the maximum $\max_m \rho_m$. The parameter $\beta$ thus controls the degree of selectivity, interpolating smoothly between averaging and maximum selection.

For compact notation one sets
\[
S(\vecsym{x})=\sum_{m=1}^n e^{\beta \rho_m(\vecsym{x};\vecsym{z}^{(m)})}, 
\qquad 
w_m=\frac{e^{\beta \rho_m}}{S},
\]
so that the weights $w_m$ form a convex partition of unity. Differentiation with respect to any parameter $z_i^{(a)}$ yields
\[
\pp{A}{z_i^{(a)}} = \sum_{m=1}^n w_m \pp{\rho_m}{z_i^{(a)}}.
\]
Hence, the gradient of the aggregated field is a weighted combination of the individual pill gradients, with weights that concentrate on the dominant pill as $\beta$ increases.

For the second derivatives one obtains the following expression, in which the second term corresponds to the covariance of the gradients under the distribution $w_m$,
\begin{align}
\ppp{A}{z_i^{(a)}}{z_j^{(b)}} &=
\sum_{m=1}^n w_m \ppp{\rho_m}{z_i^{(a)}}{z_j^{(b)}} \notag \\[0.3em]
&\quad + \beta \sum_{m=1}^n w_m 
\Bigl(\pp{\rho_m}{z_i^{(a)}}-\sum_{r=1}^n w_r \pp{\rho_r}{z_i^{(a)}}\Bigr)
\Bigl(\pp{\rho_m}{z_j^{(b)}}-\sum_{s=1}^n w_s \pp{\rho_s}{z_j^{(b)}}\Bigr).
\label{eq:softmax_hessian}
\end{align}
The first term aggregates the individual pill Hessians with weights $w_m$, while the second term introduces smooth cross–feature couplings scaled by $\beta$. In the limit $\beta\to\infty$ the weights collapse onto the pill with largest $\rho_m$ and the operator becomes a differentiable approximation of the maximum.
\noindent
Since $\rho_m$ depends only on its own parameter block $\vecsym{z}^{(m)}$, $\partial\rho_m/\partial z_i^{(a)}=0$ for $m\neq a$. Accordingly, the sums in the gradient and in the first term of the Hessian collapse to the contributing features.

\subsection*{Sum–Softcap aggregation}

A further variant aims at controlling accumulation by introducing smooth saturation above a prescribed threshold $\tau>0$. Each pseudo–density $\rho_m(\vecsym{x};\vecsym{z}^{(m)})$ is first mapped by a monotone function $\varphi_\beta:[0,1]\to\R_{\ge0}$, which may be chosen as the identity or as a weakly saturating transformation. The transformed values are then summed and passed through a concave cap function,
\[
A(\vecsym{x}) := \operatorname{cap}_\tau\!\Biggl(\sum_{m=1}^n 
\varphi_\beta\bigl(\rho_m(\vecsym{x};\vecsym{z}^{(m)})\bigr)\Biggr).
\]

A convenient choice of cap function is
\[
\operatorname{cap}_\tau(s)=
\tau-\tfrac{1}{\beta_c}\log\bigl(1+e^{\beta_c(\tau-s)}\bigr),
\qquad \beta_c>0,
\]
which is smooth, strictly increasing, concave and asymptotically saturates at $\tau$. Its derivatives are
\[
\operatorname{cap}_\tau'(s)=\frac{1}{1+e^{\beta_c(\tau-s)}}\in(0,1),
\qquad
\operatorname{cap}_\tau''(s)=-\beta_c\,\operatorname{cap}_\tau'(s)
\bigl(1-\operatorname{cap}_\tau'(s)\bigr)\le 0.
\]

Denoting by $s(\vecsym{x})=\sum_{m=1}^n \varphi_\beta(\rho_m(\vecsym{x};\vecsym{z}^{(m)}))$, the sensitivities with respect to any parameters $z_i^{(a)}$ and $z_j^{(b)}$ are
\[
\pp{A}{z_i^{(a)}}=
\operatorname{cap}_\tau'(s)\,
\sum_{m=1}^n \varphi_\beta'(\rho_m)\,\pp{\rho_m}{z_i^{(a)}},
\]
\begin{align}
\ppp{A}{z_i^{(a)}}{z_j^{(b)}} &=
\operatorname{cap}_\tau'(s)\,
\sum_{m=1}^n\Bigl[
\varphi_\beta'(\rho_m)\,\ppp{\rho_m}{z_i^{(a)}}{z_j^{(b)}}
+\varphi_\beta''(\rho_m)\,\pp{\rho_m}{z_i^{(a)}}\pp{\rho_m}{z_j^{(b)}}
\Bigr] \notag \\[0.4em]
&\quad +\operatorname{cap}_\tau''(s)\,
\Biggl(\sum_{m=1}^n \varphi_\beta'(\rho_m)\,\pp{\rho_m}{z_i^{(a)}}\Biggr)
\Biggl(\sum_{r=1}^n \varphi_\beta'(\rho_r)\,\pp{\rho_r}{z_j^{(b)}}\Biggr).
\label{eq:sumsoftcap_hessian}
\end{align}

The first term scales the individual pill Hessians by $\operatorname{cap}_\tau'(s)$ and incorporates additional curvature through $\varphi_\beta''$. The second term introduces global coupling between features, weighted by the concavity of $\operatorname{cap}_\tau$. Below the threshold ($s\lesssim\tau$) the cap function acts nearly linearly and the second term is negligible, so accumulation is essentially additive. Above the threshold, the negative curvature of $\operatorname{cap}_\tau$ suppresses marginal gains and enforces smooth saturation.
Note: with the softplus–based cap, $A$ can be slightly negative for very small $s$ (by at most $\log 2/\beta_c$).

\paragraph{On clipping and normalization.}
A natural postprocessing step is the pointwise clipping $A\mapsto \min\{A,1\}$. The mapping $s\mapsto \min\{s,1\}$ is, however, only $C^0$ with a kink at $s=1$. As a consequence, the composite field $x\mapsto \min\{A(x),1\}$ is merely continuous and typically exhibits gradient discontinuities along the moving level set $\{x: A(x)=1\}$. In this way, second-order smoothness is lost, which may impair gradient-based optimization.

Likewise, global normalizations such as $A \mapsto n^{-1/p} A$ for the $p$–norm reduce peak values below one. Since a single pill is expected to attain a density of one by itself, such a scaling forces multiple features to overlap in order to achieve the same effect, thereby introducing a bias in the design. Normalizations based on a data-dependent ``active count'' (e.g., the number of features above a threshold) reintroduce discontinuities whenever features enter or leave the active set.

For the same reason, pointwise extrema introduce non–smoothness: the hard maximum $A(x)=\max_m \rho_m(x)$ (limit $p\to\infty$ or $\beta\to\infty$) and the hard minimum (e.g., for intersections) are only $C^0$ in general, with kinks where the active argmax/argmin switches.
In contrast, softly saturating constructions such as the Sum–Softcap (Sec.~\ref{sec:aggregation_methods}) achieve boundedness through smooth, concave capping and preserve differentiability.

\begin{figure}[ht]
\centering
\includegraphics[width=0.75\textwidth]{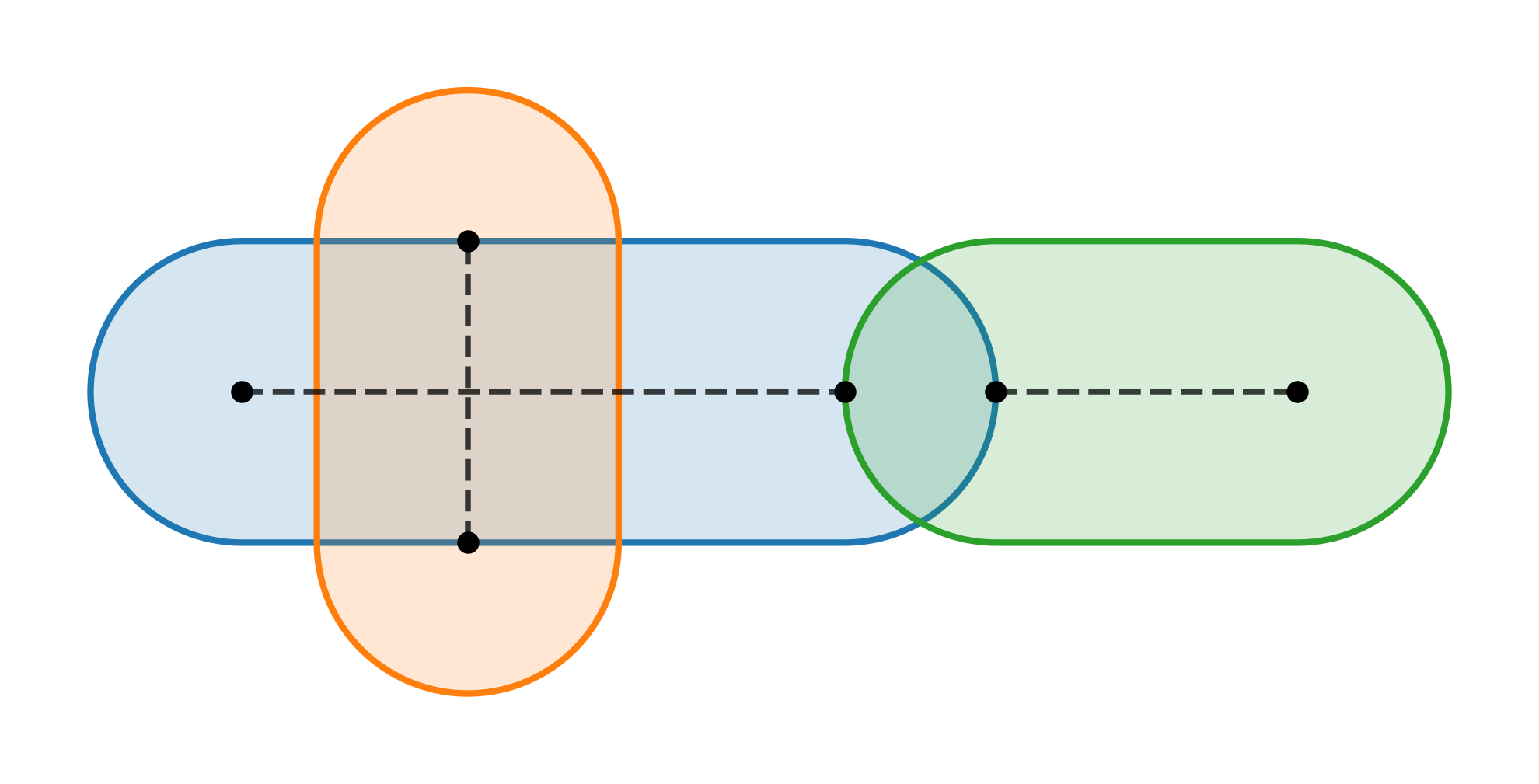}
\caption{Arrangement of three features and scanline $y=\tfrac12$.}
\label{fig:features_aggregation}
\end{figure}

\begin{figure}[ht]
\centering
\includegraphics[width=0.9\textwidth]{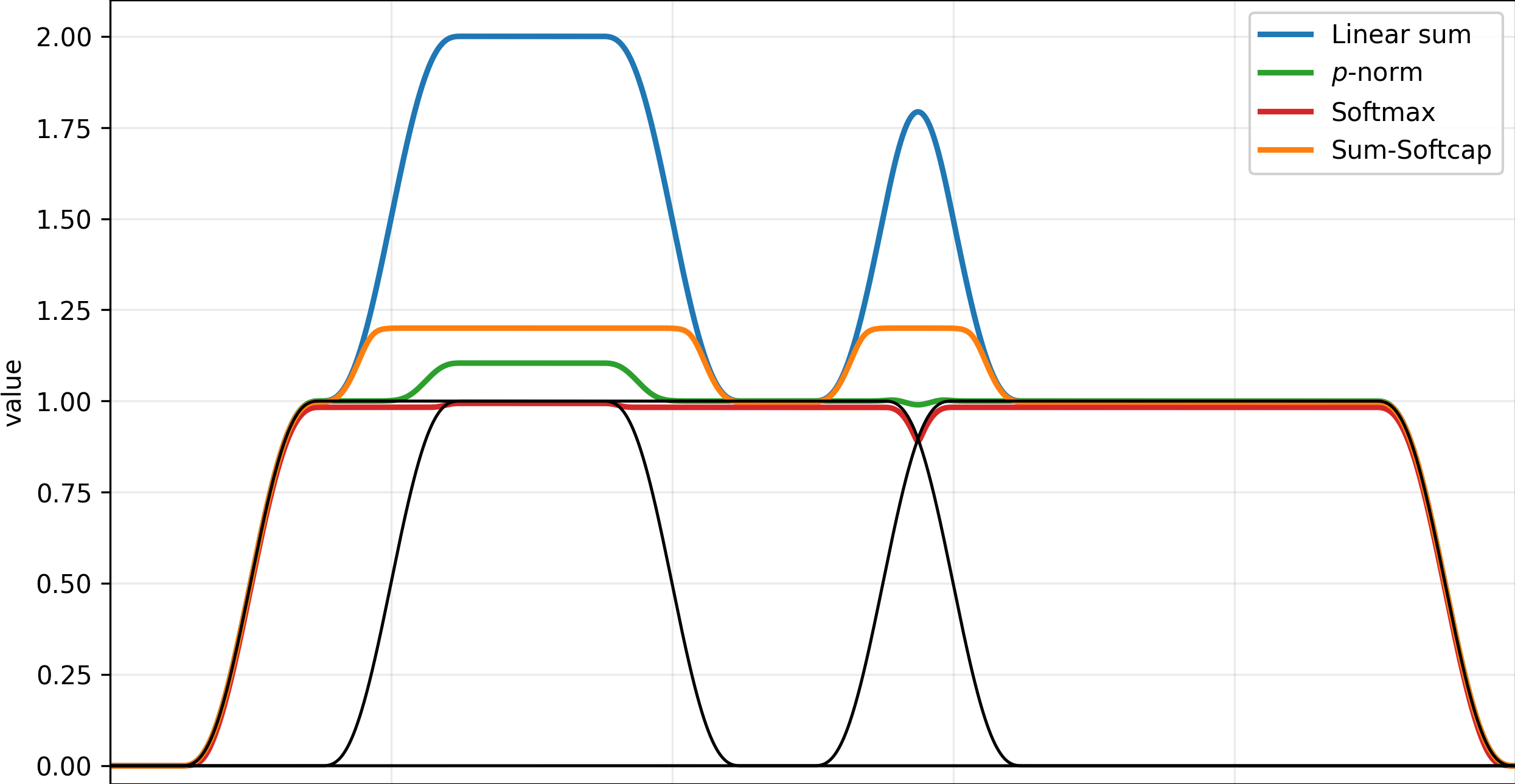}
\caption{Aggregated densities along $y=\tfrac12$ for different operators.}
\label{fig:aggregation_comparison}
\end{figure}

To illustrate the behavior of the aggregation operators, three pill features are arranged as shown in \cref{fig:features_aggregation}. Along the scanline $y=\tfrac12$, the profile alternates between regions with no pill, single coverage, pairwise overlap and full separation. The aggregated densities in \cref{fig:aggregation_comparison} reflect these situations consistently. The linear sum grows additively and can exceed unity under full overlap. The $p$–norm with $p=7$ and the softmax with $\beta=64$ both interpolate smoothly toward the maximum, suppressing secondary contributions and emphasizing dominant features. The Sum–Softcap with threshold $\tau=1.2$ behaves additively in low-density regions but saturates once the cap is exceeded, thereby avoiding excessive accumulation. All operators remain differentiable and preserve locality, but they differ in their balance between additivity, selectivity and saturation.

With these constructions, the analytic representation of pill fields is complete: signed distances are converted to differentiable pseudo–densities and consistently aggregated into a global aggregation field. The next step is the projection of this continuous field onto a finite element mesh, which enables numerical analysis and optimization.
\chapter{Mapping to a Fixed Grid\markboth{Mapping to a Fixed Grid}{}}
\label{sec:grid_mapping}

For numerical analysis and optimization, the continuous pseudo–density field \(\rho:\Omega\to[0,1]\) is transferred to a finite element discretization. The design domain is represented by a fixed Cartesian mesh and each element is characterized by a constant pseudo–density together with its first and second derivatives with respect to the design parameters. This mapping is performed through element–wise quadrature, providing a consistent interface between analytic shape functions and discrete finite element representations.

\section{Domain discretization}
\label{sec:grid_discretization}

Let \(\Omega\subset\R^2\) be partitioned into a uniform Cartesian mesh of rectangular finite elements \(\{\Omega_e\}_{e\in\mathcal E}\) with side lengths \(h_x,h_y>0\). For an element with lower-left corner \((x_e,y_e)\) one has
\[
\Omega_e=[x_e,x_e+h_x]\times[y_e,y_e+h_y].
\]

\begin{figure}[H]
 \centering
 \includegraphics[width=.78\textwidth]{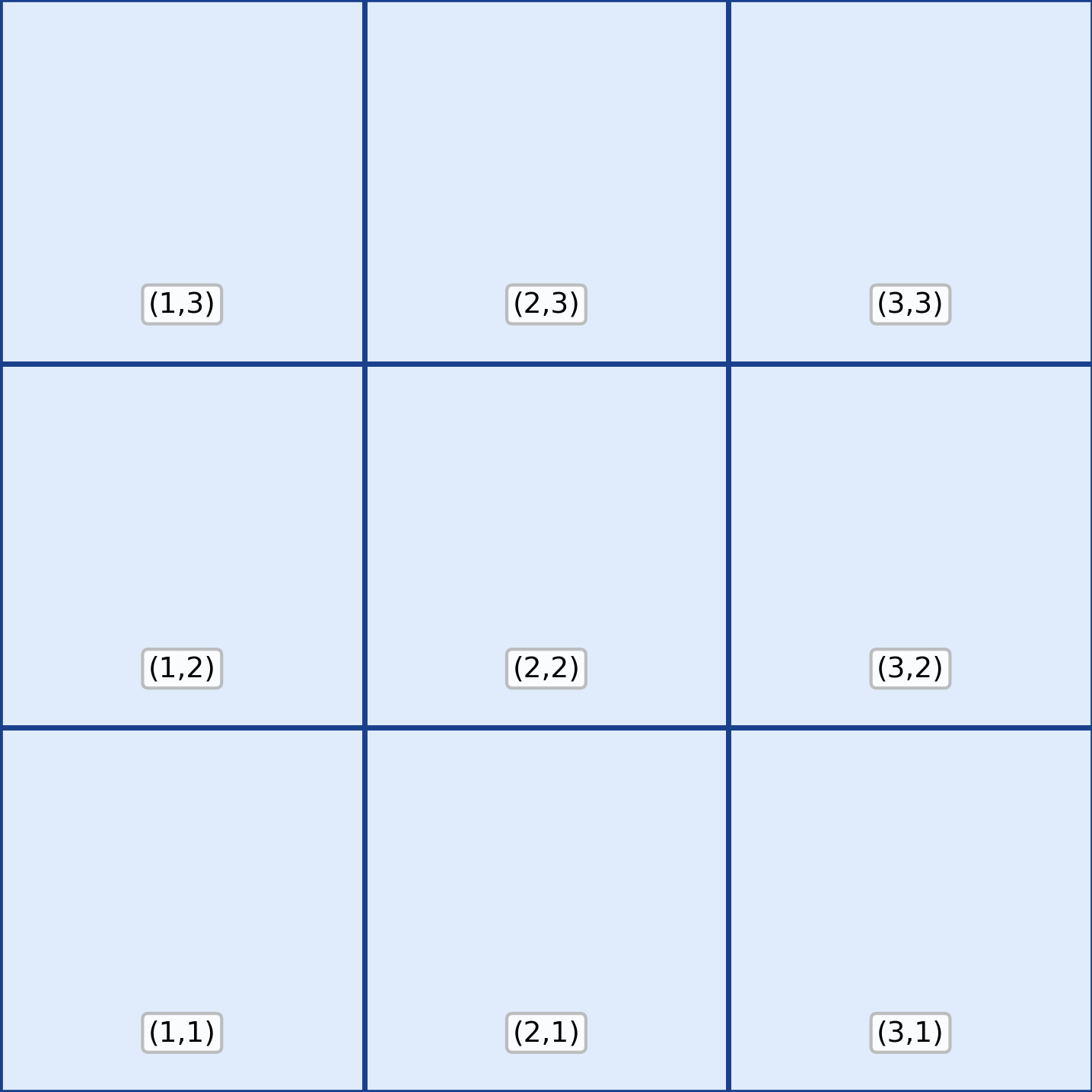}
 \caption{Cartesian mesh subdivided into \(3 \times 3\) elements.}
 \label{fig:integration_domain}
\end{figure}

In this representation, each element \(\Omega_e\) is regarded as an independent integration domain. The element pseudo–density is obtained by evaluating the contributions of all shapes \(\rho_m\) at a set of integration points \(x_e^i \in \Omega_e\). At each point, the shape values are combined by an aggregation function \(A(\cdot)\) and the element density results from the arithmetic mean of these aggregated values. Depending on the chosen quadrature order, the integration points may coincide with the element center, the vertices or a denser tensor grid, which controls the accuracy with which sharp transitions in \(\rho\) are captured.

\section{Element quadrature rules}
\label{sec:quadrature_rules}

The element density is obtained as the arithmetic mean of the aggregated pill contributions evaluated at a set of quadrature points. For a rectangular element, the points are placed by a uniform tensor–product subdivision into a $q \times q$ grid, resulting in $N_{\mathrm{ip}} = q^2$ sampling nodes. The element average is then given by
\[
\bar{\rho}^e \;=\; \frac{1}{N_{\mathrm{ip}}} \sum_{i=1}^{N_{\mathrm{ip}}} 
A\!\Bigl(\{\rho_m(x_e^i)\}_{m}\Bigr),
\]
where $\rho_m(x)$ denote the contributions of the individual shapes and $A(\cdot)$ the chosen aggregation function. Each integration point contributes with equal weight, so that accuracy is solely determined by the subdivision parameter $q$. For $q=1$ the construction reduces to the midpoint rule, which is second–order accurate. With $q=2$ the trapezoidal rule is recovered by averaging the four corner values. For $q=3$ a $3 \times 3$ tensor grid is obtained, which incorporates corners, edge midpoints and the element center. Increasing $q$ provides higher resolution and permits sharp gradients in $\rho$ to be represented with controllable precision.

\noindent
As the projection is given by an arithmetic mean, differentiation commutes with the summation. For the gradient with respect to a design parameter \(z_i\) one obtains
\[
\frac{\partial \bar{\rho}^e}{\partial z_i} 
= \frac{1}{N_{\mathrm{ip}}} \sum_{k=1}^{N_{\mathrm{ip}}} 
\sum_{a} \frac{\partial A}{\partial \rho_a}\!\Bigl(\{\rho_s(x_e^k)\}_{s}\Bigr)\;
\frac{\partial \rho_a(x_e^k)}{\partial z_i},
\]
where the outer sum runs over integration points \(k\) and the inner sum over shapes \(a\).

Similarly, for the Hessian with respect to parameters \(z_i,z_j\) one obtains
\begin{align}
\frac{\partial^2 \bar{\rho}^e}{\partial z_i \,\partial z_j} 
&= \frac{1}{N_{\mathrm{ip}}} \sum_{k=1}^{N_{\mathrm{ip}}} 
\Biggl[
\sum_{a,b} 
\frac{\partial^2 A}{\partial \rho_a \,\partial \rho_b}
  \!\Bigl(\{\rho_s(x_e^k)\}_{s}\Bigr)\,
\frac{\partial \rho_a(x_e^k)}{\partial z_i}\,
\frac{\partial \rho_b(x_e^k)}{\partial z_j}
\\[0.5em]
&\qquad\qquad
+ \sum_{a} 
\frac{\partial A}{\partial \rho_a}
  \!\Bigl(\{\rho_s(x_e^k)\}_{s}\Bigr)\,
\frac{\partial^2 \rho_a(x_e^k)}{\partial z_i \,\partial z_j}
\Biggr].
\end{align}
 
Thus, element averages, gradients and Hessians are obtained by direct averaging of the analytic pointwise quantities in full consistency with the aggregation rules.

In practice, the quadrature order is chosen relative to the ratio between element size and the transition half–width of the shape functions. Low orders already yield accurate results in regions where the aggregated field varies smoothly, while higher orders are required when narrow transition zones intersect element boundaries. Since the projection is expressed as an arithmetic mean of aggregated shape contributions, the gradients and Hessians follow from the same averaging process and thus inherit the smoothness of the underlying analytic expressions without additional effort. This completes the mapping step: element densities and their sensitivities are consistently transferred to the analysis mesh, thereby forming the basis for the optimization formulations developed in the subsequent chapter.

\begin{figure}[H]
 \centering
 \includegraphics[width=0.32\textwidth]{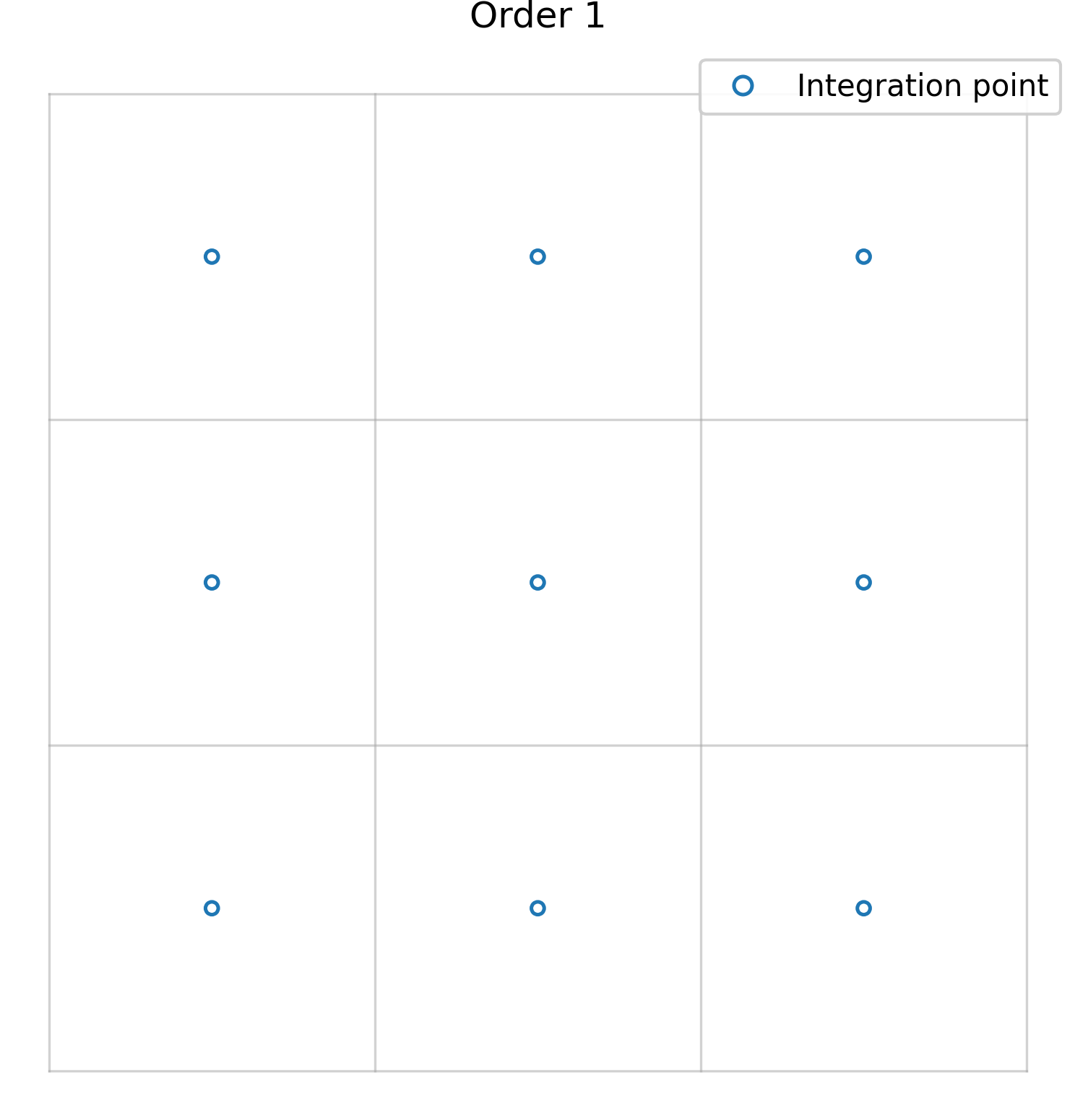}
 \includegraphics[width=0.32\textwidth]{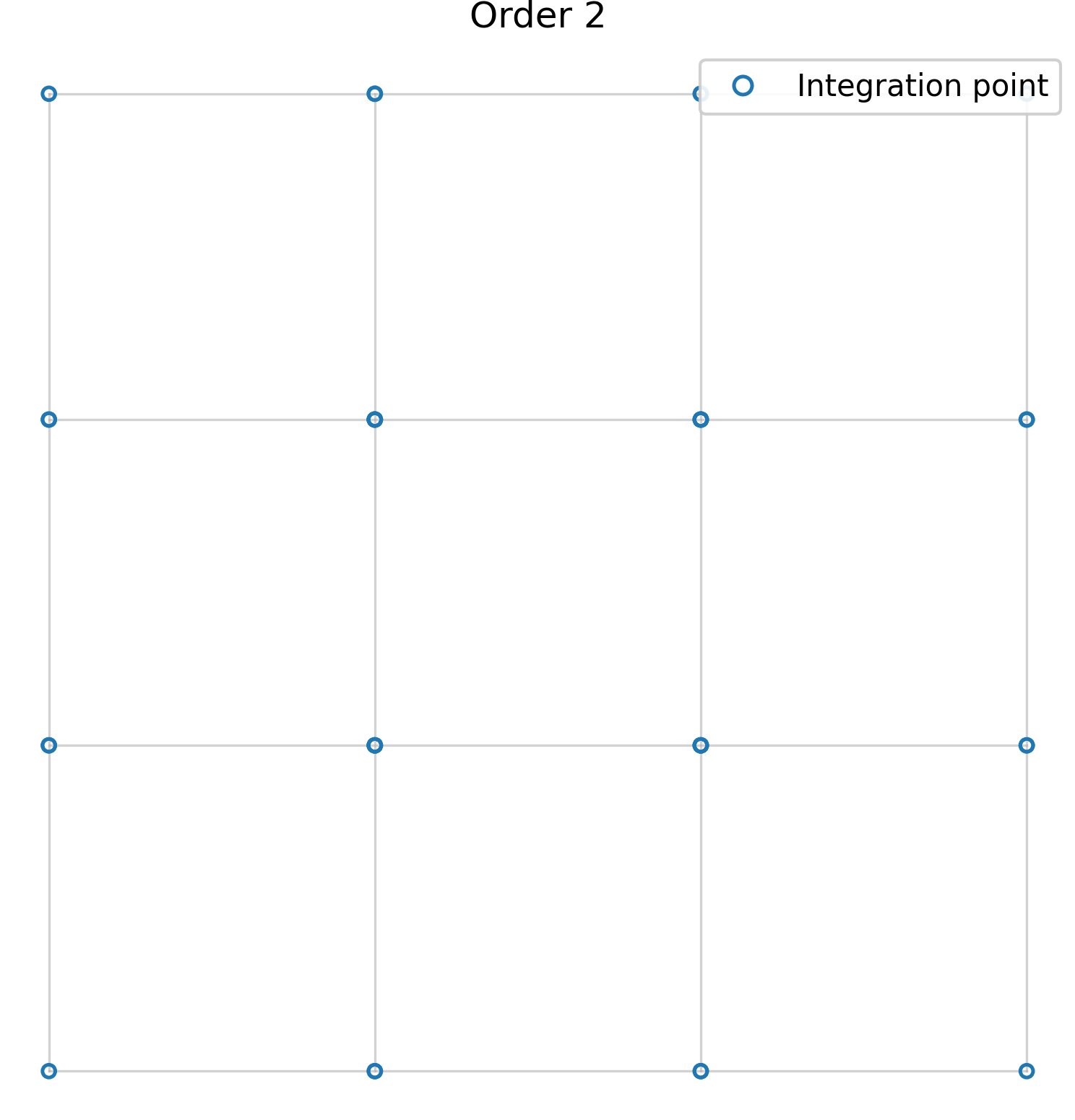}
 \includegraphics[width=0.32\textwidth]{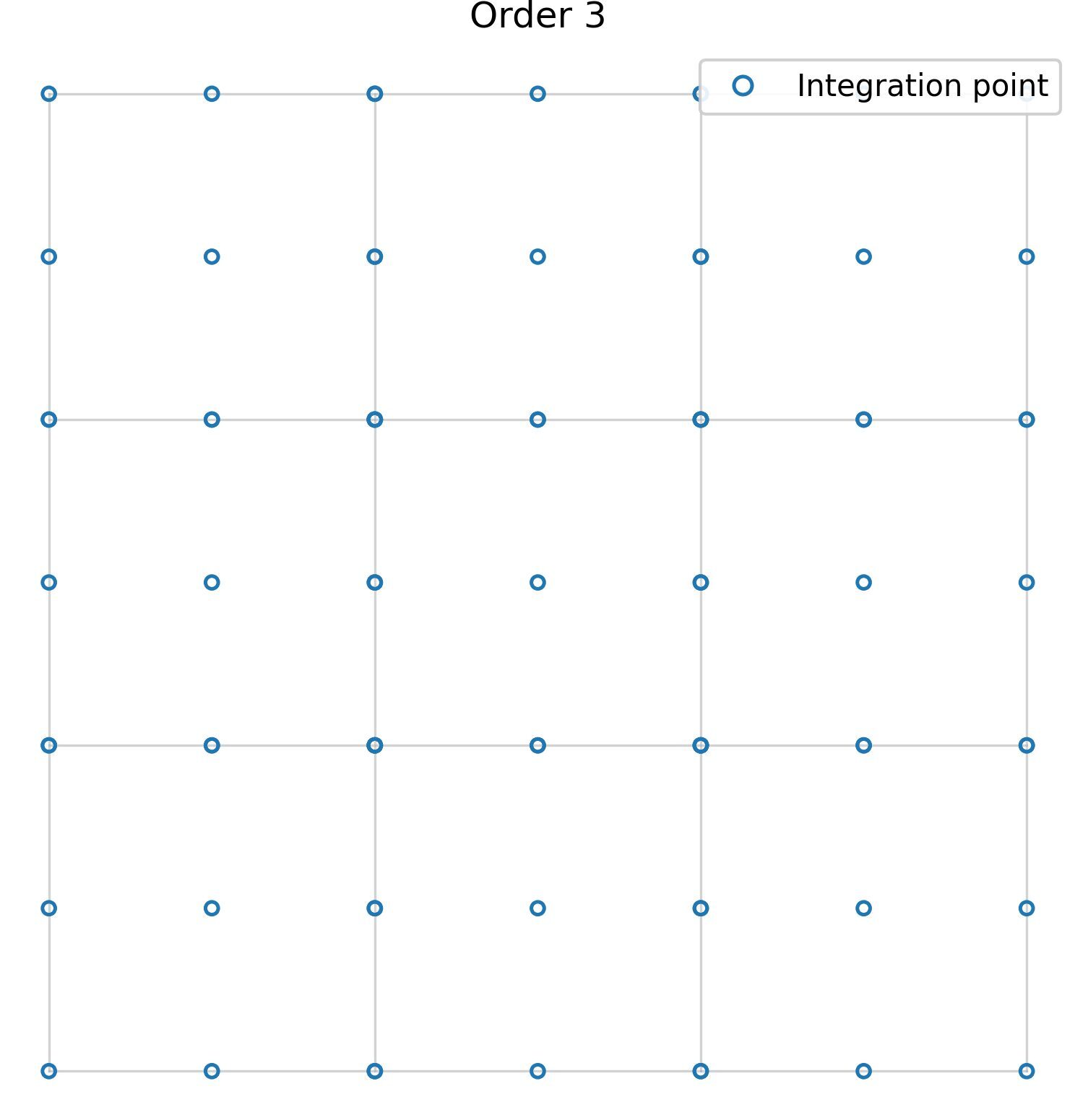}
 \caption{Tensor–product quadrature rules of order $q=1,2,3$. 
 For $q=1$ the midpoint rule is obtained with one integration node per element. 
 For $q=2$ the trapezoidal rule is recovered with nodes located at the vertices. 
 For $q=3$ a $3\times 3$ tensor grid of integration nodes is obtained.}
 \label{fig:quadrature_orders}
\end{figure}

\chapter{Formulation of the Optimization Problem\markboth{Formulation of the Optimization Problem}{}}
\label{sec:optimization_problem}

With the aggregated element densities and their sensitivities available from
\cref{sec:grid_mapping}, the design task can now be posed as an inverse problem.
The global parameter vector
\[
\vecsym\zeta=\bigl[(\vecsym z^{(1)})^\top \;\cdots\; (\vecsym z^{(n)})^\top\bigr]^\top
\in\R^{5n}
\]
is to be determined such that the aggregated pseudo–density field
\[
x \;\mapsto\; A\!\bigl(\rho_1(x;\vecsym z^{(1)}),\dots,\rho_n(x;\vecsym z^{(n)})\bigr)
\]
approximates a prescribed target field \(\rho^*:\Omega\to[0,1]\).
Evaluation is carried out on the finite element grid, where each element
\(\Omega_e\) is represented by its averaged density value
\(\bar{\rho}^e(\vecsym\zeta)\) as defined in \cref{sec:quadrature_rules}.
These element-wise quantities provide the basis for the objective functionals,
their analytic derivatives with respect to \(\vecsym\zeta\) and the geometric
constraints that delimit the feasible design space.

\section{Tracking objective}
\label{sec:tracking_objective_sec}

A least–squares functional is employed to penalize deviations from the target. The discrete objective is defined as
\begin{equation}
 F_{\mathrm{track}}(\vecsym\zeta)
 = \sum_{e\in\mathcal E} \bigl(\rho^*_e - \bar{\rho}^e(\vecsym\zeta)\bigr)^2 ,
 \label{eq:tracking_objective_sec}
\end{equation}
where $\rho^*_e$ denotes the prescribed target density and $\bar{\rho}^e(\vecsym\zeta)$ the averaged element density.
The functional rewards accurate coverage of the target while penalizing both voids and overshoot. Differentiation with respect to the pill parameters proceeds component–wise. For pill indices $a,b\in\{1,\dots,n\}$ and
component indices $i,j\in\{1,\dots,5\}$ one obtains
\begin{align}
\pp{F_{\mathrm{track}}}{z_i^{(a)}}
&= -2 \sum_{e\in\mathcal E}
\bigl(\rho^*_e-\bar\rho_e(\vecsym\zeta)\bigr)\,
\pp{\bar\rho_e}{z_i^{(a)}}, 
\label{eq:track_grad}\\[0.5em]
\ppp{F_{\mathrm{track}}}{z_i^{(a)}}{z_j^{(b)}}
&= 2 \sum_{e\in\mathcal E} \Bigl[
\pp{\bar\rho_e}{z_i^{(a)}}\,
\pp{\bar\rho_e}{z_j^{(b)}}
- \bigl(\rho^*_e-\bar\rho_e(\vecsym\zeta)\bigr)\,
\ppp{\bar\rho_e}{z_i^{(a)}}{z_j^{(b)}}
\Bigr].
\label{eq:track_hess}
\end{align}
The Hessian separates into two parts: the first term in
\eqref{eq:track_hess} is the Gauss–Newton contribution, always positive
semidefinite, while the second accounts for curvature effects that arise from the dependence of $\bar{\rho}^e$ on the design variables.

\begin{figure}[H]
 \centering
 \includegraphics[width=0.55\textwidth]{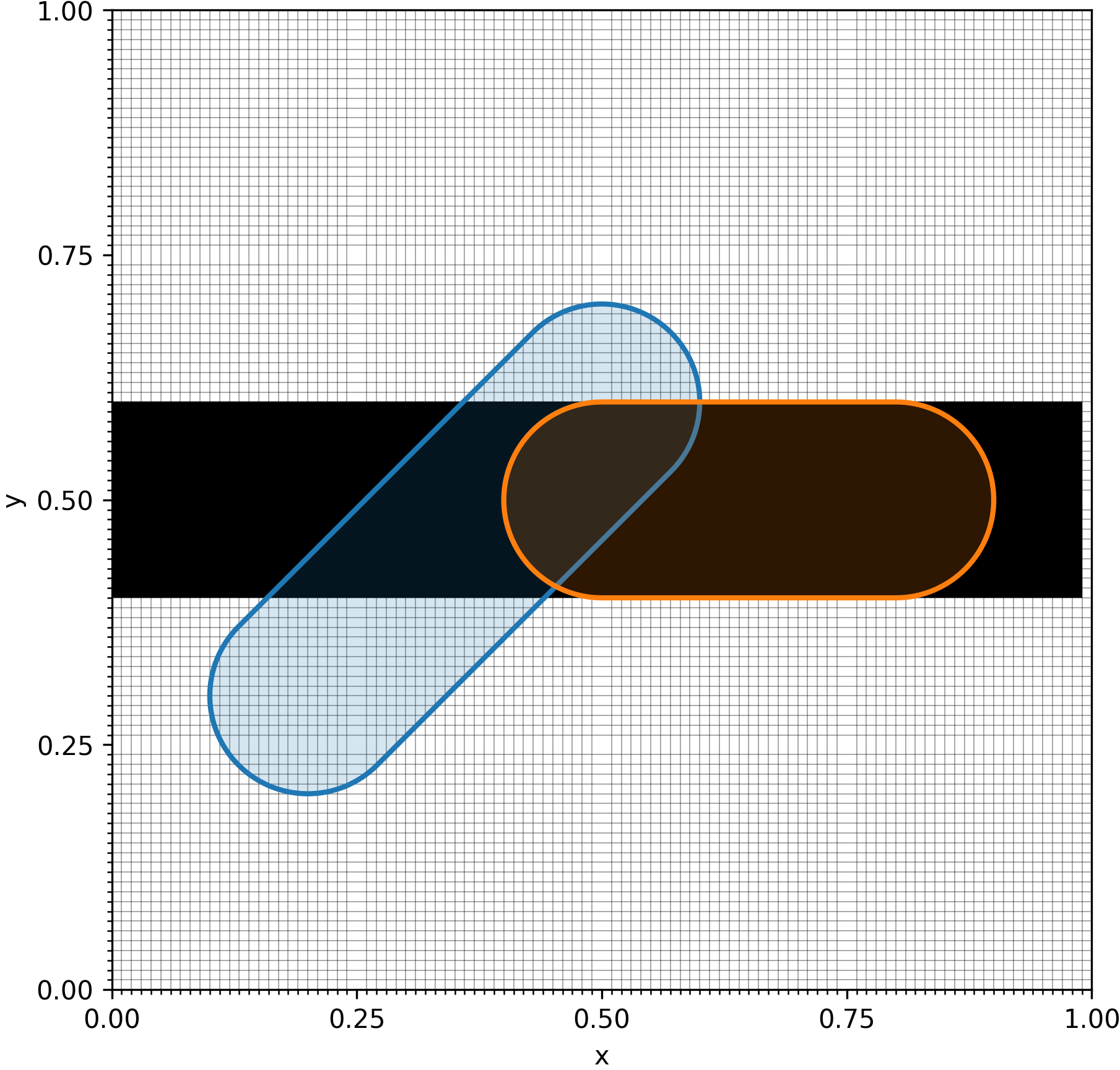}\\[1ex]
 \includegraphics[width=0.32\textwidth]{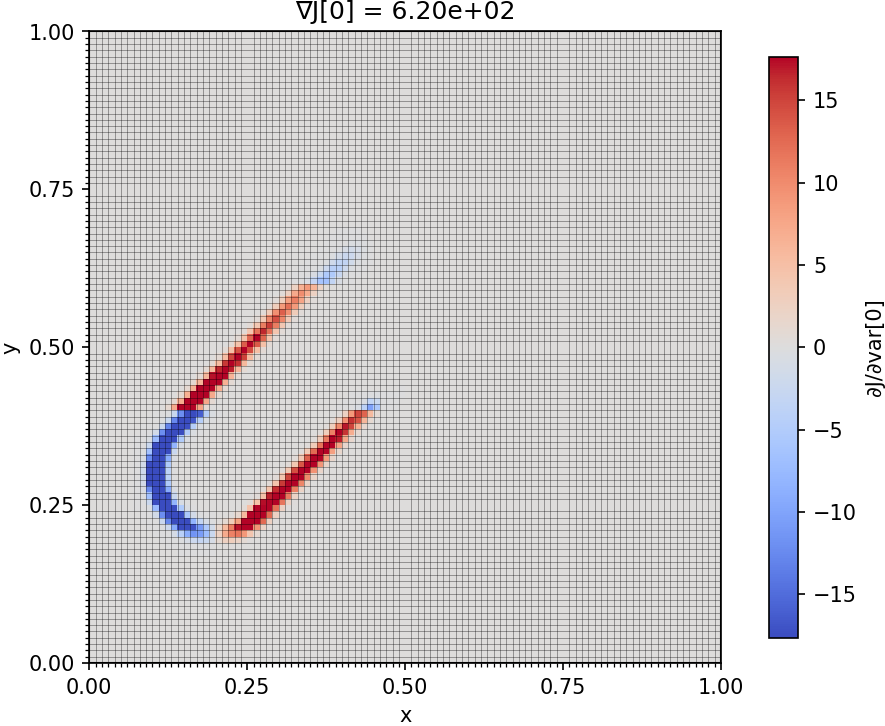}
 \includegraphics[width=0.32\textwidth]{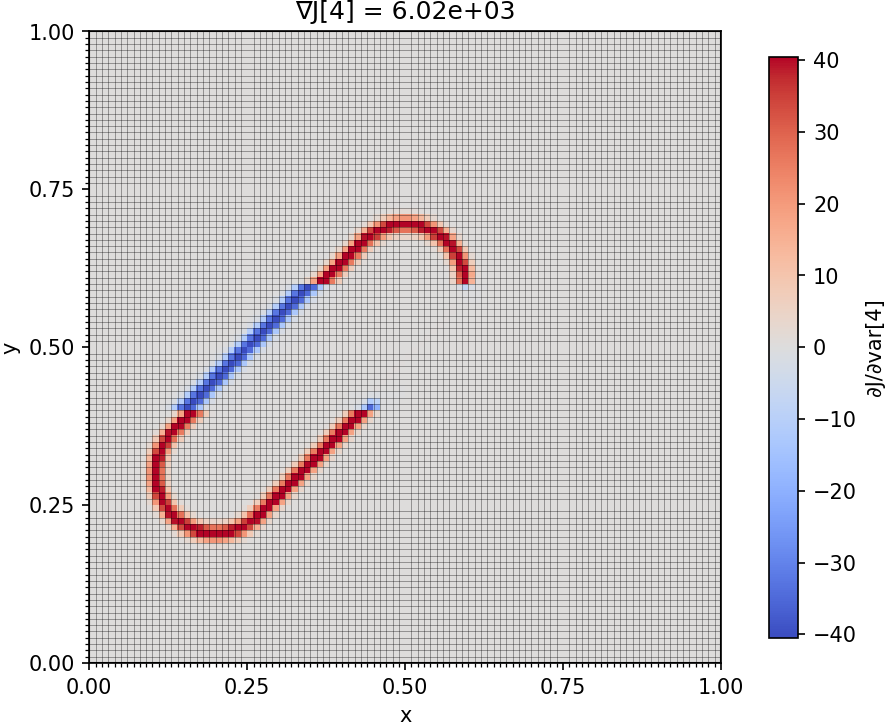}
 \includegraphics[width=0.32\textwidth]{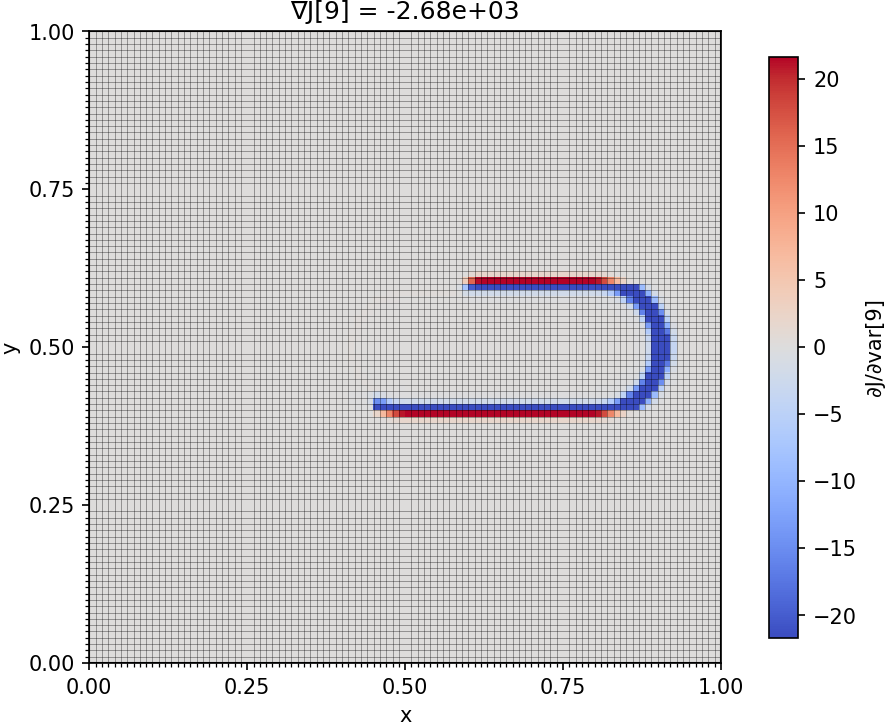}
 \caption{Tracking objective. Top: setup with two overlapping features (blue and orange) and a bar-shaped target (black). 
 Bottom: gradient fields with respect to the position parameter $p_x$ of the blue pill (left), 
 the radius of the blue pill (middle) and the radius of the orange pill (right).}
 \label{fig:tracking_obj_derivatives}
\end{figure}

The sensitivity structure in \cref{fig:tracking_obj_derivatives} confirms the locality of the gradients: nonzero contributions appear only in regions where the reconstructed density overlaps the target. Radius derivatives are strongly negative inside the target, driving features to shrink until the prescribed thickness is matched. Endpoint derivatives influence the orientation and placement of the segment, while overlap with neighboring features is reduced through redistribution across the domain. The functional therefore provides consistent and localized corrective signals that guide the design toward the target distribution.

\section{Reward-only objective}
\label{sec:reward_objective_sec}

As an alternative to the least–squares tracking functional, a reward–only formulation is introduced. It promotes overlap with the target without penalizing overshoot into void regions. A minus sign converts the maximization into a minimization problem, leading to the definition
\begin{equation}
F_{\mathrm{reward}}(\vecsym\zeta)
=\sum_{e\in\mathcal E} \bar{\rho}_e(\vecsym\zeta)*\rho^*_e.
\label{eq:reward_objective_sec}
\end{equation}
where $\rho^*_e$ denotes the prescribed target value in element $e$ and $\bar\rho_e$ the aggregated pseudo–density introduced in \cref{sec:grid_mapping}. 

Differentiation with respect to the design parameters yields only direct contributions, since the residual weighting present in \eqref{eq:tracking_objective_sec} is absent. For pill indices $a,b\in\{1,\dots,n\}$ and parameter components $i,j\in\{1,\dots,5\}$ one obtains
\begin{align}
\pp{F_{\mathrm{reward}}}{z_i^{(a)}} 
&= - \sum_{e\in\mathcal E} \rho^*_e\, \pp{\bar\rho_e}{z_i^{(a)}}, 
\label{eq:reward_grad}\\[0.5em]
\ppp{F_{\mathrm{reward}}}{z_i^{(a)}}{z_j^{(b)}}
&= - \sum_{e\in\mathcal E} \rho^*_e\, \ppp{\bar\rho_e}{z_i^{(a)}}{z_j^{(b)}}.
\label{eq:reward_hess}
\end{align}
The Gauss–Newton term that appears in the tracking case vanishes, so curvature originates exclusively from the second–order sensitivities of the element densities. This structural simplification has a distinct impact on the shape of the gradients. 

\begin{figure}[H]
 \centering
 \includegraphics[width=0.32\textwidth]{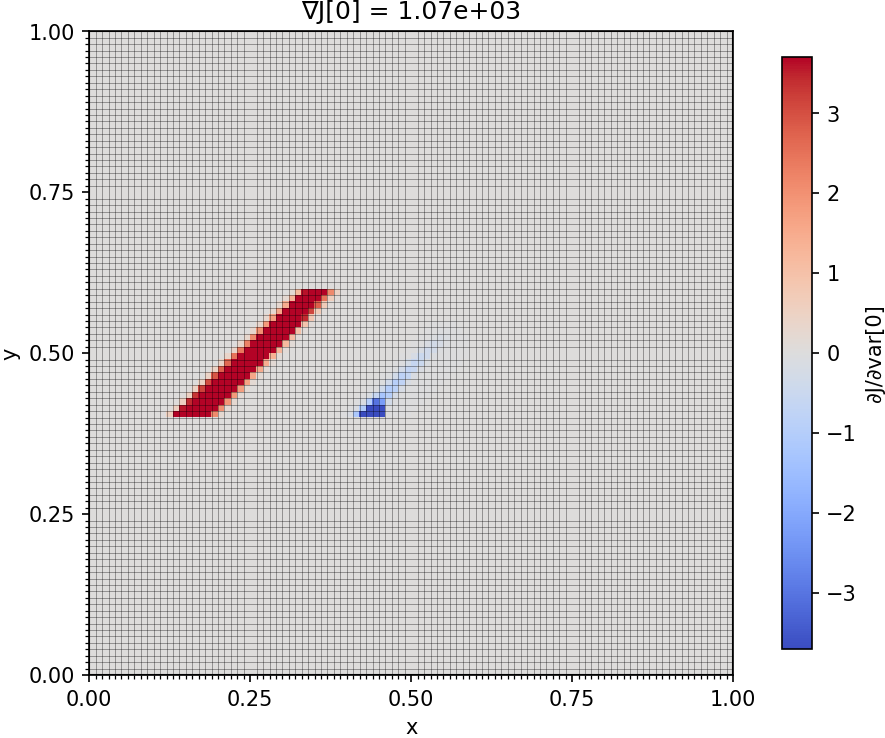}
 \includegraphics[width=0.32\textwidth]{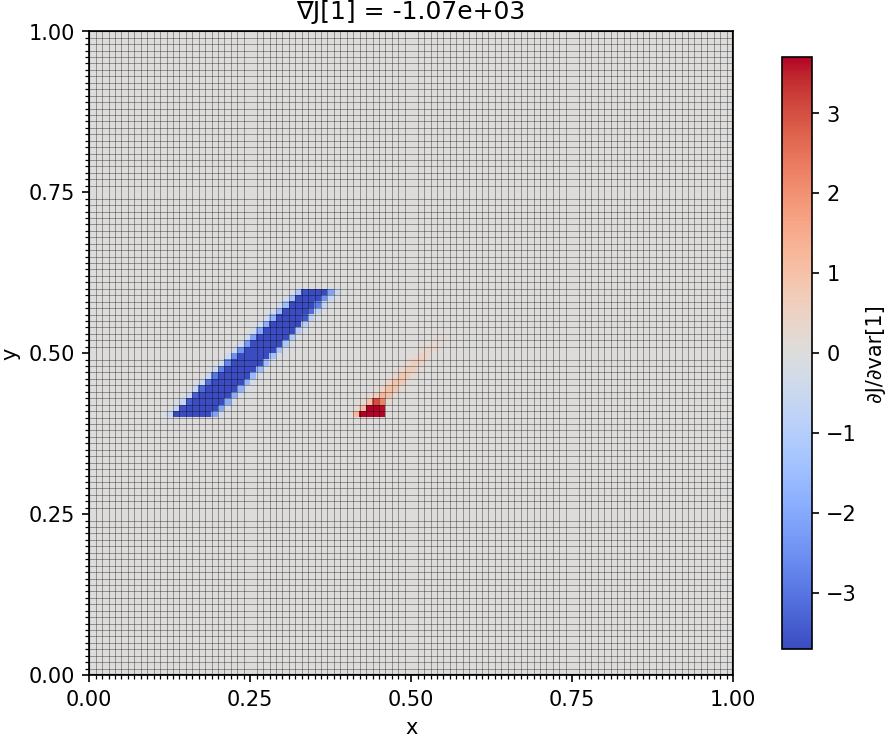}
 \includegraphics[width=0.32\textwidth]{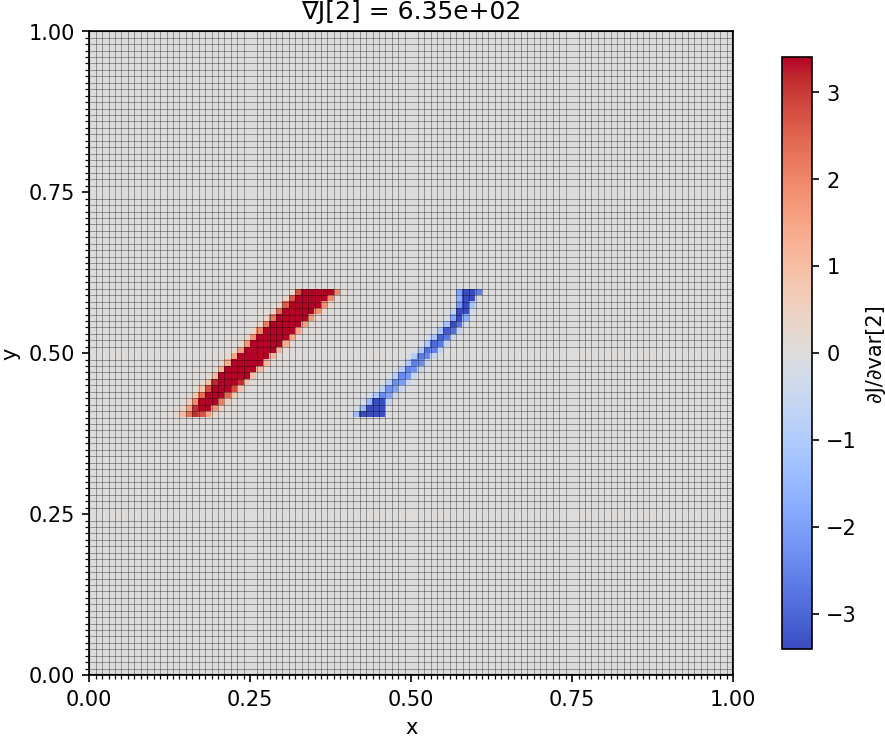}
 \includegraphics[width=0.32\textwidth]{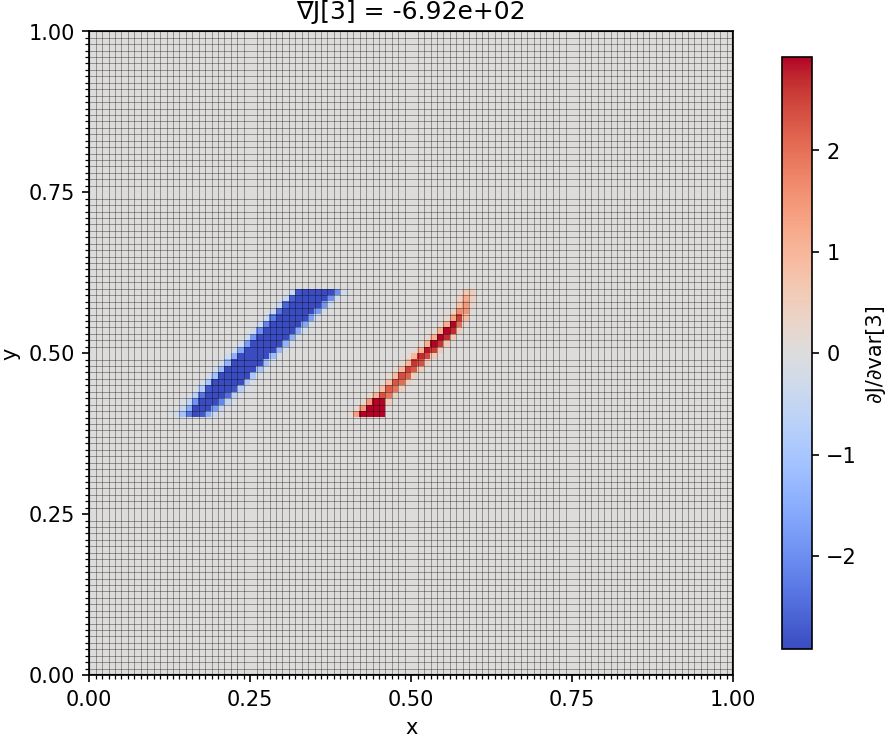}
 \includegraphics[width=0.32\textwidth]{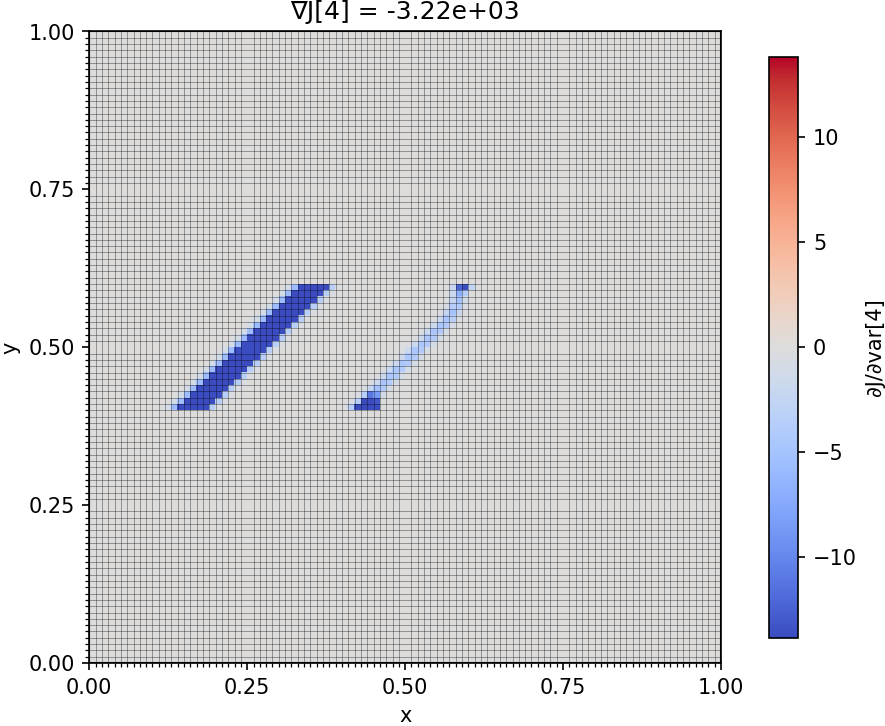}
 \includegraphics[width=0.32\textwidth]{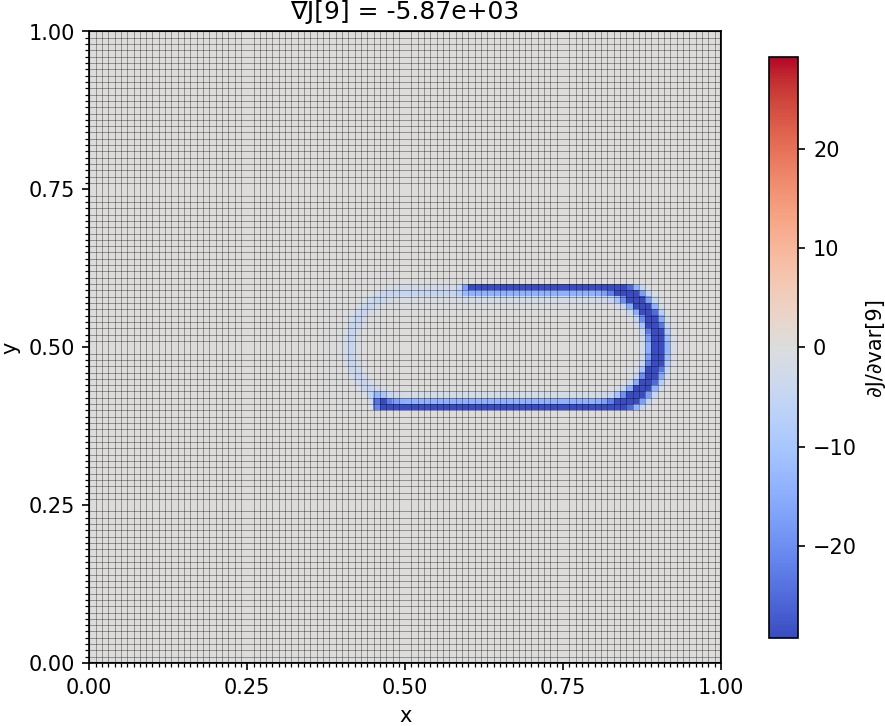}
 \caption{Reward-only objective. Gradient plots for the first pill; sensitivities remain active throughout the target region and consistently push features to expand.}
 \label{fig:reward_obj_derivatives}
\end{figure}

The qualitative behavior of the reward functional is illustrated in \cref{fig:reward_obj_derivatives}. Within the target region, sensitivities remain uniformly active, in contrast to the tracking case where activity diminishes once overlap is achieved. Derivatives with respect to the radius are consistently negative, thereby exerting an outward pressure on the pill boundary. As a consequence, features are expanded until the prescribed area fraction is saturated. This expansion is not gradual but characterized by sharp transitions: as soon as the contour of the pill intersects the target boundary, the gradient magnitude increases abruptly, producing a distinct driving force at the interface. 
The endpoint derivatives exhibit a complementary behavior. They generate alignment forces that rotate and translate the segment towards the target band, irrespective of whether parts of the pill already extend beyond the boundary. In this sense the reward functional does not penalize overshoot but continues to reinforce movement along the target direction. Compared to the tracking objective, the response is therefore more exploratory and less conservative: while tracking tends to stall once an approximate fit is reached, the reward objective maintains nonzero sensitivities that actively probe the entire target region. This property explains the consistently hard gradient transitions visible in the plots, which form a robust mechanism to prevent premature convergence but may also lead to overshoot if not balanced by additional constraints.

\section{Bounds and boundary bias}
\label{sec:bounds}

In principle, features are expected to remain inside the unit square $\Omega=[0,1]^2$, since target densities are defined only within this region. Nevertheless, due to finite step sizes in the optimization, control points may occasionally leave the domain. If no bounds are enforced, such excursions can persist, because the boundary bias described below prevents a natural return to the interior. For this reason, simple box constraints are imposed on the pill parameters,
\[
(p_x,p_y),(q_x,q_y)\in[0,1]^2,\qquad r\ge\delta,
\]
where $\delta>0$ denotes the transition half–width introduced in \cref{sec:transition_function}. An optional upper bound $r\le r_{\max}$ can be applied to restrict the pill radius.

Even with parameters formally constrained to the design domain, a systematic bias arises from the evaluation of the pseudo–density field. Since densities are computed only inside $\Omega$, any portion of a transition band extending beyond the boundary is truncated and does not contribute to the objective. As a result, features close to the boundary lose effective volume: in a corner only a quarter of the annulus remains, along an edge roughly half, whereas in the interior the full support contributes. This mechanism constitutes the boundary bias.

\begin{figure}[H]
  \centering
  \includegraphics[width=0.55\textwidth]{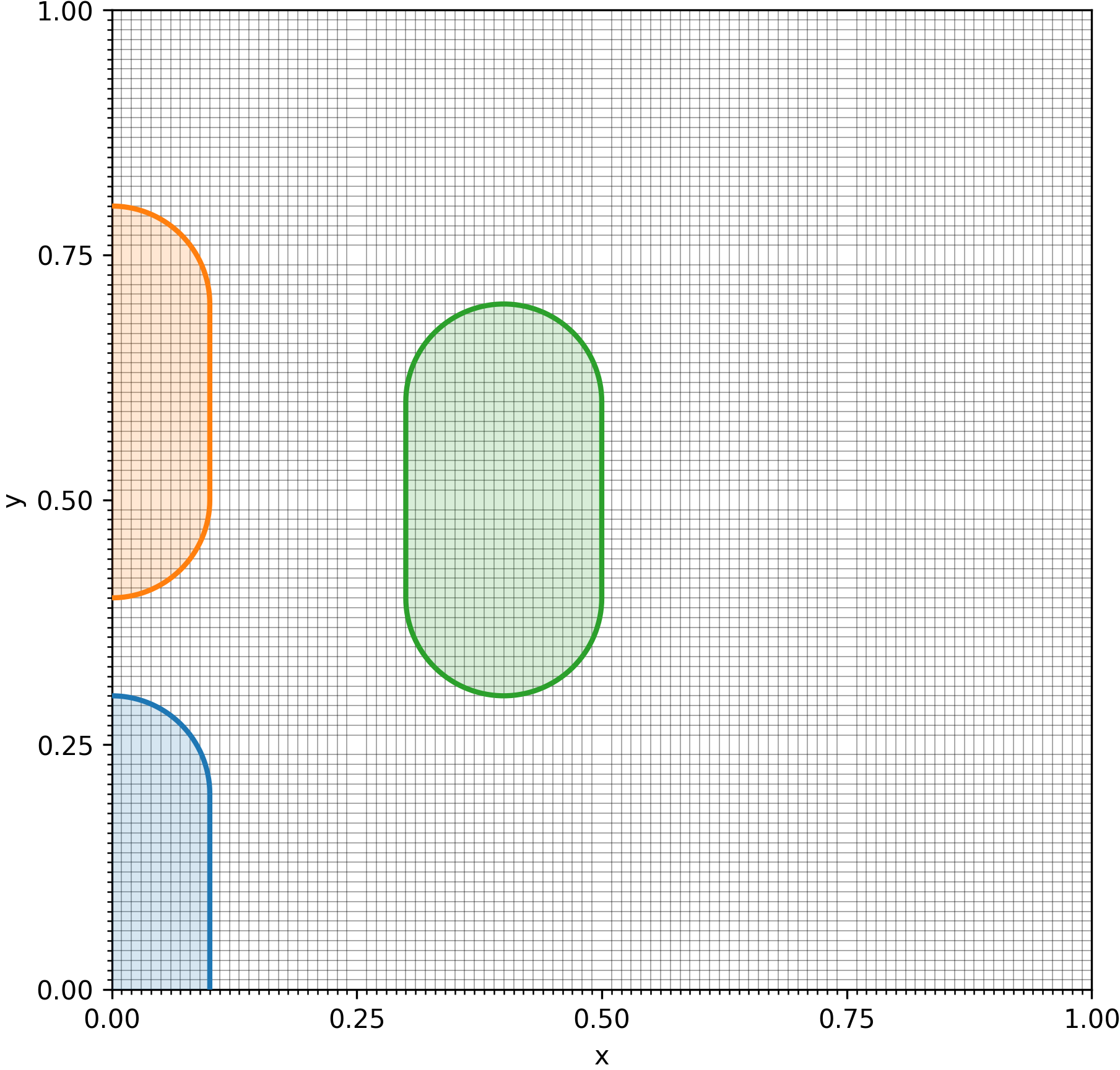}\\[1ex]
  \includegraphics[width=0.32\textwidth]{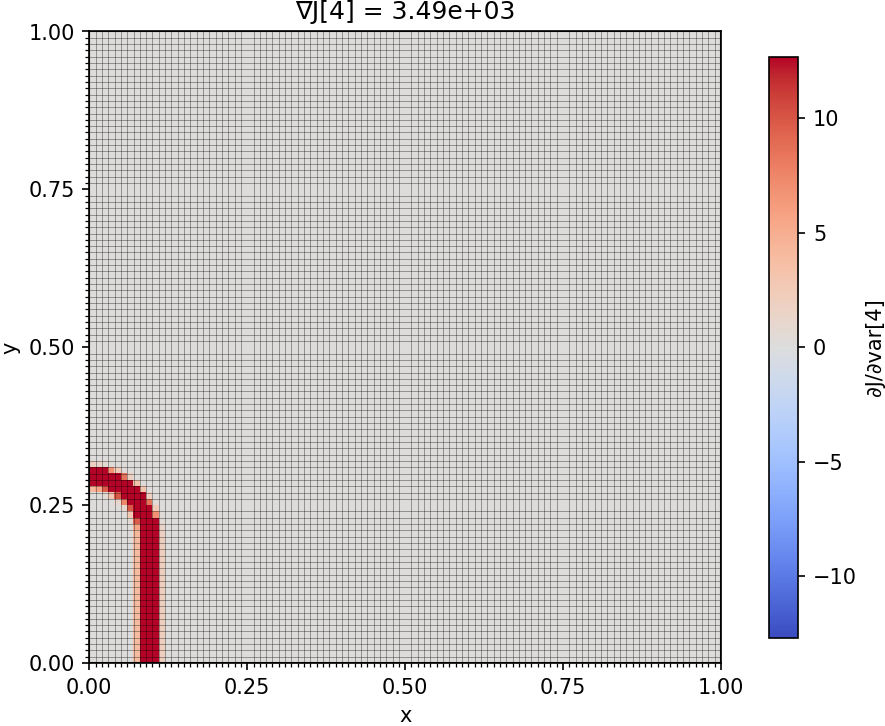}
  \includegraphics[width=0.32\textwidth]{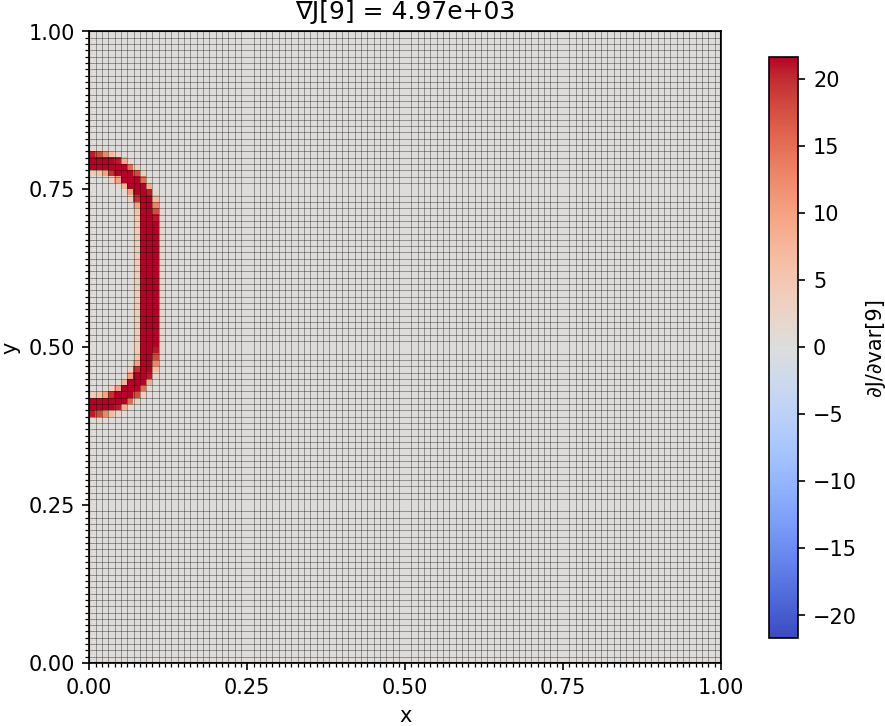}
  \includegraphics[width=0.32\textwidth]{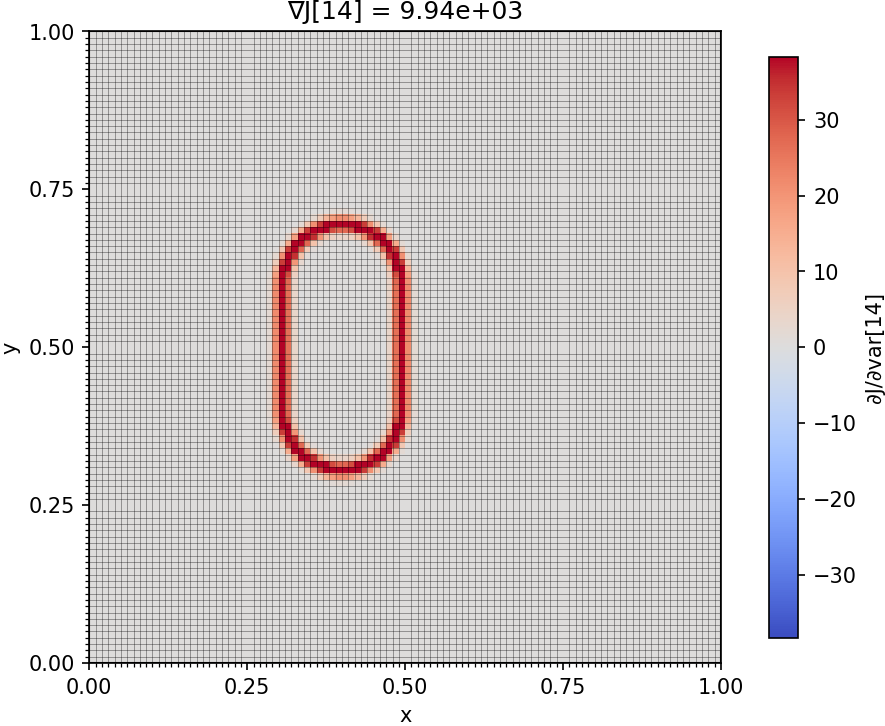}
  \caption{Boundary bias for capsules placed at the corner, edge and interior. 
  Top: setup with three identical features at different positions. 
  Bottom: gradients with respect to the radius parameter demonstrate the reduction of effective contributions near the boundary.}
  \label{fig:boundary_bias}
\end{figure}

The effect has two major implications for the tracking objective. First, the functional consistently undervalues boundary features, since truncated support lowers the effective density inside void regions. A pill can thereby decrease its objective contribution by shifting towards the boundary, even though its overlap with the target deteriorates. Second, sensitivities vanish outside the design domain, so no repulsive forces exist that could drive a boundary pill back into the interior. Together, these effects create a persistent attraction to corners and edges, where the apparent objective cost is artificially reduced. This behavior is illustrated in \cref{fig:boundary_bias}. The setup (top) shows three identical capsules positioned at a corner, along an edge and in the interior. The corresponding gradients with respect to the radius parameter (bottom) reveal the bias directly: near the boundary, truncation of the annulus reduces the gradient magnitude, thereby stabilizing the biased configurations. Only the interior pill experiences balanced sensitivities that fully reflect its geometric support. 

It is important to note that this mechanism is specific to the tracking objective. For the reward functional introduced in \cref{sec:reward_objective_sec}, density in void regions is not penalized, so a pill cannot lower its objective by moving towards the boundary. While sensitivities outside the domain remain truncated in the same way, the incentive to exploit this artifact disappears and the reward formulation is free of the boundary bias that affects the tracking case.

\section{Geometric constraints on segment length}
\label{sec:constraints}

To avoid degeneracy and to ensure numerical stability, a single global lower bound
on the segment length is imposed uniformly for all features. With $\ell^{(m)}=\|Q^{(m)}-P^{(m)}\|$, the admissible range reads
\[
\ell_{\min}\;\le\;\ell^{(m)}\;\;(\le\;\ell_{\max})\qquad\text{for all }m,
\]
where the upper bound $\ell_{\max}$ is optional.

The motivation is twofold. First, the degenerate configuration $P^{(m)}=Q^{(m)}$ must be excluded. Second, the analytic sensitivities of the segment distance contain factors of the form
$D^{-3}$ and $D^{-5}$ (with $D=\|P^{(m)}-Q^{(m)}\|$), see \cref{eq:hessian_dseg_general,eq:dseg_pp,eq:dseg_ppy,eq:dseg_pq}. As $D\to0$, these terms become large and render the gradient and, in particular, the Hessian ill-conditioned. A uniform bound $\ell_{\min}>0$ prevents this instability and stabilizes the early iterations.

In practice, the bound is enforced as a smooth inequality constraint on each pill using the common threshold $\ell_{\min}$; the corresponding derivatives follow directly from the expressions already derived and are therefore not repeated here.

\section{Constrained minimization problem}
\label{sec:constrained_problem}

Combining objectives and constraints yields the constrained optimization problem
\begin{equation}
\begin{aligned}
\min_{\vecsym{\zeta}\in\R^{5n}} \quad & \Phi(\vecsym{\zeta}) \\[0.3em]
\text{s.t.}\quad 
& (p_x^{(m)},p_y^{(m)}),(q_x^{(m)},q_y^{(m)})\in[0,1]^2, 
  \qquad m=1,\dots,n, \\[0.3em]
& r^{(m)} \ge \delta \quad (\text{optionally } r^{(m)} \le r_{\max}), 
  \qquad m=1,\dots,n, \\[0.3em]
& \ell^{(m)} \;\ge\; \ell_{\min} \quad \bigl(\text{optionally } \ell^{(m)} \le \ell_{\max}\bigr), 
  \qquad m=1,\dots,n,
\end{aligned}
\label{eq:final_constrained_form}
\end{equation}
where $\ell^{(m)}=\|Q^{(m)}-P^{(m)}\|$ denotes the segment length of pill $m$.
Here, $\Phi$ denotes either the tracking functional \eqref{eq:tracking_objective_sec} or the reward functional \eqref{eq:reward_objective_sec}. The box constraints enforce domain bounds for the endpoints and lower (optionally upper) bounds for the radii. A length threshold $\ell_{\min}>0$ is applied uniformly to all features (optionally $\ell_{\max}$); see \cref{sec:constraints} for the motivation and for derivative expressions implied by the sensitivity formulas developed earlier (e.g., \eqref{eq:hessian_dseg_general}–\eqref{eq:dseg_pq}). 

All components of the formulation are thus available: analytic objectives, box bounds and smooth length inequalities with their derivatives. Together they define a well-posed constrained minimization problem in closed analytic form. The subsequent chapter turns from the mathematical model to its algorithmic treatment, introducing initialization strategies, staged continuation and the solver framework that exploits the block-sparse structure of the gradients and Hessians.

\chapter{Optimization strategy\markboth{Optimization strategy}{}}
\label{sec:optimization_strategy}

With the constrained formulation \eqref{eq:final_constrained_form} and closed-form sensitivities in place, attention shifts from modeling to computation. A practical solution procedure is organized around three components that act on the admissible set $\mathcal A$: a feasible initialization, a staged optimization schedule and lightweight heuristics for model size control. The design vector is denoted by $\vecsym{\zeta}$ throughout for consistency with \cref{sec:constrained_problem}.

Initialization (\cref{sec:init}) provides admissible starting designs by projecting cross-seeded (and optionally randomized) pill constellations onto $\mathcal A$, thereby ensuring domain bounds and the global segment-length threshold while delivering uniform coverage of $\Omega=[0,1]^2$. These states supply stable warm starts for subsequent solver runs.

A staged procedure (\cref{sec:stages}) then advances the optimization in three phases. An exterior radius inflation $\varepsilon_{\mathrm{ext}}\ge0$ is used to widen transition bands during early exploration and is reduced to zero as resolution is increased. The exploration stage employs the reward objective to promote target overlap and robust orientation; the bridging stage switches to the tracking objective under reduced extension to convert overlap into residual reduction; the convergence stage removes extension entirely and sharpens transition and aggregation parameters $(\delta,p,\beta)$ under tight tolerances. Phases are warm-started from their predecessors and acceptance criteria enforce feasibility and monotone decrease of the active objective.

To control redundancy and stabilize multi-feature layouts, heuristic editing is applied before the final stage. A pruning and merging pass (\cref{subsec:heuristics}) removes weak or non-unique contributors and consolidates aligned neighbors based on area and exclusivity ratios computed from the analytic density fields. When residual structure remains, an additive refinement loop (\cref{subsec:additive}) augments the design by orienting a candidate with the reward objective on a residual mask and admitting it only upon demonstrable improvement under the tracking objective. All subproblems use the same admissible set and sensitivity chain as the main formulation, preserving sparsity and block structure.

The remainder of the chapter details these components and the associated parameter schedules and solver settings, emphasizing the exploitation of locality and block sparsity in gradient and Hessian evaluations to enable large-scale feature optimization.

\section{Initialization}
\label{sec:init}

Initialization is crucial for the stability and efficiency of the subsequent
optimization. A feasible design vector $\vecsym{\zeta}^{(0)}\in\R^{5n}$ must satisfy
the box and length constraints from \cref{sec:bounds,sec:constraints} while
providing sufficiently uniform coverage of the design domain
$\Omega=[0,1]^2$.

A practical approach is \emph{cross seeding}, where $\Omega$ is partitioned into
an $R\times C$ grid of cells with side lengths $d_x=1/C$ and $d_y=1/R$. In each
cell, up to two diagonal segments are placed. With cell center
$\vecsym c=(c_x,c_y)$ and admissible segment length
\[
L=\min\!\bigl(0.95\,\sqrt{d_x^2+d_y^2},\;\ell_{\max}\bigr),
\]
the orthogonal unit vectors
\[
\hat{\vecsym u}=\tfrac{(d_x,d_y)}{\sqrt{d_x^2+d_y^2}},\qquad
\hat{\vecsym v}=\tfrac{(d_x,-d_y)}{\sqrt{d_x^2+d_y^2}}
\]
define the endpoint pairs
\[
(P,Q)=\bigl(\vecsym c\mp\tfrac{L}{2}\hat{\vecsym u}\bigr),\qquad
(P,Q)=\bigl(\vecsym c\mp\tfrac{L}{2}\hat{\vecsym v}\bigr).
\]
Each pill is assigned a radius $r^{(m)}\ge\delta$ and projected onto the admissible set $\mathcal A$ to enforce box bounds. Placement is truncated once the prescribed number $n$ of pills is reached. The resulting configuration provides a symmetric and uniformly distributed initialization.

To mitigate artifacts caused by symmetry, a \emph{randomized cross seeding} variant is considered. It follows the same construction but perturbs each diagonal by a small random rotation angle $\theta\sim\mathcal U(-\theta_{\max}, \theta_{\max})$. In this way, exact orthogonality and uniform alignment of the segments are avoided, while coverage and feasibility are preserved. Consequently, the risk of convergence to symmetric but suboptimal stationary points is reduced.

\section{Staged optimization}
\label{sec:stages}

\begin{figure}[t]
 \centering
 \includegraphics[width=.65\linewidth]{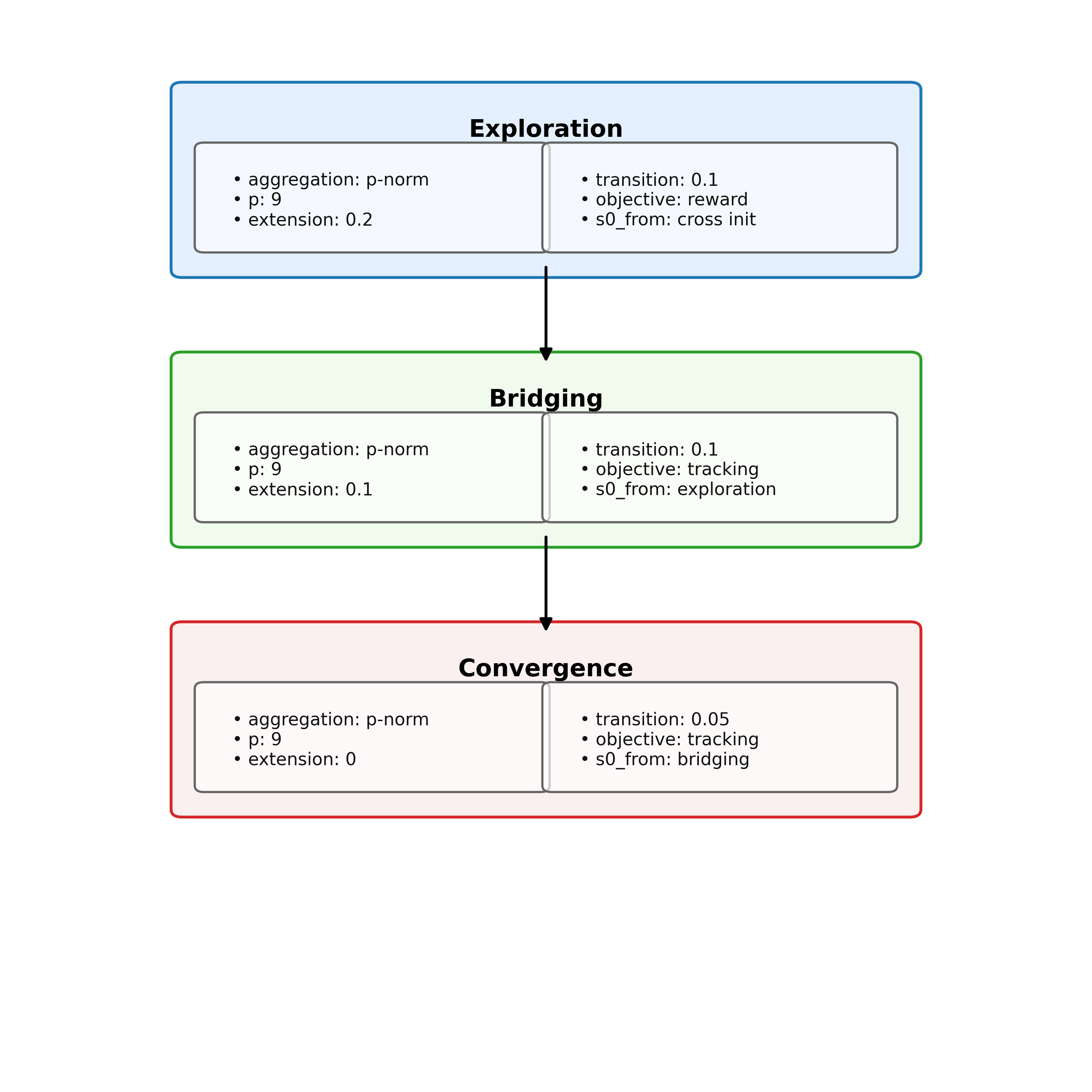}
 \caption{Three-stage pipeline: exploration with inflated radii $\to$ bridging with reduced inflation $\to$ convergence at full resolution.}
 \label{fig:pipeline_flow}
\end{figure}

The optimization is carried out in three successive stages that share the admissible set $\mathcal A$ but differ in objective, continuation parameters and solver tolerances. A temporary radius extension $\varepsilon_{\mathrm{ext}}\ge0$ is applied in the early iterations, inflating each pill radius and thereby broadening transition bands and enlarging attraction basins. This extension is gradually reduced to zero, while transition and aggregation parameters $(\delta,p,\beta)$ are sharpened and tolerances tightened, so that resolution increases progressively. In the exploration stage the reward functional with inflated radii guides pills toward dominant regions of the target, prioritizing feasibility and coarse alignment over accuracy. The bridging stage replaces the objective by the tracking functional under reduced extension, converting overlap into residual reduction while avoiding premature collapse. Finally, the convergence stage removes the extension entirely and employs refined parameters under strict tolerances, achieving full resolution. The stages are executed sequentially with warm starts, each initialized by the final iterate of its predecessor. Acceptance requires strict feasibility together with monotone decrease of the active objective: the reward functional in the first stage and the tracking functional thereafter. The parameter schedules underlying this progression are summarized in \cref{fig:pipeline_flow}.

\subsection{Heuristics}
\label{subsec:heuristics}

To reduce redundancy in the feature set before the final optimization stage, 
heuristic selection and consolidation steps are applied. 
Each pill is evaluated by two normalized measures derived from its soft area footprint 
\(\chi_m(\cdot;\vecsym{\zeta}^{(m)}):\Omega\to[0,1]\), induced through transition and aggregation 
(\cref{sec:transition_function,sec:aggregation}). 
The soft area and the total coverage are
\[
A_m = \int_\Omega \chi_m(x;\vecsym{\zeta}^{(m)})\,dx,
\qquad 
A=\sum_{j=1}^n A_j,
\]
from which the area ratio \(\mathrm{AR}_m=A_m/A\) and the unique–region ratio
\[
u_m(x)=\chi_m(x)\!\prod_{j\neq m}(1-\chi_j(x)),\qquad
\mathrm{UR}_m=\frac{\int_\Omega u_m(x)\,dx}{A_m}
\]
are obtained. While \(\mathrm{AR}_m\) expresses relative size, \(\mathrm{UR}_m\) quantifies exclusive support. 
Features with \(\mathrm{AR}_m<\alpha\) or \(\mathrm{UR}_m<\upsilon\) are removed, 
using default thresholds \(\alpha=0.15\) and \(\upsilon=10^{-4}\). 
The surviving set is retained for possible consolidation.

Grouping is then performed among nearly parallel and adjacent pills. 
For each survivor, the orientation vector and angle are
\[
\hat u_m=\frac{Q^{(m)}-P^{(m)}}{\|Q^{(m)}-P^{(m)}\|},\qquad 
\theta_m=\operatorname{atan2}(\hat u_{m,y},\hat u_{m,x})\in(-\pi,\pi],
\]
and proximity is assessed either by center distance 
\(d_{mn}=\|\tfrac{P^{(m)}+Q^{(m)}}{2}-\tfrac{P^{(n)}+Q^{(n)}}{2}\|_2\) 
or by minimal segment distance \(d_{mn}^{\mathrm{seg}}\). 
Two pills are linked if their angular difference is within \(\theta_{\lim}\) 
and their distance below \(d_{\min}\). 
Connected components of this adjacency relation form groups. 
Within each group, the representative is chosen as the longest segment, 
and its radius is assigned as the minimum group radius truncated at zero:
\[
m^\star \in \arg\max_{m\in\mathcal G}\|Q^{(m)}-P^{(m)}\|,
\qquad
(\widehat P,\widehat Q)=(P^{(m^\star)},Q^{(m^\star)}),
\qquad
\widehat r=\max\!\Bigl(0,\min_{m\in\mathcal G}r^{(m)}\Bigr).
\]
The merged pill is then \(\widehat{\vecsym{z}}_{\mathcal G}=(\widehat P_x,\widehat P_y,\widehat Q_x,\widehat Q_y,\widehat r)\), 
while ungrouped survivors are copied unchanged. 
The procedure is deterministic for fixed parameters, preserves feasibility, 
and is controlled by the thresholds \((\alpha,\upsilon,\theta_{\lim},d_{\min})\).

\subsection{Iterative refinement}
\label{subsec:additive}

In addition to pruning and merging, the feature set can also be augmented by systematically adding new primitives. This additive refinement proceeds iteratively and incorporates a candidate only if it demonstrably decreases the tracking functional. Given an initial configuration $\vecsym{\zeta}^{(0)}\in\R^{5n}$ and target distribution $S^\star$, a residual mask $R$ is computed from the current synthesized field $S_{\mathrm{cur}}$ by marking all regions where the deficit $S^\star-S_{\mathrm{cur}}$ exceeds a prescribed threshold $\tau_{\mathrm{res}}$,
\[
R(x) =
\begin{cases}
1, & S^\star(x)-S_{\mathrm{cur}}(x) > \tau_{\mathrm{res}},\\[0.3em]
0, & \text{otherwise}.
\end{cases}
\]
If $R\equiv 0$, the procedure terminates. Otherwise, a new pill is initialized at the center with short diagonal endpoints and seed radius $r_{\mathrm{seed}}$, optionally constrained to a fixed value $\bar r$. This primitive is first oriented against the residual by minimizing the reward functional subject to admissible bounds,
\[
\vecsym{\zeta}_{\mathrm{new}}^{\mathrm{ori}}
=\arg\min_{t\in\mathcal A_{\mathrm{add}}} F_{\mathrm{reward}}(t;R),
\]
and then refined by minimizing the tracking functional on $R$ without the cap, yielding $\vecsym{\zeta}_{\mathrm{new}}^{\mathrm{conv}}$. If no prior design exists, this pill becomes the initial configuration; otherwise it is concatenated with the current best $\vecsym{\zeta}_{\mathrm{best}}$,
\[
\vecsym{\zeta}_{\mathrm{cand}}=\bigl[\,\vecsym{\zeta}_{\mathrm{best}};\ 
\vecsym{\zeta}_{\mathrm{new}}^{\mathrm{conv}}\,\bigr]\in\R^{5(n+1)},
\]
and subjected to a full convergence against the true target $S^\star$, producing $\vecsym{\zeta}_{\mathrm{full}}$.

Acceptance is based on improvements of the mean–squared error $J$,
\[
\Delta_{\mathrm{abs}}=J(\vecsym{\zeta}_{\mathrm{best}})-J(\vecsym{\zeta}_{\mathrm{full}}),\qquad
\Delta_{\mathrm{rel}}=\frac{J(\vecsym{\zeta}_{\mathrm{best}})-J(\vecsym{\zeta}_{\mathrm{full}})}
{\max\{J(\vecsym{\zeta}_{\mathrm{best}}),10^{-12}\}},
\]
and the candidate is retained if either $\Delta_{\mathrm{abs}}>\varepsilon_{\mathrm{abs}}$ or $\Delta_{\mathrm{rel}}>\varepsilon_{\mathrm{rel}}$. In that case the current best is updated, a new residual mask is computed, and the loop continues until $K_{\max}$ additions are reached, $R$ becomes empty or acceptance fails. All optimization substeps are carried out with IPOPT under the same bound and radius constraints as in the main problem.

Together, pruning, grouping, and additive refinement provide heuristics that streamline the feature set, close coverage gaps, and stabilize the multi–stage optimization pipeline. Their deterministic formulation and compatibility with the analytic sensitivity framework ensure seamless integration with the solver, and these preparations conclude the methodological part of the work.

\chapter{Numerical Results and Evaluation\markboth{Numerical Results and Evaluation}{}}
\label{sec:numerical_results}

Having established the optimization strategy, the subsequent analysis is concerned with the numerical behavior of the framework. Rather than re-deriving analytic formulas, the focus lies on their practical effect: the conditions under which reliable convergence is obtained, the failure modes that emerge in the absence of safeguards and the influence of initialization, geometric constraints, transition width and second–order information on the overall robustness of the method.

To make these aspects transparent, the numerical experiments are organized in a progressive sequence. The analysis begins with single--pill reconstructions, where sensitivity patterns and failure modes can be isolated most clearly. On this basis, the effect of geometric constraints and transition parameters is examined, demonstrating how safeguards prevent collapse or overextension. The study then extends to multi--feature problems, highlighting aggregation effects, overlap handling and the role of continuation strategies. Each stage builds on the previous one, so that challenges are revealed in simple settings before solutions are generalized to more complex cases.

The chapter therefore provides a structured evaluation: it identifies limitations, illustrates corrective mechanisms and motivates the final configuration adopted for large--scale reconstructions. In this way the experiments establish the link between analytic design and numerical performance, preparing the ground for the application scenarios studied in the following chapters.

\section{Single-Feature Optimization}
\label{sec:single_feature_optimization}

As a first diagnostic step, the framework is evaluated on targets that are almost exactly representable by a single pill. Computations are performed on a uniform Cartesian grid $\Omega_h\subset[0,1]^2$ with $100\times100$ cells. Pseudo–densities are obtained by composing the signed distance with a polynomial smoothstep transition with $k=3$, ensuring $\mathcal C^2$–continuity across feature boundaries as described in \cref{sec:transition_function}.

A batch of fifty targets is generated by sampling endpoints $P^\star,Q^\star\in[0,1]^2$ with a minimal separation of $0.2$ and radii $r^\star\in[0.05,0.25]$. For each case the initial parameters $\vecsym{\zeta}_0=(P,Q,r)$ are drawn independently of the target, so that the initial overlap ranges from substantial to none. Optimization proceeds with the tracking functional from \cref{sec:tracking_objective_sec}, using analytic first- and second–order derivatives and enforcing box bounds $(p_x,p_y),(q_x,q_y)\in[0,1]^2$ together with $r\ge\delta$. For all experiments the transition half–width was fixed at $\delta=0.05$. A run is considered successful if it terminates due to satisfaction of the tolerance criterion, i.e.\ an objective reduction below $10^{-8}$, rather than by reaching the maximum iteration limit of 100.

\begin{figure}[H]
  \centering
  \includegraphics[width=\textwidth]{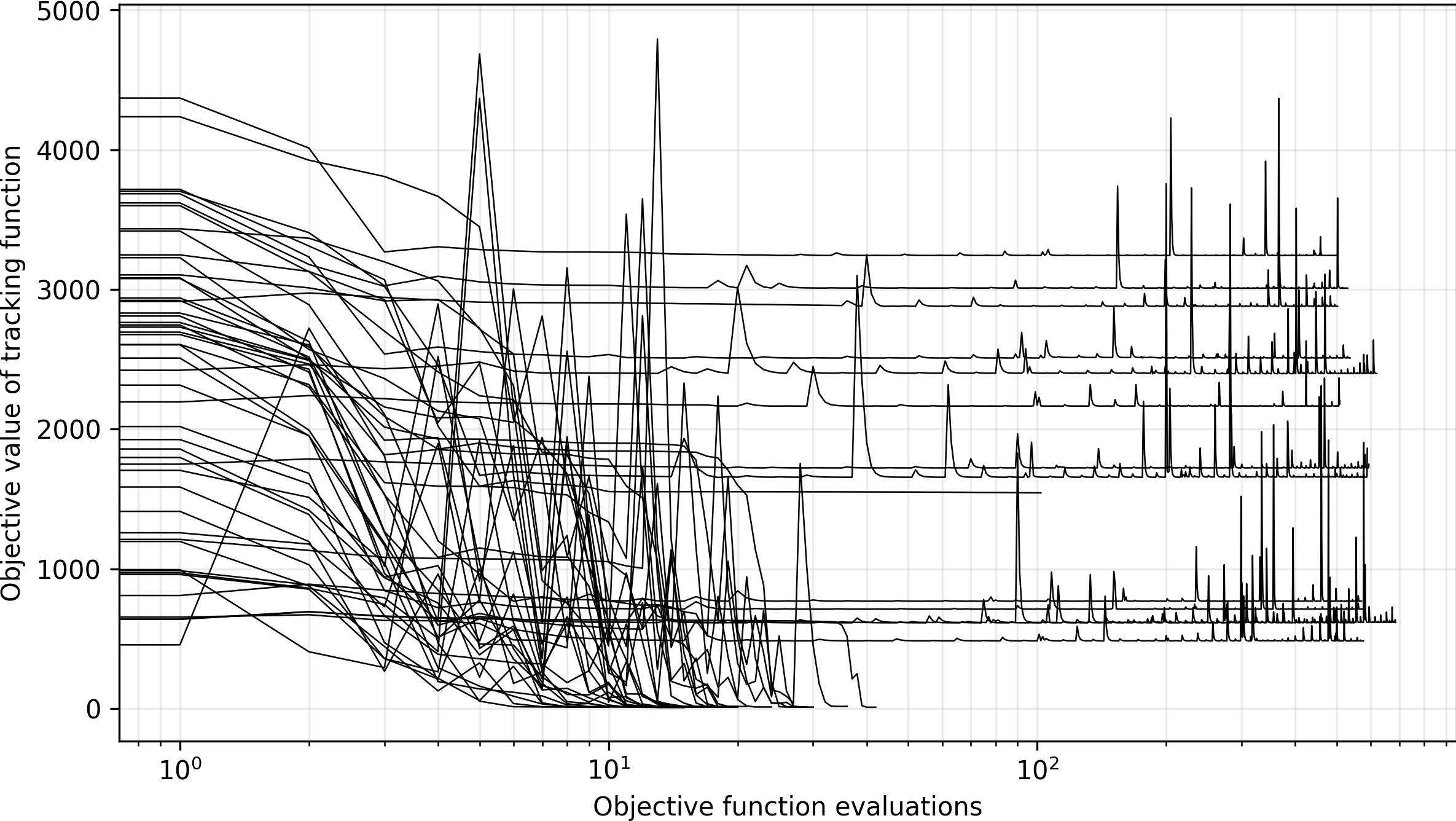}
  \caption{Objective trajectories for 50 single-feature reconstructions on a logarithmic scale. Convergence, when attained, occurs within the first fifty iterations, whereas 13 runs stagnate with irregular jumps of the objective.}
  \label{fig:all_runs_lineplot}
\end{figure}

The aggregate behavior across all runs is displayed in \cref{fig:all_runs_lineplot}. The trajectories show highly irregular progress, with many exhibiting sudden jumps in the objective value. Nevertheless, in all convergent cases the decrease occurs within the first fifty iterations, after which the objective stabilizes close to zero. Out of the fifty runs, 37 converge reliably, while 13 stagnate despite continued oscillations of the objective. The trajectories therefore reveal two distinct behaviors: in the majority of runs the objective decays rapidly and stabilizes near zero, while in a considerable number of cases the progress halts and the values fluctuate around elevated plateaus. The underlying reasons for this stagnation will be examined in the subsequent analysis.

\begin{figure}[H]
  \centering
  \includegraphics[width=0.32\textwidth]{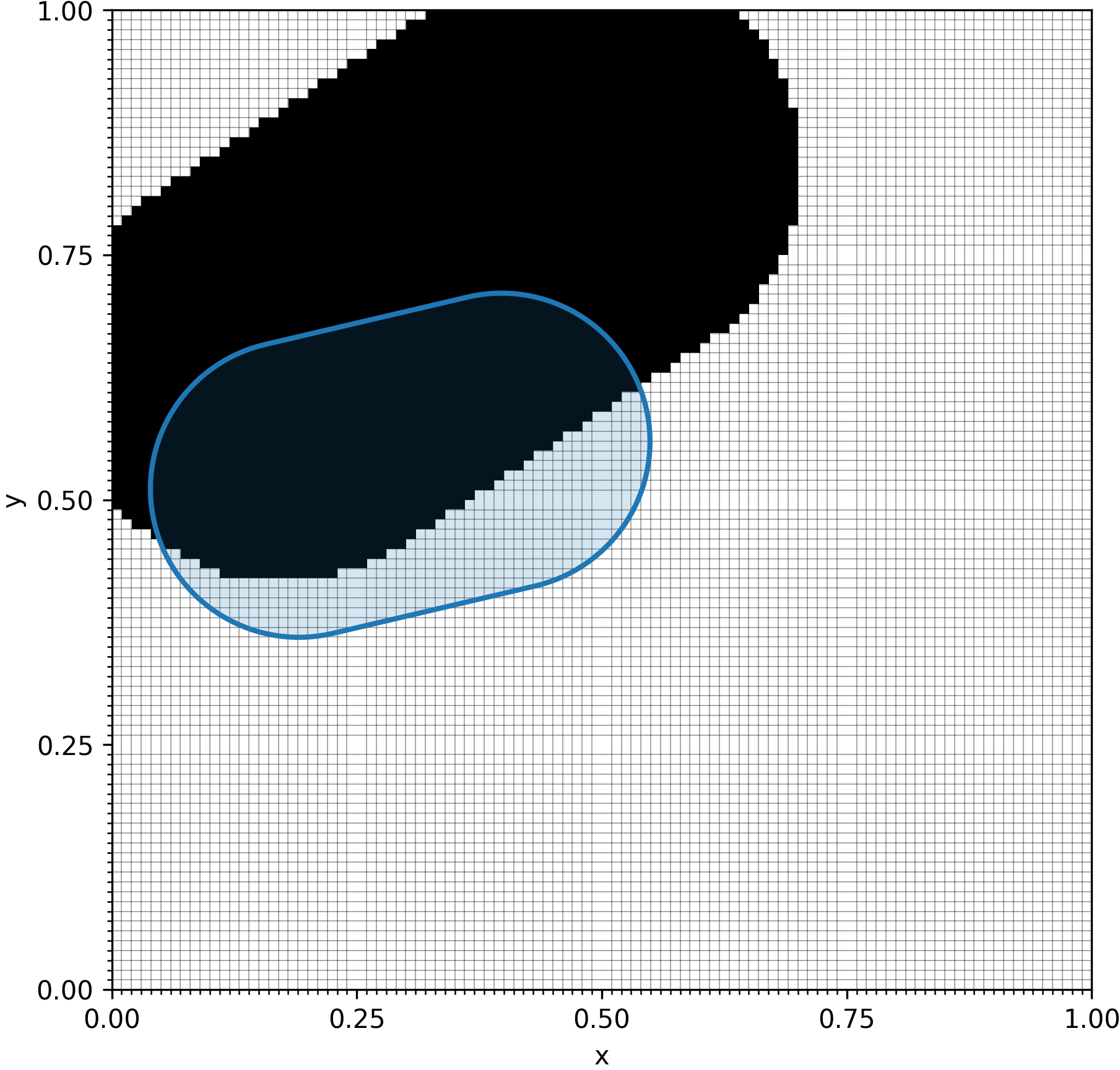}\hfill
  \includegraphics[width=0.32\textwidth]{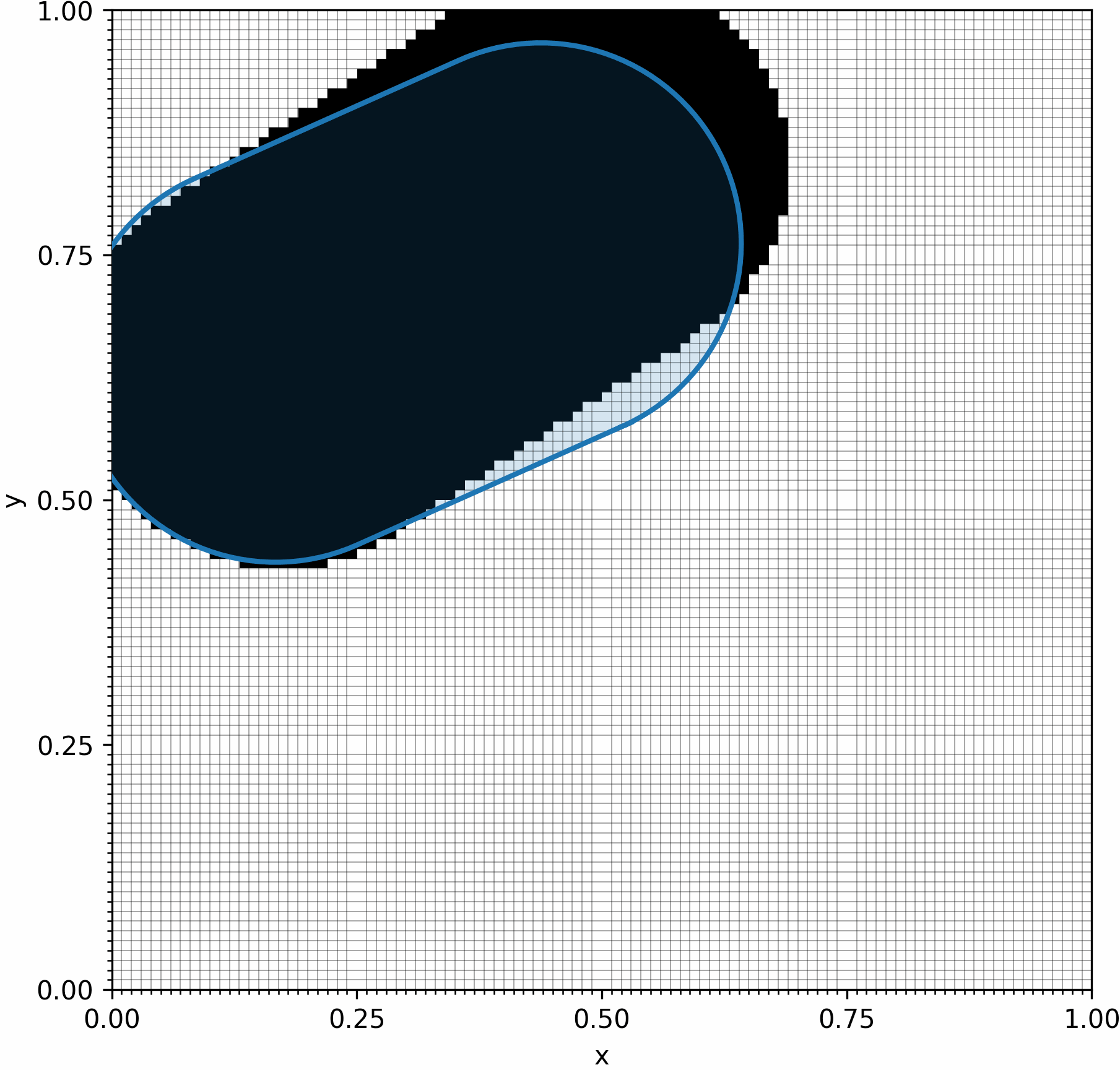}\hfill
  \includegraphics[width=0.32\textwidth]{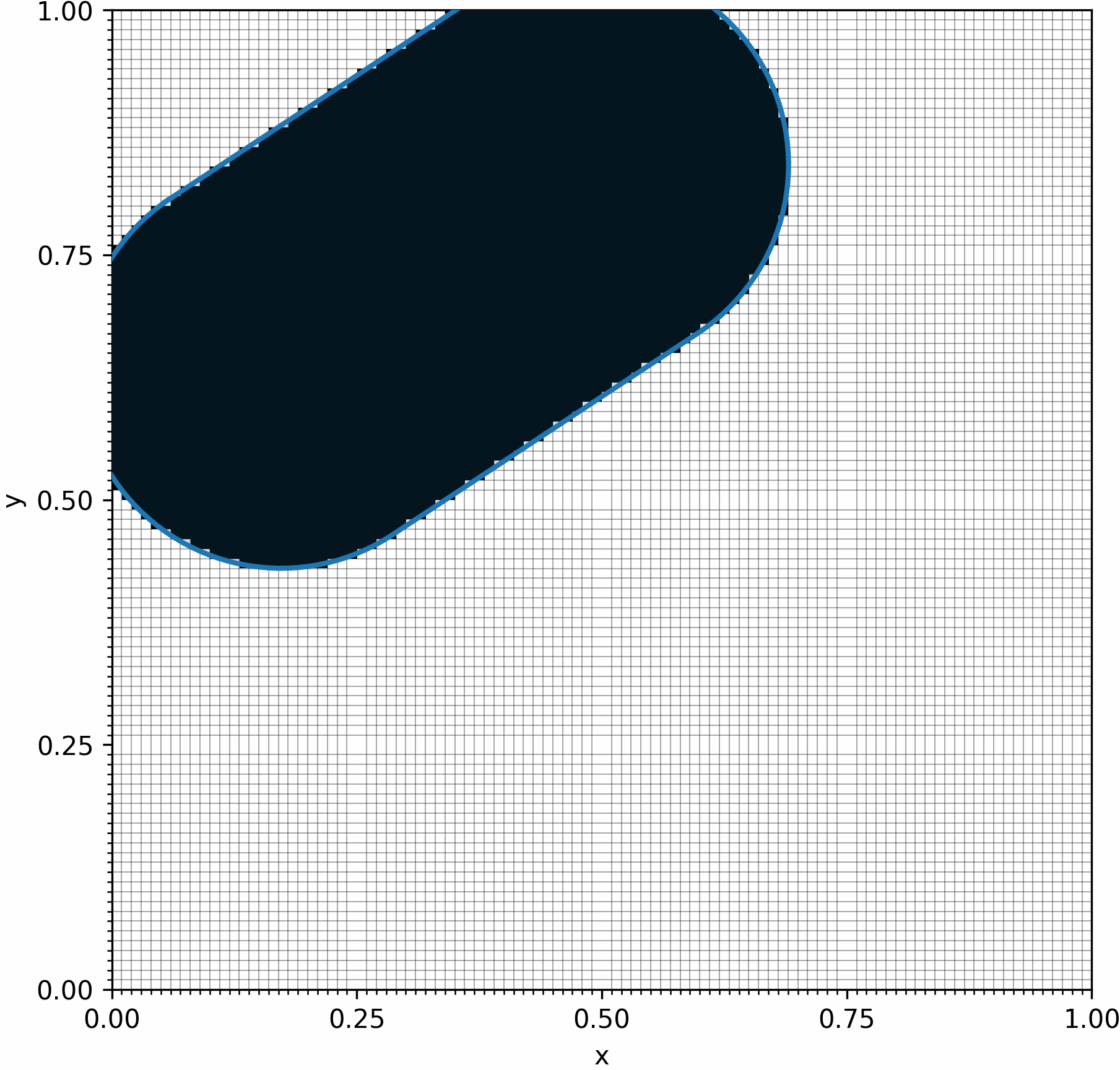}
  \caption{Representative success from the single-feature batch. Starting from a favorable initialization (left), the pill adapts within a few iterations (middle) and reproduces the target accurately by iteration~10 (right).}
  \label{fig:success_example}
\end{figure}

A representative success is illustrated in \cref{fig:success_example}. Here favorable initialization yields substantial initial overlap with the target, so that only moderate adjustments of endpoints and radius are required. The pill aligns rapidly and by iteration~10 the target is reproduced with high accuracy.

\begin{figure}[H]
  \centering
  \includegraphics[width=0.32\textwidth]{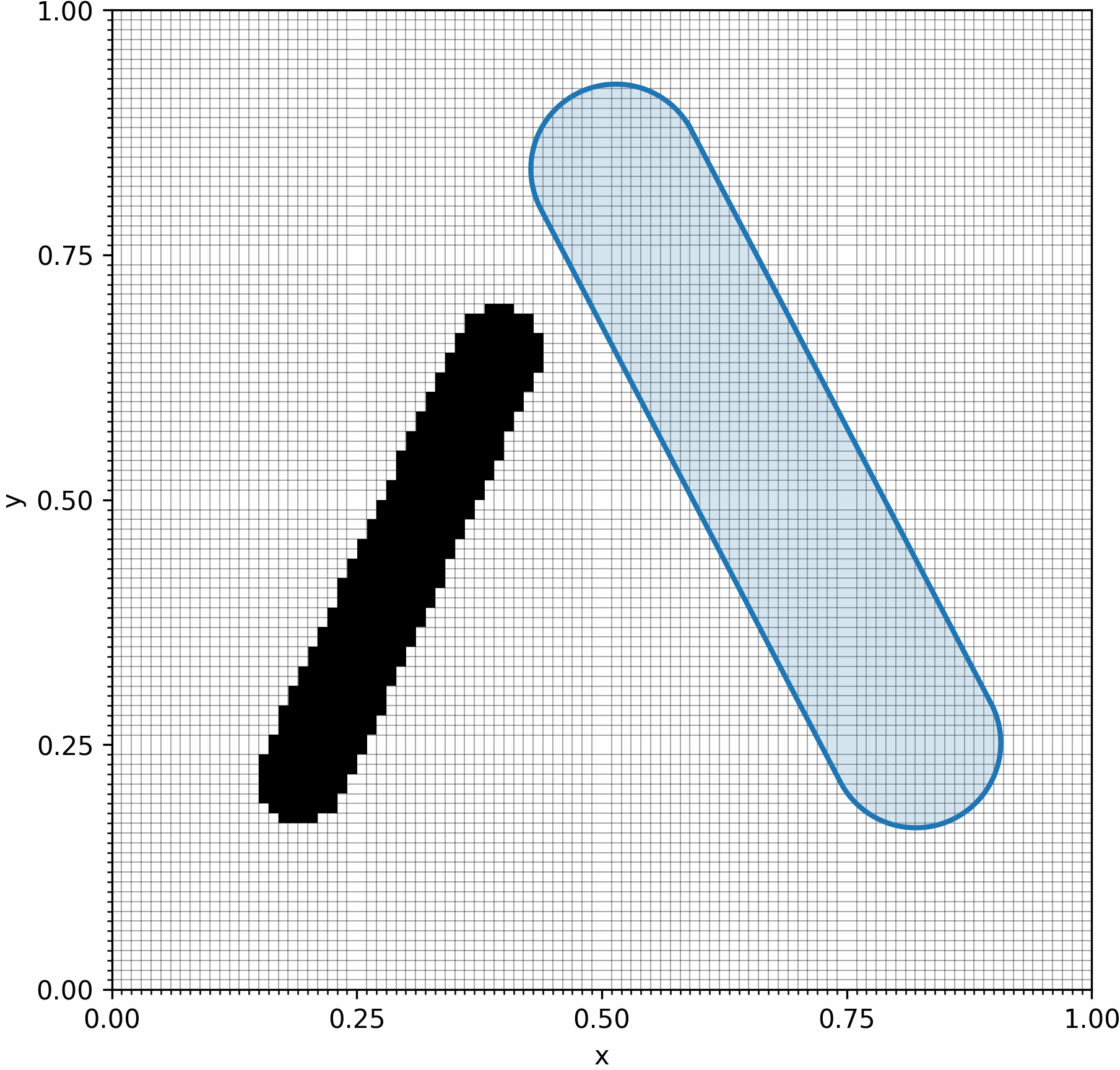}\hfill
  \includegraphics[width=0.32\textwidth]{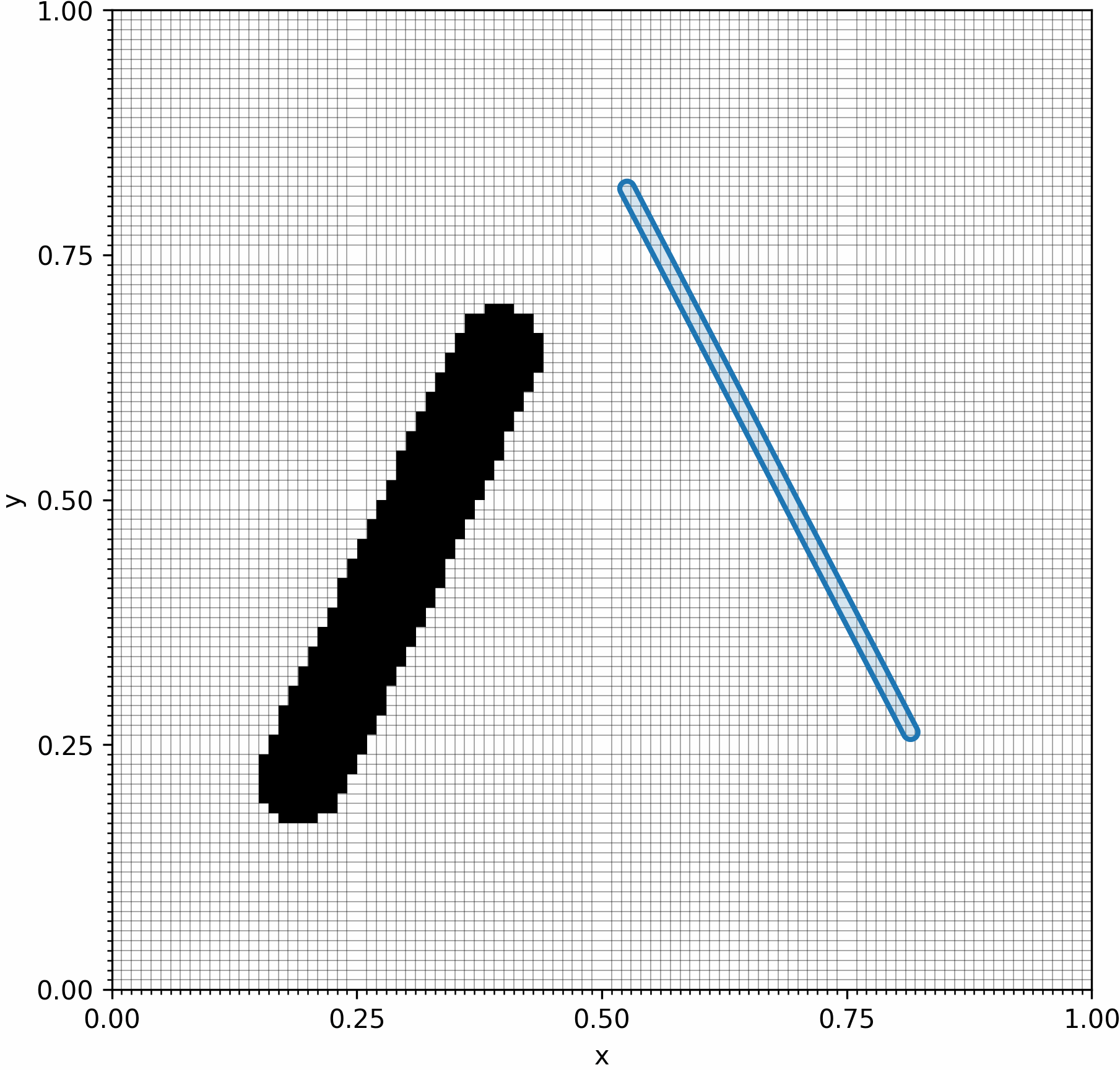}\hfill
  \includegraphics[width=0.32\textwidth]{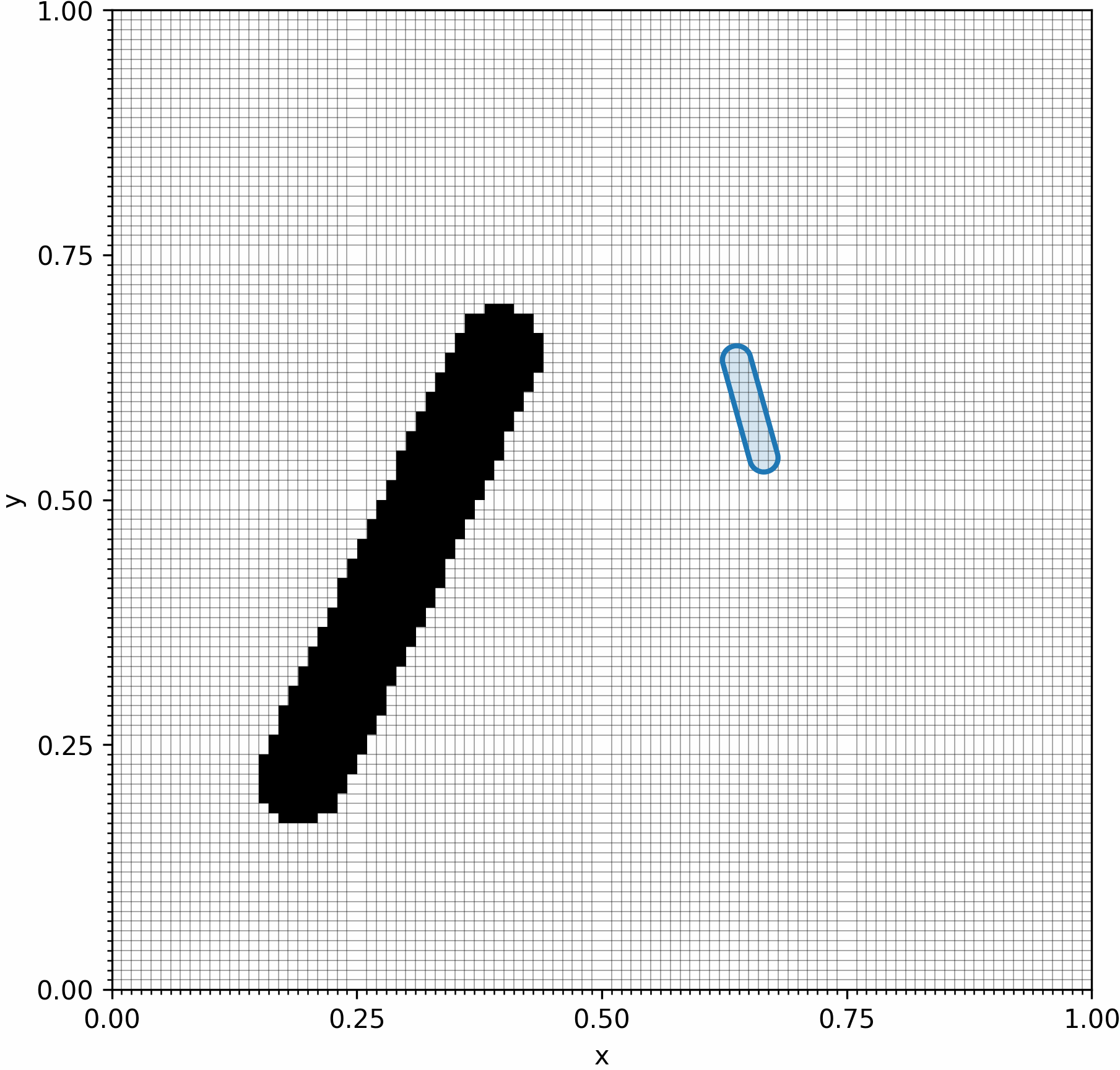}\\[0.5ex]
  \includegraphics[width=0.32\textwidth]{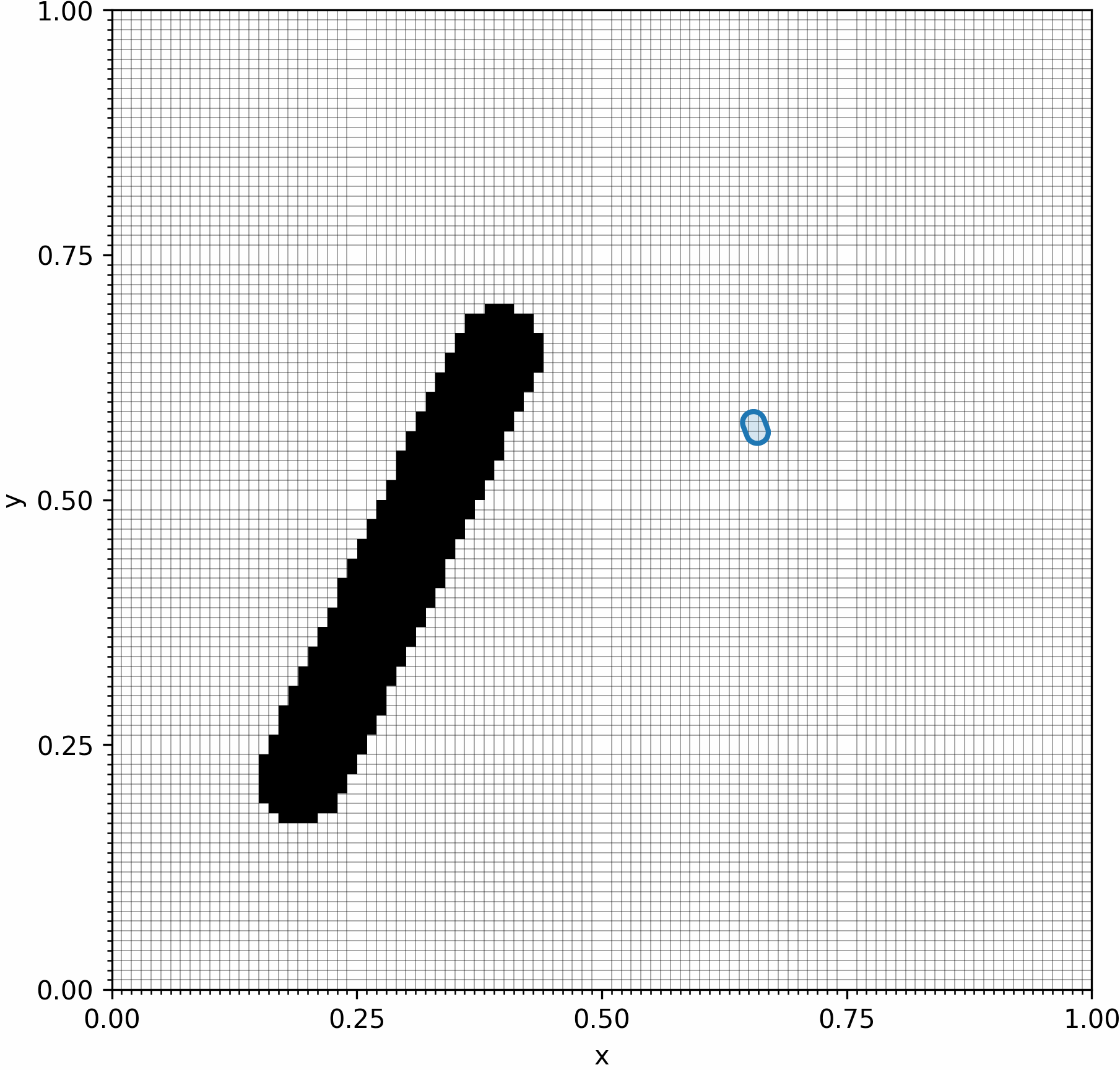}\hfill
  \includegraphics[width=0.32\textwidth]{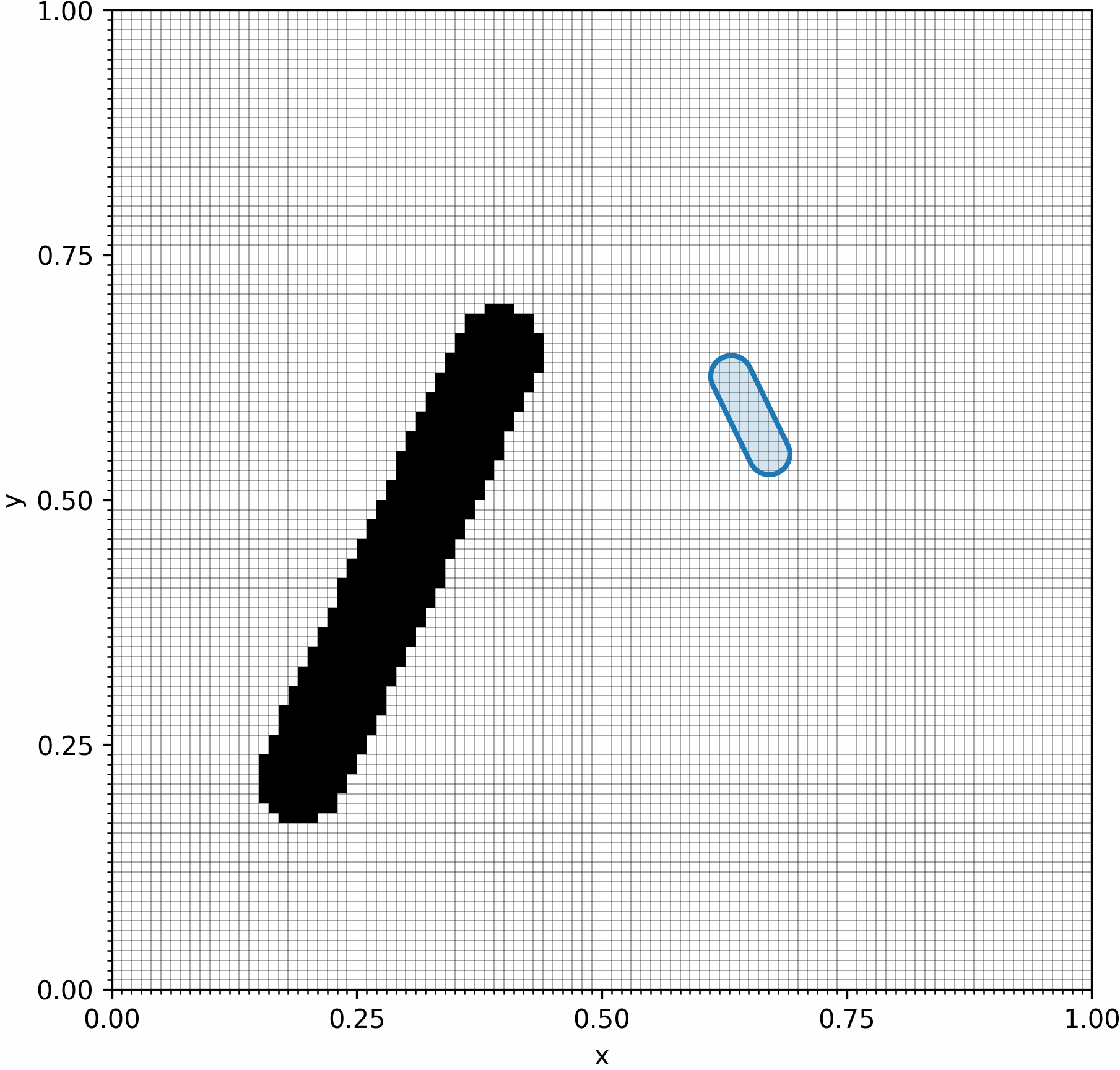}\hfill
  \includegraphics[width=0.32\textwidth]{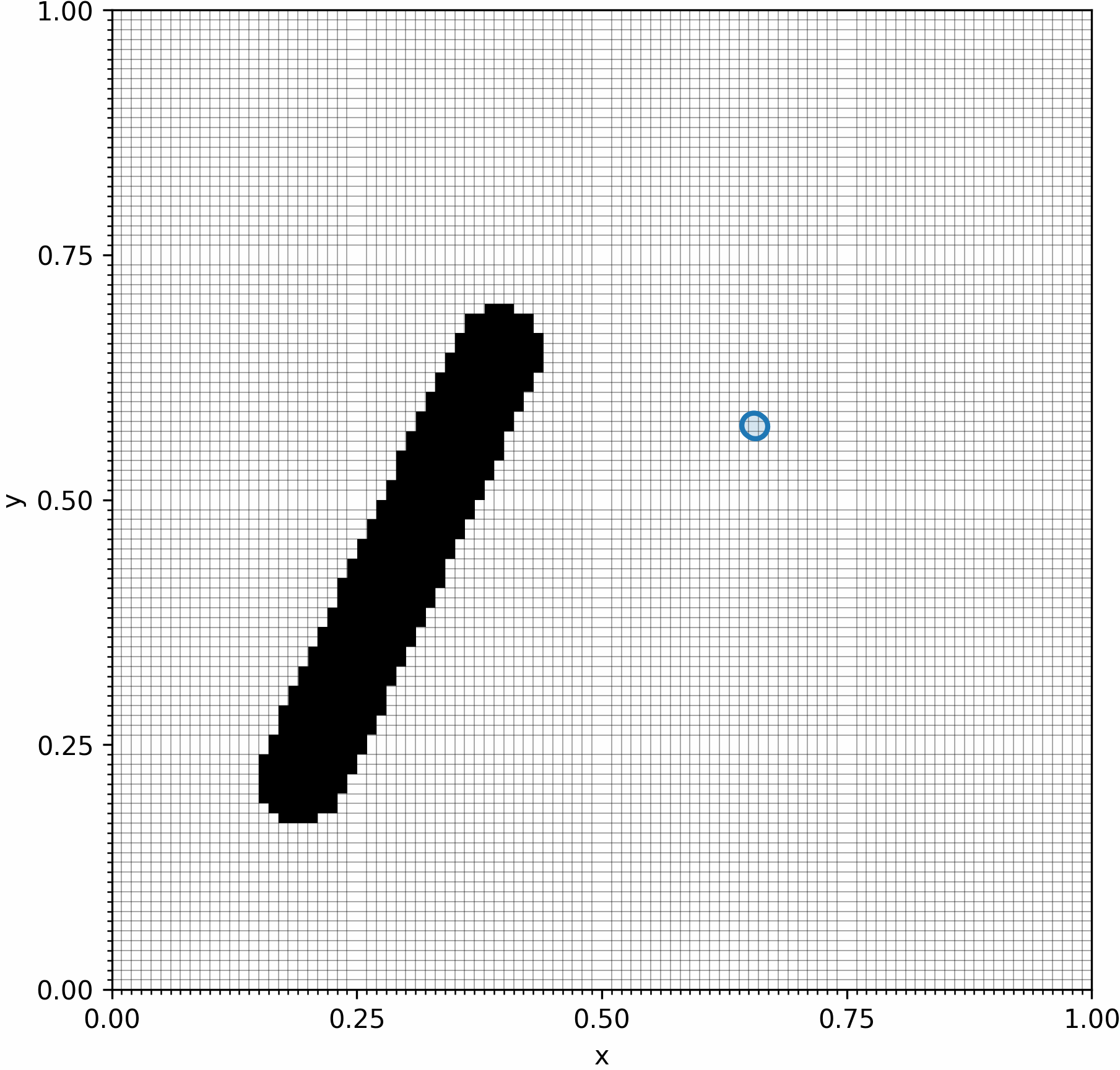}
  \caption{Representative stagnation trajectory from a run without initial overlap. The pill contracts until nearly collapsing (top row), then re-elongates (bottom row, left) and subsequently contracts again (bottom row, right). This oscillatory behavior corresponds to the jumps visible in the objective curves and prevents convergence.}
  \label{fig:failure_example}
\end{figure}

In contrast, several runs without initial overlap display a qualitatively different behavior. Instead of establishing contact with the target, the pill undergoes repeated cycles of contraction and re-elongation. During the first iterations the radius shrinks and the segment endpoints move closer together until the pill nearly collapses. Subsequently the radius grows again and the segment re-elongates, only to contract once more in later iterations. This alternating pattern is reflected in the objective traces as sharp spikes and plateaus and prevents the optimization from reaching the tolerance criterion. Representative snapshots of such a trajectory are shown in \cref{fig:failure_example}. The figure illustrates how the pill oscillates between nearly point-like and elongated states without generating overlap with the target. A more detailed examination of this mechanism in the absence of target signal is given in the following section.

\subsection{Feature behavior without target intersection}
\label{subsec:feature_no_overlap}

The stagnation observed in several runs of \cref{sec:single_feature_optimization} is explained by the intrinsic behavior of a single feature in the absence of target overlap. This regime is represented by the extreme case $\rho^\star\equiv 0$ on the evaluation domain, for which the tracking objective reduces to
\begin{equation}
  \Phi_{\mathrm{track}}(\vecsym \zeta)
  = \sum_{(a,b)\in\Omega_h} \bigl(0 - \bar\rho_{a,b}(\vecsym{\zeta})\bigr)^{2}
  = \sum_{(a,b)\in\Omega_h} \bar\rho_{a,b}(\vecsym{\zeta})^{2}
  \label{eq:no_overlap_J}
\end{equation}

Any nonzero density increases $\Phi$, so that descent directions arise solely from the feature’s own sensitivities.

A tilted segment with endpoints $P=(0.3,0.4)$ and $Q=(0.6,0.7)$ is used to visualize these directions. The gradients with respect to $P$ concentrate on the transition band near the cap around $P$ and along the adjacent segment. The dominant sign is negative for both $\partial \Phi_{\mathrm{track}}(\vecsym{\zeta})/ \partial p_x$ and $\partial \Phi_{\mathrm{track}}(\vecsym{\zeta})/ \partial p_y$, indicating that increments of $p_x$ or $p_y$ lower the objective. Small positive lobes occur on the opposite side of the band, but they largely cancel in the net displacement. The effective motion is governed by the cap and pushes $P$ toward $Q$.
\begin{figure}[H]
  \centering
  \includegraphics[width=0.45\textwidth]{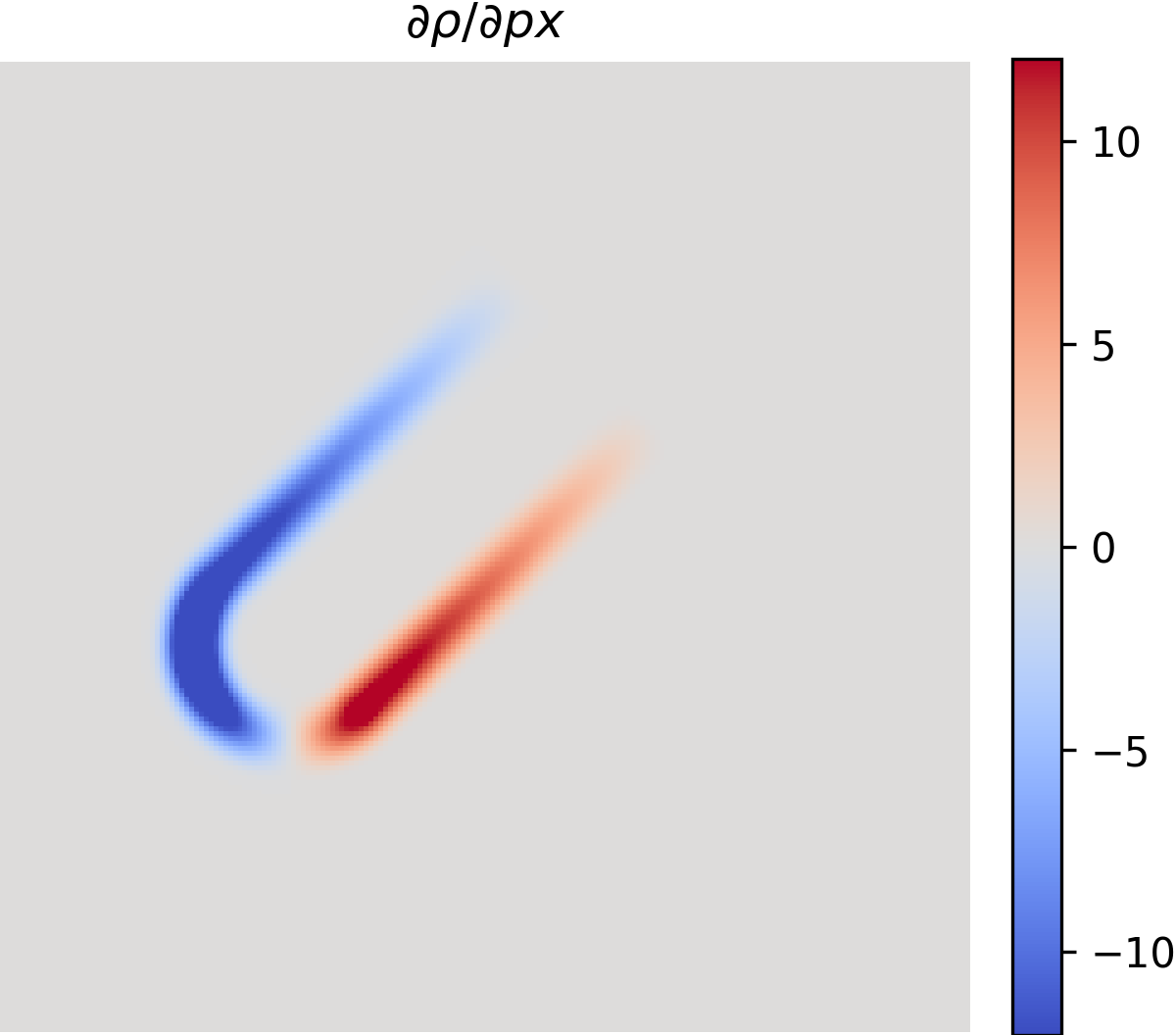}\hfill
  \includegraphics[width=0.45\textwidth]{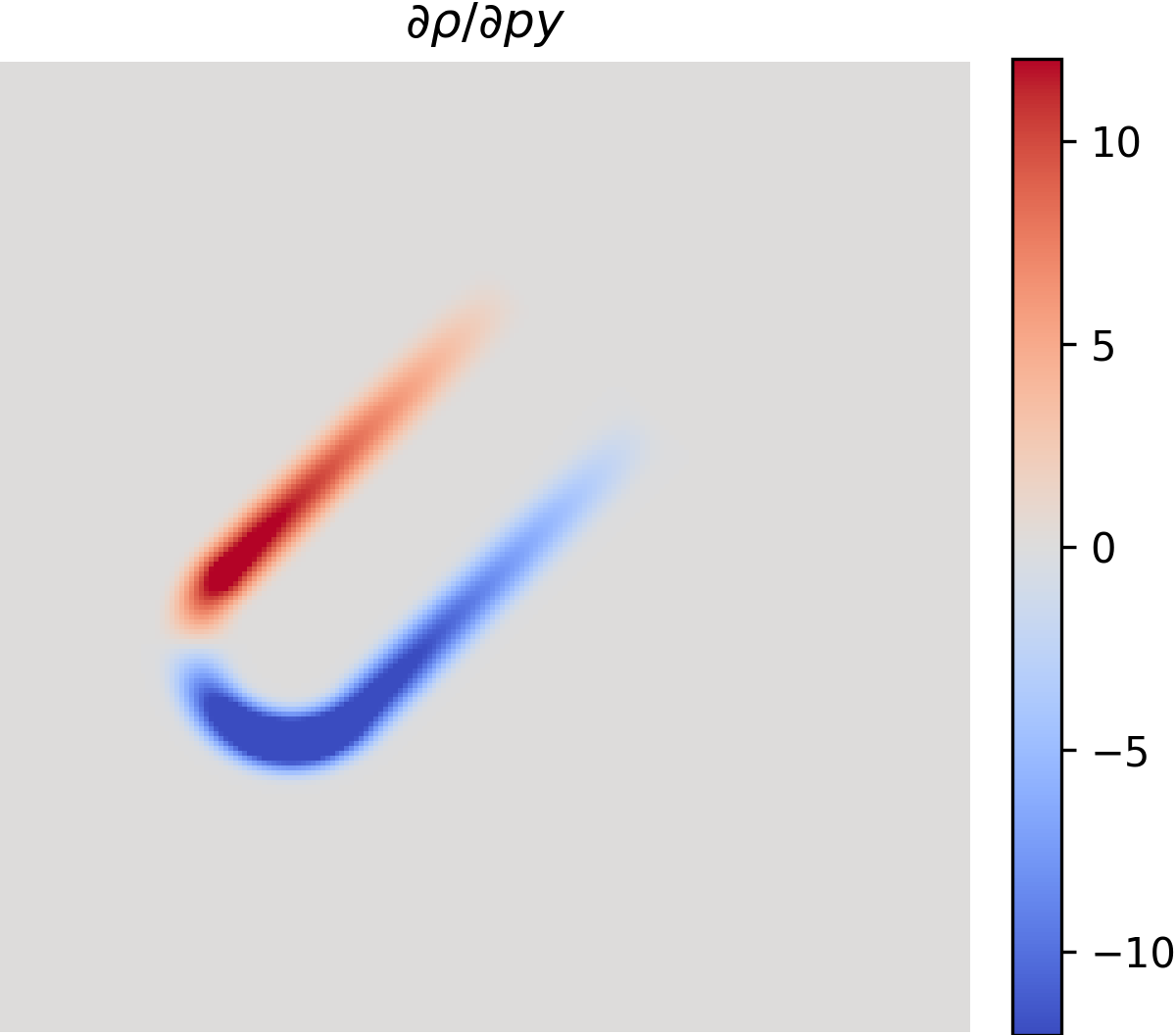}
  \caption{Endpoint $P$: sensitivities w.r.t.\ $p_x$ (left) and $p_y$ (right).}
  \label{fig:grad_P_tilted}
\end{figure}

The complementary structure is observed for $Q$: around its cap, gradients are predominantly positive, which, in a minimization setting, induces decreases in both $q_x$ and $q_y$. Negative lobes along the adjacent segment cancel with positive values across the band and contribute little to the net motion. The result is a pull of $Q$ toward $P$, i.e., consistent shortening of the segment.
\begin{figure}[H]
  \centering
  \includegraphics[width=0.45\textwidth]{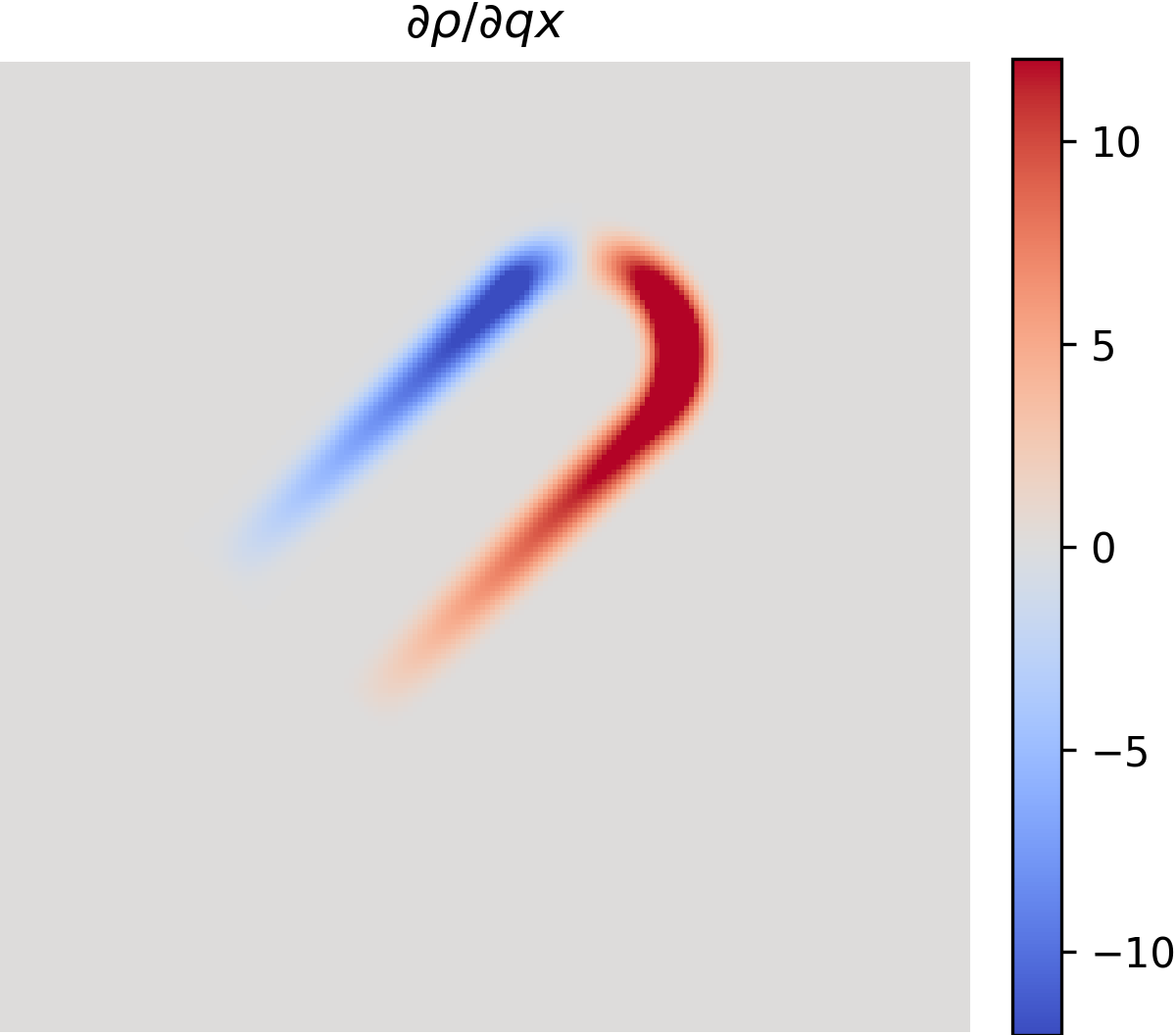}\hfill
  \includegraphics[width=0.45\textwidth]{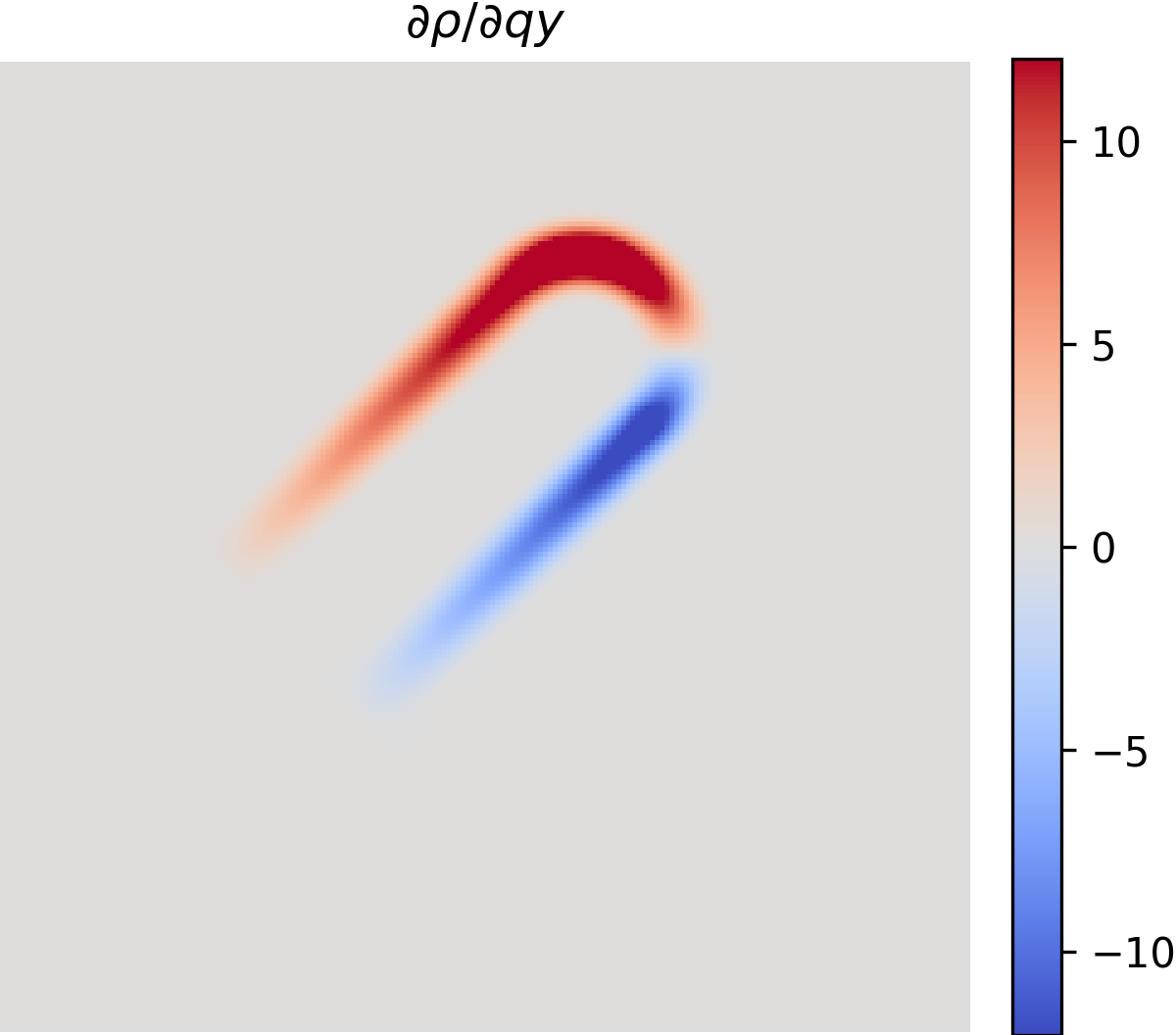}
  \caption{Endpoint $Q$: sensitivities w.r.t.\ $q_x$ (left) and $q_y$ (right).}
  \label{fig:grad_Q_tilted}
\end{figure}

The derivative with respect to the radius is strictly positive on the transition
band, hence $r$ is driven to its lower bound:
\begin{figure}[H]
  \centering
  \includegraphics[width=0.45\textwidth]{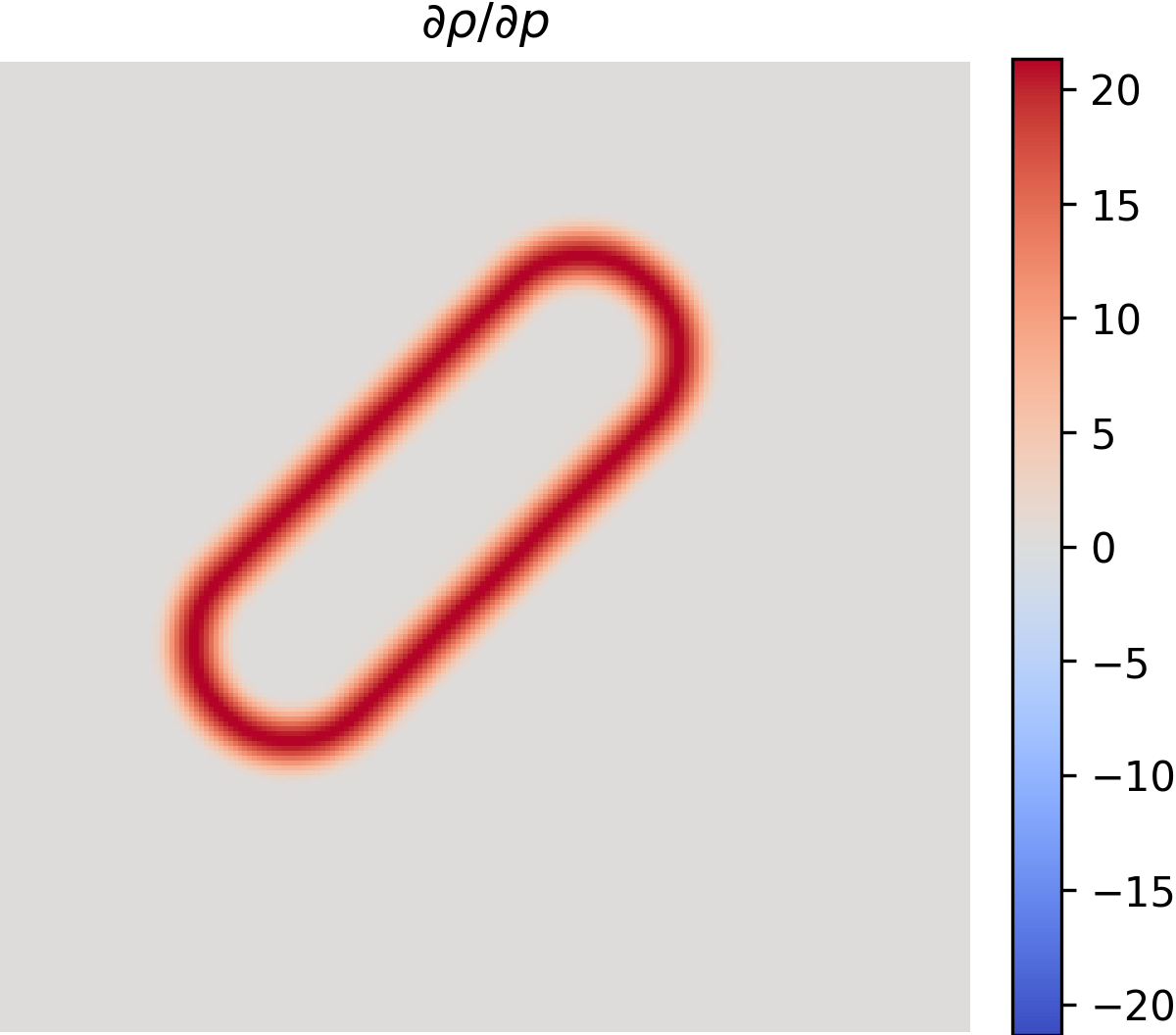}
  \caption{Sensitivity w.r.t.\ radius $r$. Positive values enforce shrinkage.}
  \label{fig:grad_r_tilted}
\end{figure}

\begin{figure}[H]
  \centering
  \includegraphics[width=0.32\textwidth]{images/hess_px_px.jpg}\hfill
  \includegraphics[width=0.32\textwidth]{images/hess_px_py.jpg}\hfill
  \includegraphics[width=0.32\textwidth]{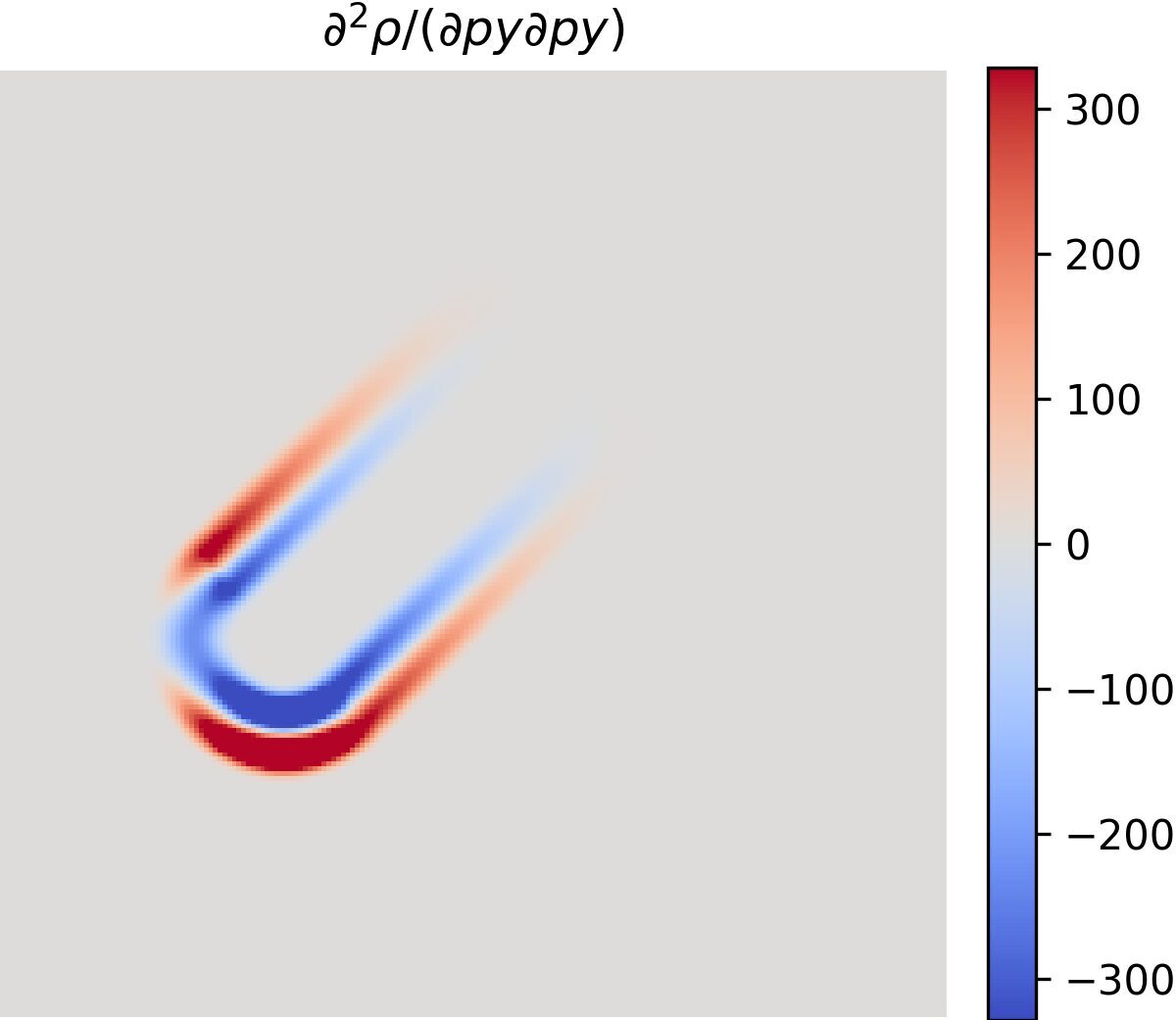}
  \caption{Representative Hessian entries for endpoint $P$. Diagonal terms are
  mostly positive; the mixed term alternates in sign.}
  \label{fig:hess_P}
\end{figure}

The second-order structure reinforces this contraction. Around $P$, diagonal
Hessian entries are predominantly positive, while the mixed entry
$\partial^2 \Phi/(\partial p_x\,\partial p_y)$ alternates in sign, producing
sign-changing lobes across the band. Curvature is sharply localized; exterior
parts of the band contribute positive curvature that dominates local behavior.
Representative plots are shown in \cref{fig:hess_P}.

\begin{figure}[H]
  \centering
  \includegraphics[width=0.32\textwidth]{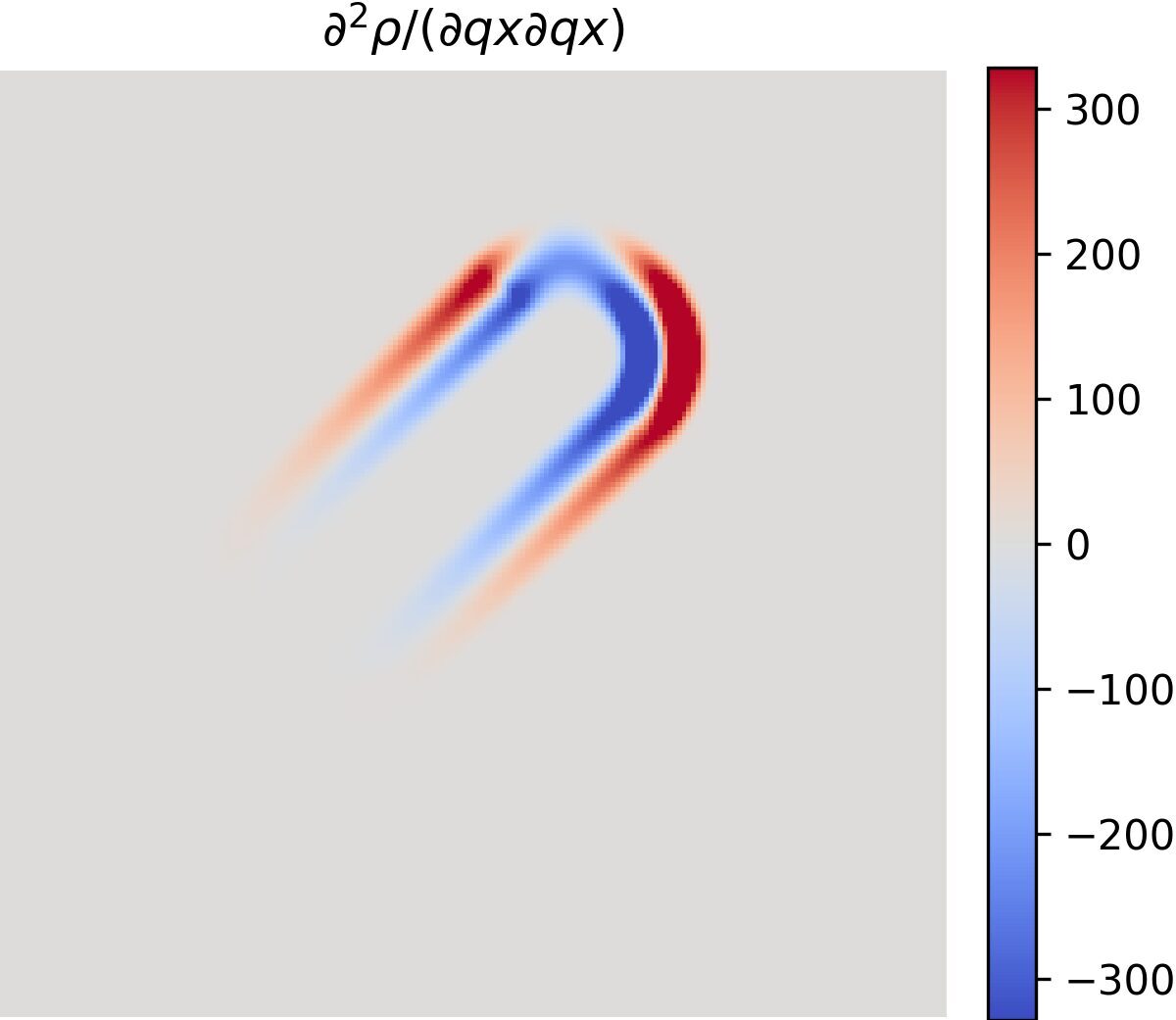}\hfill
  \includegraphics[width=0.32\textwidth]{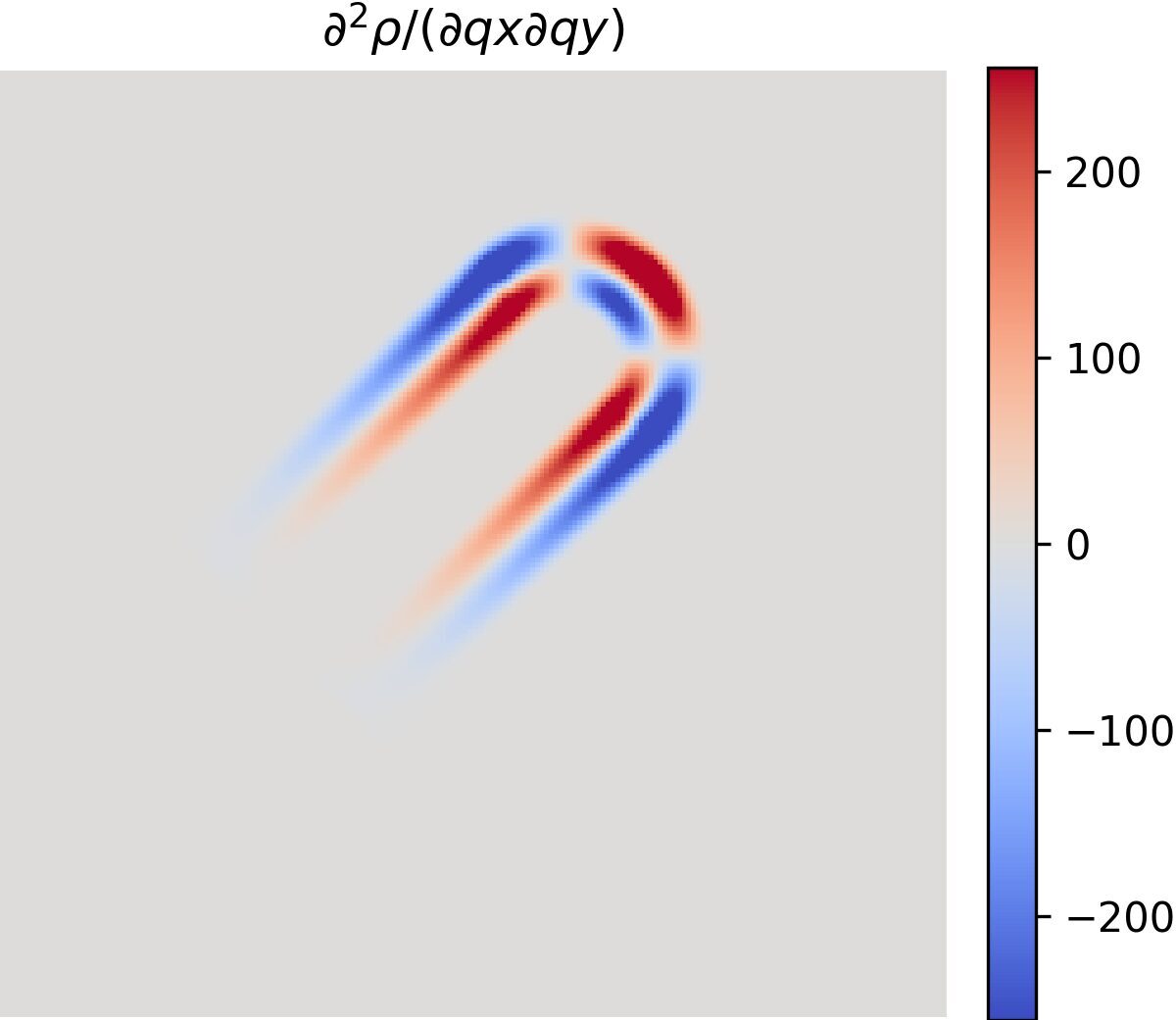}\hfill
  \includegraphics[width=0.32\textwidth]{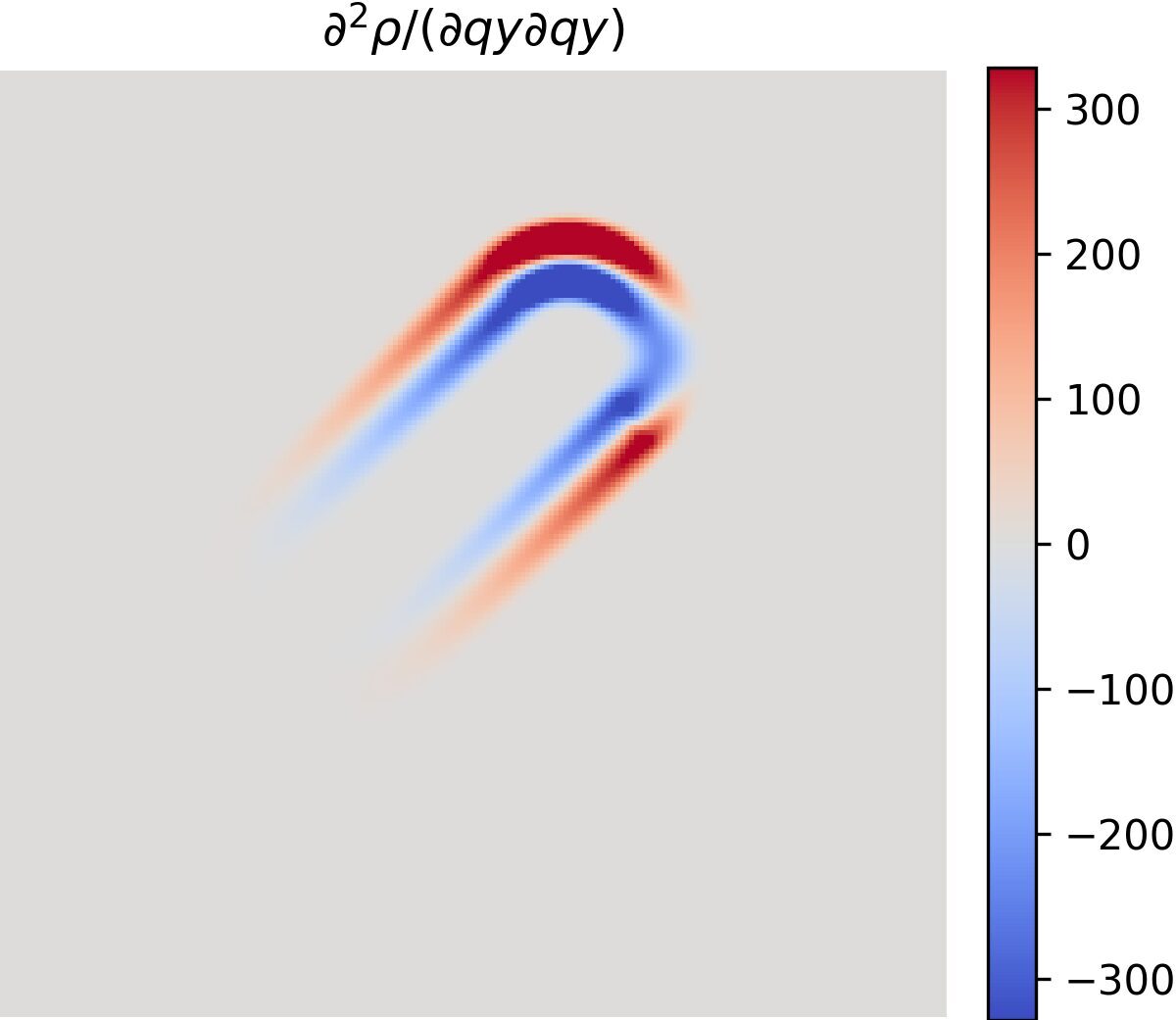}
  \caption{Representative Hessian entries for endpoint $Q$. Diagonal curvature
  is positive; the mixed entry alternates in sign.}
  \label{fig:hess_Q}
\end{figure}

\begin{figure}[H]
  \centering
  \includegraphics[width=0.85\textwidth]{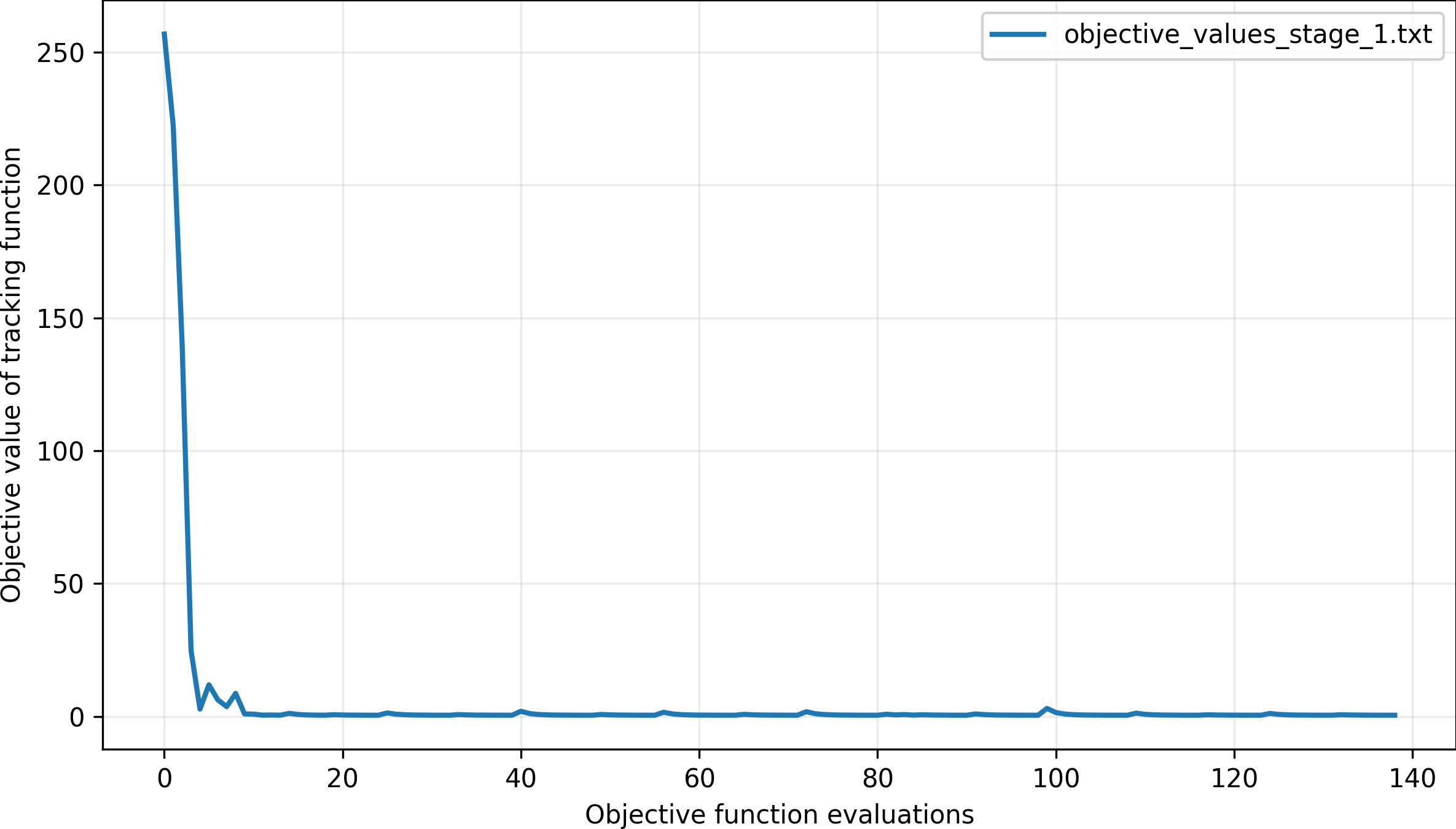}
  \caption{Objective history in the no-target case. Overall decay is interrupted by spikes caused by temporary re-elongation.}
  \label{fig:degenerate_objective}
\end{figure}

An analogous pattern appears at $Q$; see \cref{fig:hess_Q}. Importantly, the
analytic expressions of the signed-distance derivatives contain factors with
negative powers of the segment length $\ell=\|Q-P\|_2$ (e.g., $\ell^{-3}$ and
$\ell^{-5}$, cf.\ \cref{sec:feature_derivatives}). As $\ell\to 0$, these terms
grow rapidly and make the parameterization increasingly stiff. Together with the
first-order tendencies above, this explains the observed alternation between
contraction and re-elongation in the no-overlap runs and the associated spikes in
the objective traces.
To examine these effects in practice, a pill initialized at $\vecsym{\zeta}=(0.3,0.4,0.6,0.7,0.1)$ is optimized without target density. The tracking functional $\Phi_{\mathrm{track}}(\vecsym{\zeta})$ decreases overall, but the objective history in \cref{fig:degenerate_objective} exhibits intermittent upward spikes rather than smooth monotone decay. Each spike corresponds to a brief re-elongation of the segment. Between iterations~20 and~21, for instance, the objective increases by nearly ninety percent as the segment length jumps from a near-degenerate state to a slightly elongated one before contraction resumes.

\begin{figure}[H]
  \centering
  \includegraphics[width=0.3\textwidth]{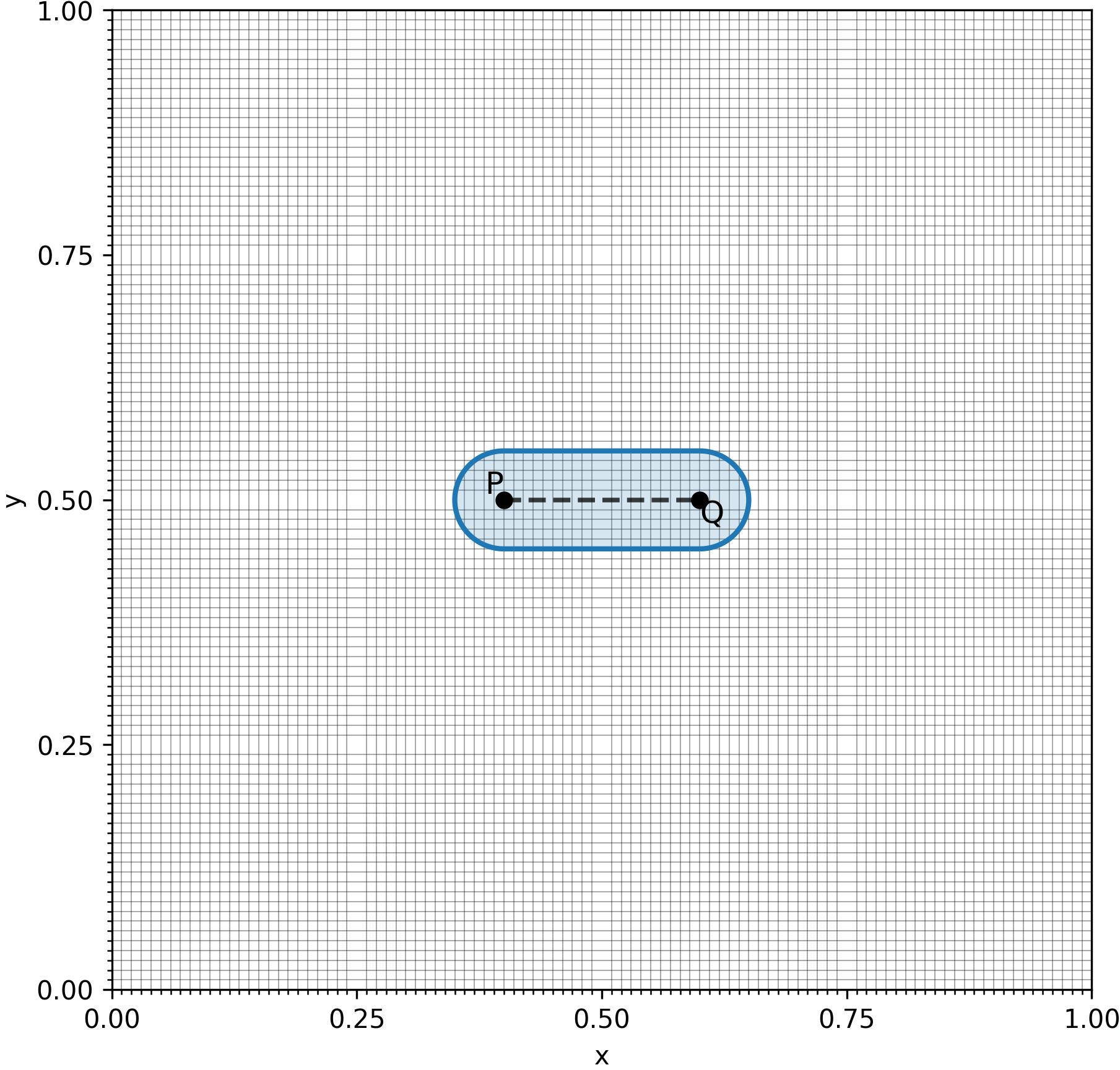}
  \includegraphics[width=0.3\textwidth]{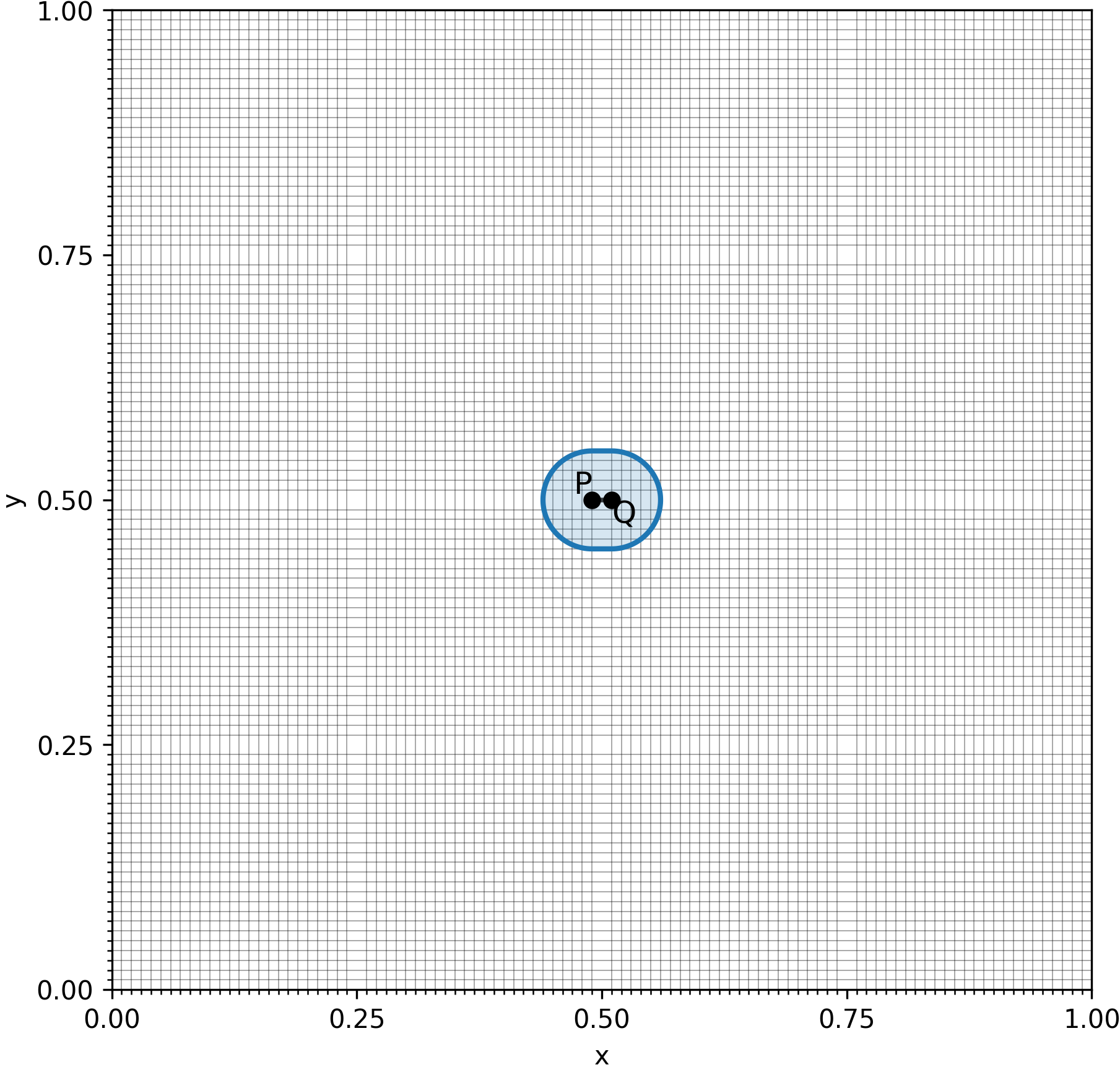}
  \includegraphics[width=0.3\textwidth]{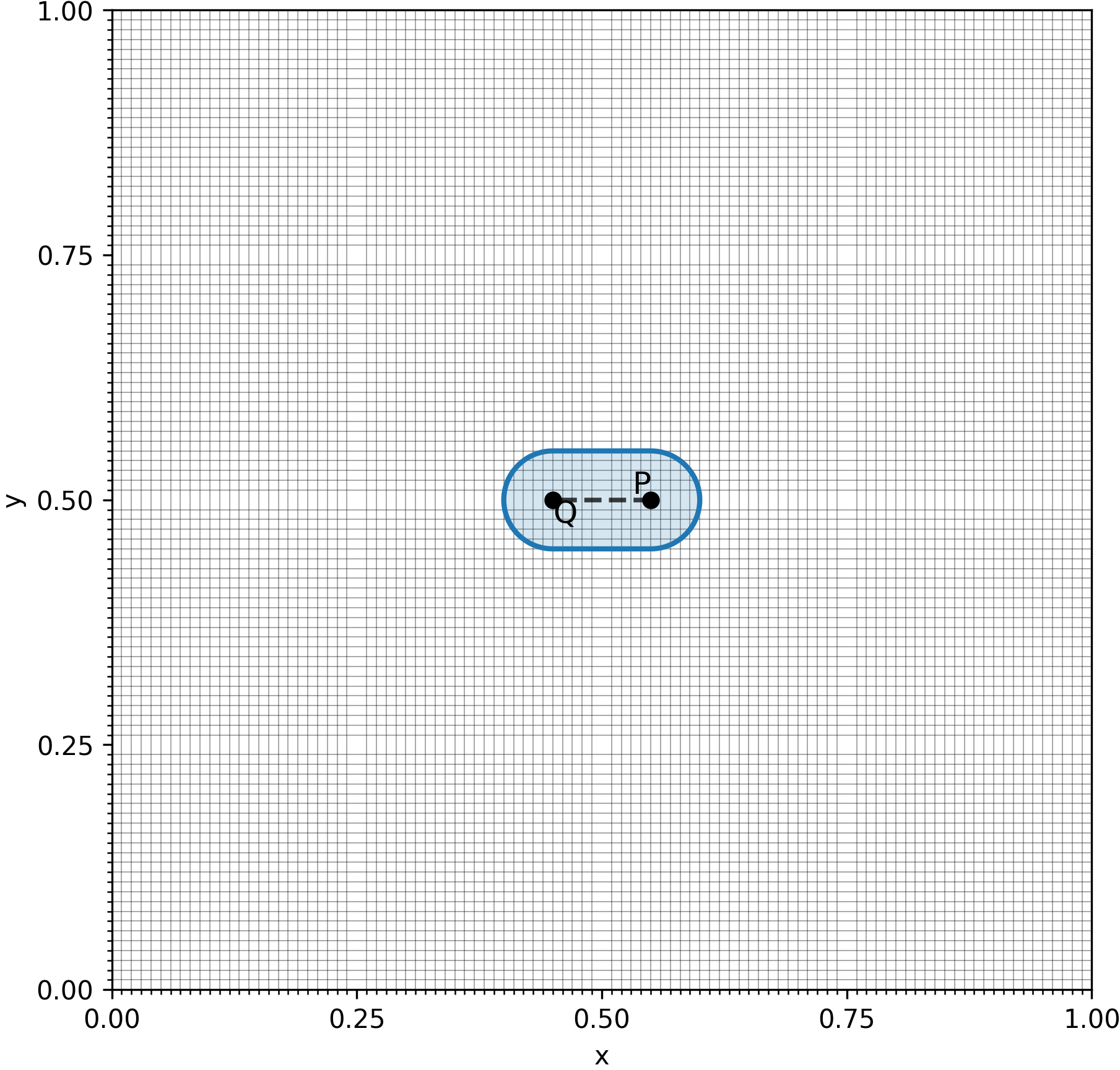}
  \caption{Geometry snapshots in the no-target case. Top: elongated initial state. Middle: near-degenerate disk. Bottom: temporary re-elongation after endpoint swap.}
  \label{fig:degenerate_geometry}
\end{figure}

The geometry snapshots in \cref{fig:degenerate_geometry} confirm this cycle. Initially, the pill contracts toward a disk-like state as $P$ and $Q$ approach each other. Near degeneracy, however, the rapidly growing Hessian terms destabilize the solver, causing the endpoints to overshoot and swap positions. This produces a short re-elongation, visible as a spike in the objective trace. The alternation between collapse and endpoint swapping repeats until termination, explaining the irregular behavior and confirming the derivative-based analysis.

\subsubsection*{Stabilization by minimum segment length}

Instabilities of the unconstrained setting can be mitigated by enforcing a lower bound on the segment length. A minimum length $\ell_{\min}>0$ prevents complete collapse and thereby avoids the alternating shrink--re-elongation cycles observed before. The pressure of $P$ and $Q$ toward each other remains present, but once $\ell_{\min}$ is reached further contraction is blocked. As a result, the likelihood of endpoints swapping sides is reduced, since such swaps would require the segment to pass through degeneracy. 

\begin{figure}[H]
  \centering
  \includegraphics[width=0.85\textwidth]{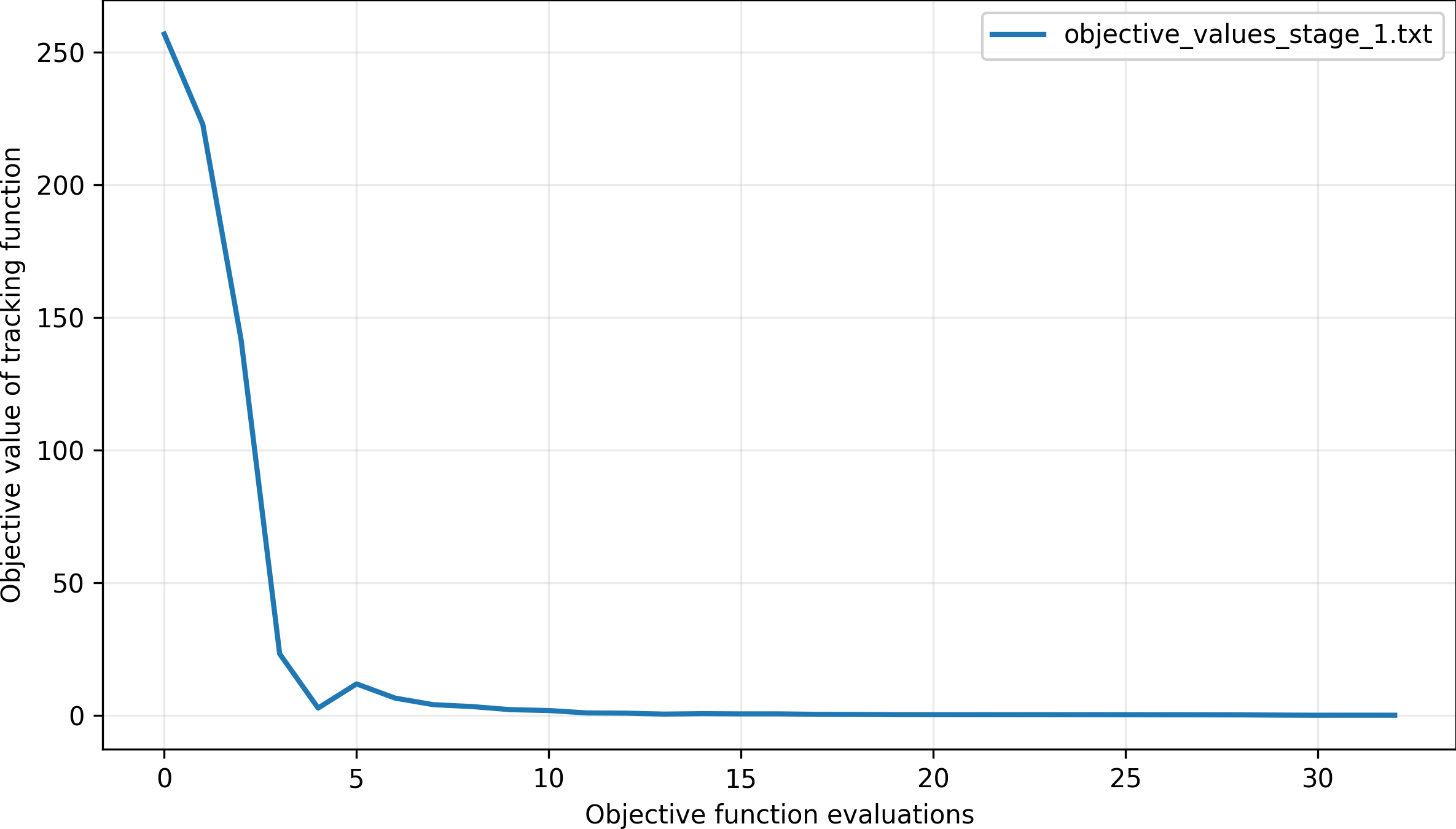}
  \caption{Effect of a minimum segment length $\ell_{\min}=0.005$. Convergence is smooth, with no spikes and far fewer function evaluations than in the unconstrained case.}
  \label{fig:seglen_constraint}
\end{figure}

The effect is illustrated in \cref{fig:seglen_constraint} for $\ell_{\min}=0.005$. The objective decays monotonically without spikes and the solver accepts nearly every Newton step. While in the unconstrained case almost 140 function evaluations were required for 30 iterations, the constrained run needs just above 30 function evaluations. This improvement reflects the removal of ill-conditioning: without the constraint the Hessian diverges as $\ell\to0$, forcing repeated backtracking, whereas with bounded curvature the iterations track smoothly. Larger values of $\ell_{\min}$ increase stability further by making endpoint swaps increasingly unlikely, but they also limit the smallest admissible pill size. The constraint thus represents a trade-off: small values suffice to remove numerical pathologies, while larger choices improve convergence speed but reduce geometric resolution.

\subsubsection*{Asymmetric extension}
\label{subsubsec:extension}

The analysis above has shown that even with a minimum length constraint, convergence may fail when pills are initialized far from the target, since localized gradients provide no effective translation or orientation cues. To mitigate this, the asymmetric transition introduced in \cref{sec:asym_scaled} is employed. By enlarging only the exterior flank while keeping the interior half--width unchanged, nonzero sensitivities are generated outside the pill, so that long–range cues are provided that can initiate motion toward the target. The price is a reduction of regularity from $C^2$ to $C^0$ at the medial axis, while compact support is retained. Small extensions preserve near–symmetric behavior and favor precise alignment, whereas larger values broaden the influence region and accelerate the establishment of overlap at the cost of diffuse density in void. The following tests illustrate both the potential and the limitations of this mechanism. Since the asymmetric profile increases the pressure on the endpoints to move toward each other, a larger minimum segment length $\ell_{\min}$ is recommended than in the symmetric case. This reduces the likelihood of endpoint swapping and suppresses re-elongation, ensuring smoother convergence when extension is active.

\begin{figure}[H]
  \centering
  \includegraphics[width=0.48\textwidth]{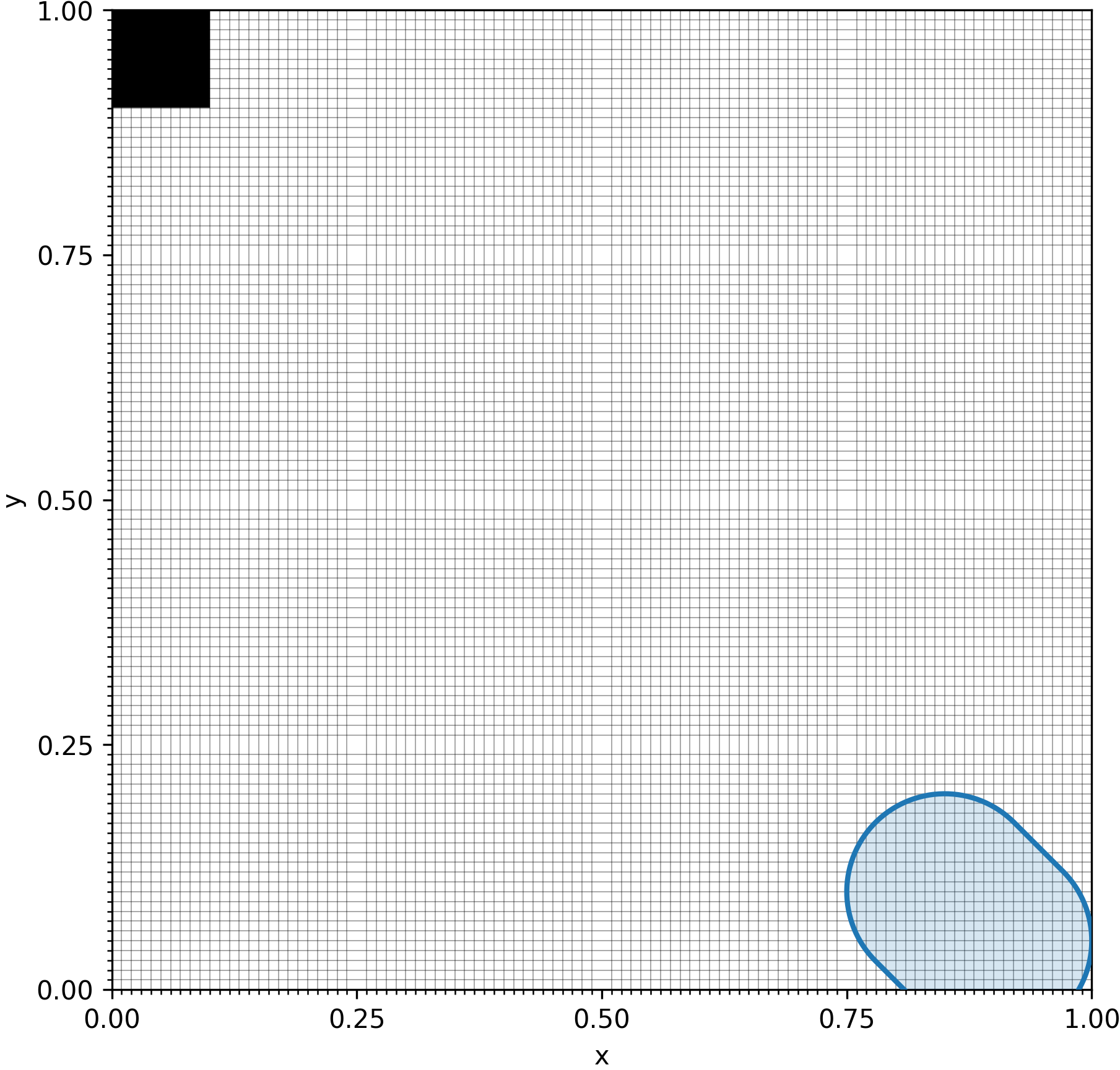}\hfill
  \includegraphics[width=0.48\textwidth]{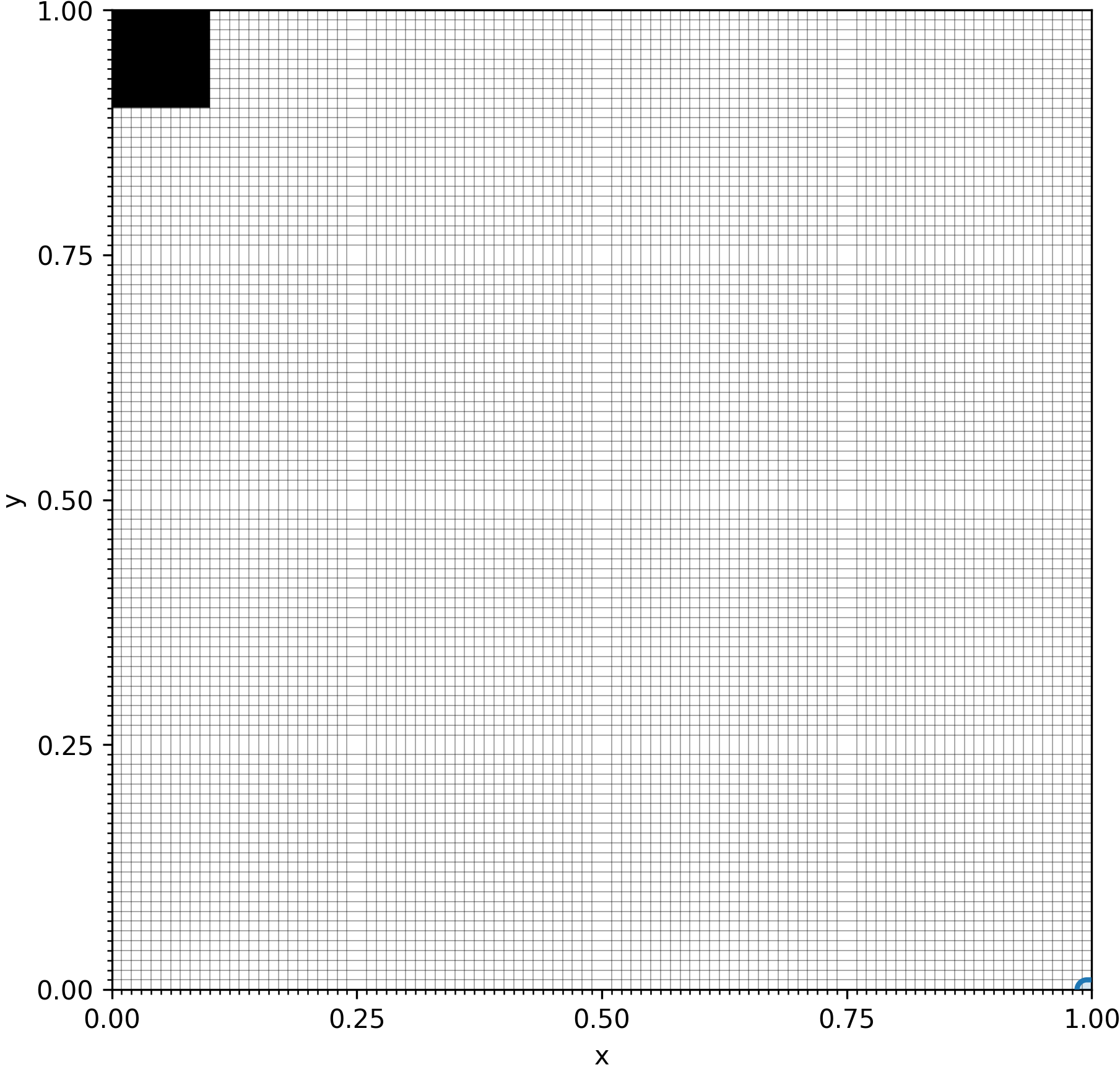}
  \caption{Corner--to--corner test with extension. Left: initialization in the lower--right and target in the upper--left. Right: final state after optimization, where the pill contracts and drifts into the corner instead of approaching the target.}
  \label{fig:corner_extreme}
\end{figure}

\begingroup\sloppy
In practice, the extension parameter is fixed at initialization: large values favor recovery from unfavorable starts, whereas small values improve accuracy in well–aligned cases.
\endgroup
An extreme stress test illustrates the limitation. The pill is placed in the lower--right corner with 
\[
(p_x,p_y,\allowbreak q_x,q_y,\allowbreak r)=(0.85,0.10,\allowbreak 0.90,0.05,\allowbreak 0.10),
\]
while the target is a thin density block in the upper--left. The symmetric transition offers no usable gradient and even with $\mathrm{ext}=1.4$, sufficient to span the diagonal of the unit square, the pill fails to approach the target. Instead, optimization proceeds as in the no--overlap case: segment length and radius decrease monotonically, local shrinkage dominates the weak far--field pull and boundary bias accelerates collapse by truncating exterior density. The final state is a nearly collapsed pill pinned in the corner, with no progress toward the target despite the nominal extension (\cref{fig:corner_extreme}).

Two coupled mechanisms explain this failure. First, the objective is more effectively reduced by shrinking the pill than by following the weak long–range cues introduced by the extension, which decay rapidly with distance. Second, boundary bias favors collapse: once the pill protrudes beyond the domain, truncated mass reduces the objective further. Large extensions make this protrusion unavoidable and even moderate values increase the risk compared to the symmetric case.

\subsubsection*{Remedies for boundary bias}
\label{subsubsec:remedies}

A first remedy is to enlarge the evaluation domain while keeping the design variables constrained to \([0,1]^2\). When quadrature is carried out on an extended window \([-\Delta,1+\Delta]^2\) with \(\Delta\) large enough to contain the entire diffuse interface, the artificial reduction of objective values near the boundary disappears. In the corner--to--corner test the run then follows the expected pattern: after initial contraction, the extended sensitivities accumulate, the pill rotates and translates across the domain and overlap with the target is established within a few tens of iterations. While effective, this approach increases computational cost in proportion to the padded area.

A complementary strategy is to modify the objective itself by removing void penalties instead of enlarging the box. Switching temporarily from the tracking functional to the reward formulation,
\[
  \Phi_{\text{reward}}(\vecsym{\zeta})
  = -\sum_{(a,b)} \bar\rho_{a,b}(\vecsym{\zeta})\,\rho^\star_{a,b},
\]
suppresses any incentive to shrink in void regions. Since gradients vanish wherever \(\rho^\star=0\), no boundary bias arises even with large extension and the optimizer receives clean directional cues from the target alone.

\begin{figure}[H]
  \centering
  \includegraphics[width=0.48\textwidth]{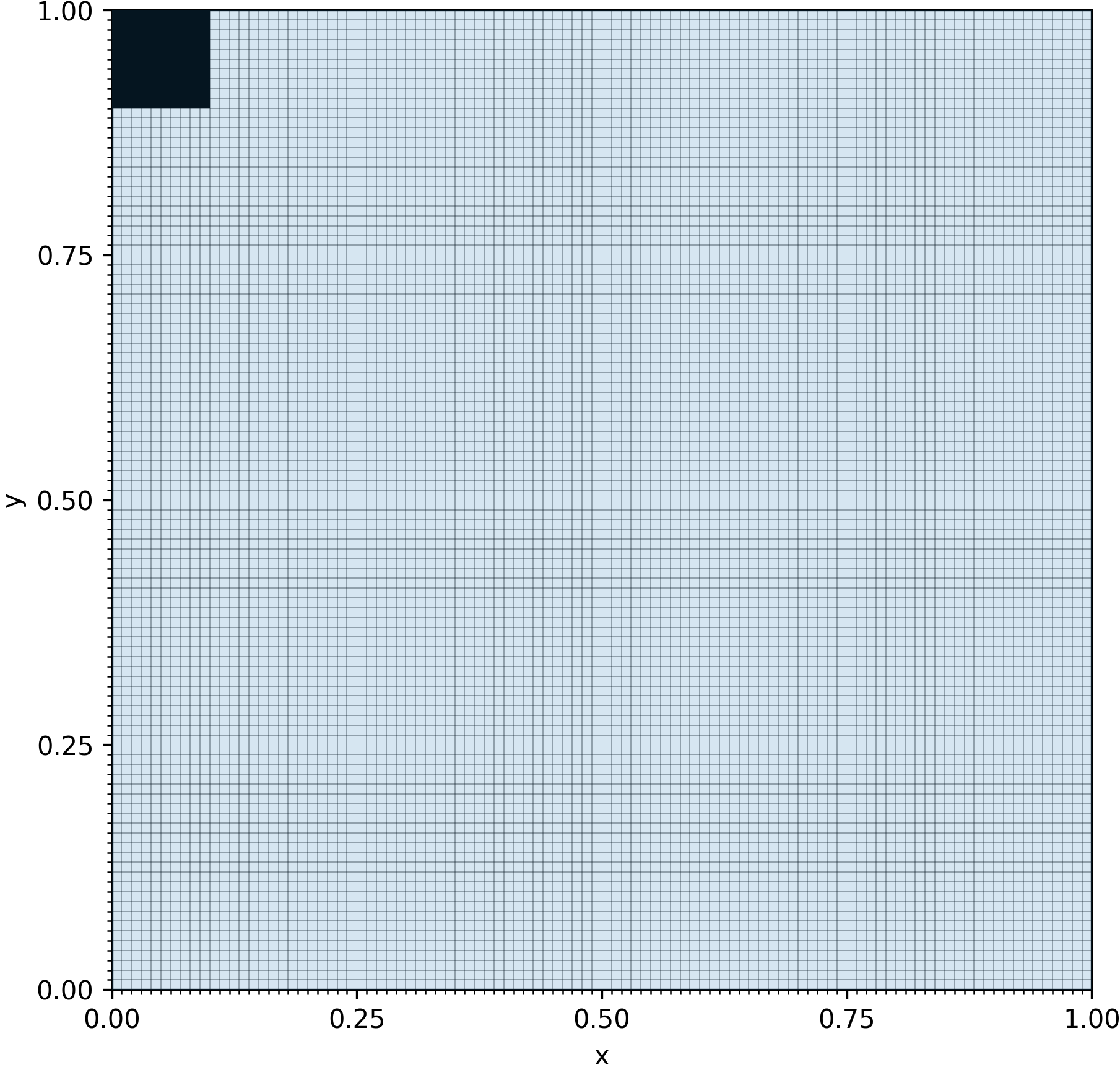}\hfill
  \includegraphics[width=0.48\textwidth]{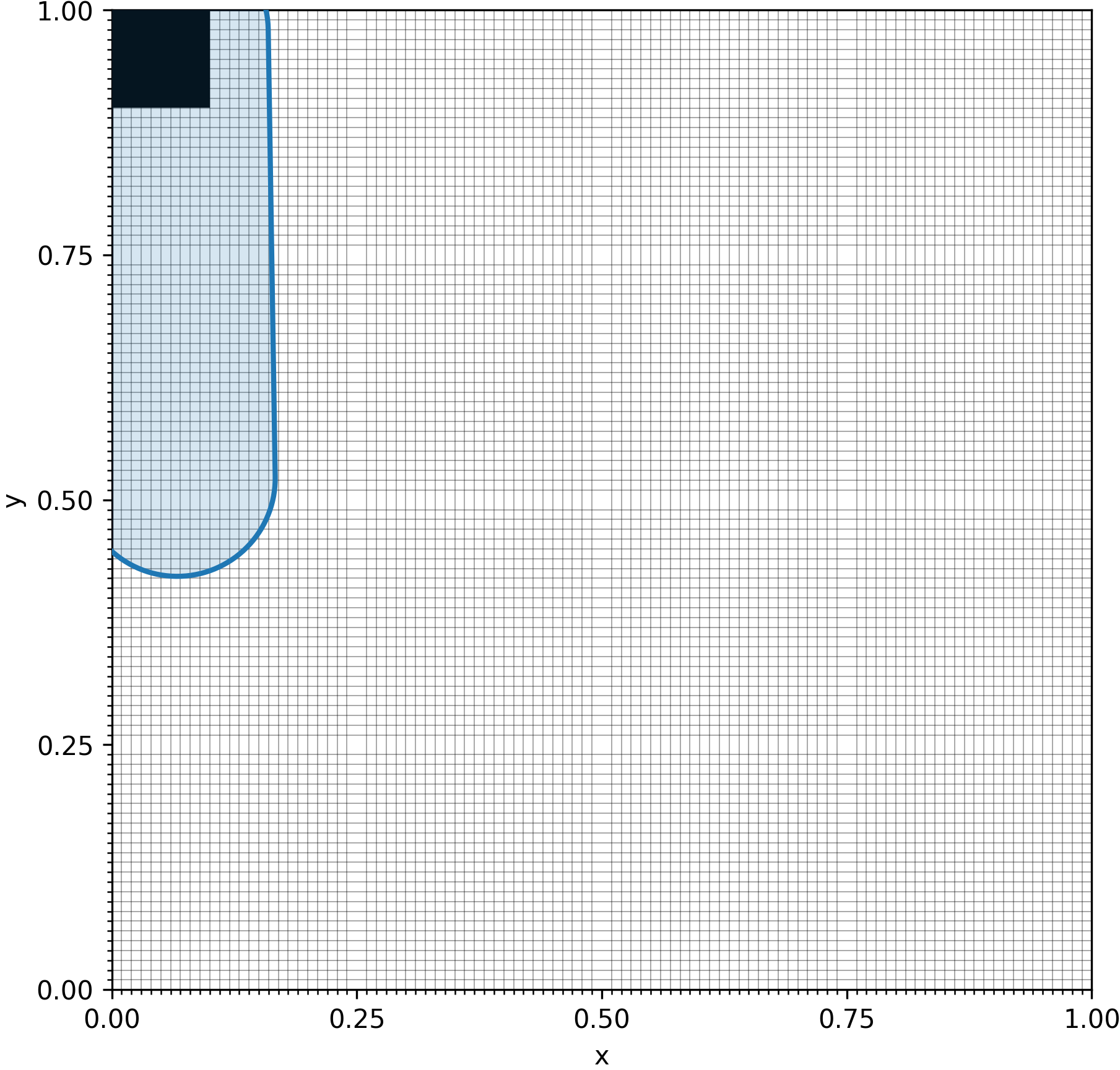}
  \caption{Reward-only exploration. Left: free \(r\) causes overexpansion. Right: fixing \(r\) yields a compact, target-aligned pose.}
  \label{fig:reward_phase}
\end{figure}

In the same corner stress test, the reward objective indeed drives the pill across the domain. Yet another side effect emerges: because density outside the target is not penalized, the optimizer tends to inflate the radius and stretch the segment, sometimes covering the entire domain. This trivial maximizer of $\Phi_{\mathrm{reward}}$ is not suitable as an initializer for subsequent tracking. The issue is avoided by fixing the radius during the reward phase, which restricts the optimizer to translation and orientation updates. The difference is illustrated in \cref{fig:reward_phase}: , the pill expands aggressively and blankets the domain, whereas with fixed radius the endpoints move decisively toward the target and settle into a meaningful pose.

\subsubsection*{Exploration-to-tracking workflow}
\label{subsubsec:explore_to_track}

The observations above motivate a two–stage procedure that is used throughout the experiments. In the first (\emph{exploration}) stage, long–range cues are created by activating the asymmetric extension from \cref{sec:asym_scaled} and temporarily switching to the reward objective $\Phi_{\mathrm{reward}}$ from \cref{sec:reward_objective_sec}. The radius is held fixed and a minimum segment length is enforced. Fixing $r$ eliminates the trivial incentive to inflate the pill in void regions, while the length bound suppresses jump–like endpoint swaps that otherwise arise under strong extension. Because $\Phi_{\mathrm{reward}}$ carries no penalty in void, boundary bias is neutralized even when the exterior flank is wide and the resulting gradients provide translation and orientation signals until nontrivial overlap with the target is attained.

Once overlap has been established and the endpoints $P,Q$ have migrated into an informative region, the second (\emph{tracking}) stage restores the symmetric transition (extension off), switches back to the least–squares objective $\Phi_{\mathrm{track}}$ from \cref{sec:tracking_objective_sec} and releases the radius variable. The minimum segment–length constraint remains active. In this stage the pose obtained during exploration is refined, the radius adapts to the prescribed thickness and residual misalignment is removed. In  practice, this yields stable convergence without the computational overhead of an enlarged integration window; it also avoids the corner–collapse pathology of direct tracking from unfavorable initializations.

\begin{figure}[H]
  \centering
  \includegraphics[width=.62\textwidth]{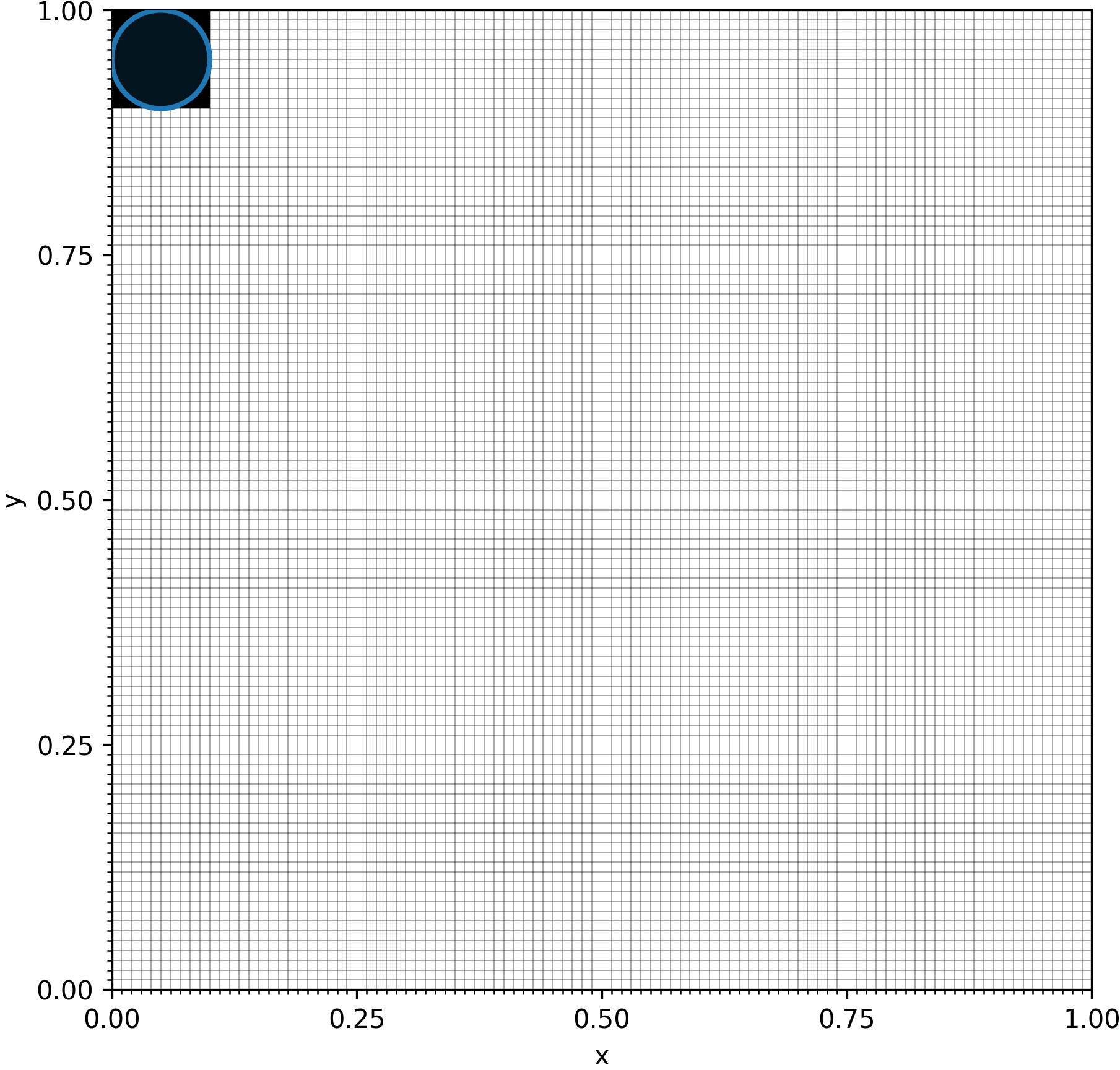}
  \caption{Corner stress test. Reference run with an enlarged evaluation window
  (for visualization) shows the intended cross–domain motion toward the target.
  Reward–only exploration with fixed radius achieves a comparable pose on the
  unit box and provides a robust initializer for tracking.}
  \label{fig:big_box_reference}
\end{figure}

\section{Multiple features}
\label{sec:multiple_feature}

In the preceding chapter, optimization was performed using a single pill only 
(cf.~Section~\ref{sec:single_feature_optimization}). 
In that setting, initialization was straightforward: 
placing one pill at the domain center with sufficiently large extension 
ensured complete coverage of the target structure. 
This allowed the optimizer to approximate the desired density field 
with relatively low complexity. 

The situation changes fundamentally once multiple pills are introduced. 
If large extensions are retained, extensive overlaps occur, 
which smear the resulting density field and impair approximation quality. 
At the same time, a central initialization of all pills is not feasible. 
A structured initialization scheme is therefore required. 
The cross initialization strategy introduced in Section~\ref{sec:init} 
is adopted here, distributing pills symmetrically and homogeneously. 
This guarantees adequate domain coverage 
and avoids unfavorable bias at the start of optimization. 

A further distinction from the single-pill case arises from the need for an aggregation operator. 
The global density field results from combining the contributions of all pills, 
and the aggregation rule has a direct impact on approximation accuracy. 
In this study, several operators and parameterizations are examined: 
a linear sum, $p$-norms with $p \in \{3,5,7,9\}$, 
softmax functions with $\beta \in \{6,10,14,18\}$, 
and a capped softmax variant with $\beta=18$ and cap $1.1$. 
Each operator is tested in optimization runs with 
five pills (minimal coverage), 
eight pills (symmetric cross initialization), 
and thirteen pills (near-complete coverage). 

Following the aggregation tests, the influence of initialization is investigated 
by comparing cross initialization with random placement. 
The effect of discretization quality is then studied 
by increasing the finite element mesh resolution 
and by varying the quadrature order. 
The role of second-order information is subsequently assessed 
by contrasting exact Hessians with BFGS approximations. 
In addition, the impact of pill count is analyzed 
for values between three and eighteen, 
highlighting phenomena of undercoverage and overcoverage. 
To alleviate undercoverage, an iterative refinement procedure is introduced 
that adaptively inserts additional pills during optimization. 
Conversely, redundancy is mitigated by heuristic strategies 
that remove superfluous pills or consolidate overlapping ones. 
Finally, the methodology is applied to a more challenging target, 
the cantilever, to demonstrate both the strengths and the limitations 
of the proposed framework.

\subsection{Global optimization settings}
\label{sec:global-config}

\begin{figure}[H]
    \centering
    \includegraphics[height=0.3\textwidth,width=0.6\textwidth]{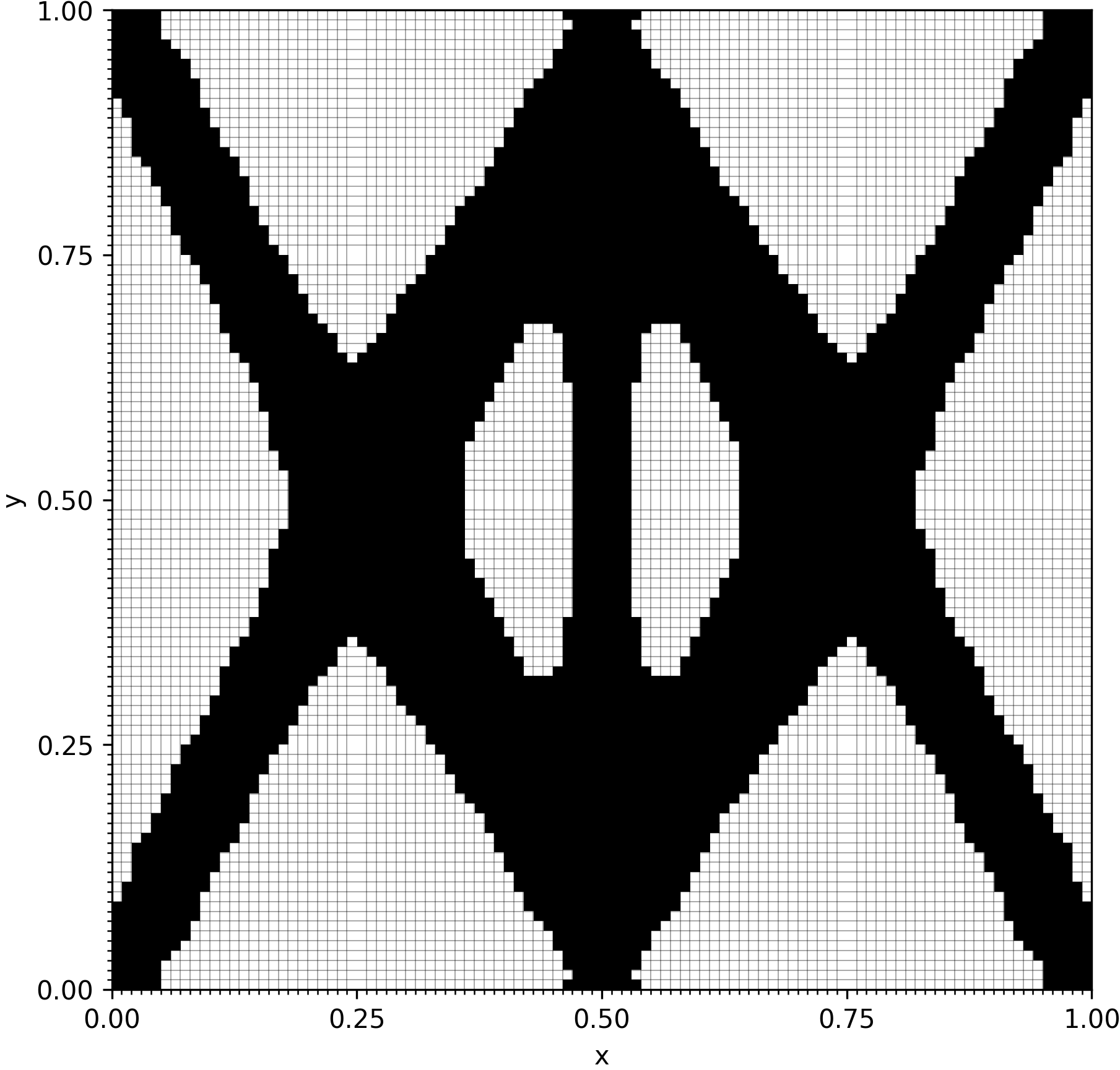}
    \caption{Target density field (five-bar structure) used for multi-pill optimization. 
    The image is stretched to reflect the $2{:}1$ aspect ratio of the design domain, 
    which improves comparability but slightly distorts the visual proportions.}
    \label{fig:target-5bar}
\end{figure}

Unless stated otherwise, all optimizations reported in this chapter were carried out under the following settings:

\begin{itemize}
  \item Target density field: five-bar structure (\cref{fig:target-5bar}).
  \item Finite element discretization: $120 \times 60$ elements.
  \item Integration order: $3$.
  \item Transition function: smoothstep with $C^k,\, k=3$.
  \item Minimum segment length: $0.05$.
  \item IPOPT settings: maximum number of iterations set to $100$.
  \item Tolerances: $10^{-2}$ (exploration), $10^{-3}$ (bridging), $10^{-7}$ (convergence).
  \item Transition width: $0.05$.
  \item Pill extension: $0.2$ (exploration), $0.1$ (bridging), $0.0$ (convergence).
  \item Pill radius: fixed $r = 0.05$ in exploration, 
        variable $0.005 < r < 0.5$ otherwise.
\end{itemize}

The target density field employed in all experiments, apart from the cantilever benchmark, is shown in \cref{fig:target-5bar}. 
The structure consists of five dominant load-bearing members that must be reconstructed by the pill-based representation. 

\begin{figure}[H]
    \centering
    \begin{subfigure}{0.32\textwidth}
        \includegraphics[width=\linewidth]{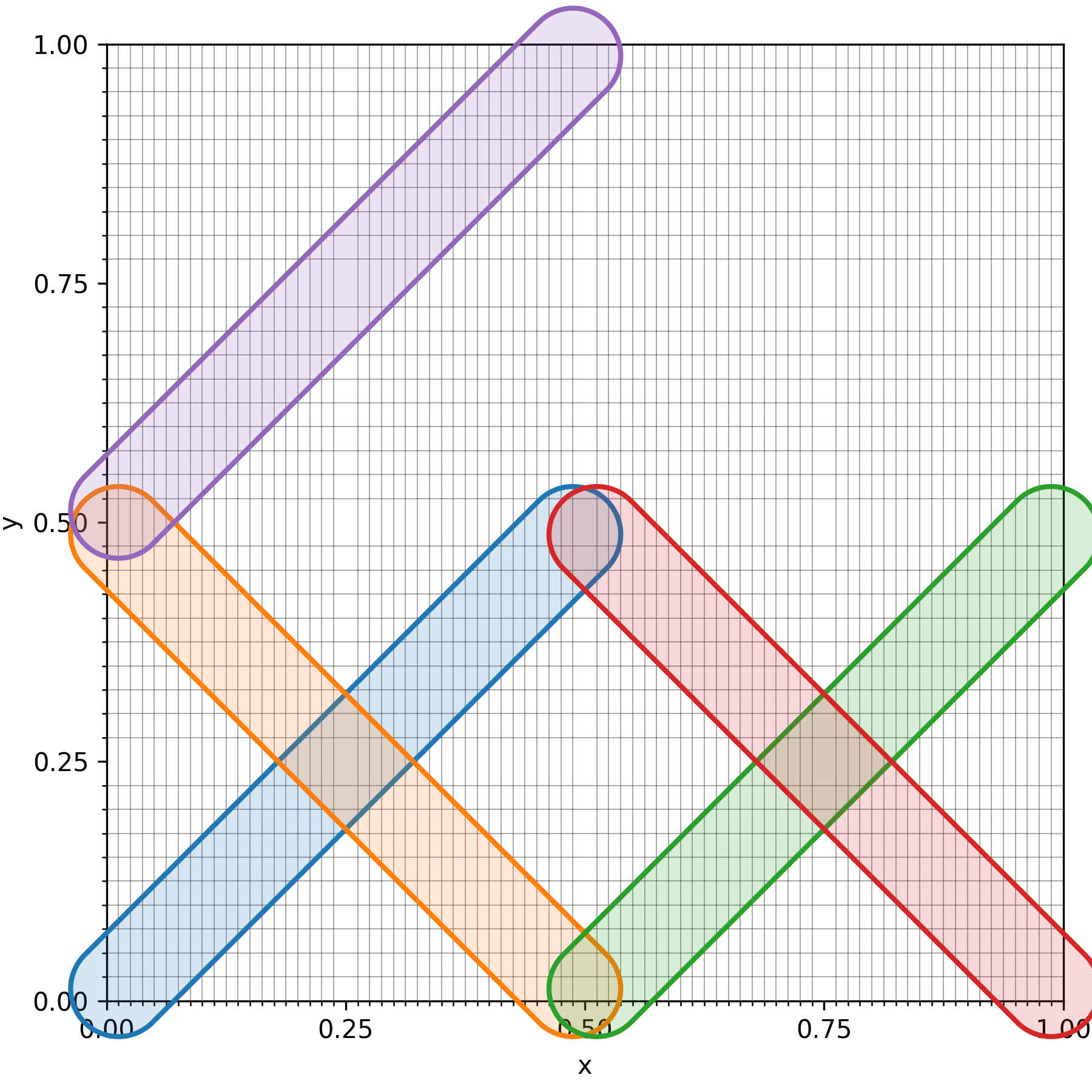}
        \caption{Five pills}
    \end{subfigure}
    \hfill
    \begin{subfigure}{0.32\textwidth}
        \includegraphics[width=\linewidth]{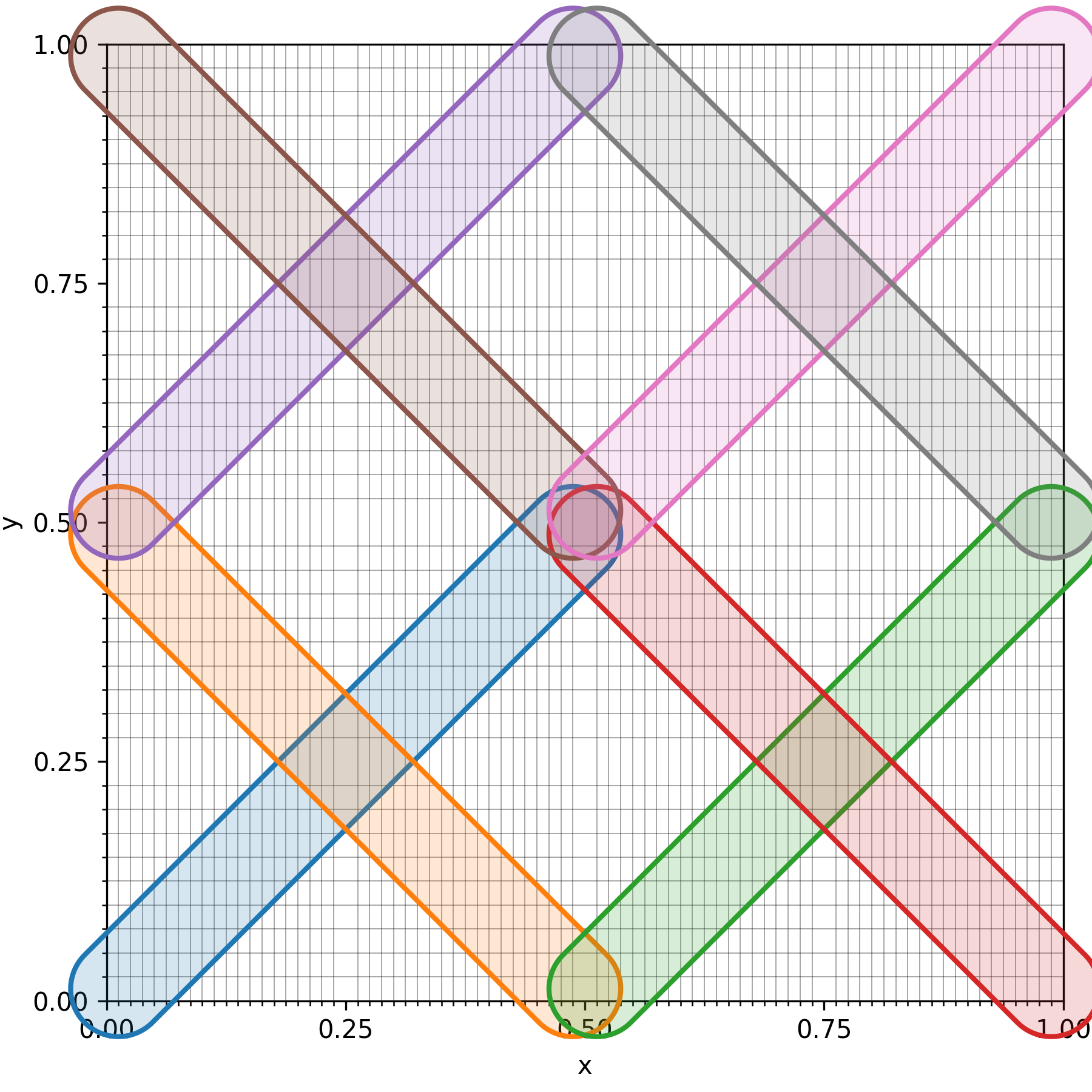}
        \caption{Eight pills}
    \end{subfigure}
    \hfill
    \begin{subfigure}{0.32\textwidth}
        \includegraphics[width=\linewidth]{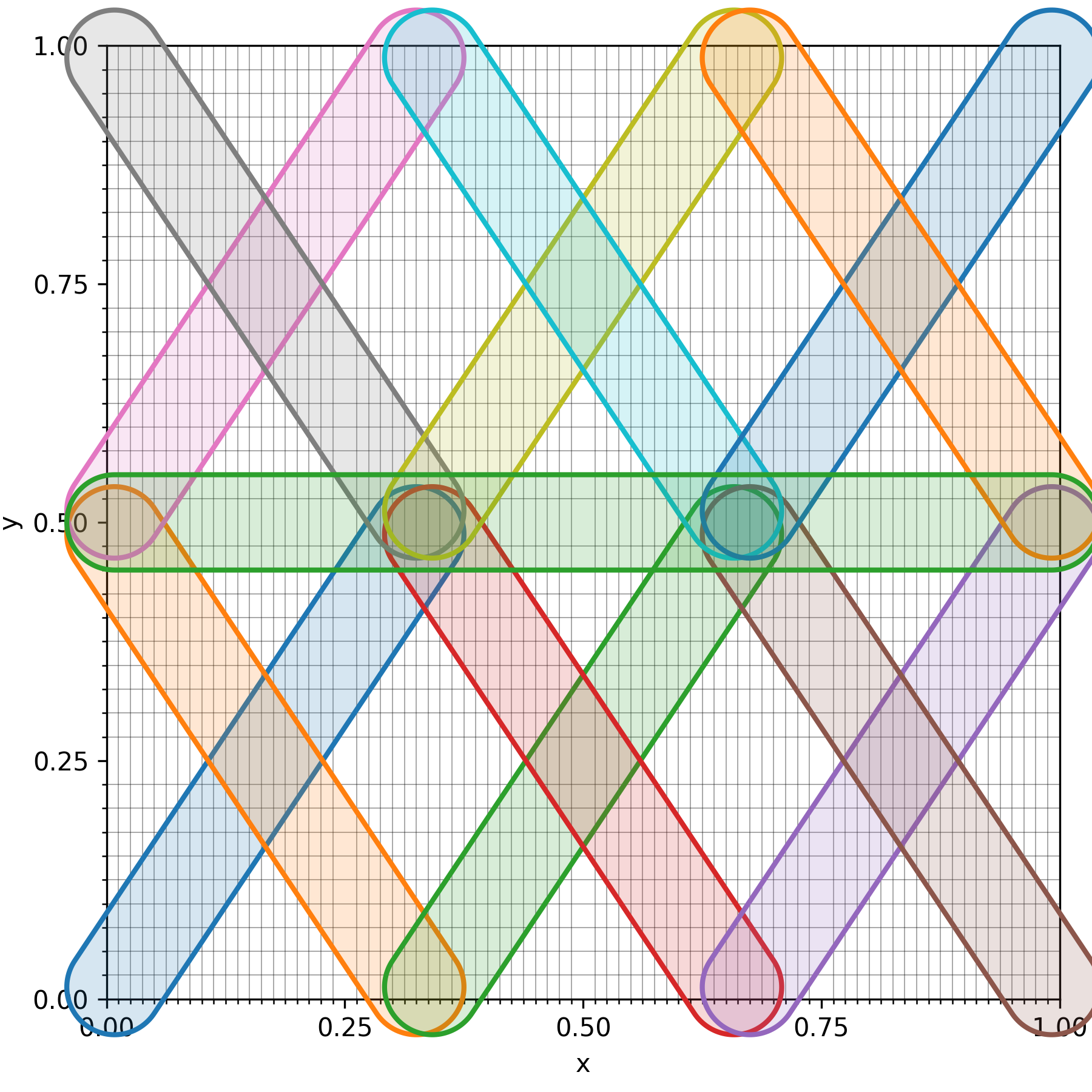}
        \caption{Thirteen pills}
    \end{subfigure}
    \caption{Cross initialization for different pill counts. 
    The design domain has an aspect ratio of $2:1$, although images are shown in square format.}
    \label{fig:initial-cross}
\end{figure}

The corresponding cross initializations are shown in \cref{fig:initial-cross}. 
With five pills, only partial coverage is achieved, leaving pronounced gaps across the domain. 
Eight pills result in a symmetric initialization that provides moderate coverage, 
whereas thirteen pills produce a denser configuration that additionally includes a vertical pill in the center.
These configurations illustrate the trade-off between computational effort and geometric coverage. 
A more detailed discussion of this trade-off is given in Section~\ref{sec:pill-count}. 

\subsection{Aggregation tests}
\label{sec:agg-tests}

Before comparing different aggregation operators, the diagnostic quantities used throughout this section are introduced. 
Objective values alone are not directly comparable, since both the tracking and reward objectives depend on the aggregated density 
$\bar\rho$, which varies with the chosen operator. 
To obtain a measure of reconstruction quality that is independent of aggregation, 
a set of auxiliary visualizations is employed. 
These complement the numerical objectives and provide insight into how the optimizer places pills 
and how closely the resulting density field matches the target.

In addition to the objective trajectories, four types of fields are considered. 
The density plot shows the aggregated field obtained from the superposition of all pills. 
It is particularly sensitive to the choice of aggregation operator and extension parameters, 
as both directly determine how overlaps and smooth transitions are expressed. 
The reward field represents what is evaluated during the exploration stage, 
highlighting regions where the current approximation coincides with the target. 
The residual field records the signed difference between target and approximation, 
distinguishing between uncovered target (positive values) and excess density in void (negative values). 
Finally, the absolute residual displays the total deviation irrespective of sign, 
avoiding cancellation and providing a direct measure of approximation quality.

Together, these visualizations provide a comprehensive assessment of the optimization process. 
The density field illustrates the achieved coverage, 
the reward field reflects the criterion promoted during exploration, 
and the residual fields quantify the deviations that remain. 
These representations are used consistently in the subsequent comparisons of aggregation rules 
and thereby establish a direct connection between numerical objectives and spatial behavior.

\begin{figure}[H]
  \centering
  \setlength{\tabcolsep}{2pt}
  \renewcommand{\arraystretch}{0.5}
  \begin{tabular}{ccc}
      \includegraphics[height= 0.16\textwidth,width=0.32\textwidth]{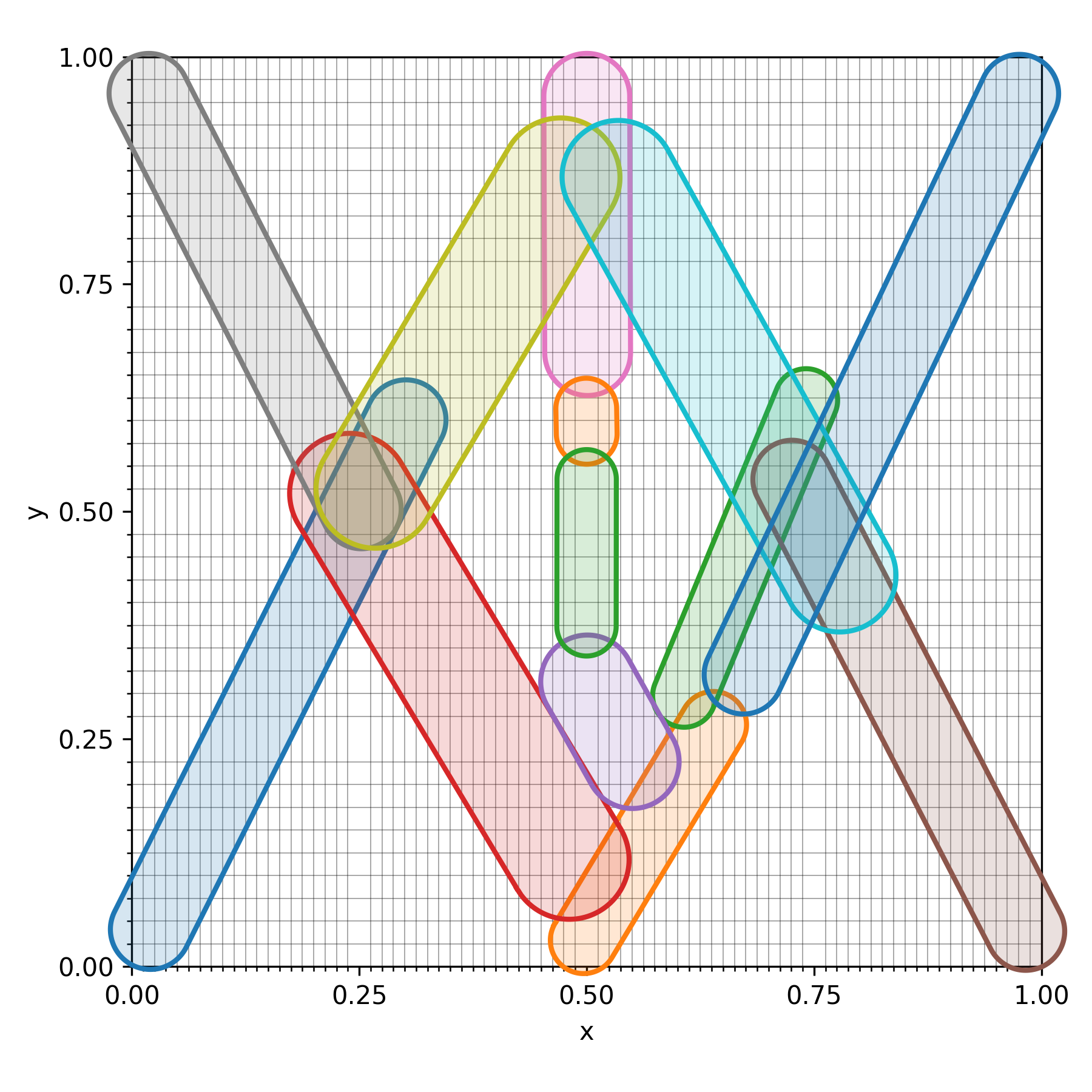} &
      \includegraphics[height= 0.16\textwidth,width=0.32\textwidth]{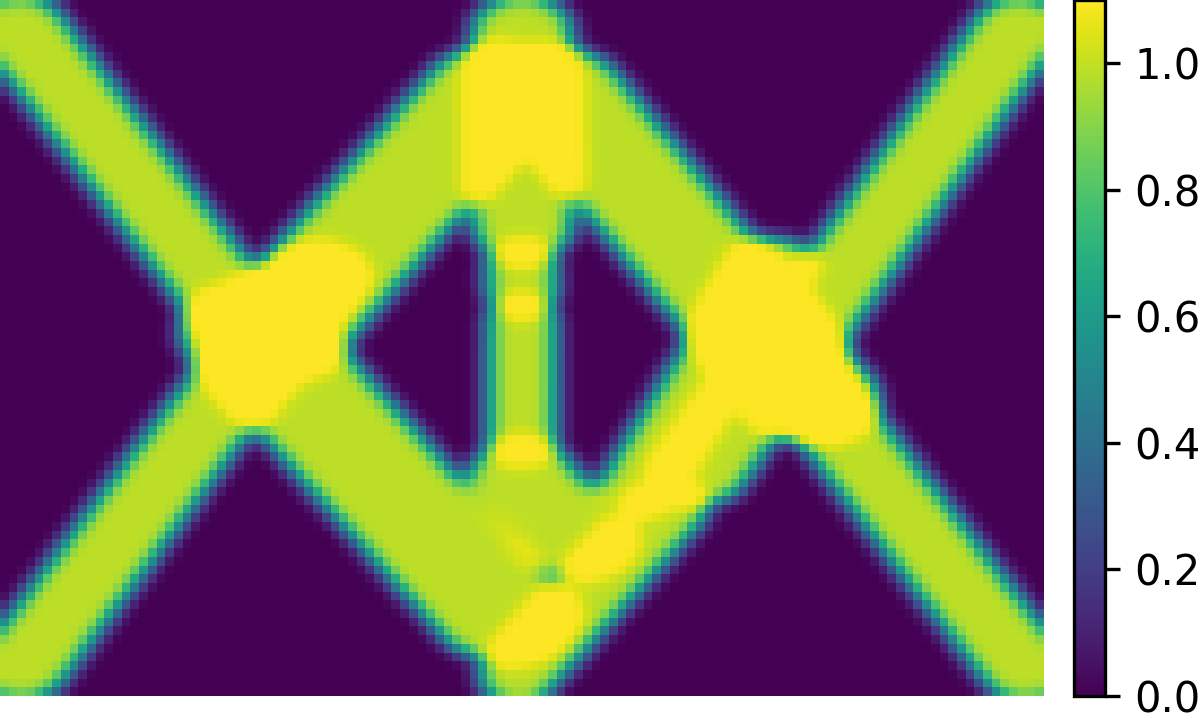} &
      \includegraphics[height= 0.16\textwidth,width=0.32\textwidth]{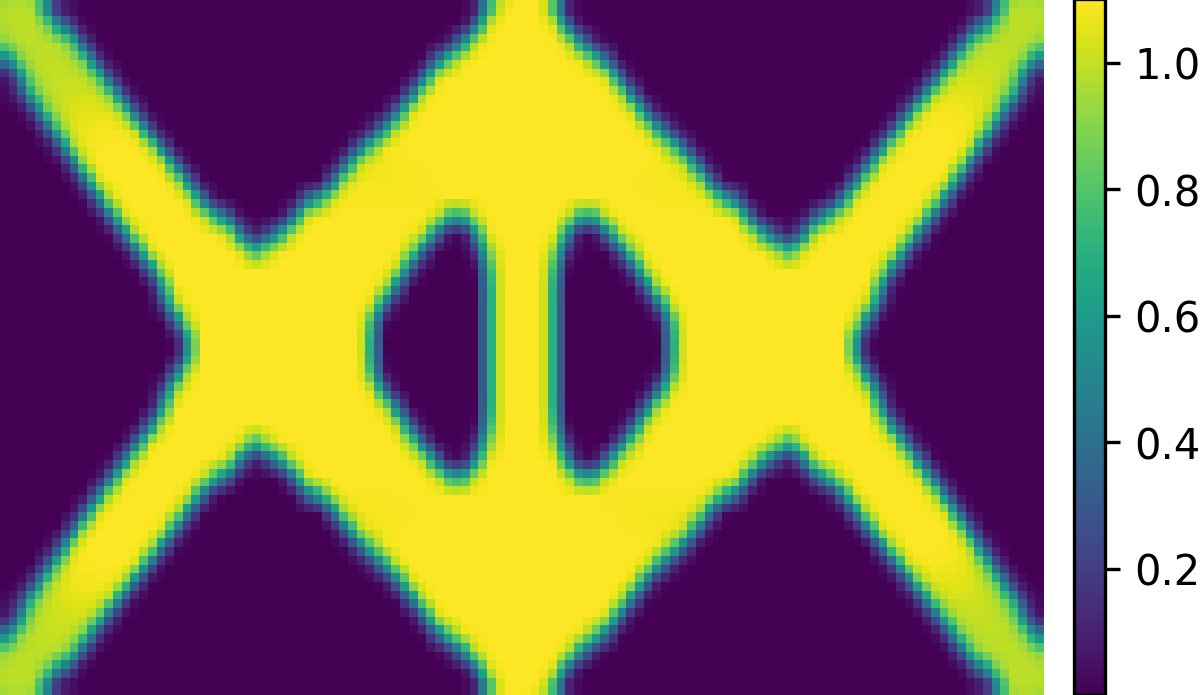} \\
      \includegraphics[height= 0.16\textwidth,width=0.32\textwidth]{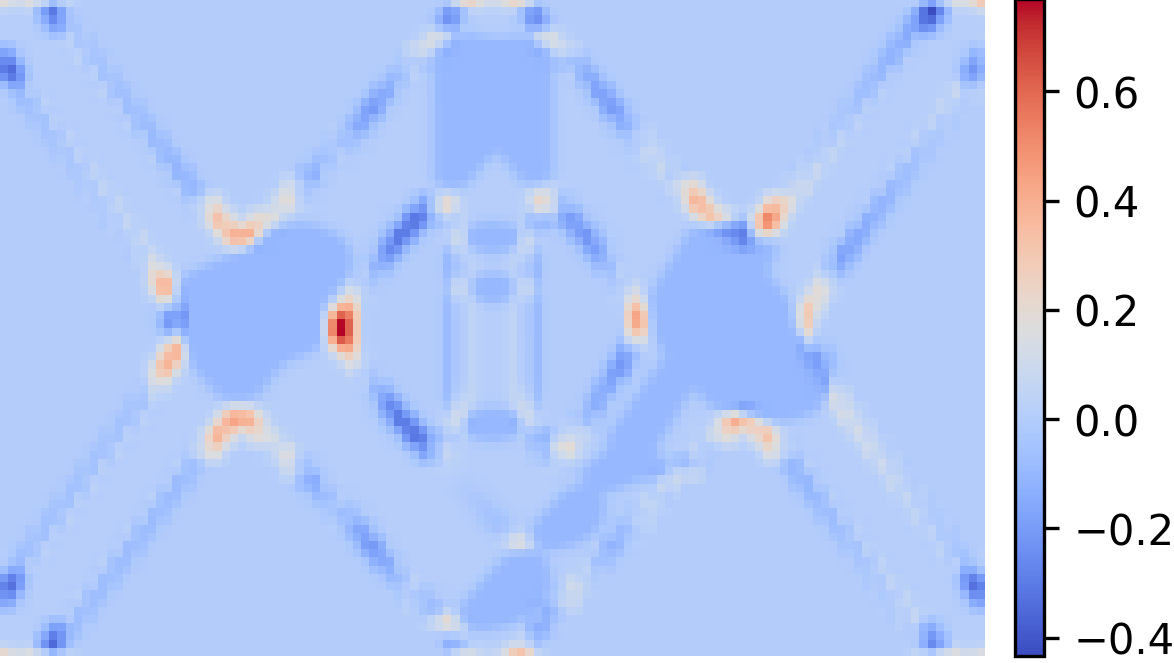} &
      \includegraphics[height= 0.16\textwidth,width=0.32\textwidth]{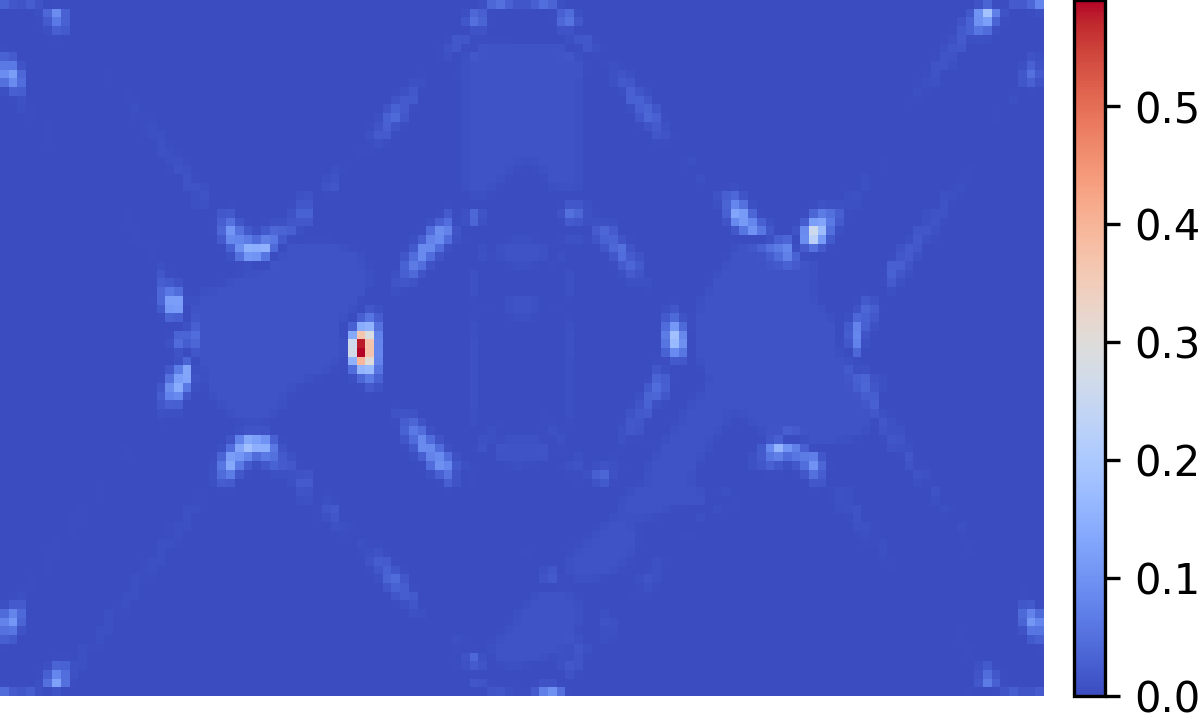} &
      % Leeres Feld für Ausrichtung
      \\
  \end{tabular}
  \caption{Diagnostic visualizations used throughout the aggregation tests. 
  Top row: pill placement (left), aggregated density field (middle) and reward field (right). 
  Bottom row: residual field (left) and absolute residual field (middle).}
  \label{fig:agg-diagnostics}
\end{figure}

\subsubsection{Exploration stage}
\label{sec:agg-exploration}

During the exploration stage, the choice of aggregation operator strongly affects how pills respond to the reward signal. 
Operators from the $p$--norm and softmax families increasingly suppress weaker contributions as their parameters grow, 
whereas the sum and sum--softcap operators treat all contributions equally. 
As a result, low--parameter $p$--norm and softmax as well as sum--based rules tend to promote overlap, 
since covering the same target region with several pills yields a higher reward. 
This effect is clearly visible in the reward fields, where values seldom exceed~1.2 except under pure sum aggregation, 
in which overlapping pills raise the aggregated density above unity. 
The tendency toward overlap is therefore intrinsic to the reward functional 
and is mitigated only when stronger nonlinearities ($p=9$ or $\beta=18$) suppress weaker contributions. 

The behavior of the sum--softcap operator illustrates this particularly well. 
With eight and thirteen pills, the pills accumulate into a single homogeneous density region at the center, 
maximizing mutual overlap rather than separating onto distinct target parts. 
The effect becomes more pronounced with larger extension values, 
indicating that overlap is favored whenever it increases the immediate reward. 
In contrast, high--parameter $p$--norm and softmax runs yield more separation 
and discourage redundant superposition.  

Despite these differences, all configurations exhibit difficulty in reconstructing the central vertical bar: 
most pills intersect it obliquely without aligning to its exploration. 
Only in the thirteen--pill runs do several configurations succeed in capturing it directly, 
notably the $p$--norm with $p=9$ and softmax with $\beta=6,10$. 
A drawback of higher parameters is, however, reduced mobility of weakly contributing pills. 
For instance, in the thirteen--pill runs, some pills under high--parameter $p$--norm or softmax 
remain close to their initial positions when their potential contribution is already covered by stronger neighbors. 
This lack of mobility suggests that max--biased aggregation emphasizes reinforcement of existing coverage 
rather than exploratory repositioning.  

Overall, exploration under high--parameter $p$--norm and softmax produces configurations 
in which overlap is minimized and each pill tends to reinforce its assigned region. 
This behavior proves advantageous in later stages: 
when the objective shifts from reward to tracking, 
the overlaps promoted by sum--type aggregation become detrimental, 
causing pills to repel each other and spread erratically. 
By contrast, max--biased aggregation already enforces a more target--oriented distribution during exploration 
and thus provides a more favorable starting point for subsequent convergence. 

\begin{figure}[H]
  \centering
  \includegraphics[width=0.79\textwidth]{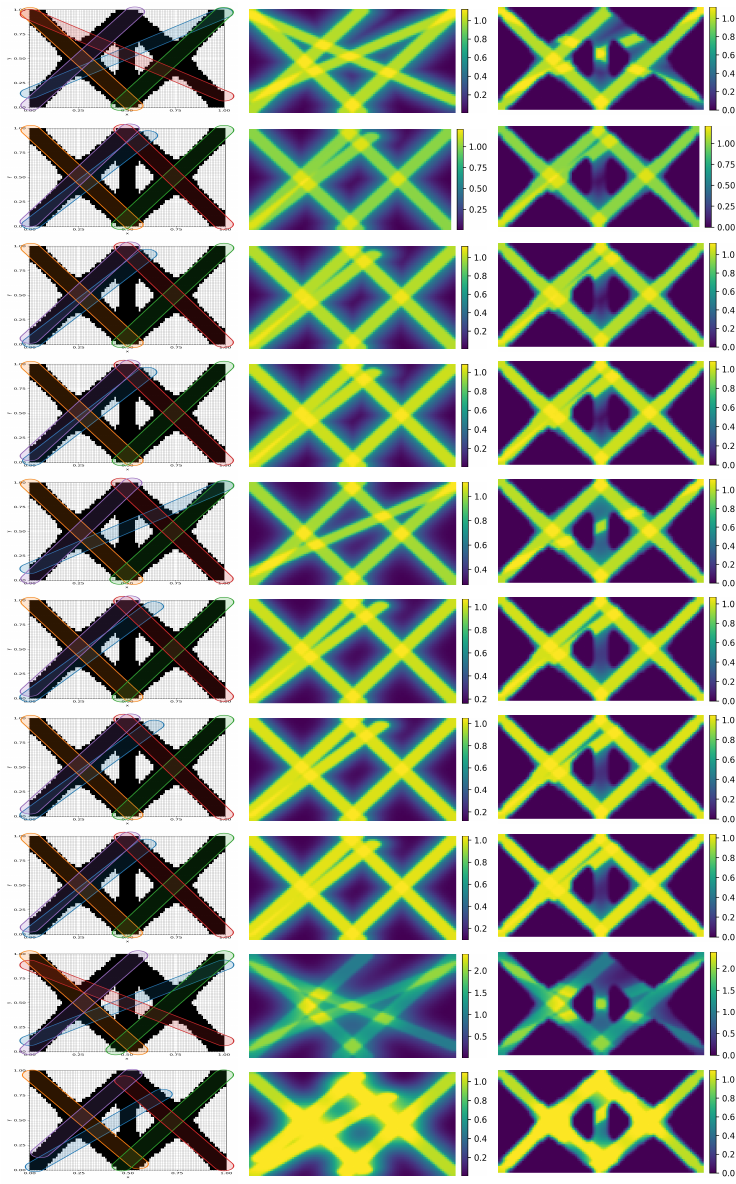}
  \caption{Exploration stage reconstructions with five pills.
  From top to bottom: $p$--norm ($p=3,5,7,9$), softmax ($\beta=6,10,14,18$),
  sum and sum--softcap ($\beta=18$, cap $=1.1$).
  Columns: pill placement, density, reward.}
  \label{fig:five-pills-pdf}
\end{figure}

\begin{figure}[H]
  \centering
  \includegraphics[width=0.79\textwidth]{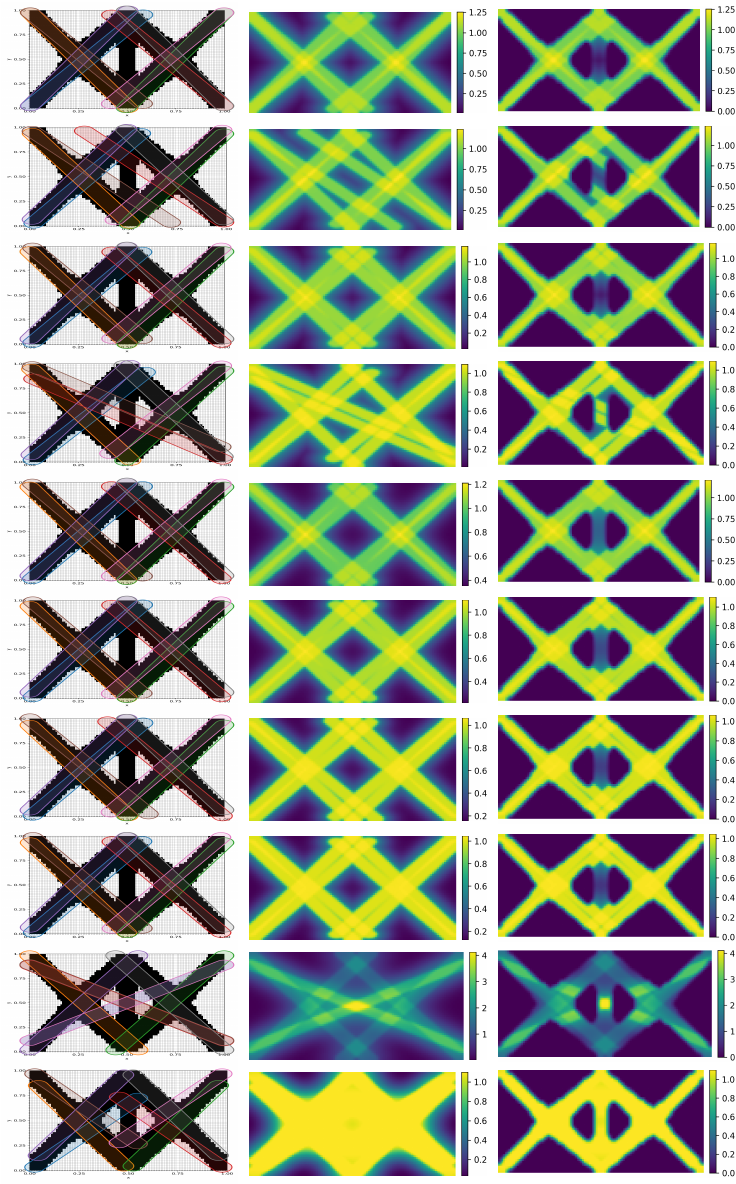}
  \caption{Exploration stage reconstructions with eight pills. 
  From top to bottom: $p$--norm ($p=3,5,7,9$), softmax ($\beta=6,10,14,18$), 
  sum and sum--softcap ($\beta=18$, cap $=1.1$). 
  Columns: pill placement, density, reward.}
  
  \label{fig:eigth-pills}
\end{figure}

% --- Feature count 13 ---
% >>> INSERT DISCUSSION TEXT HERE <<<
\begin{figure}[H]
  \centering
  \includegraphics[width=0.79\textwidth]{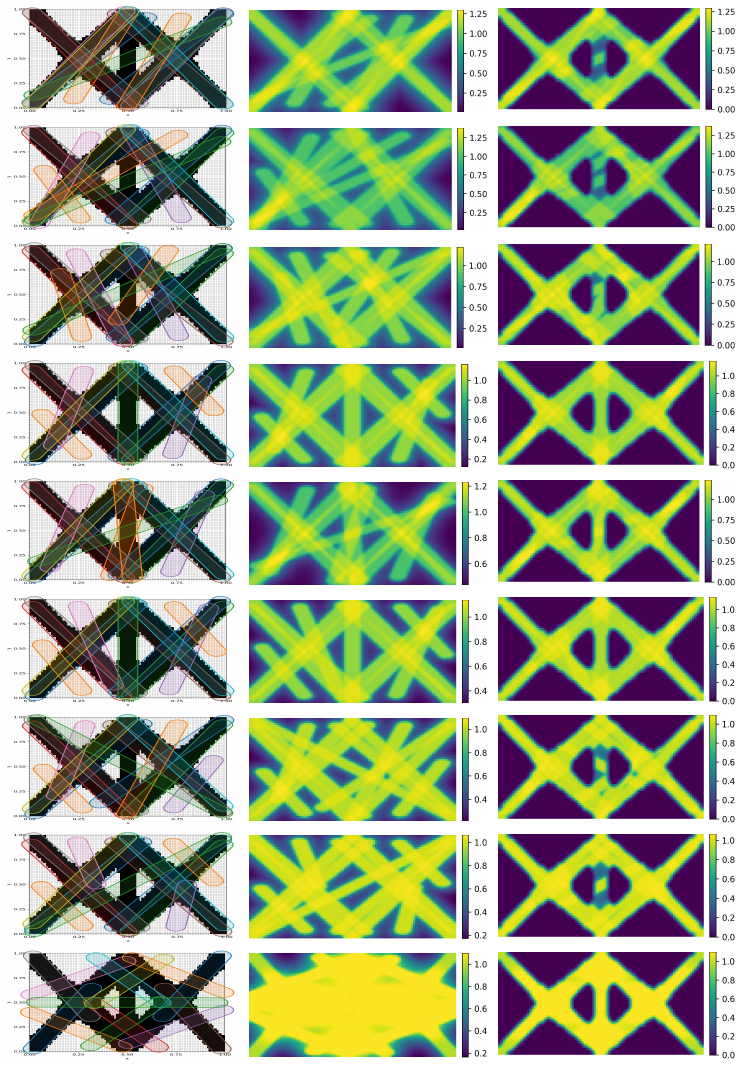}
  \caption{Exploration stage reconstructions with thirteen pills. 
  From top to bottom: $p$--norm ($p=3,5,7,9$), softmax ($\beta=6,10,14,18$) and sum--softcap ($\beta=18$, cap $=1.1$). 
  Columns: pill placement, density, reward.}
  \label{fig:thirteen-pills}
\end{figure}

\subsection{Cosine--weighted aggregation}
\label{subsec:cosine_agg}

As an additional variant, a cosine--weighted shaping of the raw sum was examined. 
The idea is to reward the first unit of mass while smoothly penalizing further overlap. 
Formally, let \(S=\sum_{m=1}^{n}\rho_m\) and define \(A:\R_{\ge0}\to\R\) by
\[
A(S)=
\begin{cases}
\sin\!\bigl(\tfrac{\pi}{2}\,S\bigr), & 0\le S \le 1,\\[0.4ex]
a + (1-a)\,\dfrac{1+\cos\!\bigl(\tfrac{\pi}{N-1}(S-1)\bigr)}{2}, & 1< S < N,\\[0.9ex]
a, & S\ge N,
\end{cases}
\qquad
a=1-\dfrac{(N-1)^2}{2}.
\]
The aggregate is then given by \(g(x,y)=A(S(x,y))\), with analytic first and second derivatives
\[
\nabla_{\!\vecsym{\zeta}} g \;=\; A'(S)\sum_{m=1}^{n} \nabla_{\!\vecsym{\zeta}} \rho_m,
\qquad\\
\nabla^2_{\!\vecsym{\zeta}\vecsym{\zeta}} g \;=\; A''(S)\Bigl(\sum_{m} \nabla_{\!\vecsym{\zeta}} \rho_m\Bigr)\!
\Bigl(\sum_{j} \nabla_{\!\vecsym{\zeta}} \rho_j\Bigr)^{\!\top}
+ A'(S)\sum_{m} \nabla^2_{\!\vecsym{\zeta}\vecsym{\zeta}}\rho_m.
\]
This formulation ensures smooth repulsion against excessive co--location 
while retaining full second--order information suitable for Newton--type optimization methods.

In practice, however, this operator proved conservative during exploration. 
Starting from the eight--pill cross initialization on the five--bar target 
with extension parameter \(\text{ext}=0.2\), the pills largely remained within 
their initial territories and migrated only sluggishly toward uncovered regions. 
As a result, coverage stagnated: several parts of the target remained only partially 
captured and residual gaps persisted throughout. 
The cosine formulation thus promoted dispersion and mutual avoidance, 
but at the expense of the active reorganization required for sharper alignment 
with the target.

\begin{figure}[H]
  \centering
  \includegraphics[width=0.78\textwidth]{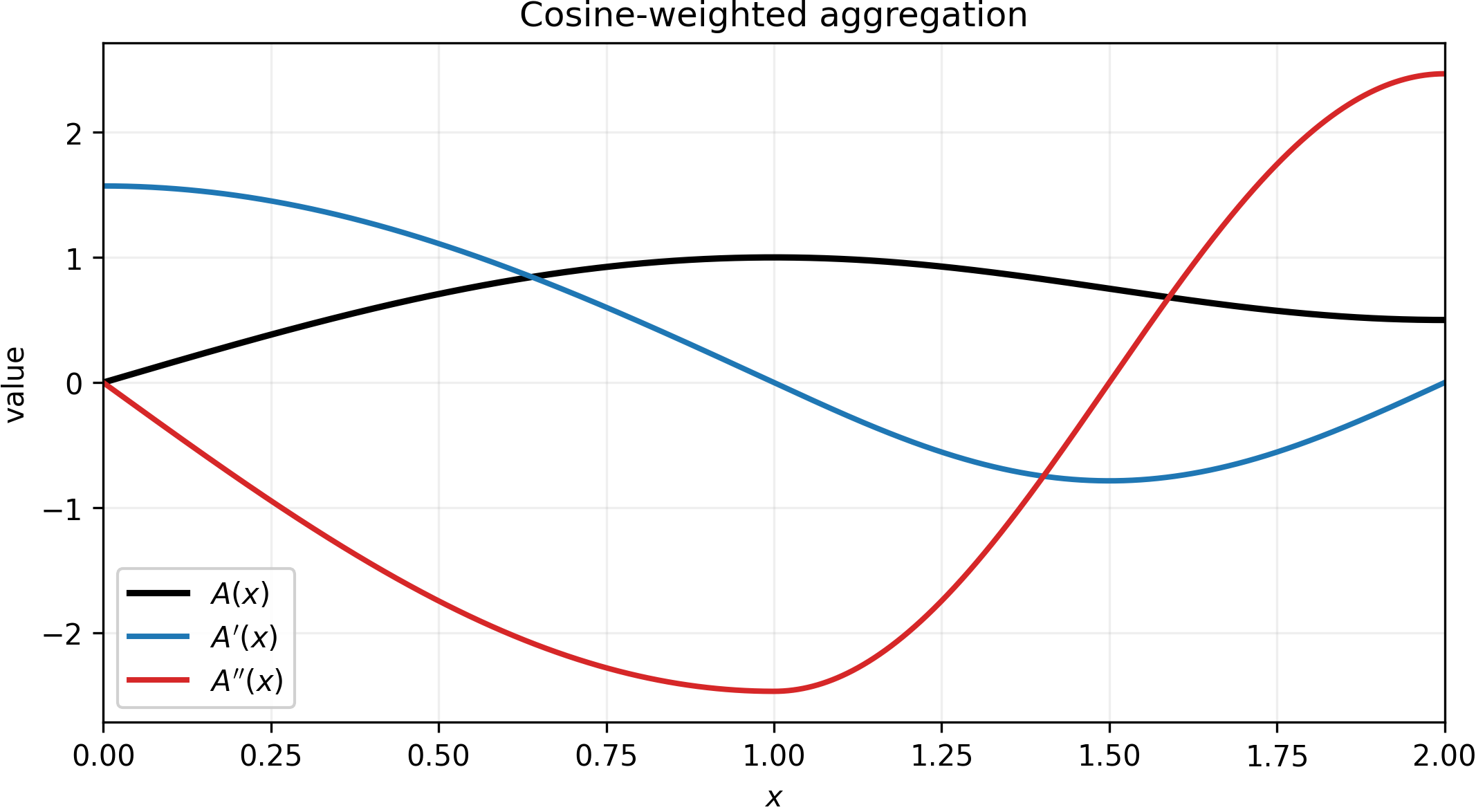}
  \caption{Cosine--weighted shaping \(A(S)\) and its first two derivatives. 
  The first unit of mass is rewarded, while additional overlap yields diminishing 
  and eventually vanishing marginal gain.}
  \label{fig:cosine_weighted}
\end{figure}

\begin{figure}[H]
  \centering
  \includegraphics[width=0.58\textwidth]{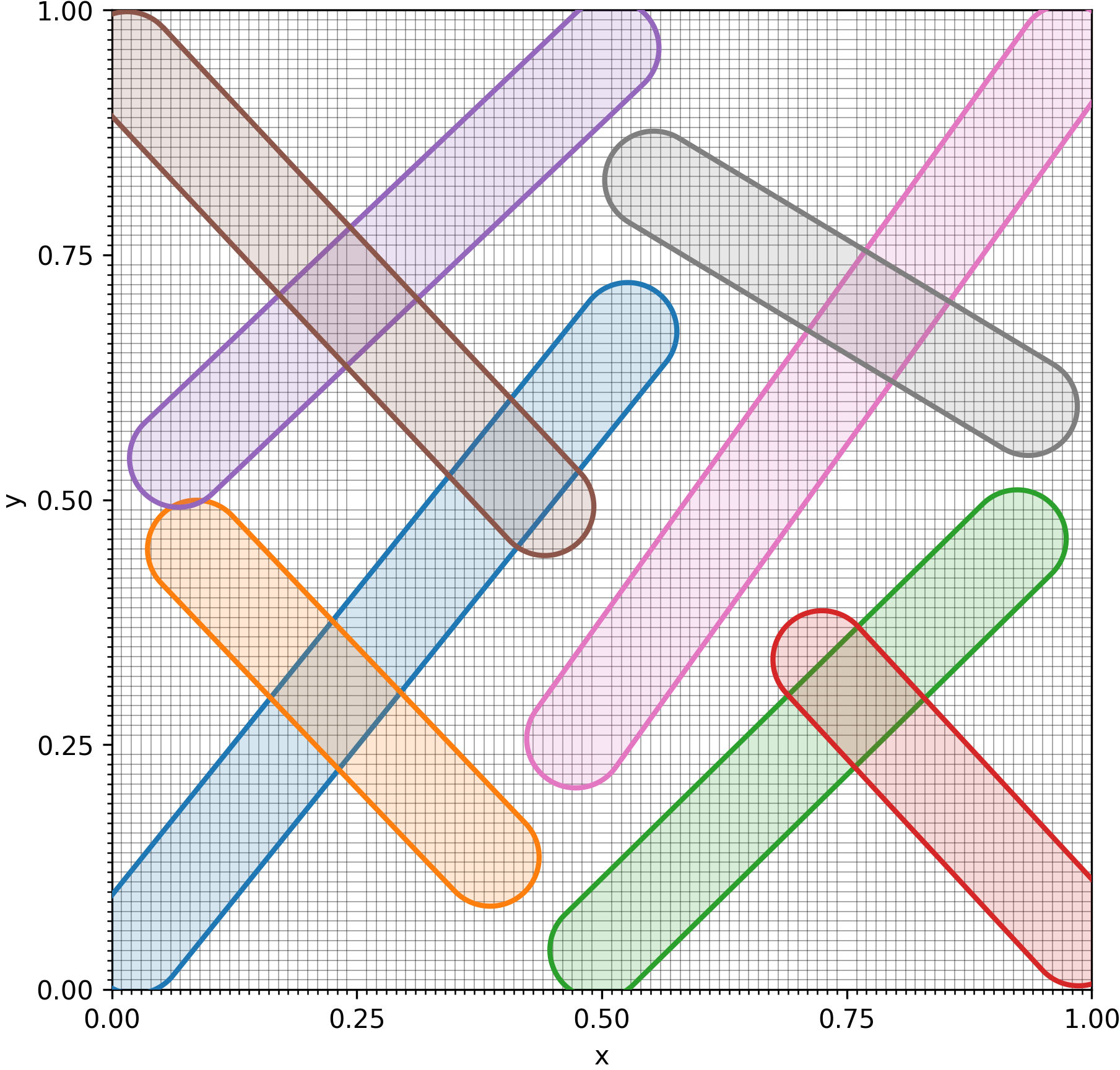}
  \caption{Exploration under cosine--weighted aggregation. 
  Features remain dispersed and reorganize only sluggishly toward uncovered bars, 
  leading to slower closure of gaps compared to $p$--norm and softmax aggregation.}
  \label{fig:cosine_layout}
\end{figure}

\subsubsection{Bridging stage}
\label{sec:agg-bridging}

The bridging stage continues directly from the exploration results: 
each run is initialized with the final pill configuration obtained 
under the respective aggregation rule at the end of exploration. 
From this starting point, the optimizer seeks to improve alignment with the target, 
while the extension parameter is reduced to promote sharper localization. 
The following analysis focuses on how well the different aggregation operators 
support reorganization at this stage, in particular regarding the recovery of 
the central vertical bar of the target, which is frequently missed during exploration.  

For five pills, bridging succeeds in covering the central region 
only if it had already been intersected during exploration. 
This is particularly evident for $p=3$ and for the sum--softcap variant, 
where a pill that previously crossed the middle realigns more directly with the bar. 
In the other five--pill runs, the central part remains uncovered. 
The residual fields confirm that alignment is only partial, 
with substantial vertical gaps persisting.

With thirteen pills, the situation changes markedly. 
Here the sum--softcap variant performs less favorably: 
with more pills present, the overlap penalty becomes dominant 
and overall coverage deteriorates. 
The optimization distributes the pills at larger mutual distances, 
reducing their ability to follow the target consistently. 
By contrast, the $p$--norm and softmax operators yield substantially better approximations. 
Softmax with $\beta=6$ and $\beta=10$ reconstructs the central bar successfully, 
while the higher settings $\beta=14$ and $\beta=18$ still fail to capture it. 
The residual fields show that $p$--norm solutions achieve comparable accuracy 
across different exponents, but their internal distributions differ: 
for higher $p$, one dominant pill tends to expand while smaller ones retract, 
whereas lower $p$ values distribute the load among several pills of similar size. 

Overall, the bridging stage confirms that $p$--norm and softmax operators 
outperform sum--based aggregation, with the most reliable reconstructions 
achieved for intermediate parameters ($p=7,9$ or $\beta=6,10$). 
Nevertheless, the results remain strongly influenced by the preceding exploration: 
if the central bar was not intersected during that phase, 
bridging rarely succeeds in establishing it. 
In general, configurations in which pills settle at smaller mutual distances 
prove more effective, as they are more likely to form coherent structures, 
while larger separations tend to leave residual gaps.

\begin{figure}[H]
  \centering
  \includegraphics[width=0.79\textwidth]{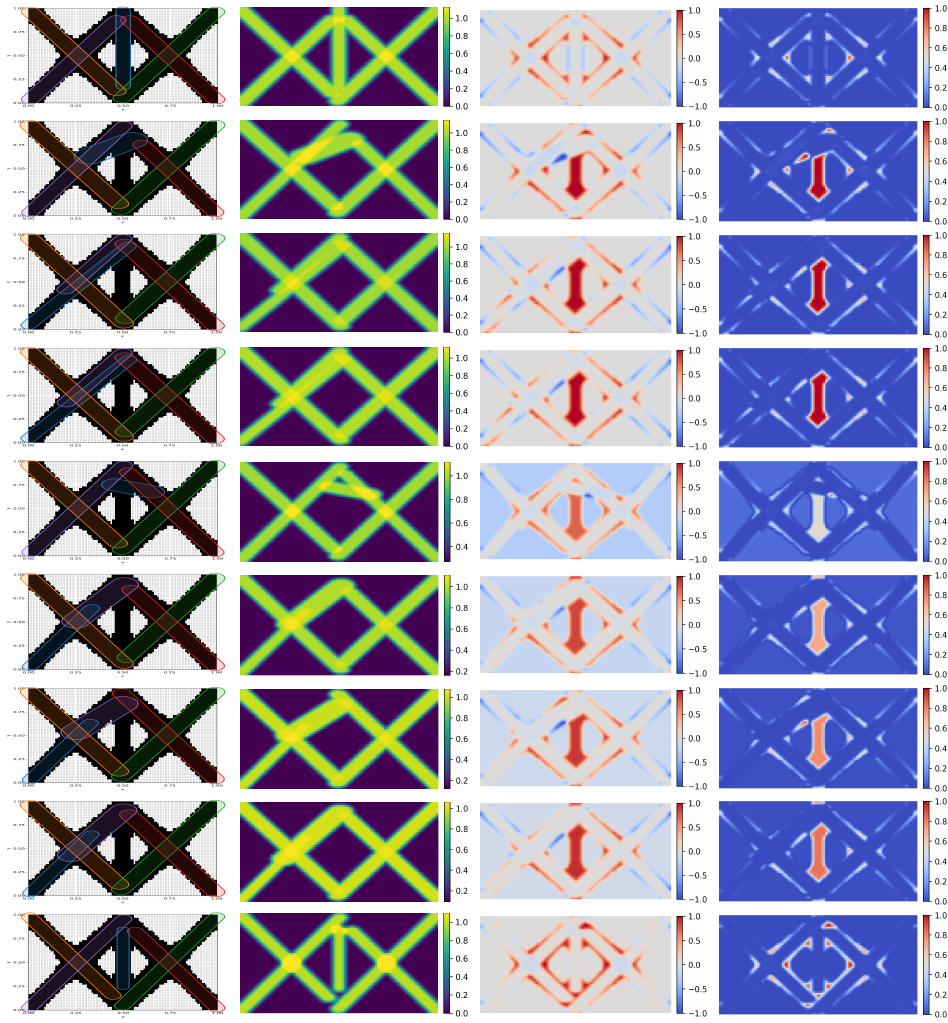}
  \caption{Bridging stage reconstructions with five pills. 
  From top to bottom: $p$--norm ($p=3,5,7,9$), softmax ($\beta=6,10,14,18$) and sum--softcap ($\beta=18$, cap $=1.1$). 
  Columns: pill placement, density, reward.}
  \label{fig:five-pills-bridging}
\end{figure}
\begin{figure}[H]
  \centering
  \includegraphics[width=0.79\textwidth]{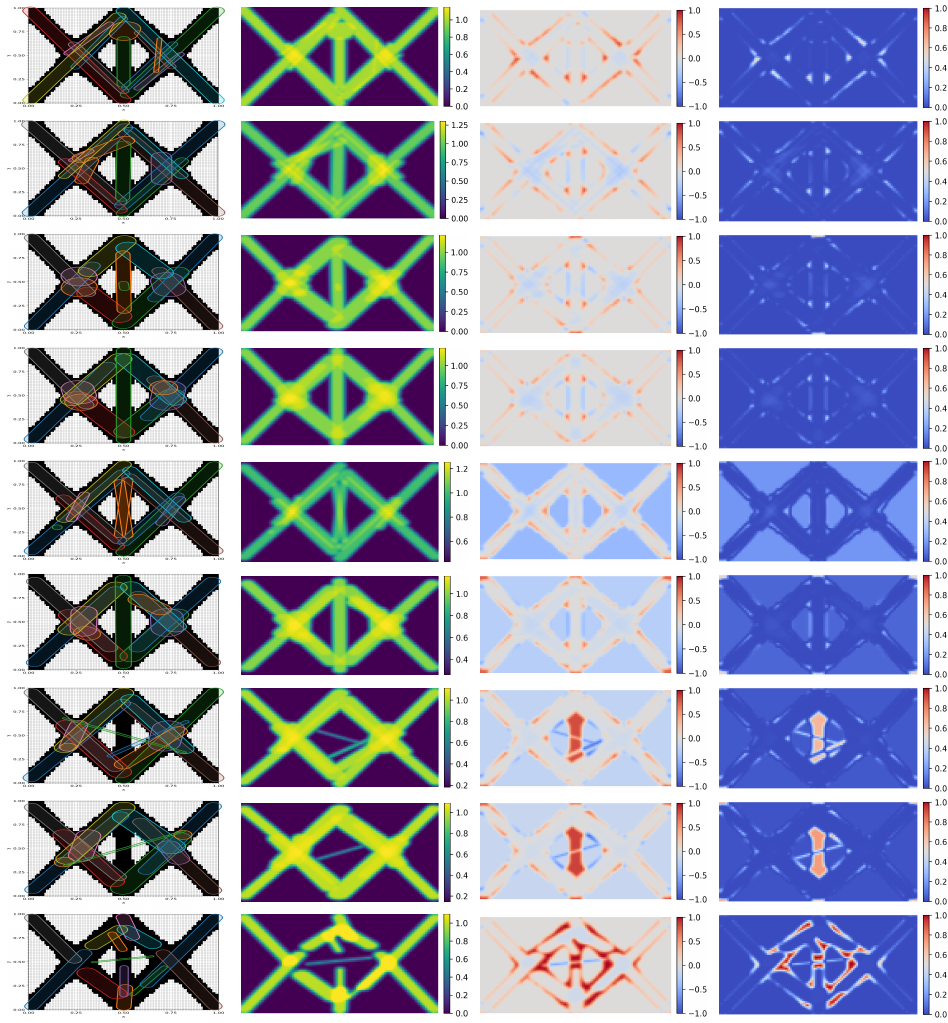}

  \caption{Bridging stage reconstructions with thirteen pills. 
  From top to bottom: $p$--norm ($p=3,5,7,9$), softmax ($\beta=6,10,14,18$) and sum--softcap ($\beta=18$, cap $=1.1$). 
  Columns: pill placement, density, reward.}
  \label{fig:thirteen-pills-bridging}
\end{figure}

\subsubsection{Convergence stage}
\label{sec:agg-convergence}

In the final convergence stage, the analysis is restricted to runs with five pills. 
At this point, the pills undergo only minor adjustments, 
primarily sharpening the approximation as the extension parameter is reduced to zero. 
The overall solution quality is largely determined by the preceding stages, 
since pill positions no longer change substantially but stabilize in the configuration inherited from bridging.  

As in earlier stages, the spacing between pills depends strongly on the aggregation rule. 
Higher values of $p$ or $\beta$ encourage tighter clustering with more dominant individual pills, 
whereas smaller values or the plain sum distribute the pills more widely 
and produce several segments of comparable strength. 
The final density and residual fields confirm these tendencies, 
demonstrating that the character of the reconstruction is essentially fixed before convergence 
and only refined in this stage.  

Across all experiments, the plain sum operator proved to be the least effective choice. 
By treating all overlaps equally, it tends to blur the density field already during exploration, 
which compromises subsequent stages. 
The capped sum variant alleviates this effect only marginally. 
By contrast, differences between $p$--norm and softmax operators are less pronounced, 
with both producing comparably accurate solutions. 
Nevertheless, the parameter setting leaves a distinct imprint: 
higher values yield reconstructions dominated by a few strong pills surrounded by smaller auxiliaries, 
while lower values distribute influence more evenly among several pills of similar size. 
These behaviors reflect the inherent properties of the respective operators 
and remain consistent across all optimization stages.

\begin{figure}[H]
  \centering
  \includegraphics[width=0.79\textwidth]{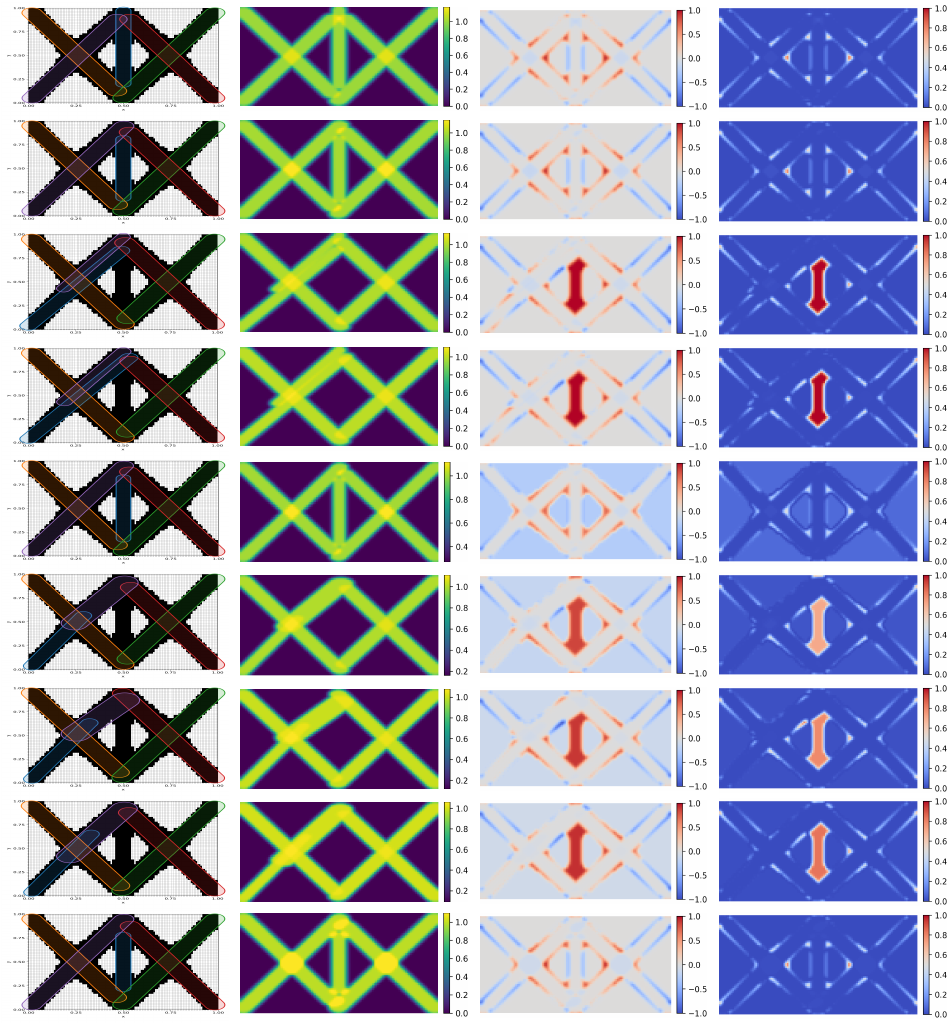}

  \caption{Convergence stage reconstructions with five pills. 
  From top to bottom: $p$--norm ($p=3,5,7,9$), softmax ($\beta=6,10,14,18$) and sum--softcap ($\beta=18$, cap $=1.1$). 
  Columns: pill placement, density, reward.}
  \label{fig:five-pills-convergence}
\end{figure}

\subsection{Random versus cross initialization}
\label{sec:init_opt}

\begin{figure}[H]
  \centering
  \setlength{\tabcolsep}{2pt}
  \renewcommand{\arraystretch}{0.5}
  \begin{tabular}{ccccc}
    \includegraphics[height= 0.095\textwidth,width=0.19\textwidth]{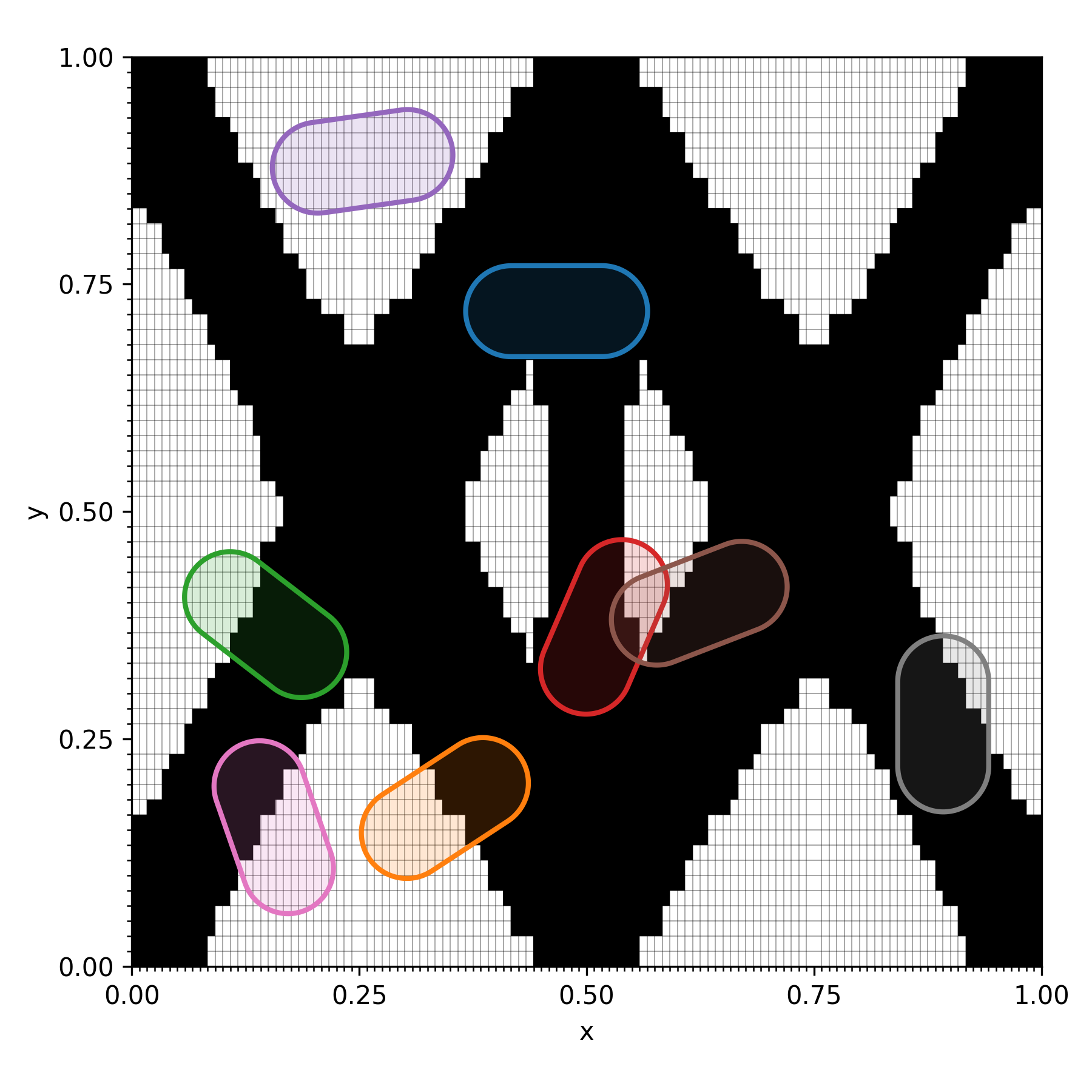} &
    \includegraphics[height= 0.095\textwidth,width=0.19\textwidth]{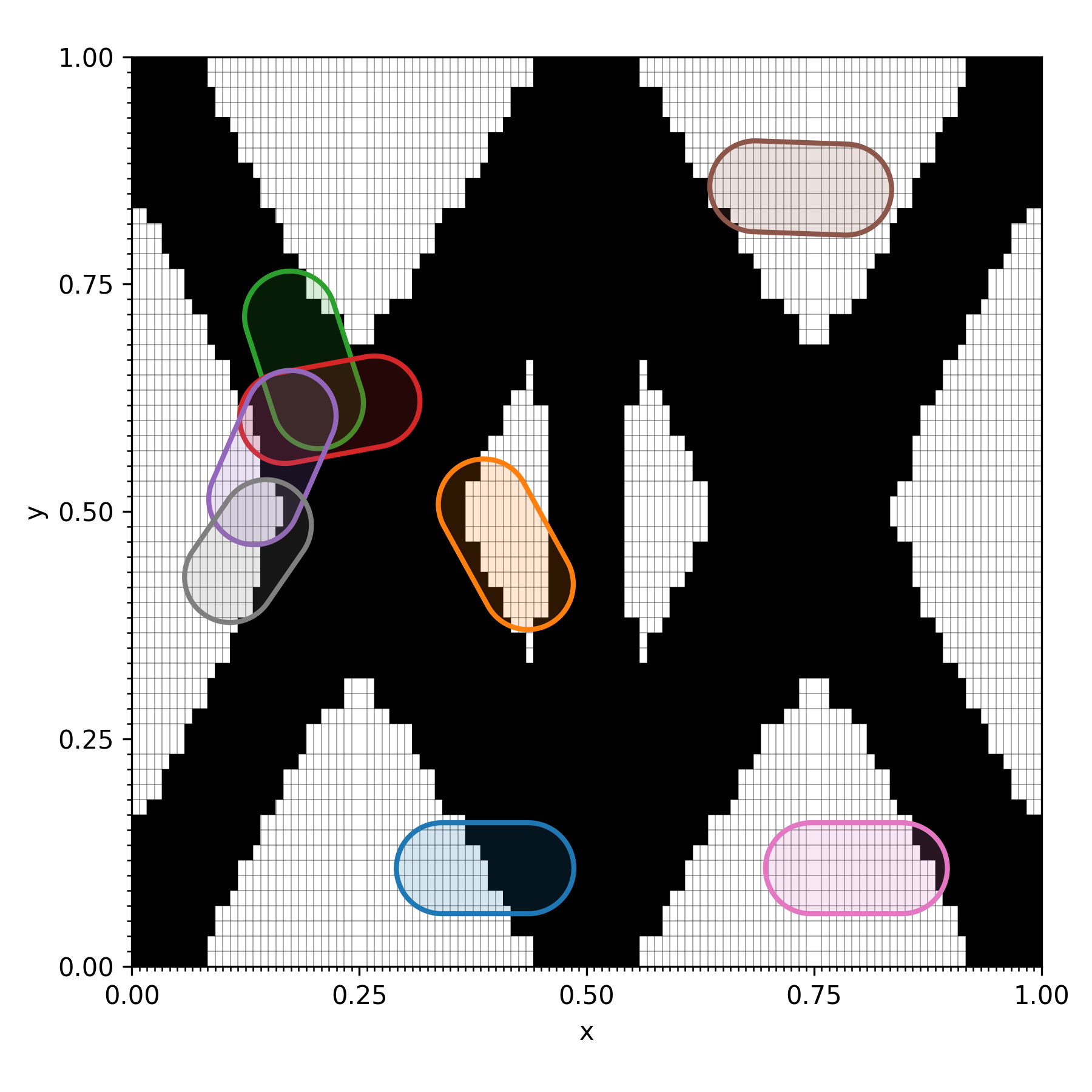} &
    \includegraphics[height= 0.095\textwidth,width=0.19\textwidth]{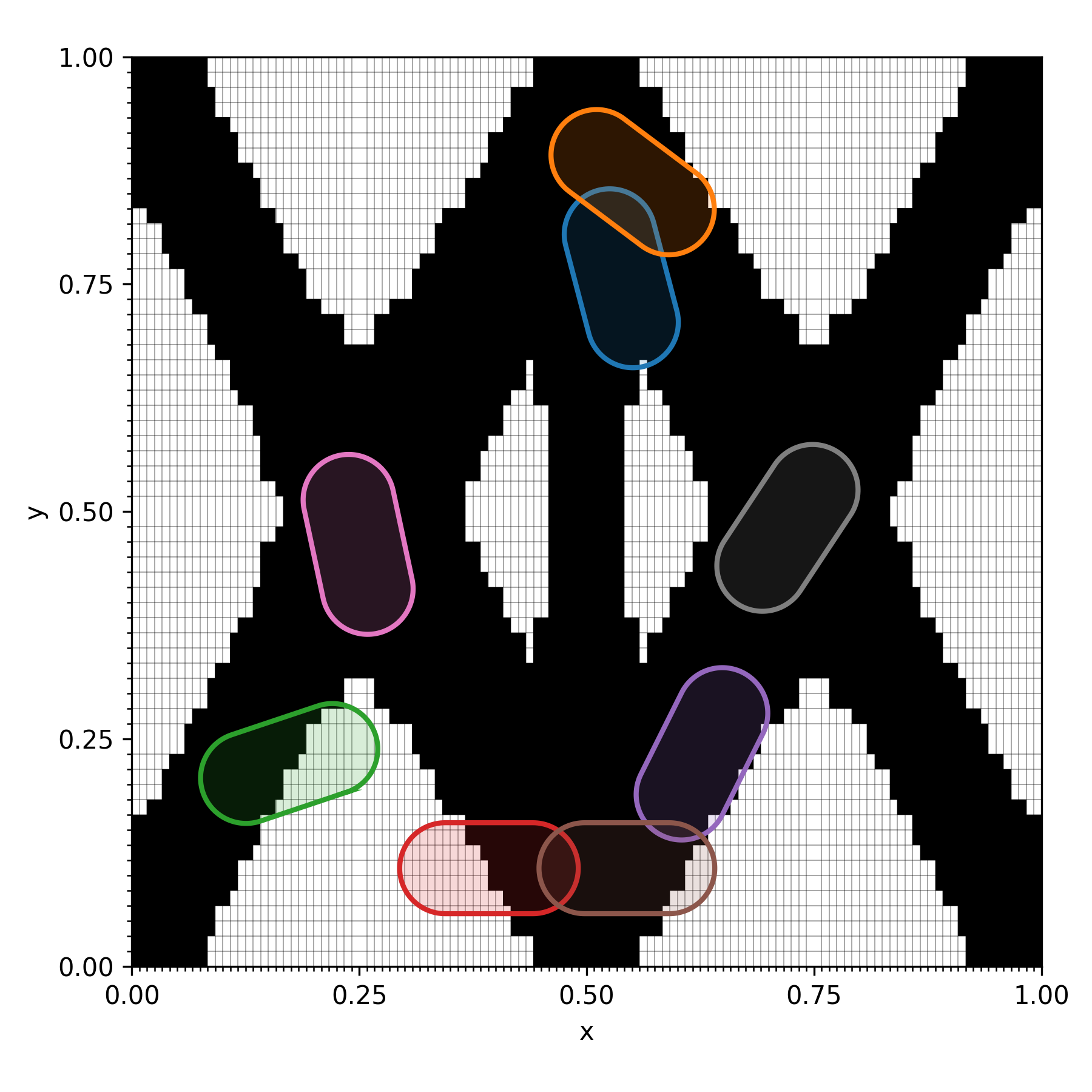} &
    \includegraphics[height= 0.095\textwidth,width=0.19\textwidth]{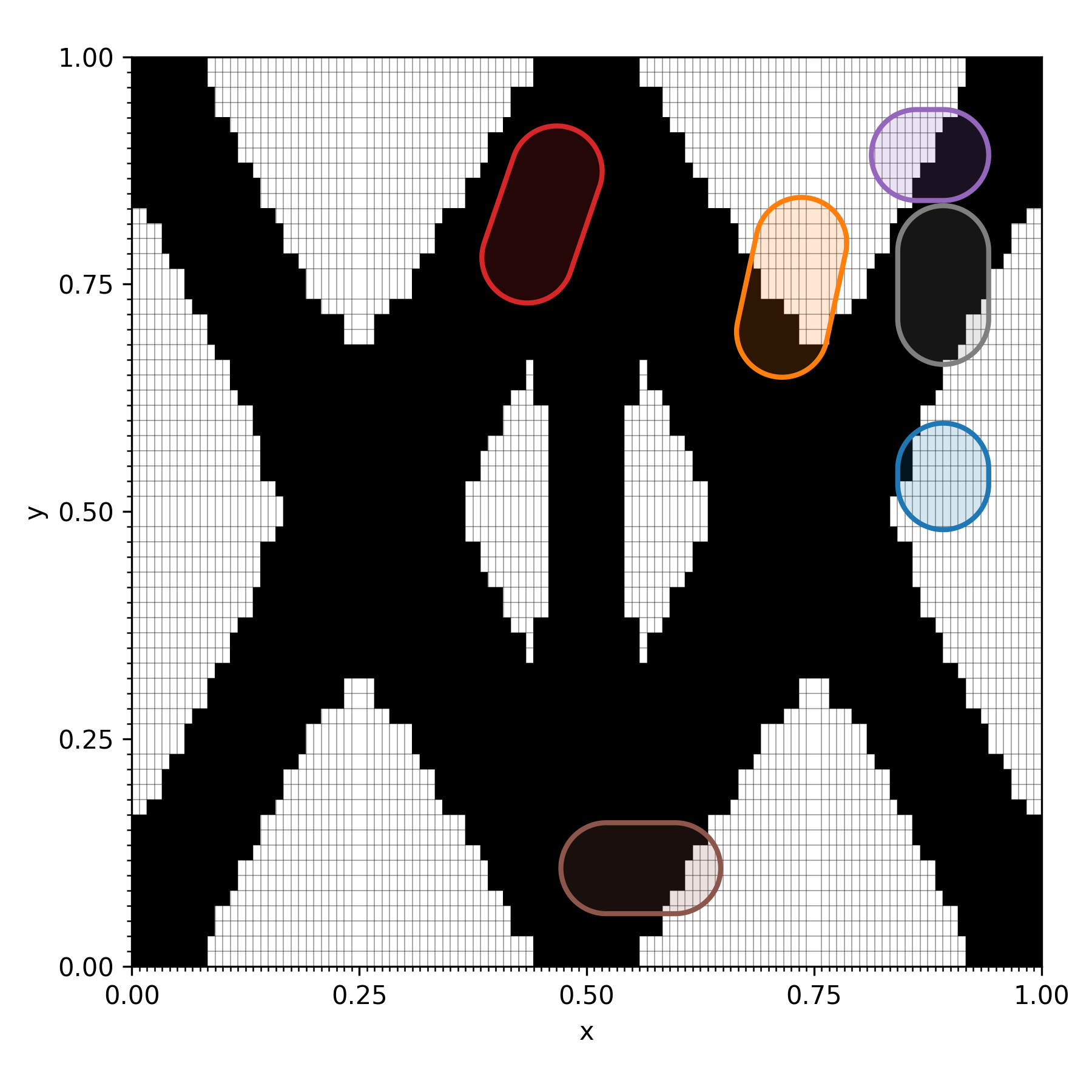} &
    \includegraphics[height= 0.095\textwidth,width=0.19\textwidth]{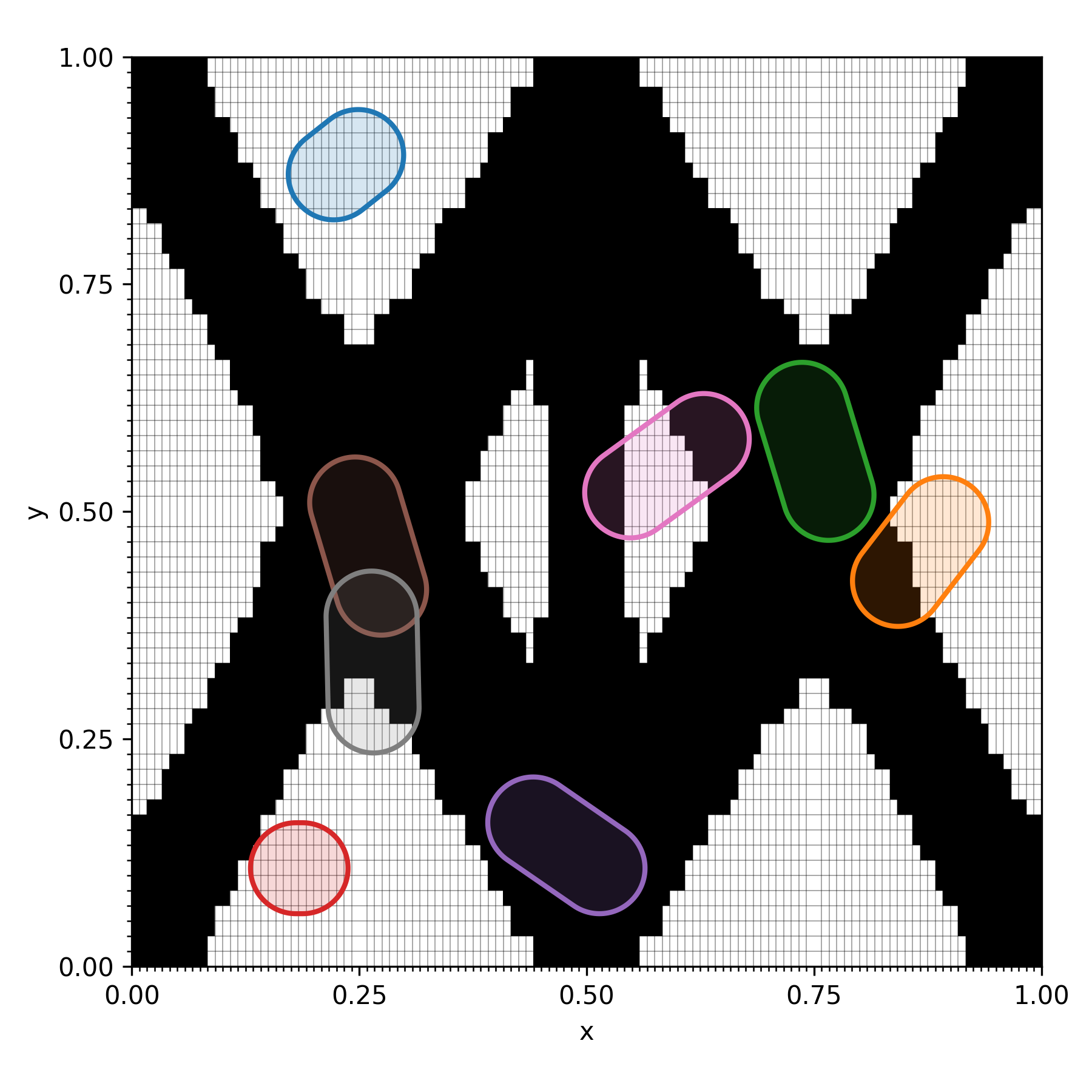} \\
    \includegraphics[height= 0.095\textwidth,width=0.19\textwidth]{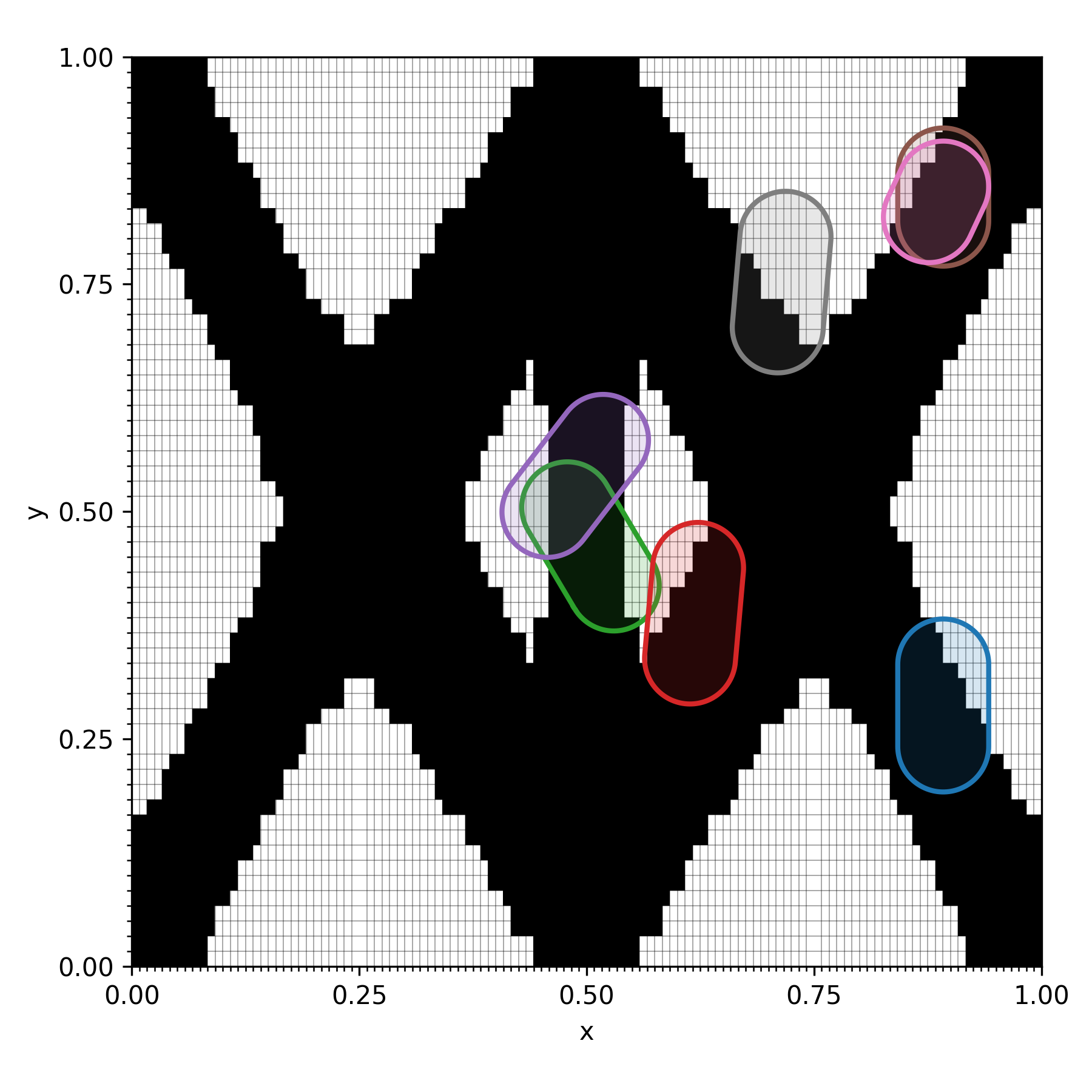} &
    \includegraphics[height= 0.095\textwidth,width=0.19\textwidth]{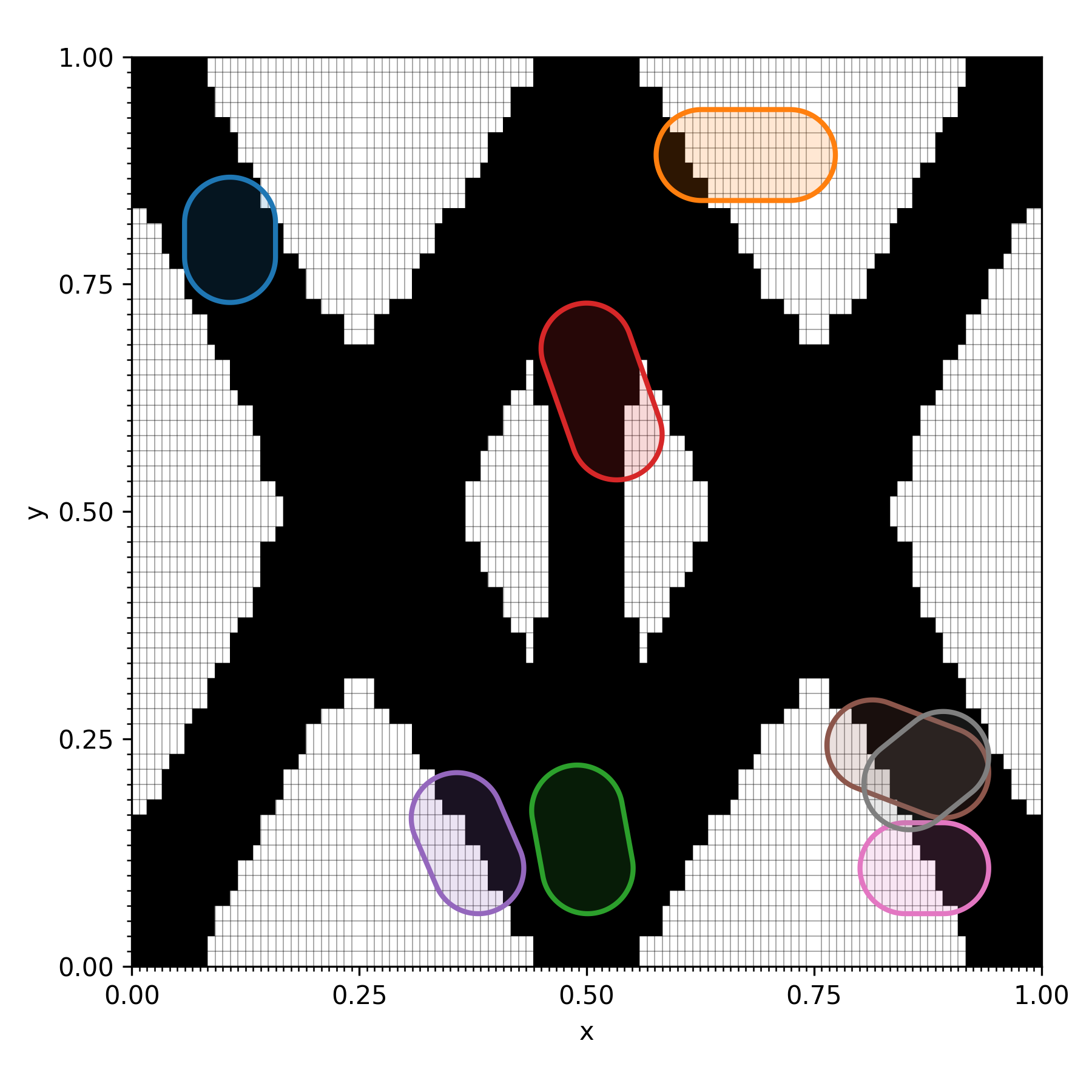} &
    \includegraphics[height= 0.095\textwidth,width=0.19\textwidth]{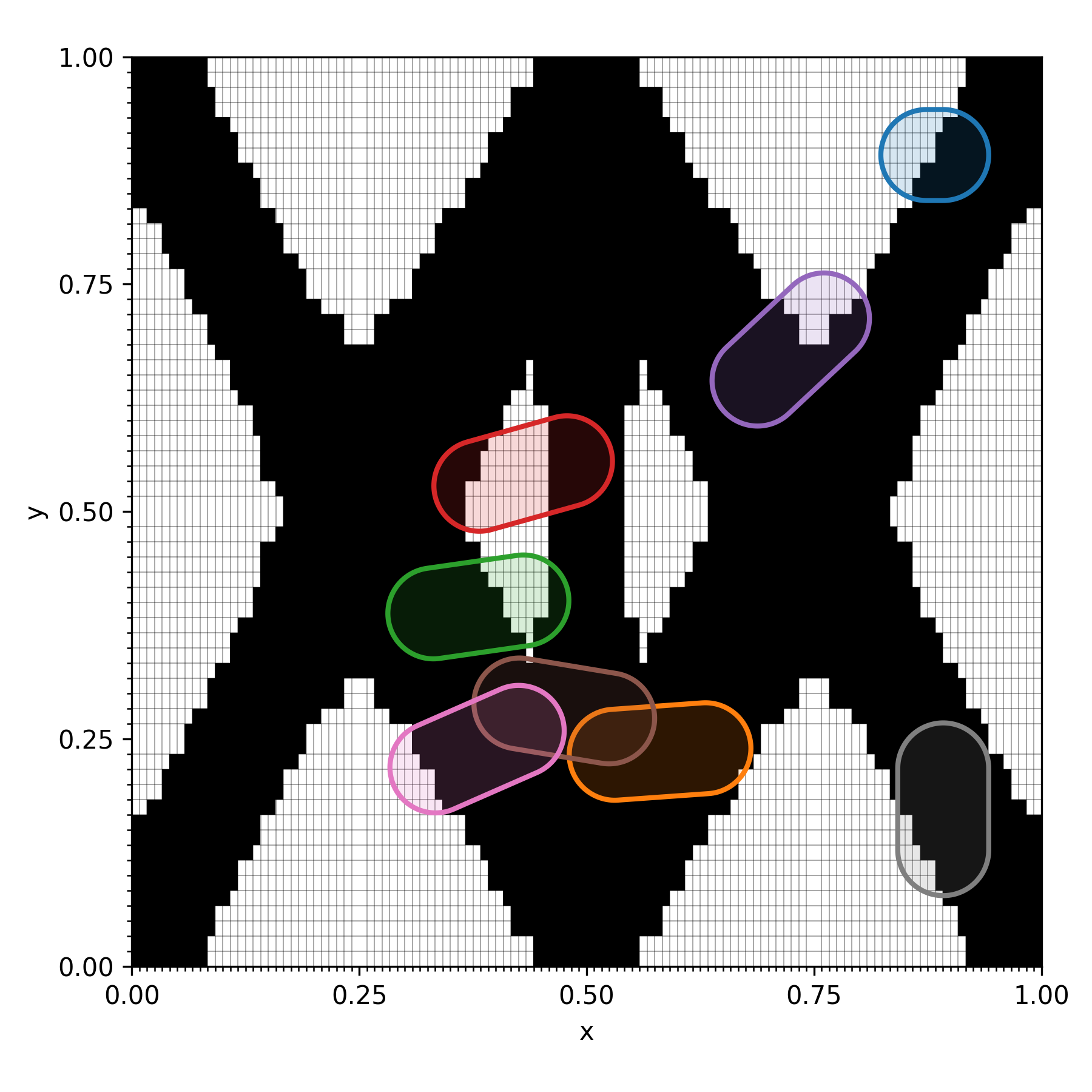} &
    \includegraphics[height= 0.095\textwidth,width=0.19\textwidth]{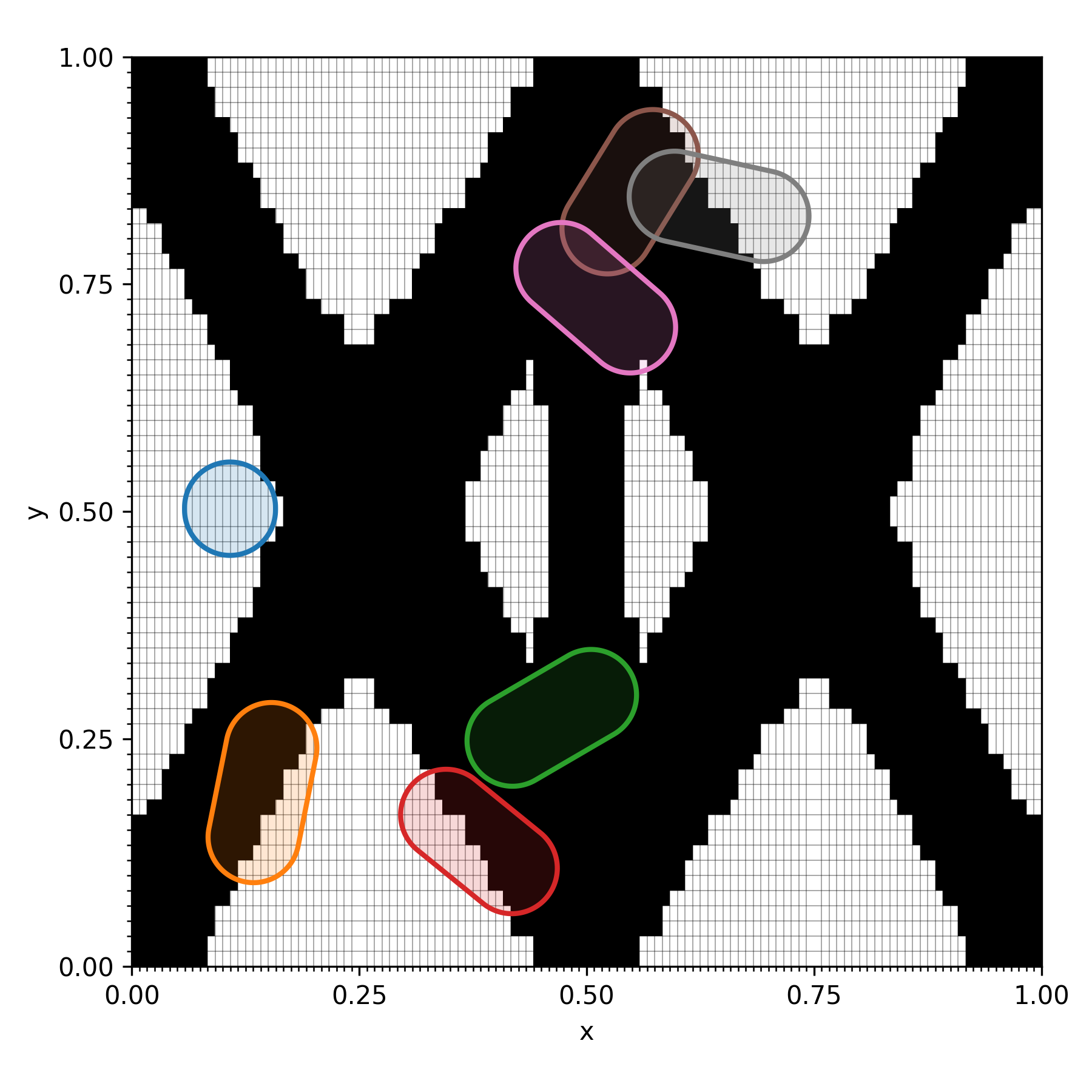} &
    \includegraphics[height= 0.095\textwidth,width=0.19\textwidth]{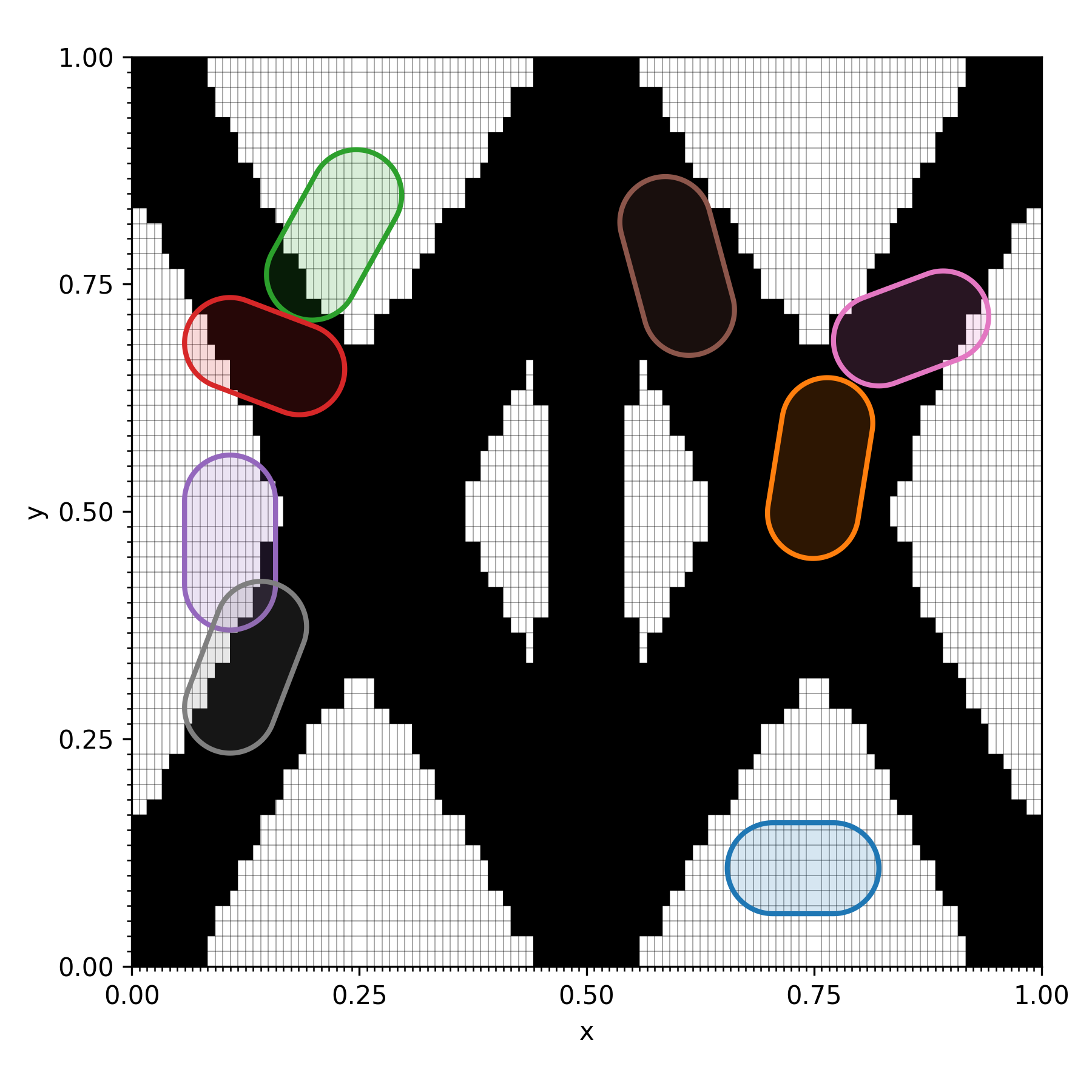} \\
  \end{tabular}
  \caption{Random initial placements used for comparison with the structured 
  cross initialization. Each panel shows the starting configuration of pills 
  before optimization.}
  \label{fig:random-inits}
\end{figure}

The influence of initialization on the reconstruction process was examined
by comparing two strategies: a structured \emph{cross} configuration and
ten independently generated \emph{random} placements. 
The cross initialization provides a deterministic and symmetric starting point, 
ensuring that all four quadrants of the domain are initially covered. 
In contrast, the random seeds represent a wide spectrum of possible pill 
arrangements and serve as a benchmark for assessing both the robustness 
and the potential limitations of the cross setup. 
By contrasting a fixed, symmetry--biased initialization with a statistical 
ensemble of randomized alternatives, the relative advantages and drawbacks 
of each approach can be identified. 
The ten random seeds used as baselines are shown in \cref{fig:random-inits}.

Qualitative inspection of the final layouts reveals a recurring pattern. 
The cross initialization reliably avoids severe undercoverage in the early stages 
and enforces a balanced distribution of pills across the domain. 
Its inherent symmetry bias, however, steers the search toward cross-shaped configurations. 
Features initialized at symmetric endpoints often interfere with one another, 
especially during the early stages, which restricts their ability to slide past each other 
and to adapt flexibly to the target. 
In some cases the resulting cross-shaped arrangement persists as a local minimum, 
so that pills remain in a crossed configuration even when this is suboptimal. 
This coupling leads to persistent overlap and reduces flexibility 
in regions where finer alignment would be required.  

By contrast, most random initializations circumvent this constraint 
at an early stage. Features are able to attach directly to distinct parts 
of the target structure, with less mutual blocking and fewer redundant overlaps. 
As a result, narrow bars and curved sections are often captured with greater 
precision and spurious extensions into void regions are reduced. 
Among the ten random trials, three runs are clearly inferior, 
as evident from the large coherent patches in both residual 
and absolute--residual fields. 
The remaining seven runs, however, achieve substantially sharper coverage 
of complex intersections than the cross baseline, 
demonstrating that random seeding can enable more flexible adaptation 
to the target geometry.

\begin{figure}[H]
  \centering
  \includegraphics[width=0.79\textwidth]{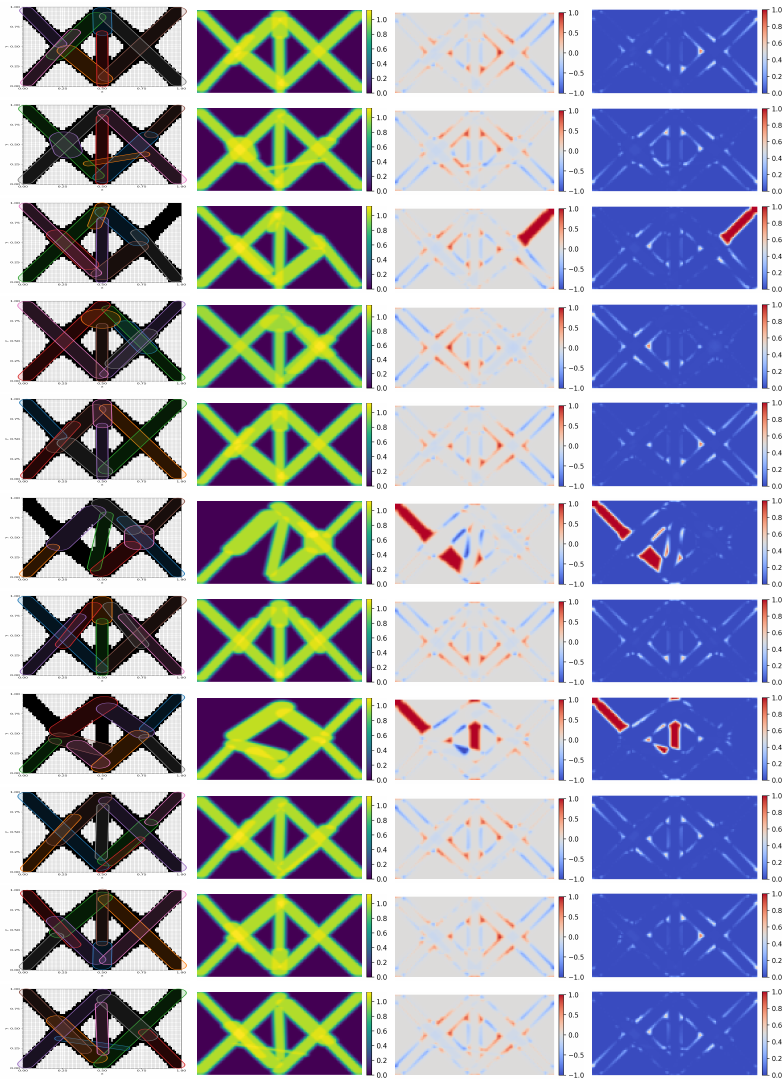}
  \caption{Final reconstructions after optimization. Each row corresponds to one initialization 
  (top ten: random seeds; bottom row: cross initialization). Columns show pill placement, 
  density field, residual and absolute residual.}
  \label{fig:rand-vs-cross}
\end{figure}
\begin{figure}[H]
  \centering
  \includegraphics[width=0.8\textwidth]{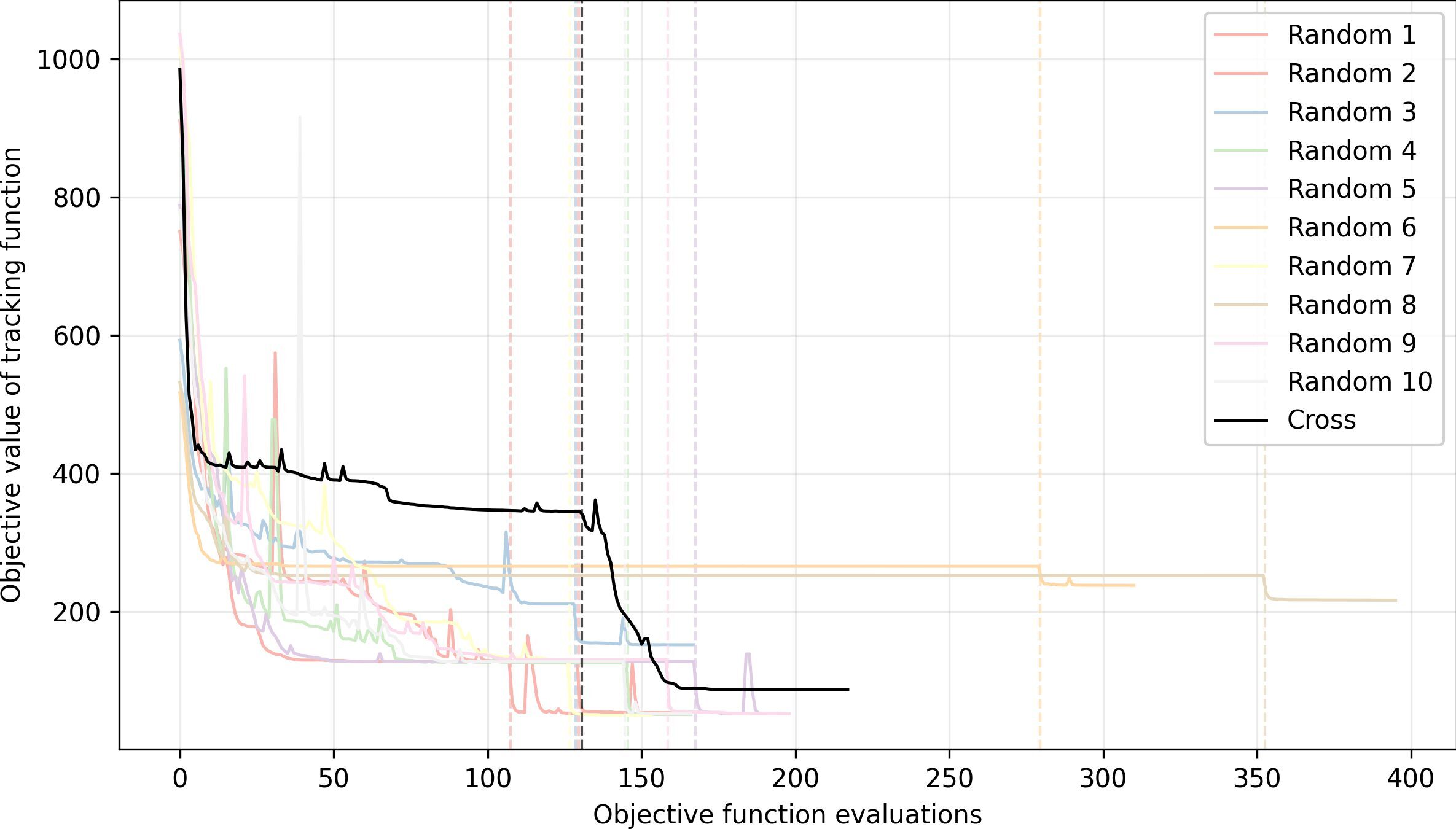}
  \caption{Objective trajectories for cross (black) and ten random
  initializations (colored). Vertical dashed lines mark the transition
  from bridging to convergence. Several random runs finish below the cross
  baseline, while three remain clearly worse, consistent with the residual 
  patterns in \cref{fig:rand-vs-cross}.}
  \label{fig:rand-vs-cross-obj}
\end{figure}

The objective trajectories substantiate these observations
(\cref{fig:rand-vs-cross-obj}). 
Several random runs converge to lower final objectives than the cross baseline, 
while three plateau at visibly higher values, reflecting the large residual patches 
observed in \cref{fig:rand-vs-cross}. 
The cross trajectory shows steady and reliable descent but does not attain the best outcomes; 
its symmetry bias appears to trap pills in mild local minima 
where overlapping endpoints hinder reorganization. 
In contrast, the better random trials exhibit marked improvements 
after the transition from bridging to convergence, 
consistent with pills having already specialized on distinct target parts during bridging. 

Overall, random seeding offers greater potential for high reconstruction quality, 
but its effectiveness depends on maintaining a reasonably balanced initial distribution. 
If pills are placed too unevenly, entire regions of the domain may remain uncovered, 
which limits reconstruction from the outset. 
A practical strategy is therefore to employ random seeding while ensuring 
that the initial placements are distributed relatively uniformly across the domain. 
This combines the higher potential of random initialization with safeguards 
against severe undercoverage. 

% ============================================================
% 3) Improve Quality
% ============================================================

\subsection{Higher resolution}
\label{sec:quality-resolution}

The influence of finite element resolution on reconstruction quality 
is assessed by performing otherwise identical optimizations on three meshes 
of different sizes: $40\times 20$, $80\times 40$ and $120\times 60$ elements. 
The target density field is the five-bar structure from \cref{fig:target-5bar}, 
and the initialization and parameter settings from \cref{sec:global-config} are used throughout. 

\begin{figure}[H]
  \centering
  \setlength{\tabcolsep}{3pt}
  \renewcommand{\arraystretch}{0.1}

  % nur innerhalb dieses Blocks: includegraphics vertikal zentrieren
  \begingroup
  \let\oldincludegraphics\includegraphics
  \renewcommand{\includegraphics}[2][]{\raisebox{-.5\height}{\oldincludegraphics[#1]{#2}}}

  \begin{tabular}{ccc}
    \includegraphics[width=0.32\textwidth,height=0.16\textwidth]{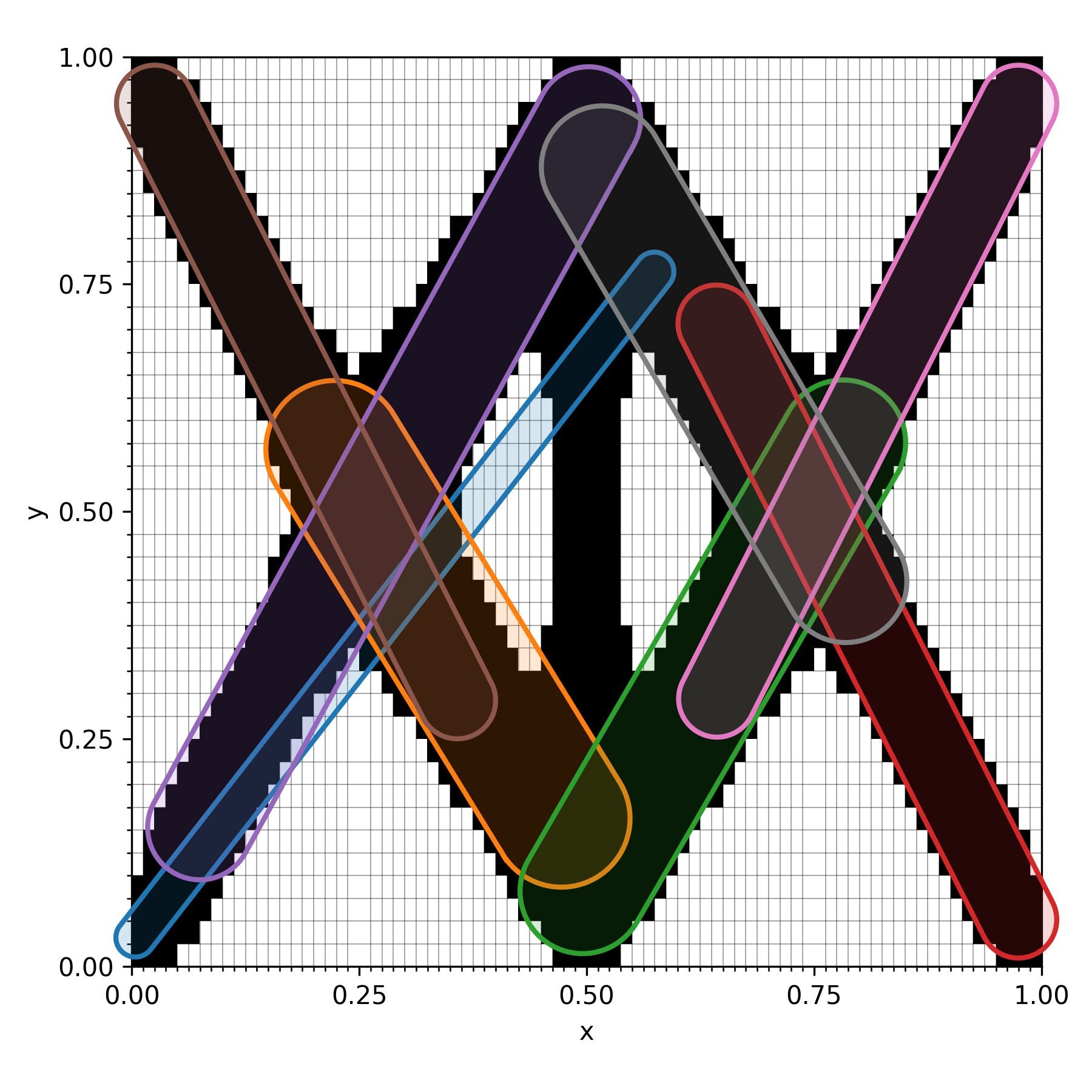} &
    \includegraphics[width=0.32\textwidth,height=0.16\textwidth]{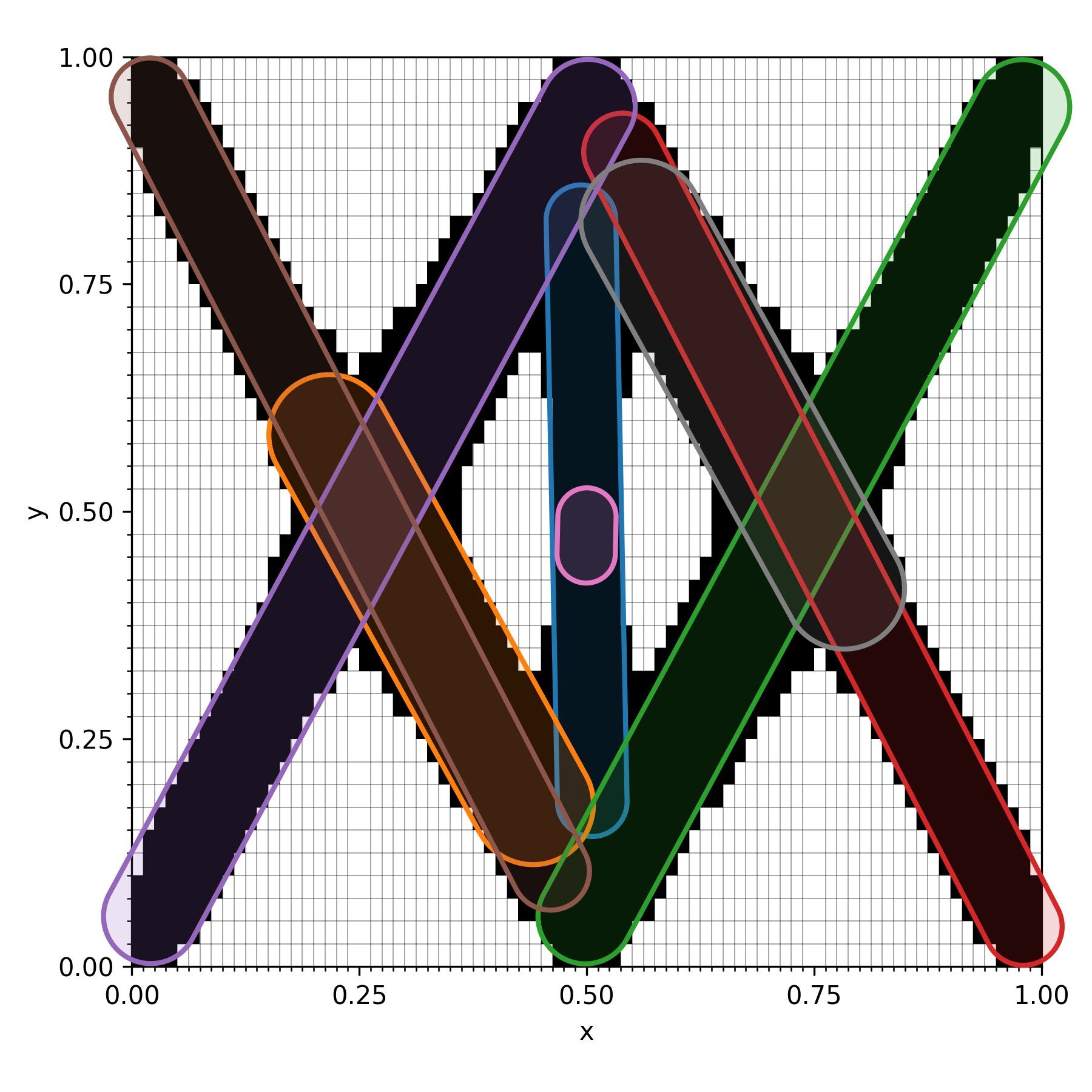} &
    \includegraphics[width=0.32\textwidth,height=0.16\textwidth]{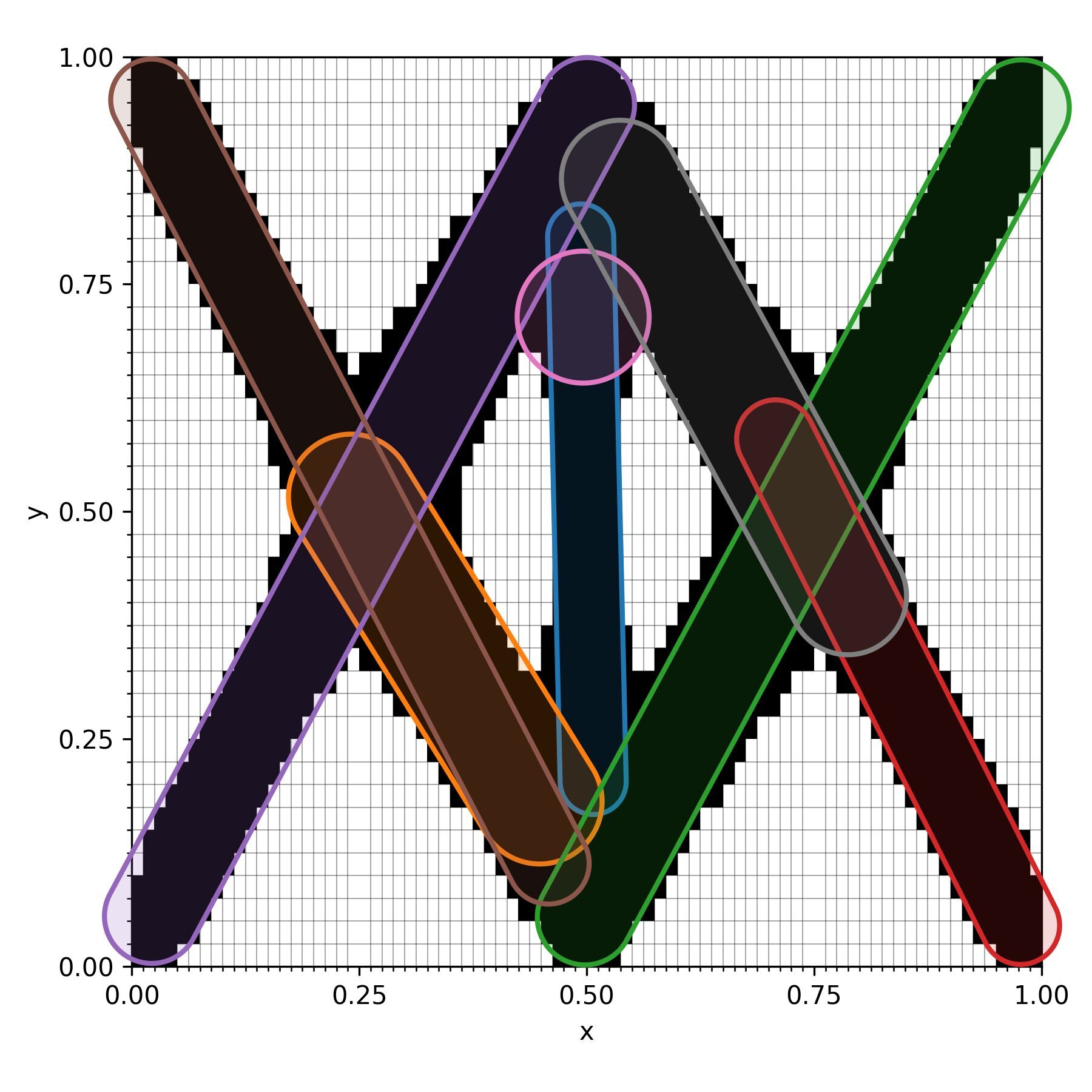} \\
    \scriptsize $40\times20$ & \scriptsize $80\times40$ & \scriptsize $120\times60$ \\
  \end{tabular}

  \endgroup

  \caption{Final pill placements at different finite element resolutions. 
  Higher resolution yields more accurate approximations of the five-bar structure.}
  \label{fig:res-pill}
\end{figure}

\cref{fig:res-pill} shows the final pill placements for each resolution. 
On the coarsest grid ($40\times 20$), the pills provide only a rough approximation of the target. 
The central vertical bar is not properly captured and several gaps remain in the coverage. 
With $80\times 40$ elements, the approximation improves substantially: 
the main bars are clearly represented, although residual deviations persist near the center. 
At the highest resolution ($120\times 60$), all five bars are well reproduced, 
and coverage is nearly complete.

\begin{figure}[H]
    \centering
    \begin{subfigure}{0.32\textwidth}
        \includegraphics[width=\linewidth]{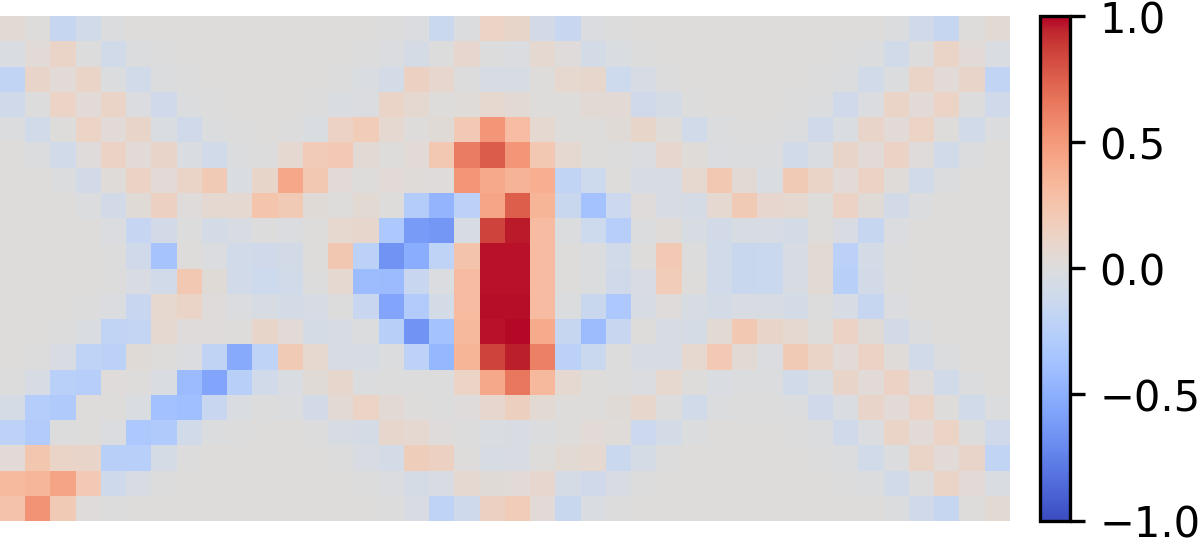}
        \caption{$40\times20$}
    \end{subfigure}
    \hfill
    \begin{subfigure}{0.32\textwidth}
        \includegraphics[width=\linewidth]{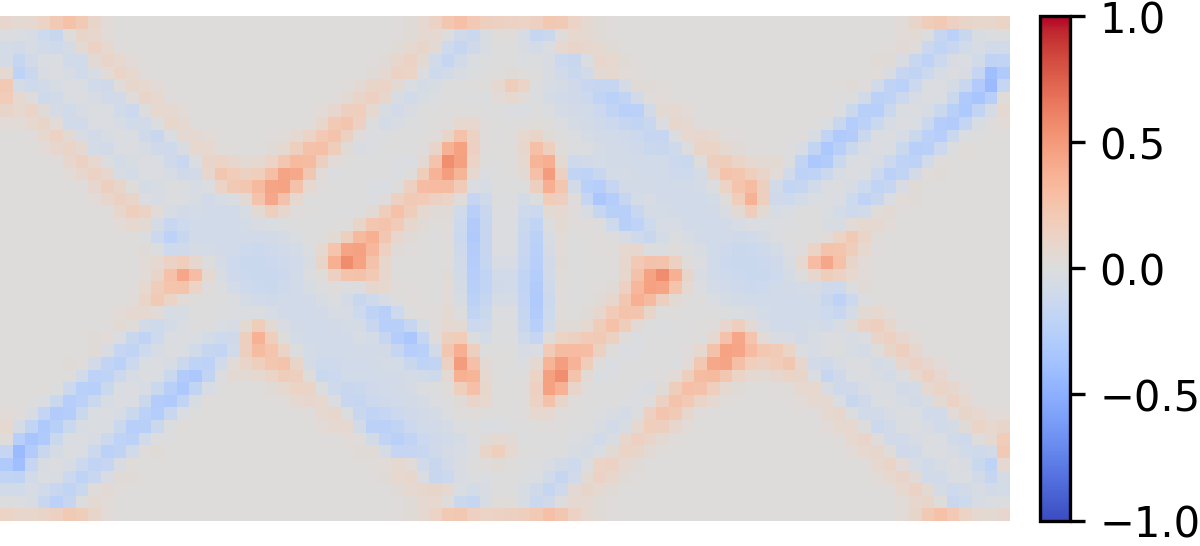}
        \caption{$80\times40$}
    \end{subfigure}
    \hfill
    \begin{subfigure}{0.32\textwidth}
        \includegraphics[width=\linewidth]{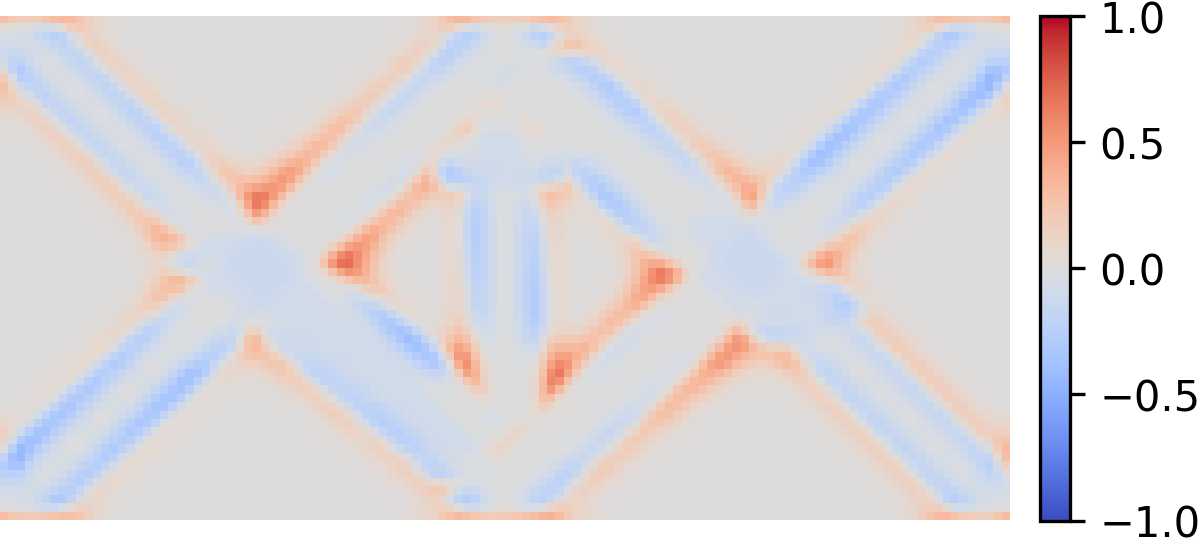}
        \caption{$120\times60$}
    \end{subfigure}
    \caption{Residual fields at different resolutions. 
    Positive values indicate uncovered target regions, negative values excess density.}
    \label{fig:res-residual}
\end{figure}

\begin{figure}[H]
    \centering
    \begin{subfigure}{0.32\textwidth}
        \includegraphics[width=\linewidth]{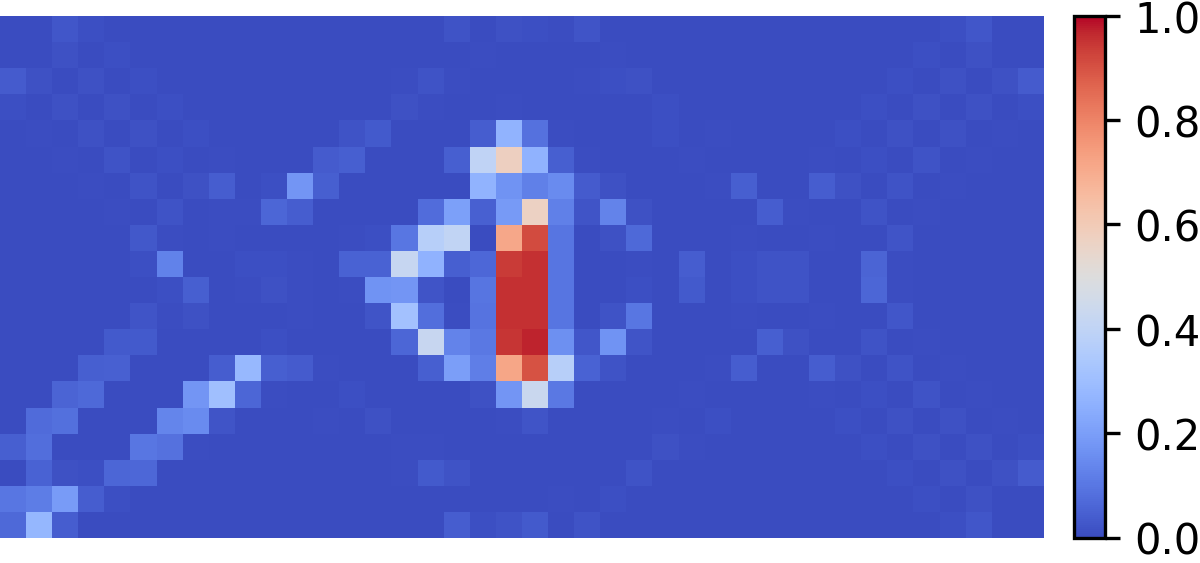}
        \caption{$40\times20$}
    \end{subfigure}
    \hfill
    \begin{subfigure}{0.32\textwidth}
        \includegraphics[width=\linewidth]{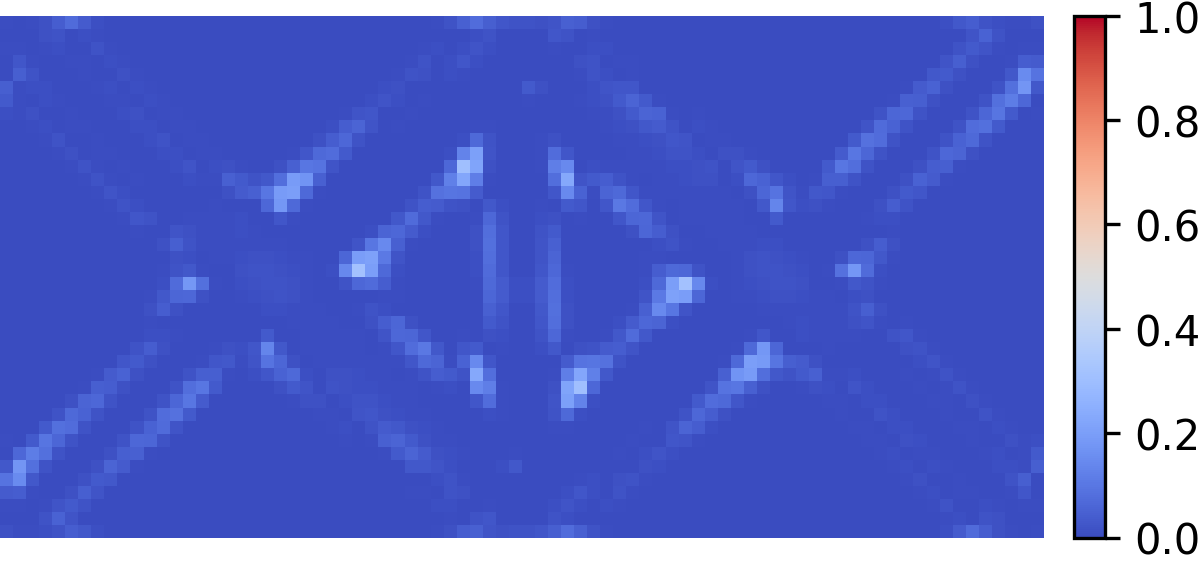}
        \caption{$80\times40$}
    \end{subfigure}
    \hfill
    \begin{subfigure}{0.32\textwidth}
        \includegraphics[width=\linewidth]{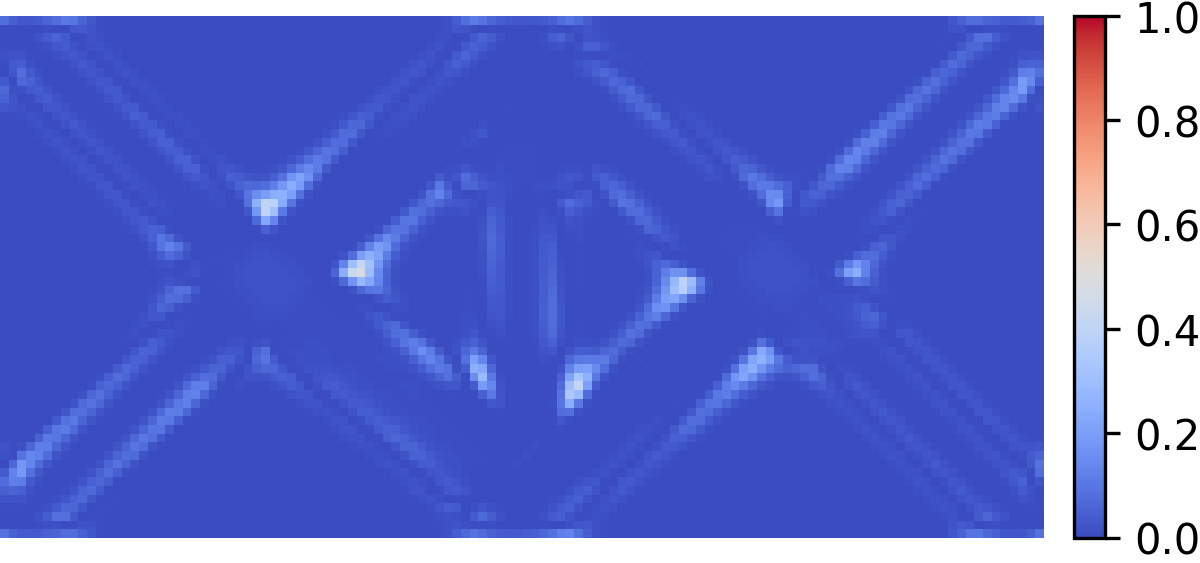}
        \caption{$120\times60$}
    \end{subfigure}
    \caption{Absolute residual fields at different resolutions. 
    Deviations diminish systematically as resolution increases.}
    \label{fig:res-residualabs}
\end{figure}

The residual fields in \cref{fig:res-residual,fig:res-residualabs} confirm this trend. 
At low resolution, large regions of missing density appear, particularly in the central bar. 
With increasing resolution, these discrepancies diminish, leaving only narrow bands of error near bar boundaries. 
The absolute residual fields make this effect especially clear: 
red regions indicating excess density and blue regions indicating uncovered target 
decrease markedly from $40\times 20$ to $120\times 60$.

\begin{figure}[H]
    \centering
    \includegraphics[width=0.75\textwidth]{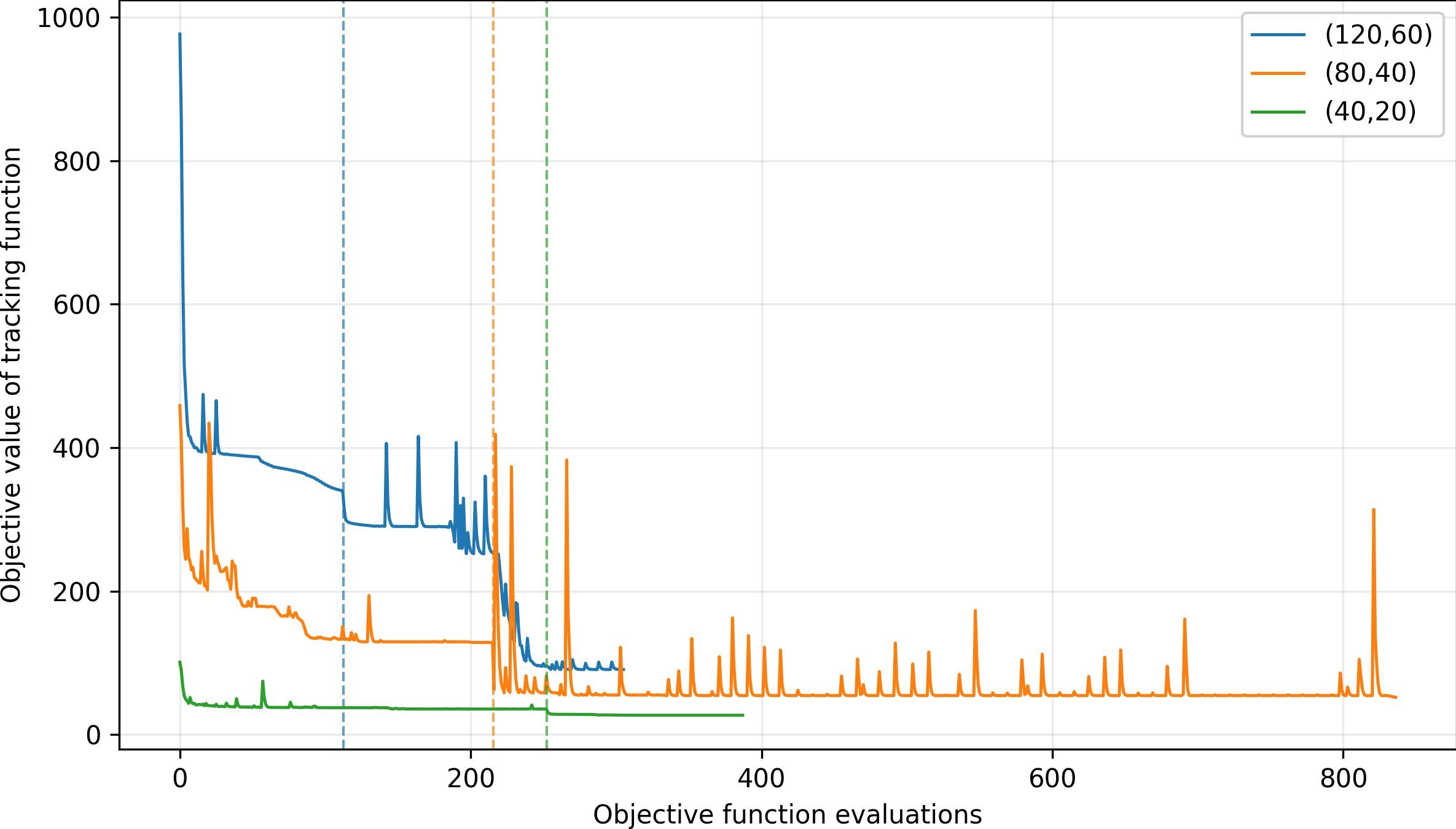}
    \caption{Objective trajectories at different resolutions. Raw objective
    values increase with element count, but normalized values reveal reduced
    average error at higher resolutions.}
    \label{fig:res-obj}
\end{figure}

The optimization traces in \cref{fig:res-obj} show that absolute objective values cannot be compared directly across resolutions, 
since the functional in \cref{sec:tracking_objective_sec} sums errors over all elements. 
To account for this, the objective $F$ is normalized by the number of elements, yielding the mean error per element,
\[
F_{\mathrm{norm}} = \frac{F}{N_xN_y}.
\]
The final values are
\[
\begin{aligned}
40\times20: &\quad \tfrac{2.7023\times10^{1}}{800} \;\approx\; 3.38\times 10^{-2},\\
80\times40: &\quad \tfrac{5.2002\times10^{1}}{3200} \;\approx\; 1.63\times 10^{-2},\\
120\times60: &\quad \tfrac{9.0709\times10^{1}}{7200} \;\approx\; 1.26\times 10^{-2}.
\end{aligned}
\]

This normalization confirms the qualitative observations from the pill placements and residual fields. 
While the raw objective values increase with the number of elements, 
the average error per element decreases, demonstrating that higher resolution enables more accurate reconstructions.

The trajectories in \cref{fig:res-obj} are plotted against the number of objective function evaluations. 
Vertical dashed lines mark the transition from the bridging to the convergence stage. 
Although all runs used the same iteration budget (100 IPOPT iterations), 
the abscissa positions of these markers differ considerably 
because one iteration may involve a variable number of evaluations due to line-search backtracking and globalization.

Among the three resolutions, the $120\times60$ case exhibits the smoothest descent, 
with long sequences of accepted steps and comparatively few backtracking phases. 
The transition to convergence therefore occurs earlier on the function-evaluation axis. 
By contrast, the $40\times20$ and $80\times40$ runs stagnate sooner, 
followed by pronounced oscillatory behavior with recurrent increases and decreases in the objective. 
This regime is characterized by frequent trial/rollback cycles in the line search, 
leading to a large number of evaluations per iteration. 
Representative IPOPT diagnostics for the $80\times40$ case report 1082 objective evaluations over 100 iterations, 
together with a steadily growing line-search counter (ls) near stagnation, 
confirming the increased globalization effort. 

Consequently, a low raw objective at a given abscissa position in \cref{fig:res-obj} 
must not be interpreted as intrinsically superior progress: 
the axis reflects evaluation effort rather than iteration count. 
Taken together with the normalized values $F_{\mathrm{norm}}=F/(N_xN_y)$, 
the traces support two conclusions. 
First, higher resolution leads to more regular and efficient convergence behavior 
(fewer backtracks and earlier stage transition in terms of function evaluations). 
Second, although absolute objectives scale with the number of elements, 
the mean error per element decreases monotonically from $40\times20$ to $120\times60$, 
in agreement with the visual improvements in the pill placements and residual fields.

\subsection{Integration schemes}
\label{sec:quality-integration}

The influence of numerical integration order on reconstruction quality 
is examined next. 
All experiments were performed on the same finite element mesh of size $80\times 40$, 
while the quadrature order was varied between one and five. 
Since the discretization is fixed, the resulting objective values 
are directly comparable across the different schemes, 
unlike in the resolution study in \cref{sec:quality-resolution}. 

\begin{figure}[H]
    \centering
    \includegraphics[width=0.32\textwidth,height=0.16\textwidth]{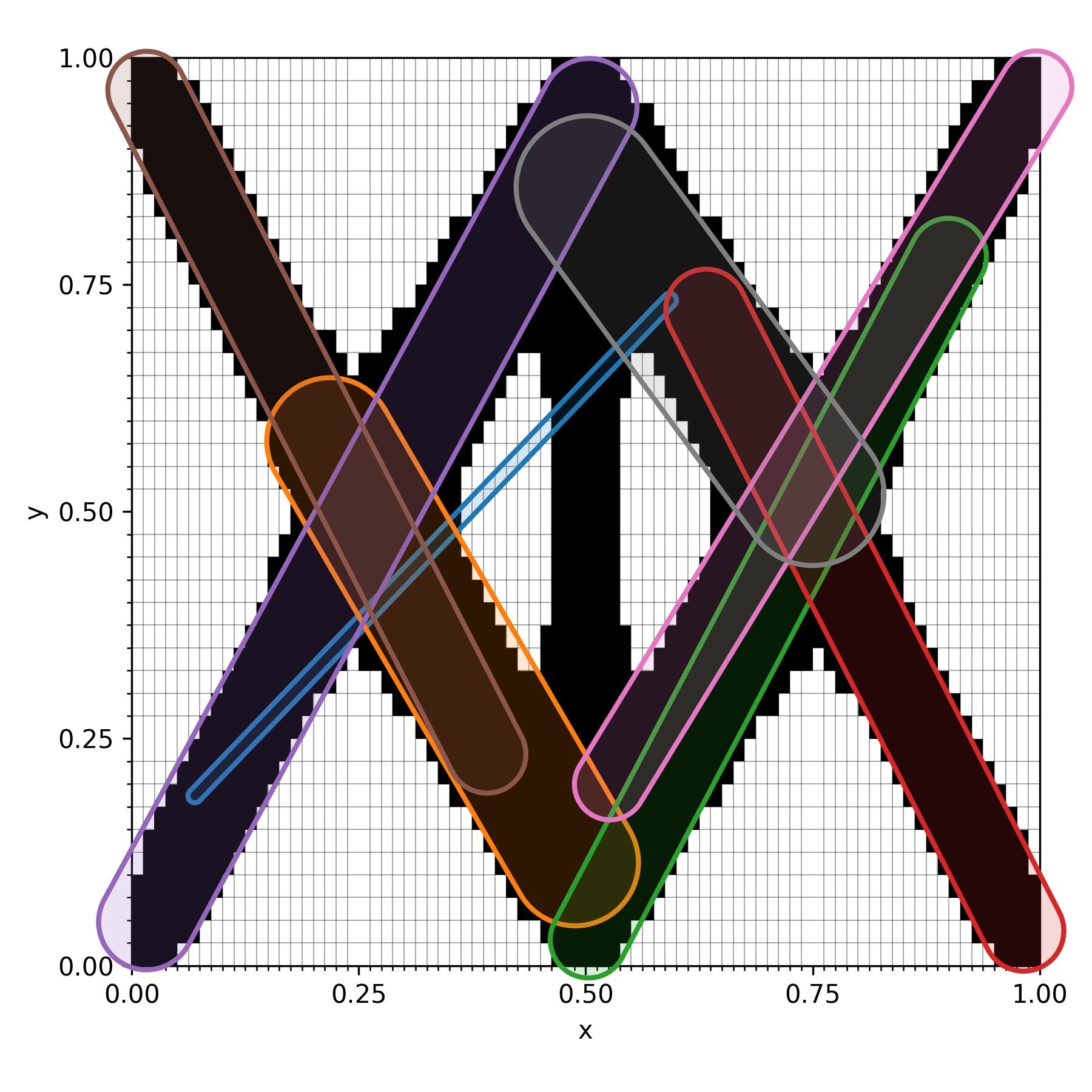}
    \includegraphics[width=0.32\textwidth,height=0.16\textwidth]{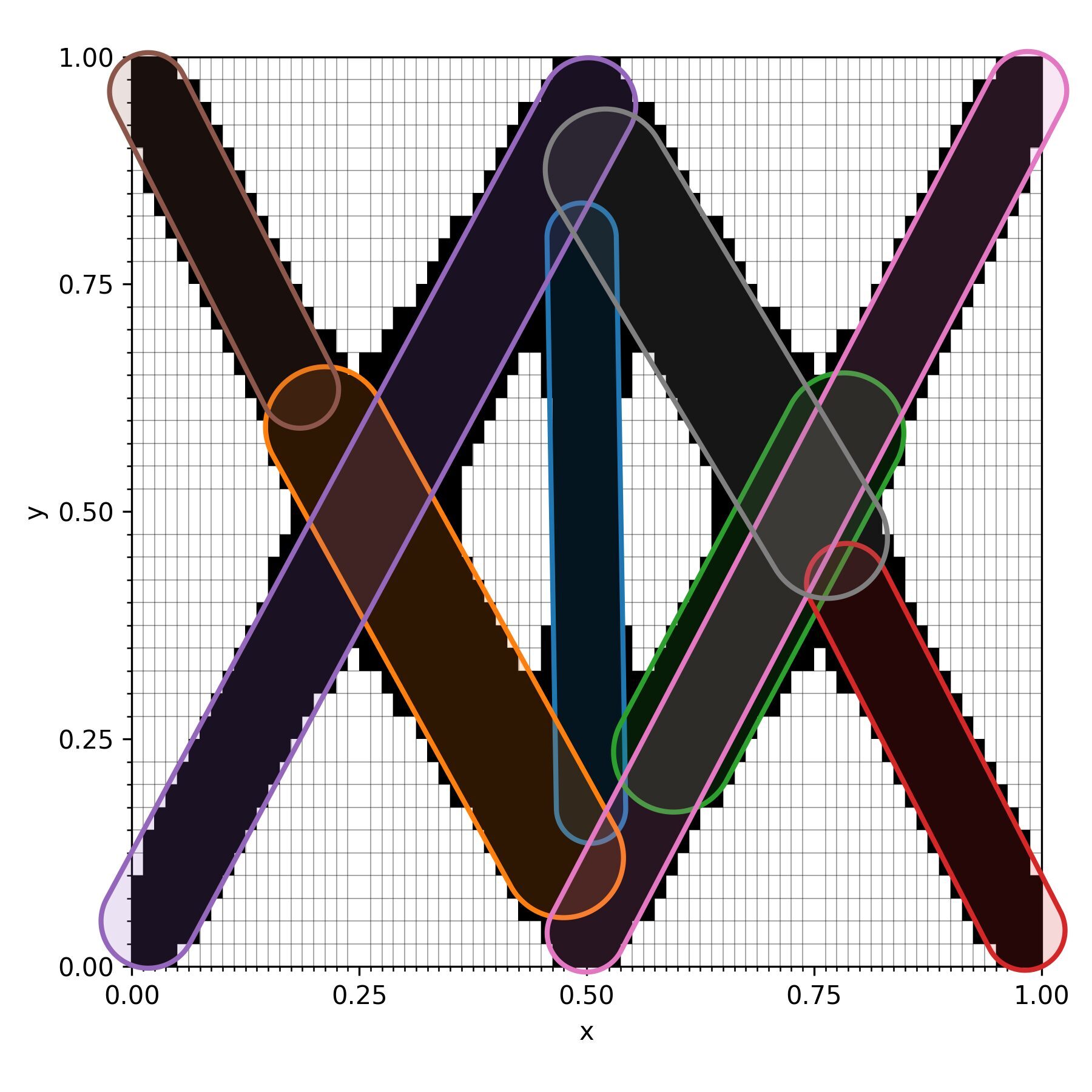}
    \includegraphics[width=0.32\textwidth,height=0.16\textwidth]{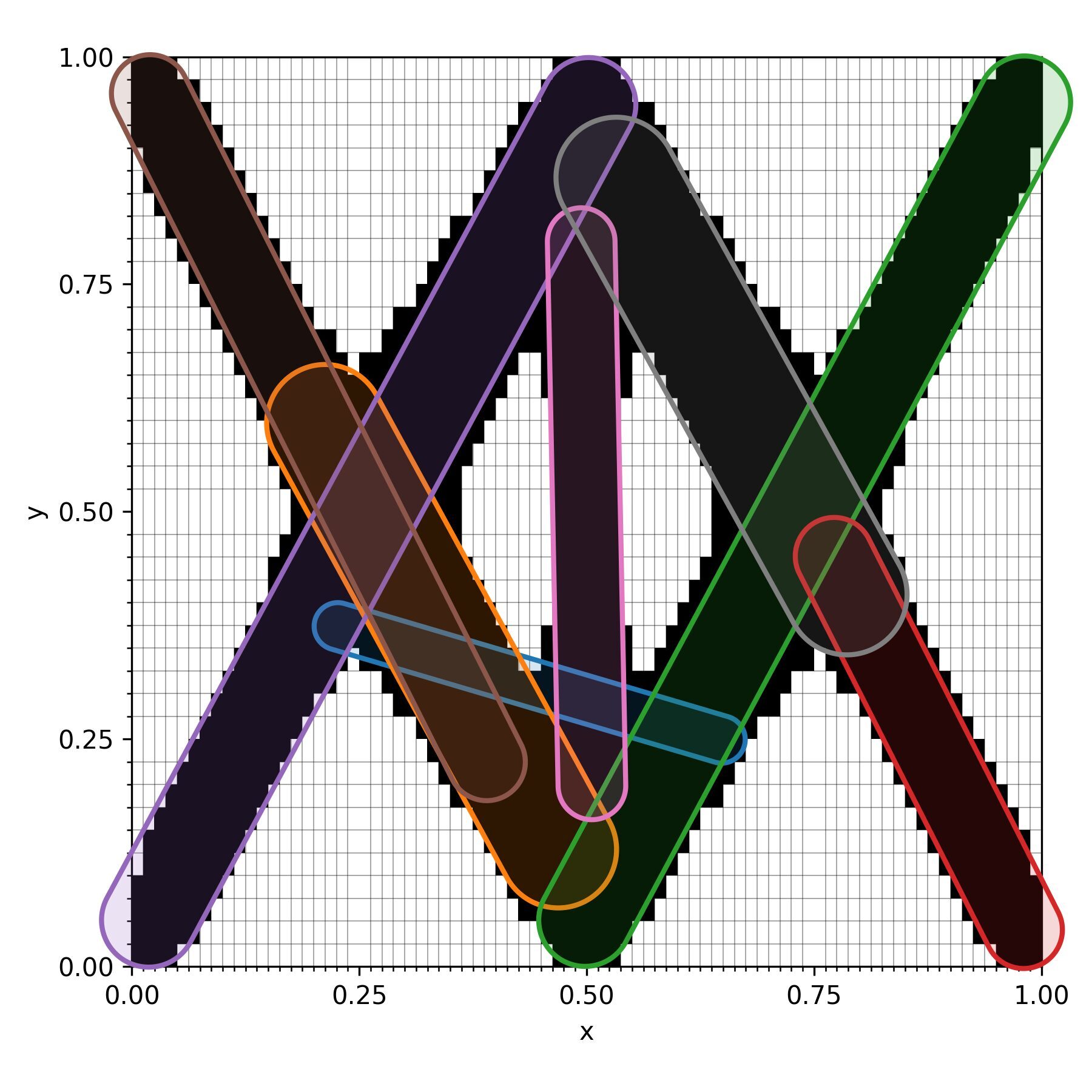}\\
    \includegraphics[width=0.32\textwidth,height=0.16\textwidth]{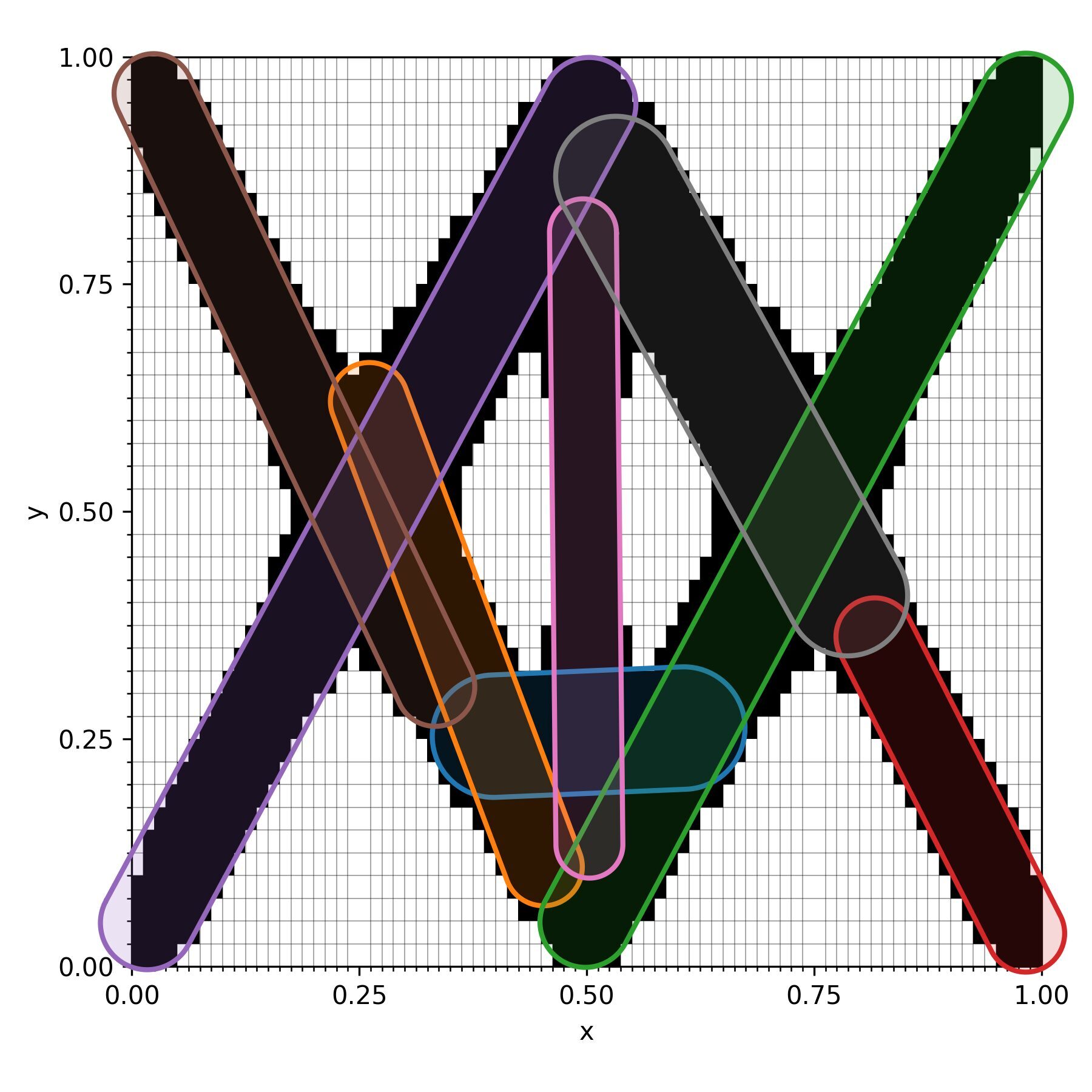}
    \includegraphics[width=0.32\textwidth,height=0.16\textwidth]{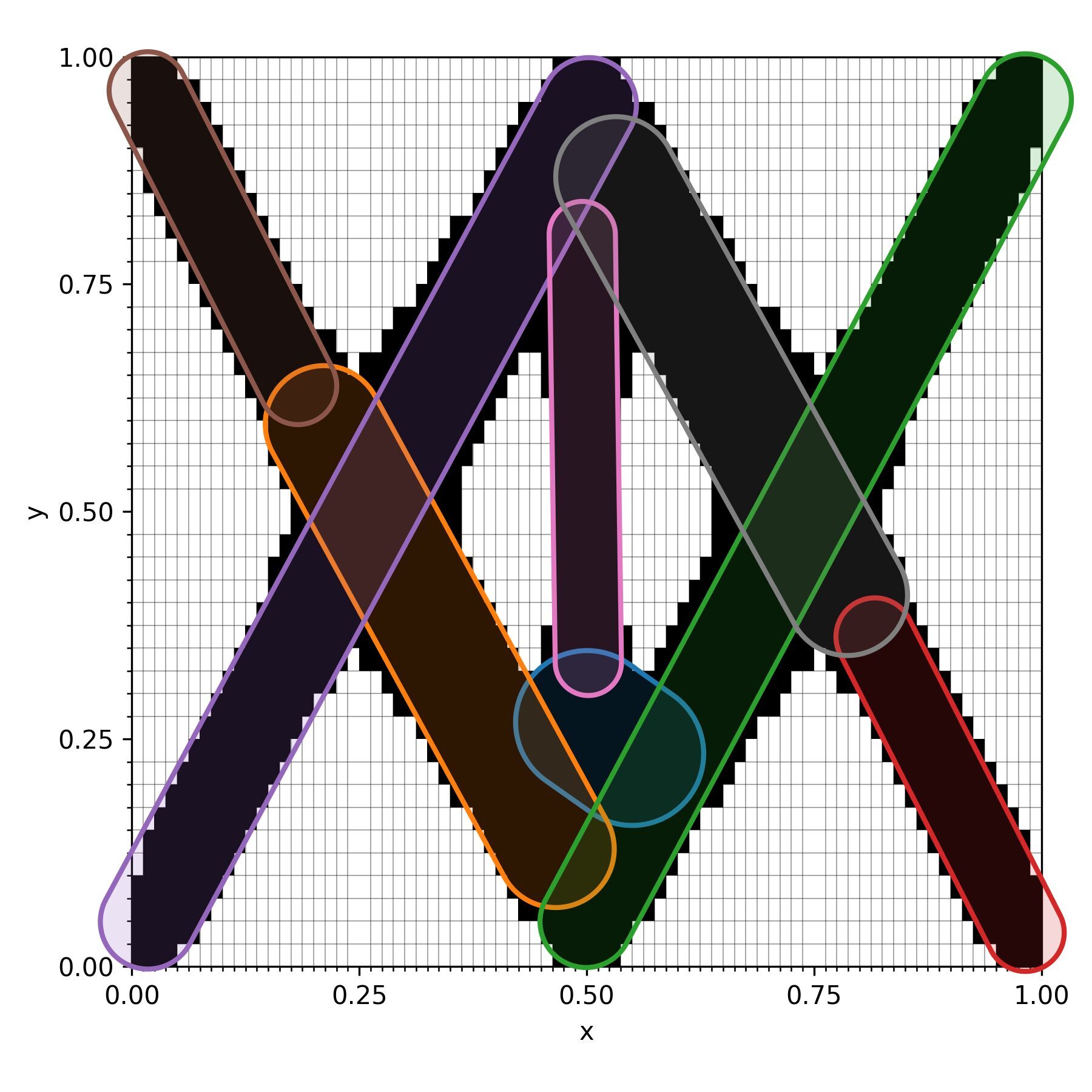}
    \caption{Pill placements for integration orders 1–5 (left to right).
    Order~1 fails to reproduce the vertical bar, 
    while orders~2–5 align well with the target structure.}
    \label{fig:order-pill}
\end{figure}

\cref{fig:order-pill} shows the final pill placements 
for all integration orders. 
With order~1, the optimizer fails to reconstruct the central vertical bar of the target structure; 
pills remain misaligned and coverage is incomplete. 
Starting from order~2, the reconstruction improves markedly: 
the vertical bar is recovered and the overall alignment of pills with the target becomes consistent. 
For orders~3–5, the pill placements are visually indistinguishable, 
indicating that quadrature accuracy beyond order~2 no longer affects the reconstructed layout in a perceptible way. 

\begin{figure}[H]
    \centering
    \includegraphics[width=0.32\textwidth]{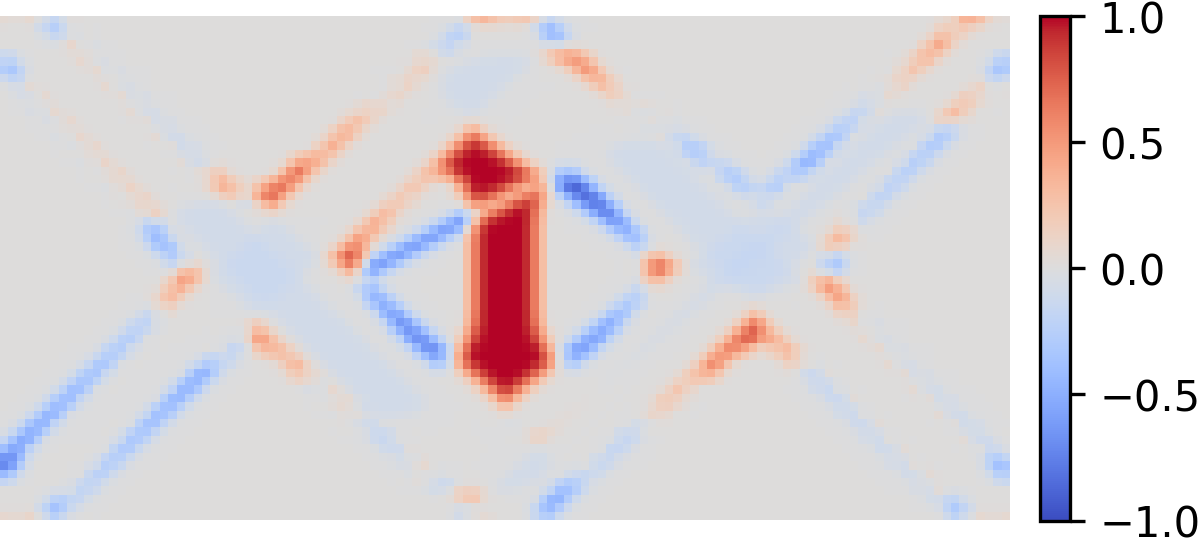}
    \includegraphics[width=0.32\textwidth]{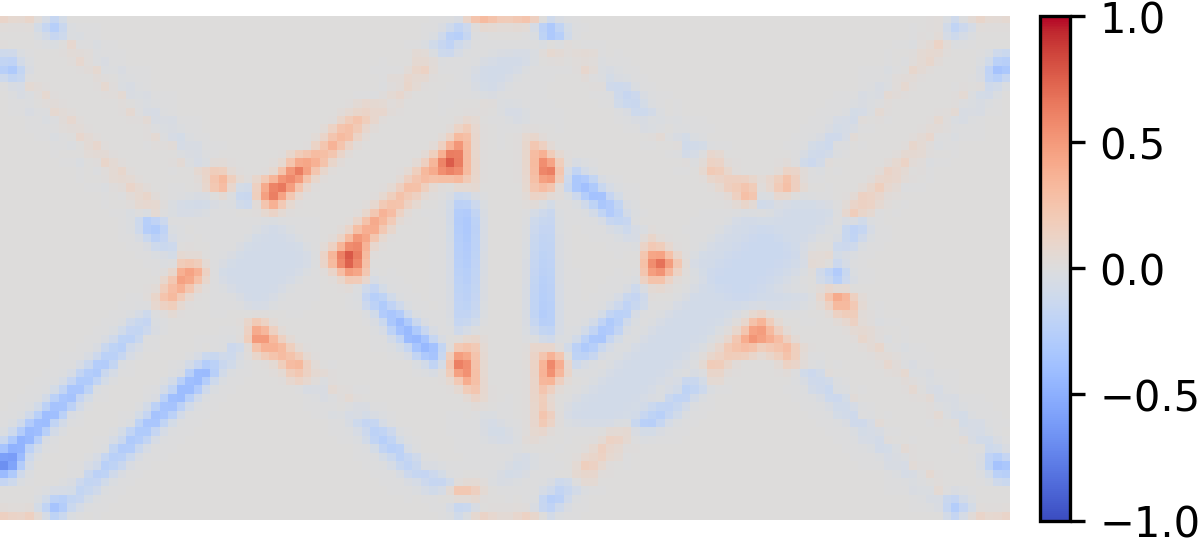}
    \includegraphics[width=0.32\textwidth]{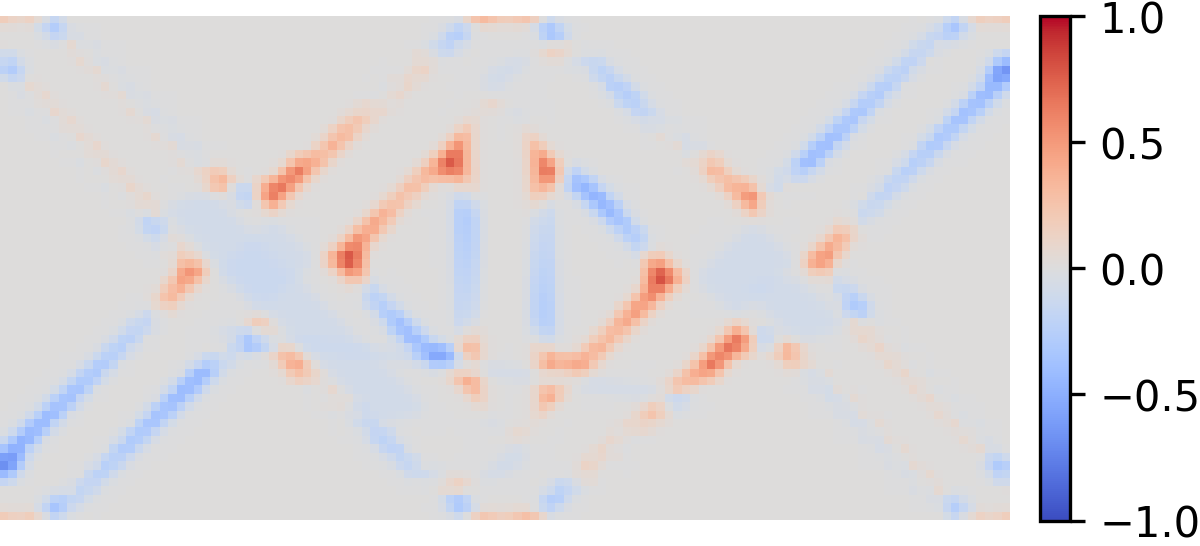}\\
    \includegraphics[width=0.32\textwidth]{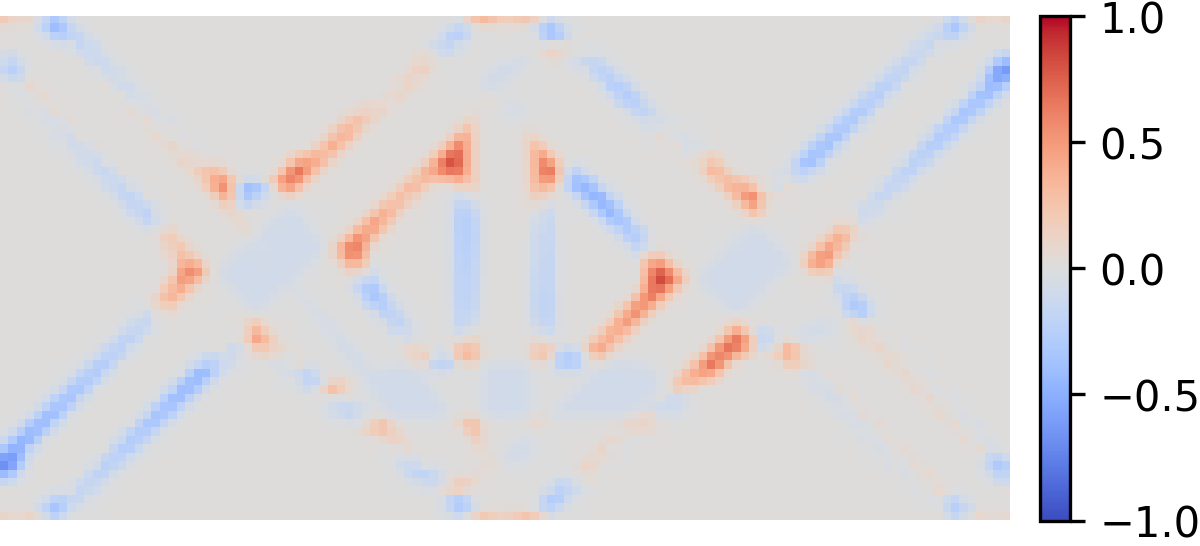}
    \includegraphics[width=0.32\textwidth]{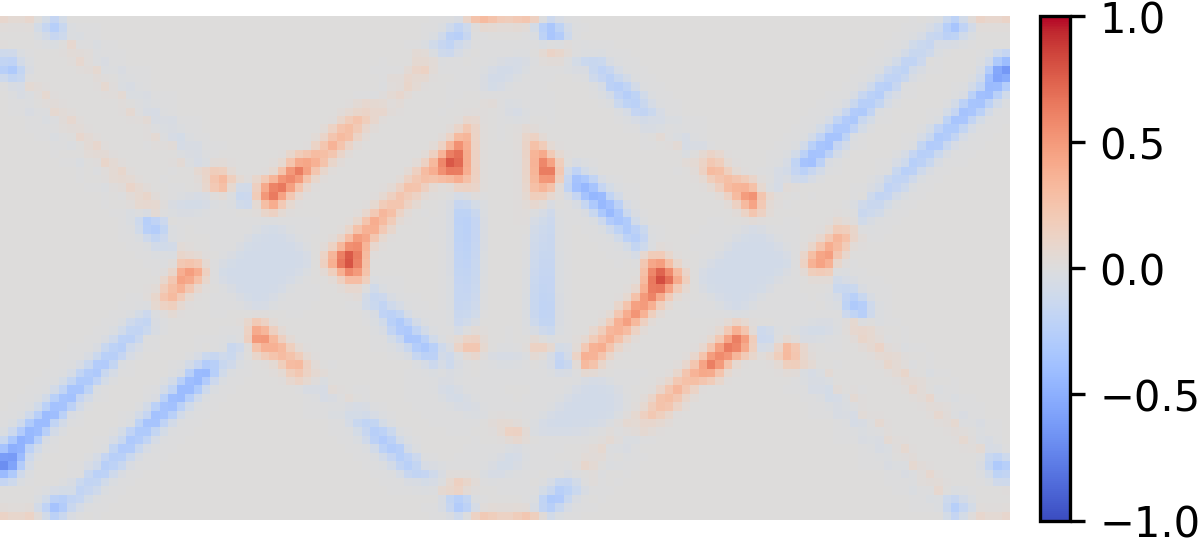}
    \caption{Residual fields for integration orders 1–5 (left to right).
    Order~1 shows large systematic errors, while orders~2–5 reduce the residuals
    to localized mismatches at intersections.}
    \label{fig:order-residual}
\end{figure}

Residual fields for the five integration orders are displayed in \cref{fig:order-residual}. 
The order~1 residual map shows large coherent deviations in the central region, 
reflecting the missing vertical bar. 
From order~2 onward, the residuals shrink substantially 
and are confined mainly to intersections and transition zones. 
The residual patterns for orders~2–5 are similar, 
with only marginal differences in magnitude, 
and no systematic improvements are observed beyond order~3. 

\begin{figure}[H]
    \centering
    \includegraphics[width=0.32\textwidth]{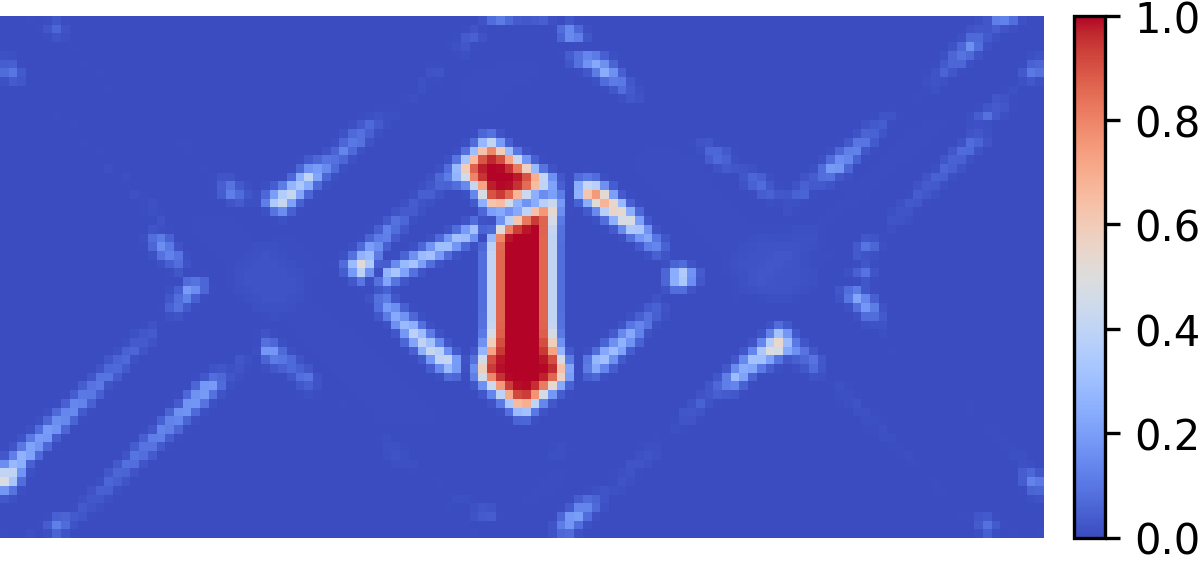}
    \includegraphics[width=0.32\textwidth]{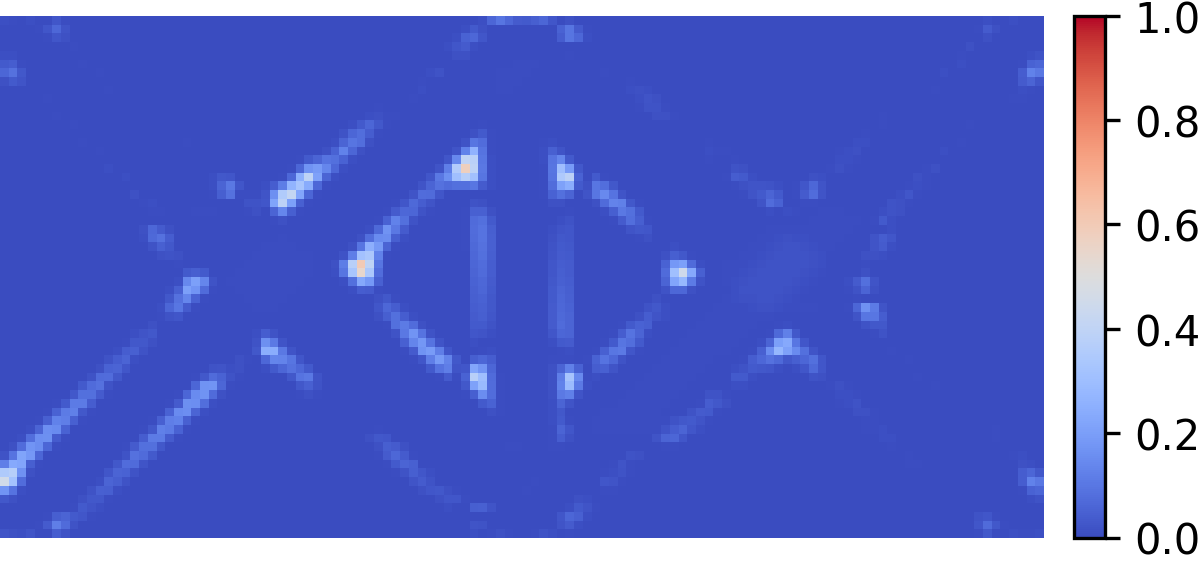}
    \includegraphics[width=0.32\textwidth]{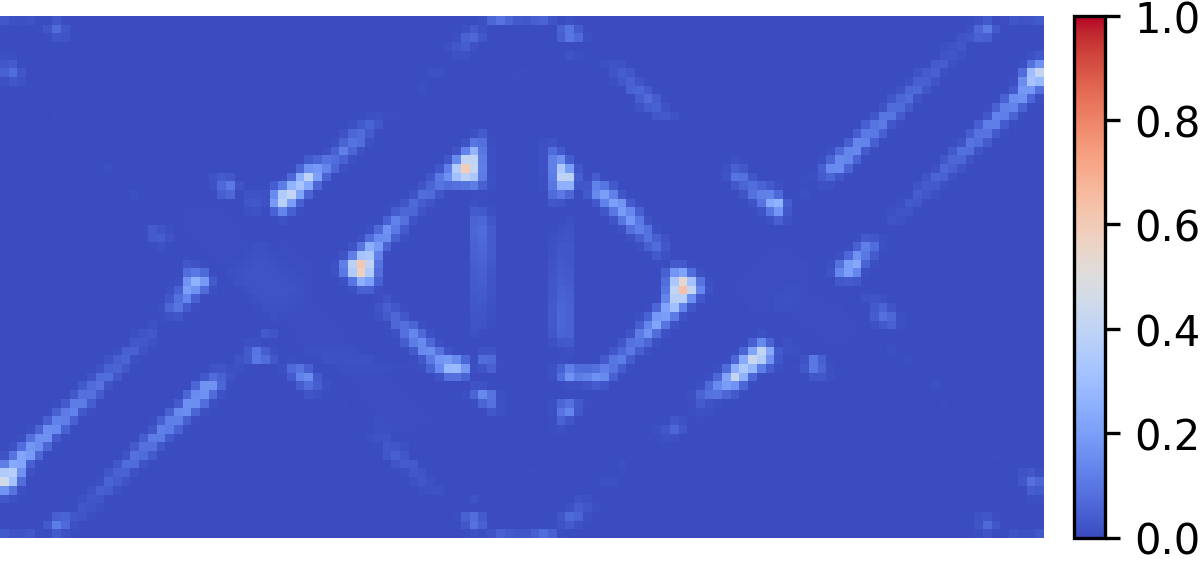}\\
    \includegraphics[width=0.32\textwidth]{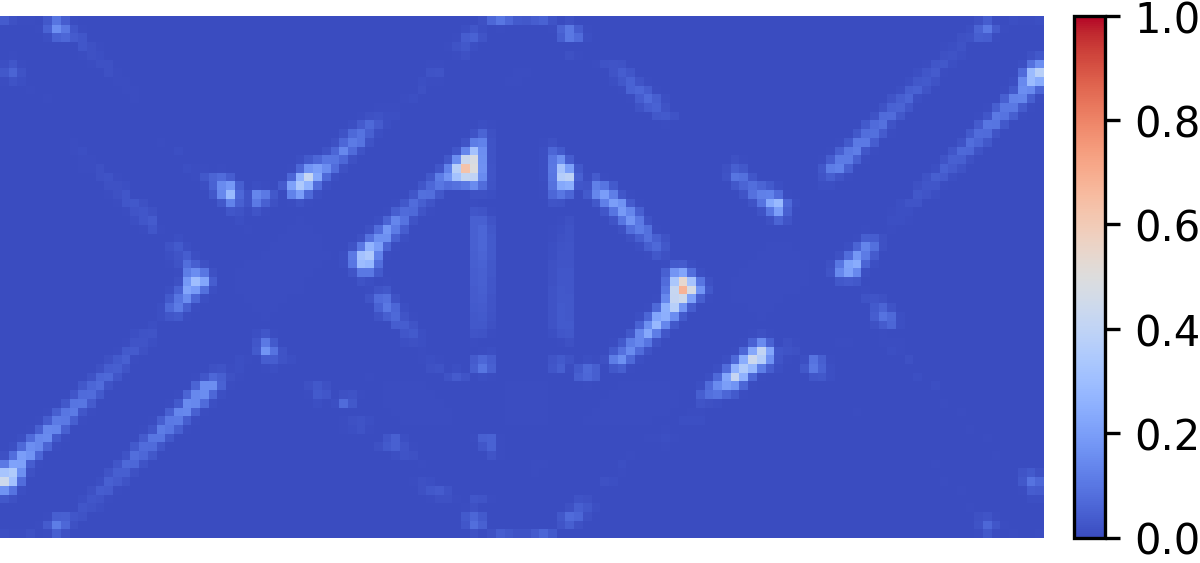}
    \includegraphics[width=0.32\textwidth]{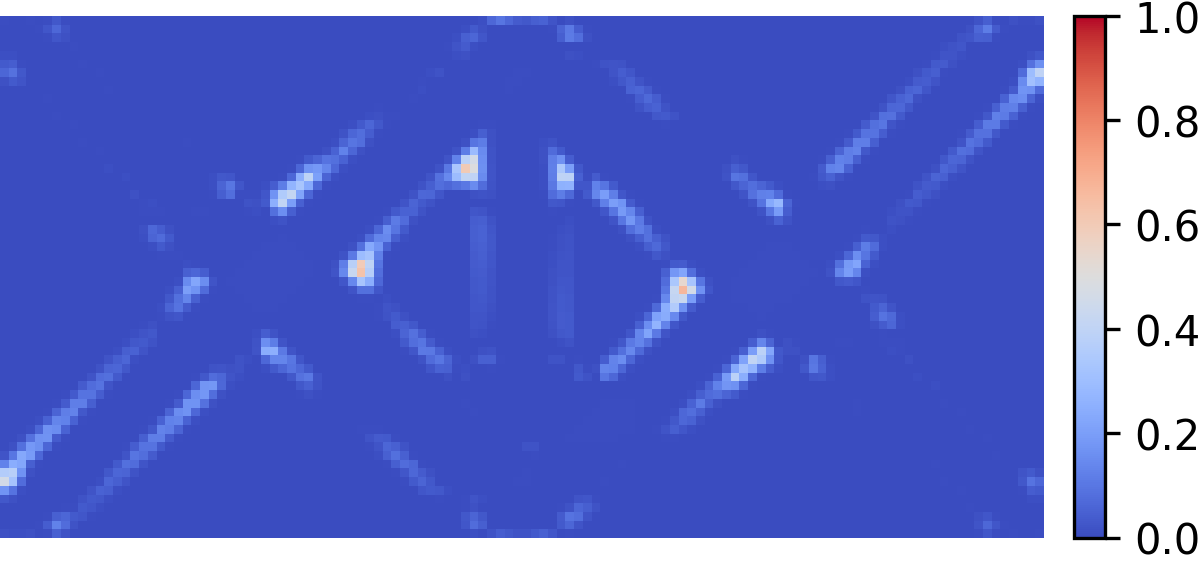}
    \caption{Absolute residual fields for integration orders 1–5 (left to right).
    Order~1 produces large errors in the central region, 
    while orders~2–5 yield consistently low error magnitudes.}
    \label{fig:order-residualabs}
\end{figure}

A direct comparison of the absolute residual fields in
\cref{fig:order-residualabs} confirms these observations. 
For order~1, high error levels dominate the central vertical region. 
Orders~2–5 achieve similar error magnitudes, 
with small deviations confined to the junctions of diagonal and vertical bars. 
The differences between orders~3, 4 and~5 are negligible, 
confirming that higher quadrature accuracy yields no additional improvement in reconstruction quality beyond order~2. 

\begin{figure}[H]
    \centering
    \includegraphics[width=0.8\textwidth]{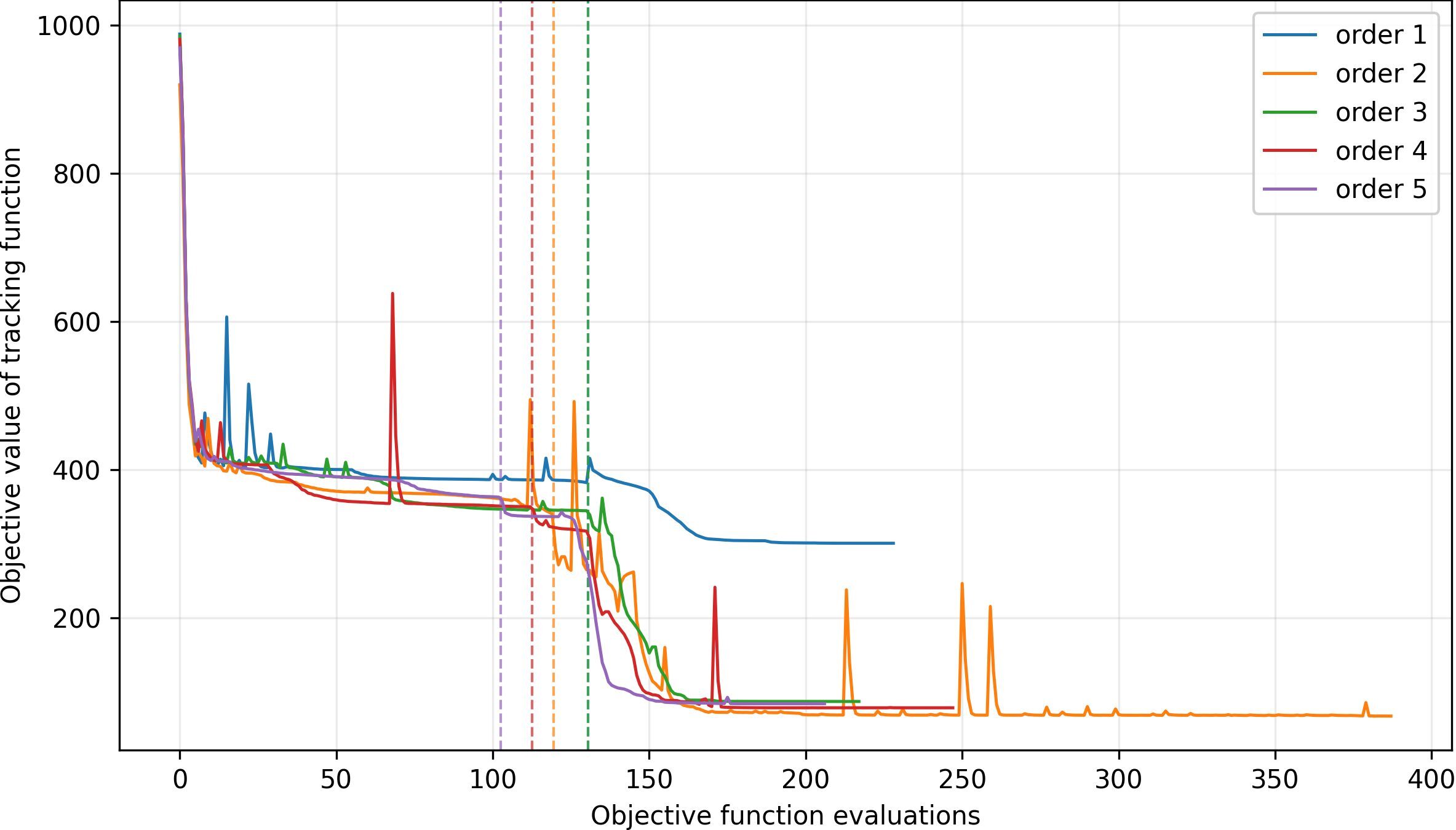}
    \caption{Objective value trajectories for integration orders~1–5. 
    Order~1 stagnates at a high value, 
    while orders~2–5 reach comparable minima. 
    Higher orders exhibit smoother convergence with fewer oscillations.}
    \label{fig:order-objval}
\end{figure}

The evolution of the objective values for all integration orders is shown in
\cref{fig:order-objval}. 
Order~1 stagnates at a significantly higher value, 
in line with its poor reconstruction quality. 
Orders~2–5 converge to comparable minima, 
but their trajectories differ: 
lower orders exhibit oscillations and abrupt jumps near the transition 
from bridging to convergence, 
while higher orders display smoother and more monotone decreases. 
The benefit of increased quadrature accuracy therefore lies primarily in the 
stability of the convergence process rather than in further reduction of the final approximation error. 

% ============================================================
% 9) Exact Hessian vs. BFGS Approximation
% ============================================================
\subsection{Exact Hessian versus BFGS approximation}
\label{sec:hessian}

The influence of second--order information was assessed by comparing optimization runs 
with exact Hessians and with a quasi--Newton approximation. 
For the latter, a limited--memory variant with history size restricted to three updates was employed, 
together with gradient--based scaling and warm starting of both primal variables and multipliers. 
These safeguards were necessary to prevent overly abrupt steps at the beginning of the convergence stage, 
which would otherwise undo progress achieved during exploration or bridging. 

\begin{figure}[H]
  \centering
  \setlength{\tabcolsep}{2pt}
  \renewcommand{\arraystretch}{0.5}
  \begin{tabular}{cccc}
    \includegraphics[width=0.24\textwidth,height=0.12\textwidth]{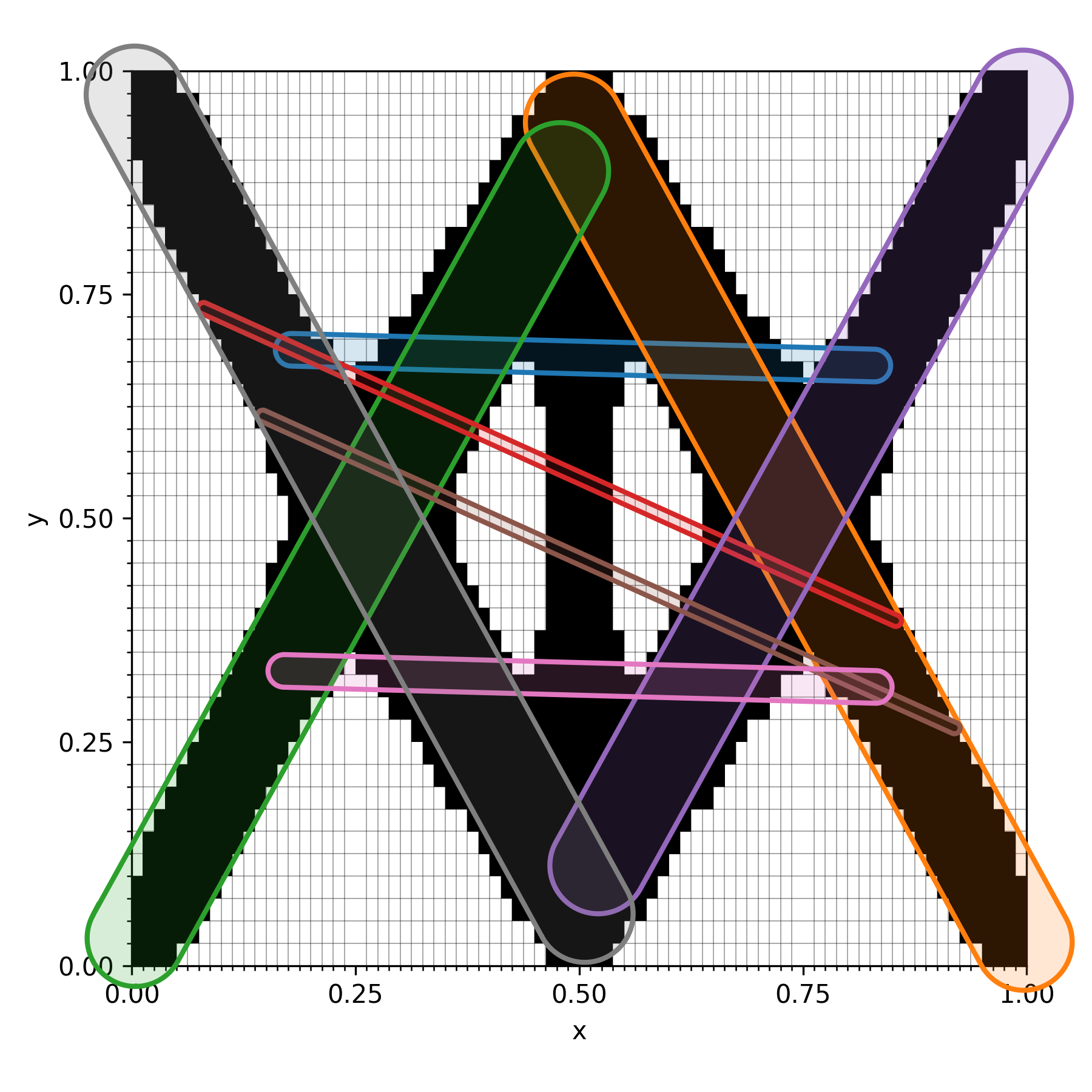} &
    \includegraphics[width=0.24\textwidth,height=0.12\textwidth]{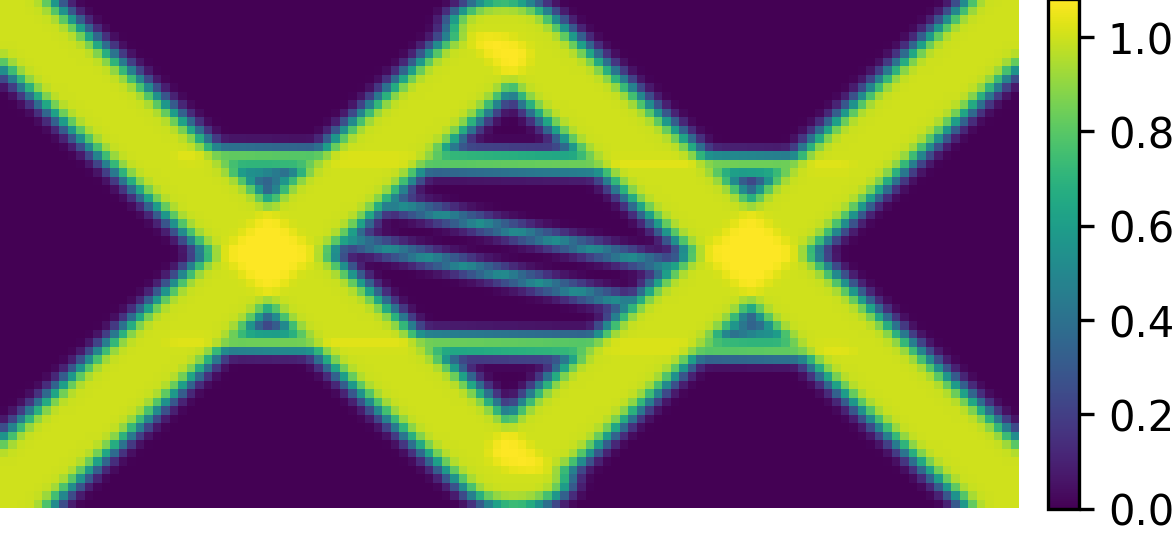} &
    \includegraphics[width=0.24\textwidth,height=0.12\textwidth]{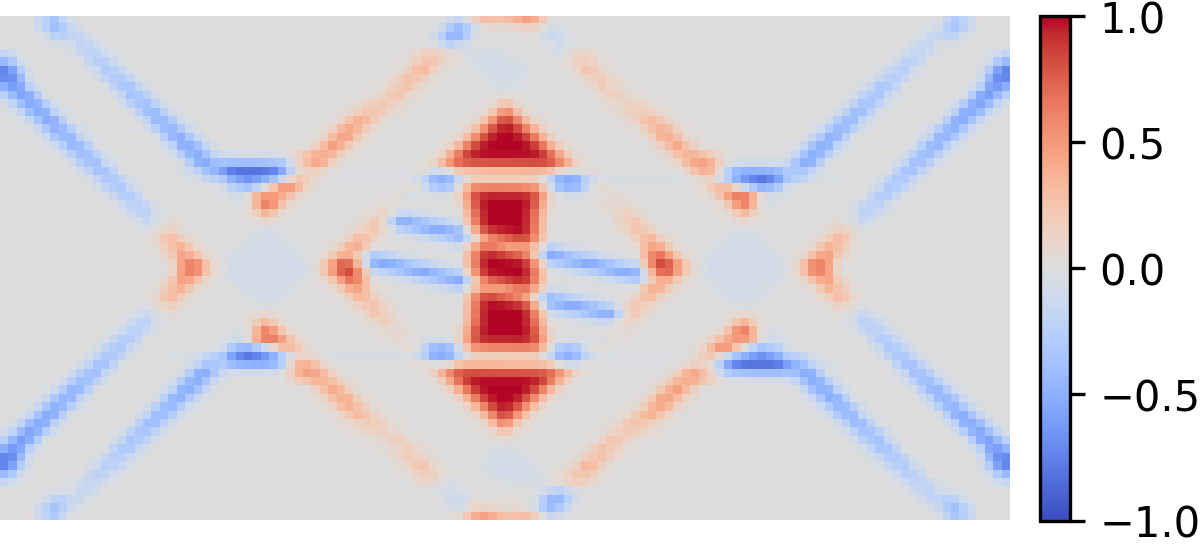} &
    \includegraphics[width=0.24\textwidth,height=0.12\textwidth]{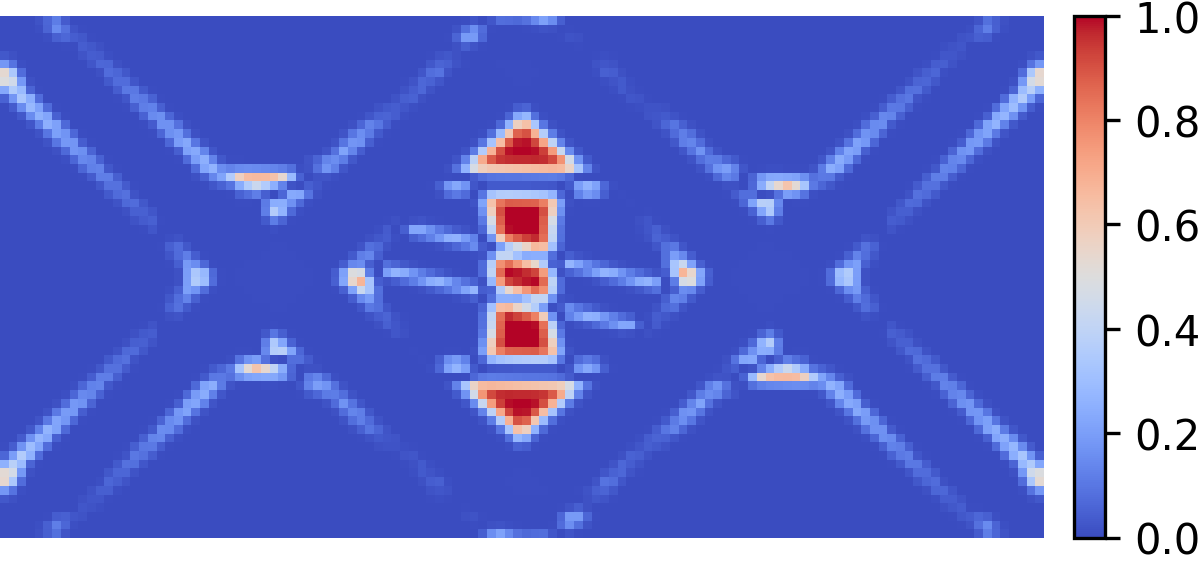} \\
    \includegraphics[width=0.24\textwidth,height=0.12\textwidth]{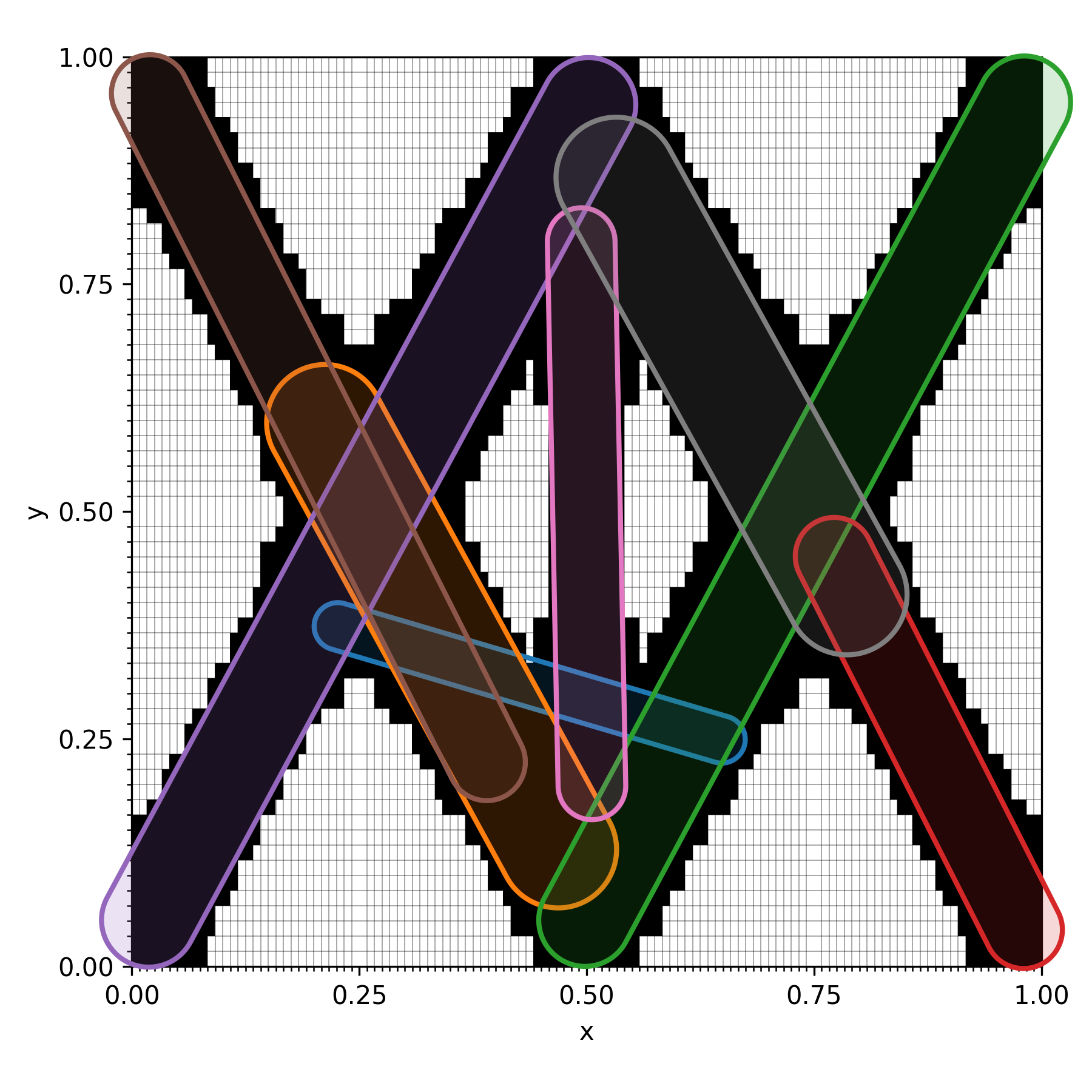} &
    \includegraphics[width=0.24\textwidth,height=0.12\textwidth]{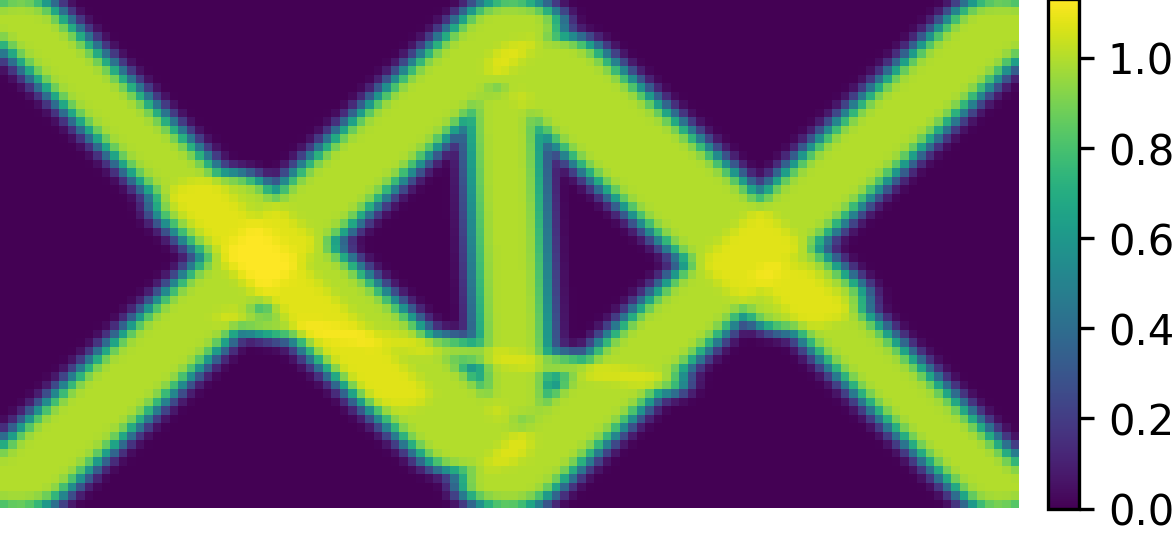} &
    \includegraphics[width=0.24\textwidth,height=0.12\textwidth]{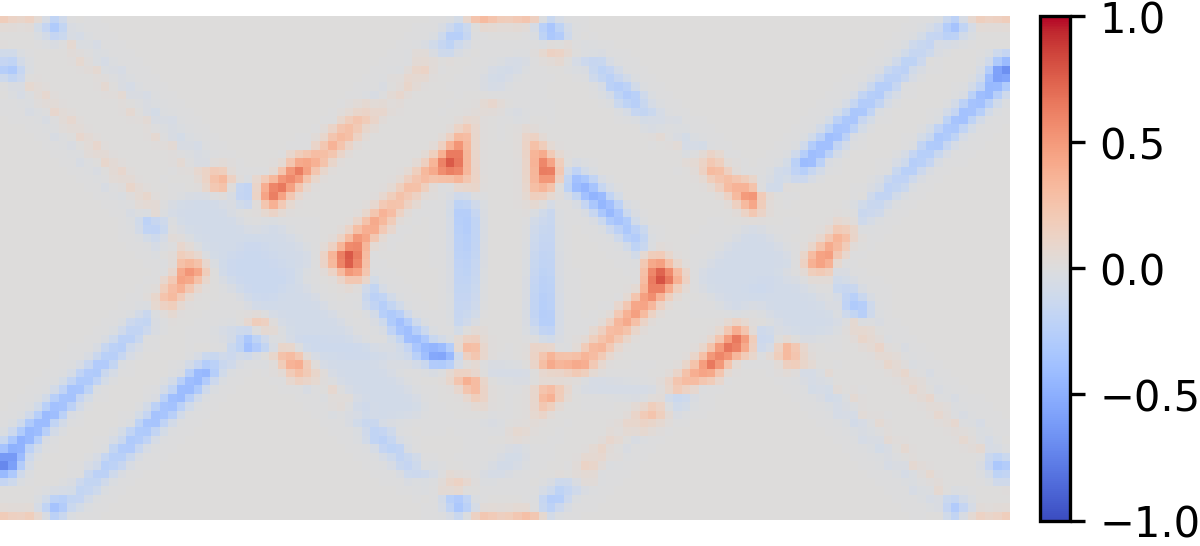} &
    \includegraphics[width=0.24\textwidth,height=0.12\textwidth]{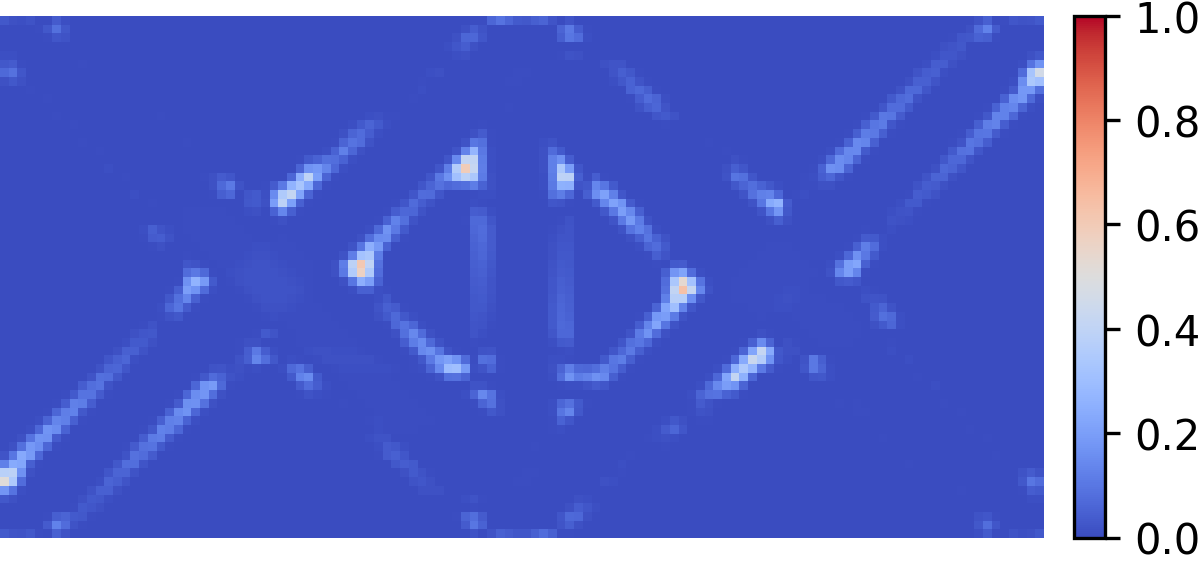} \\
  \end{tabular}
  \caption{Comparison of reconstructions with eight pills using the
  limited--memory BFGS approximation (top row) and exact Hessians (bottom row).
  Shown are pill placement, density field, residual and absolute residual.
  The approximate run leaves several regions uncovered and produces stretched
  pills across void areas, whereas the exact Hessian resolves the structure
  more faithfully.}
  \label{fig:hessian-8}
\end{figure}

\cref{fig:hessian-8} compares the reconstructions obtained with eight pills 
under the two strategies. 
The approximate variant yields a clearly inferior result: 
several pills extend diagonally across multiple disjoint target regions, 
traversing void areas without aligning to individual structures. 
Four pills in particular remain stretched across separate density pockets, 
and the central vertical bar is left largely uncovered. 
The residual maps confirm this lack of fidelity, 
showing persistent discrepancies in central and crossing regions. 
The corresponding optimization trajectory was also less stable, 
with pronounced jumps and a reduced ability of pills to escape local minima 
or to slide past one another. 

In contrast, the exact Hessian run converges to a configuration 
that aligns much more closely with the target. 
The dominant bars are reconstructed with appropriately oriented pills, 
and the central vertical region is successfully captured. 
Residuals are confined to narrow bands along bar boundaries, 
indicating a sharper and more consistent fit. 
These differences underline the role of accurate curvature information: 
while limited--memory updates can approximate descent directions, 
they fail to capture interactions between overlapping pills 
and thereby limit reconstruction quality. 

The comparison becomes even more instructive when the number of pills is
increased to thirteen. 
In the exact run, the pills distribute themselves in a coordinated manner: 
thicker regions of the target are represented by single, well--aligned pills, 
while thinner parts are covered by additional, specialized segments. 
This is particularly evident in the central vertical strip, where three pills cooperate: 
one aligned with the upper thickened section, one with the lower part, 
and one slender segment covering the narrow middle. 
The resulting placement balances local specialization with minimal redundancy, 
achieving accurate coverage with limited overlap. 

In the approximate run, by contrast, this differentiation is largely lost. 
Several pills extend across the entire height of the structure, 
stretching from top to bottom and overlapping heavily. 
In the central region, where the exact solution employs three complementary pills, 
four pills now span the same extent and are almost completely superimposed. 
The thicker target sections at the top and bottom are only coarsely approximated 
through small angular deviations between these overlapping pills. 
At the same time, one misplaced segment stretches diagonally across the void 
between two target regions, reflecting the optimizer’s difficulty in resolving 
competing coverage demands. 
Overall, the approximate update yields a less efficient allocation: 
the target density is still covered, 
but with redundant overlap, spurious extensions into void space, 
and reduced geometric fidelity in the central part. 

\begin{figure}[H]
  \centering
  \setlength{\tabcolsep}{2pt}
  \renewcommand{\arraystretch}{0.5}
  \begin{tabular}{cccc}
    \includegraphics[width=0.24\textwidth,height=0.12\textwidth]{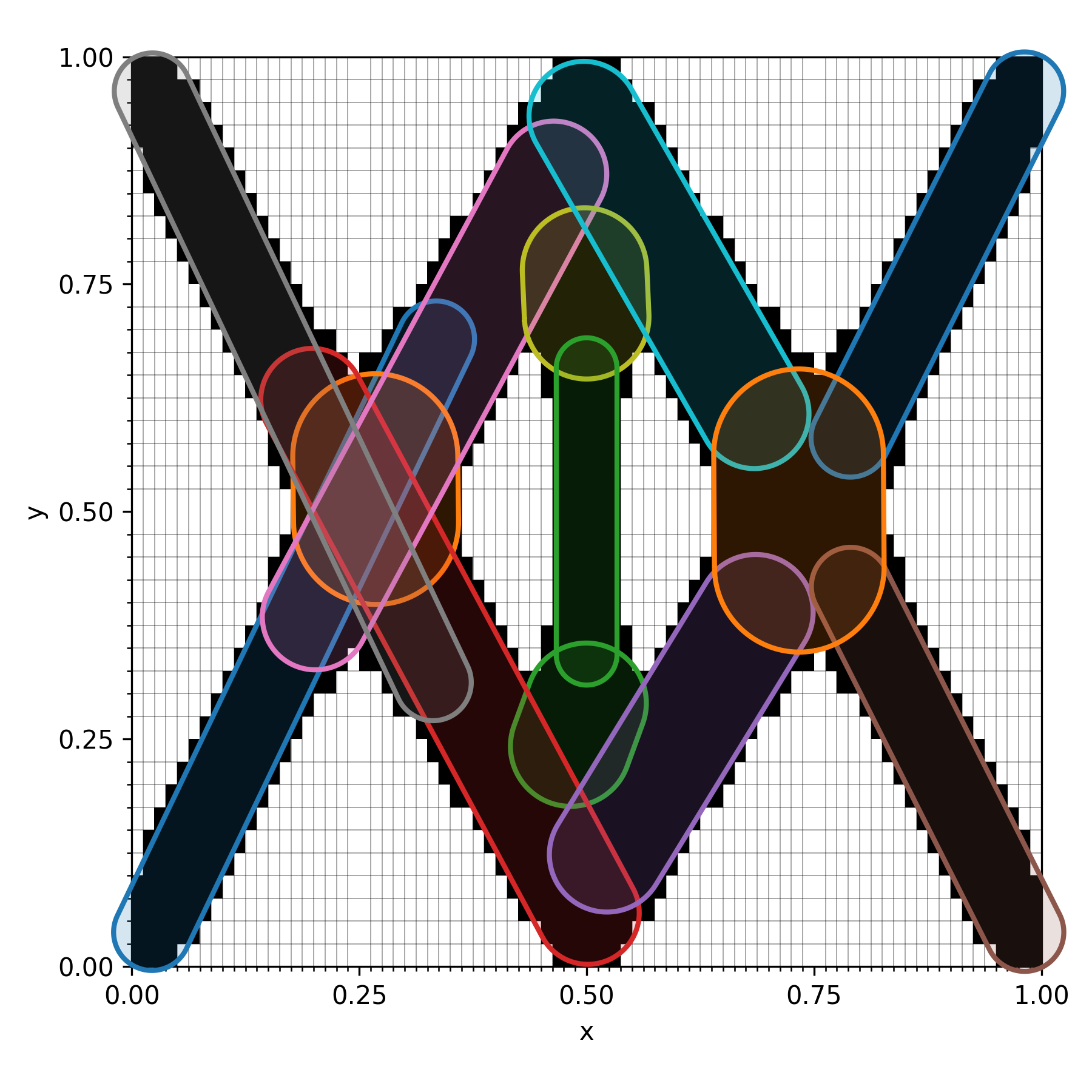} &
    \includegraphics[width=0.24\textwidth,height=0.12\textwidth]{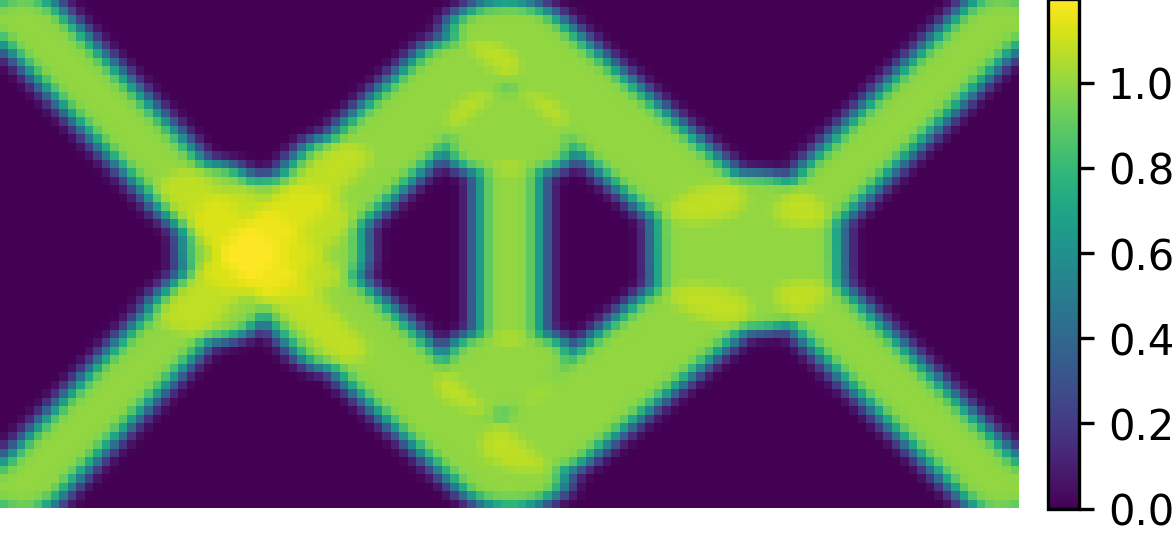} &
    \includegraphics[width=0.24\textwidth,height=0.12\textwidth]{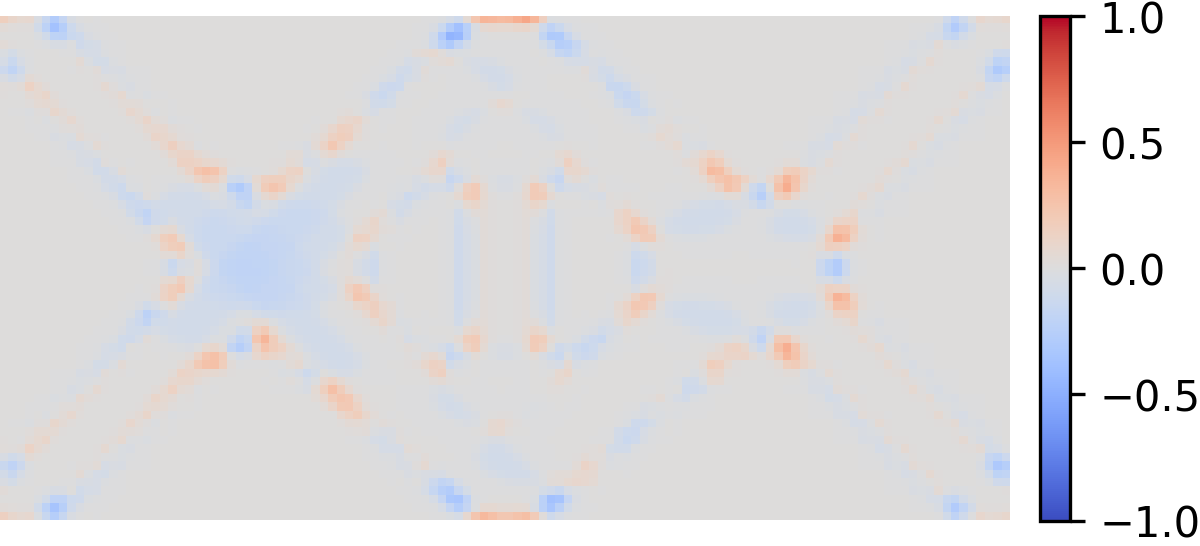} &
    \includegraphics[width=0.24\textwidth,height=0.12\textwidth]{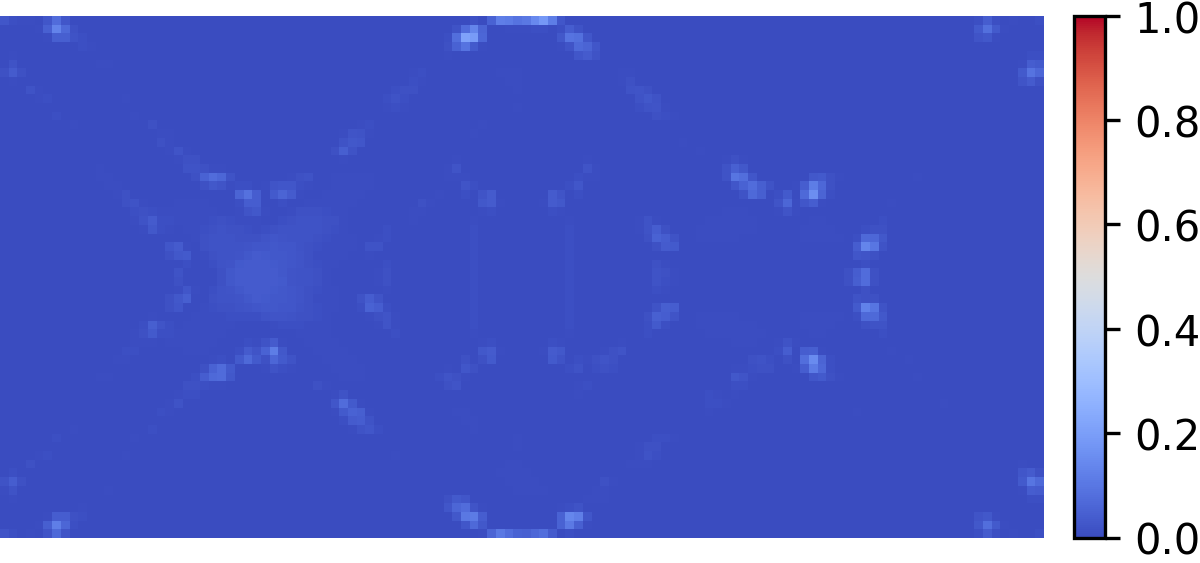} \\
    \includegraphics[width=0.24\textwidth,height=0.12\textwidth]{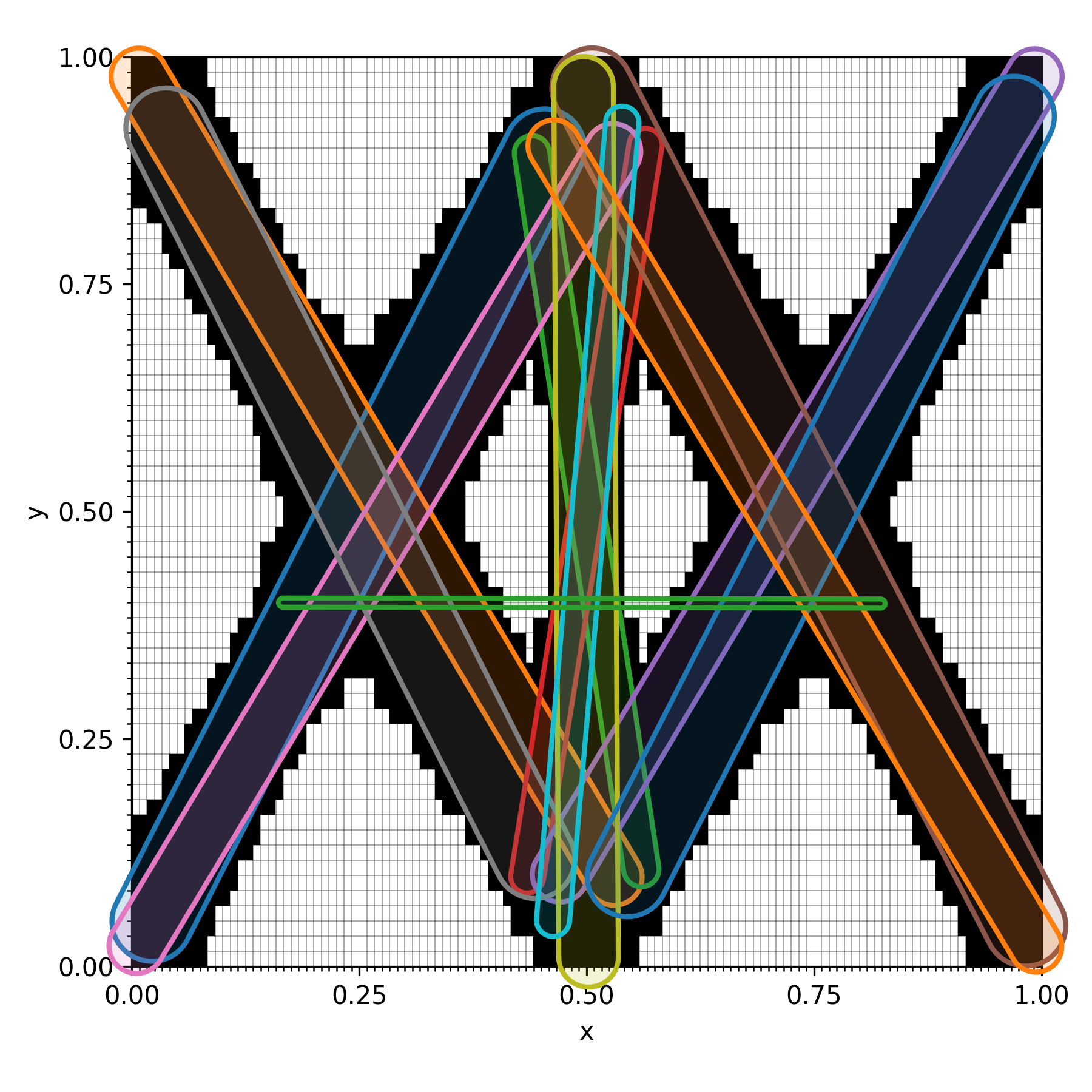} &
    \includegraphics[width=0.24\textwidth,height=0.12\textwidth]{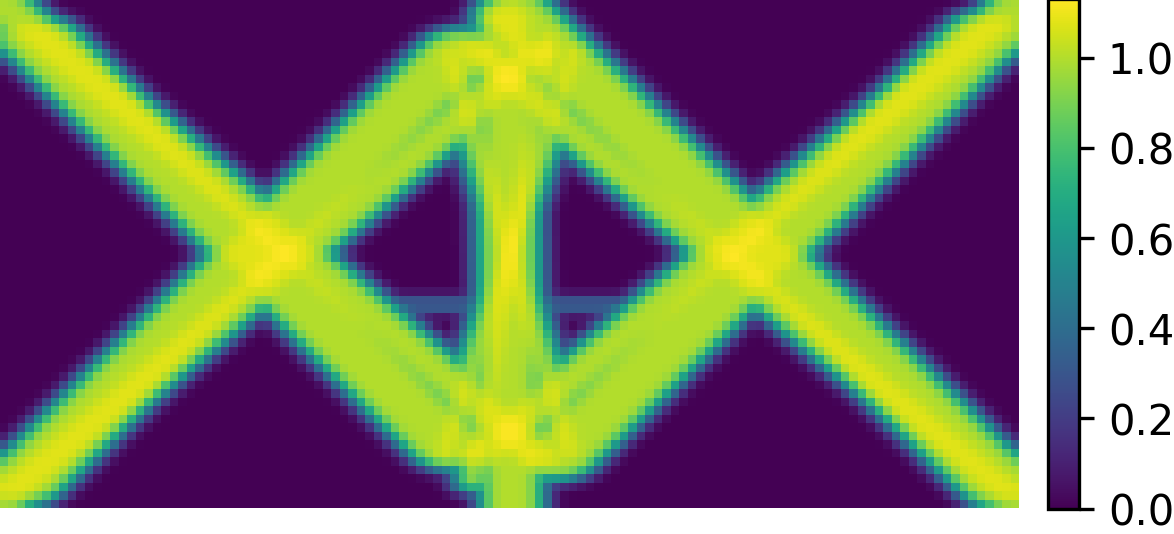} &
    \includegraphics[width=0.24\textwidth,height=0.12\textwidth]{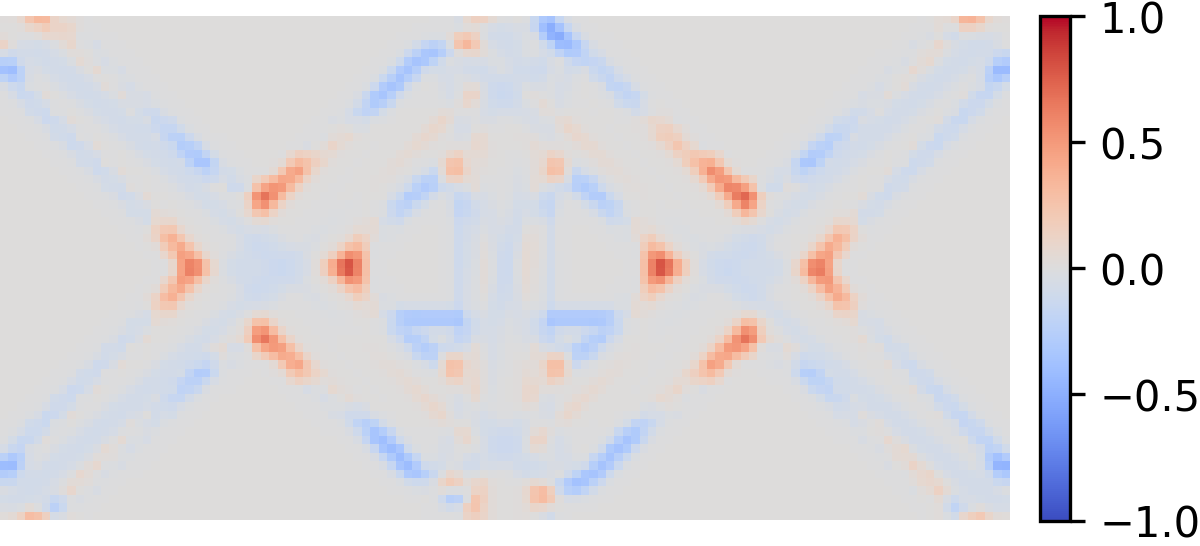} &
    \includegraphics[width=0.24\textwidth,height=0.12\textwidth]{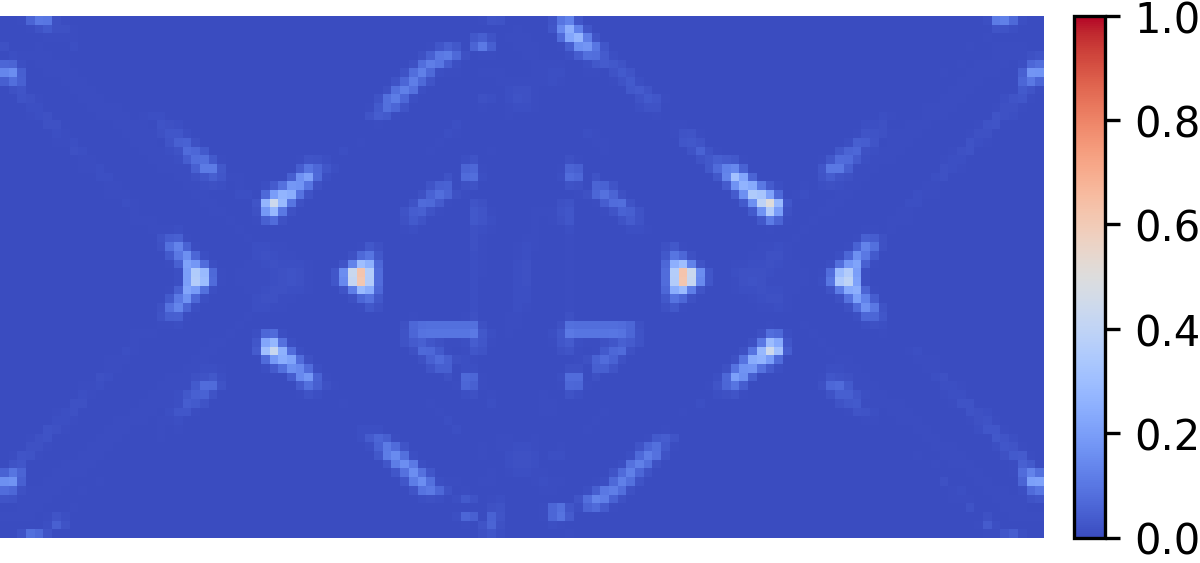} \\
  \end{tabular}
  \caption{Reconstructions with exact (top row) and approximate (bottom row) 
  Hessian updates using thirteen pills. Shown are pill placement, density, 
  residual and absolute residual. The exact variant concentrates pills in 
  complementary roles, whereas the approximate update produces strong overlap 
  and misplaced coverage across void regions.}
  \label{fig:hessian-13}
\end{figure}

\begin{figure}[H]
    \centering
    \includegraphics[width=0.75\textwidth]{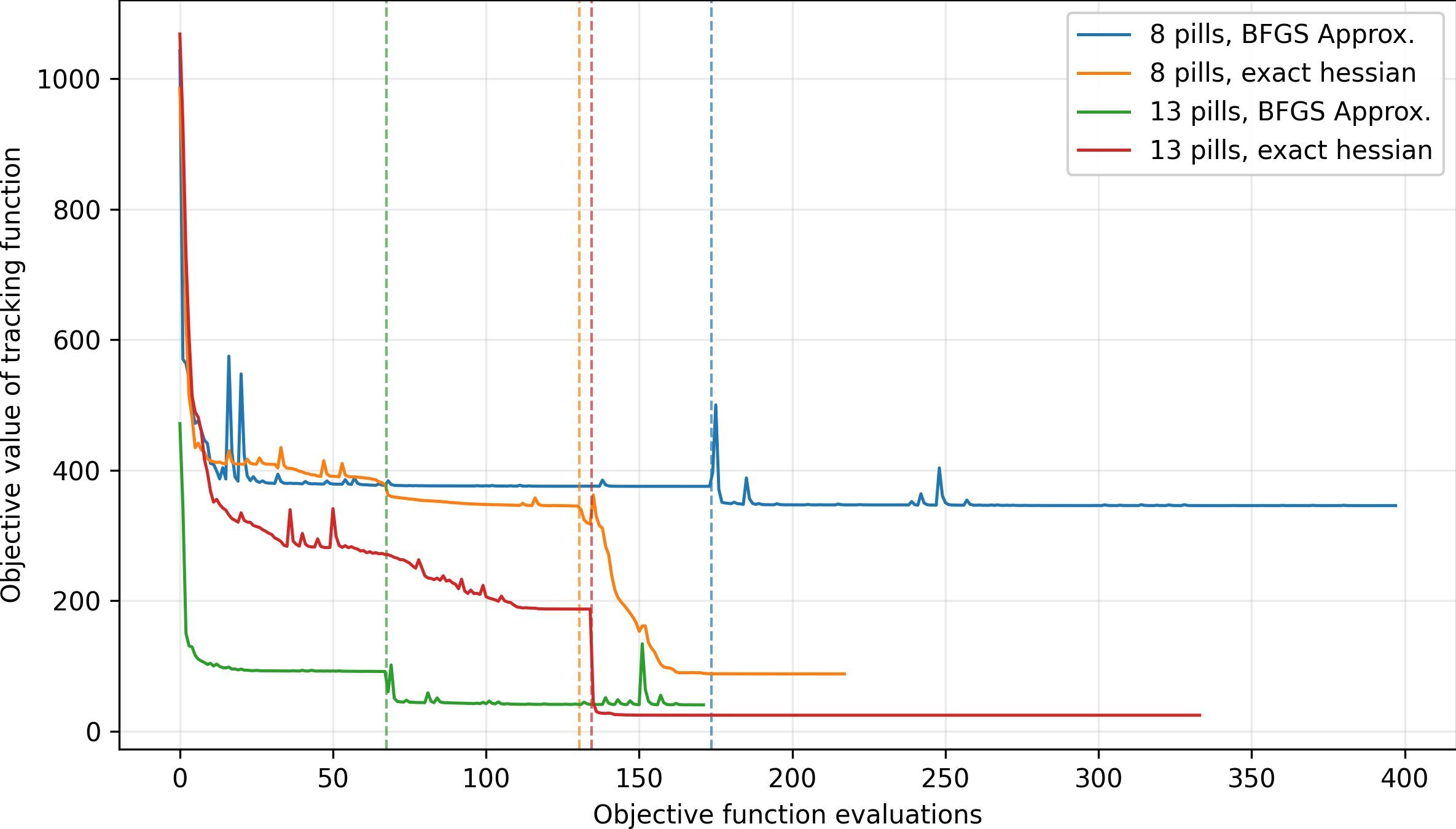}
    \caption{Objective trajectories for approximate (BFGS) and exact Hessian runs 
    with eight and thirteen pills. Vertical dashed lines mark the transition 
    from bridging to convergence. Exact Hessian updates consistently achieve 
    lower final objectives, with particularly pronounced improvement in the 
    eight--pill case.}
    \label{fig:obj-approx-exact}
\end{figure}

The numerical trajectories mirror these observations, 
\cref{fig:obj-approx-exact} shows the evolution during bridging and convergence, 
with vertical dashed lines marking the transition points. 
In both cases, the exact Hessian achieves a lower final objective than the approximate update. 
The improvement becomes particularly clear after the transition to convergence, 
where exact second--order information enables more consistent and decisive descent. 

For the eight--pill runs, the final objective is 
$3.46\times 10^{2}$ with the approximation and 
$8.76\times 10^{1}$ with the exact Hessian, 
a substantial difference that confirms the qualitative superiority of the exact placement. 
For the thirteen--pill runs, the gap in objective values is smaller 
($4.01\times 10^{1}$ vs.\ $2.43\times 10^{1}$), 
which might suggest similar quality if only numerical results were considered. 
However, the pill layouts demonstrate that the approximate variant 
suffers from extensive overlap and misplaced coverage, 
whereas the exact Hessian maintains a more balanced distribution. 
This underlines that purely numerical comparisons may underestimate 
the qualitative differences between the two approaches. 

\subsection{Feature count}
\label{sec:pill-count}

Beyond the choice of aggregation and the availability of second--order information, 
the number of pills itself has a decisive influence on reconstruction quality. 
If too few pills are available, large parts of the target remain uncovered, 
resulting in systematic underfitting. 
Conversely, an excessive number of pills increases the dimensionality of the problem 
and leads to redundancy and blocking effects between overlapping elements. 
Striking a balance is therefore essential. 
To examine this effect in detail, experiments were carried out with 
$3$, $4$, $5$, $8$, $13$ and $18$ pills under otherwise identical conditions.

\begin{figure}[H]
  \centering
  \includegraphics[width=0.79\textwidth]{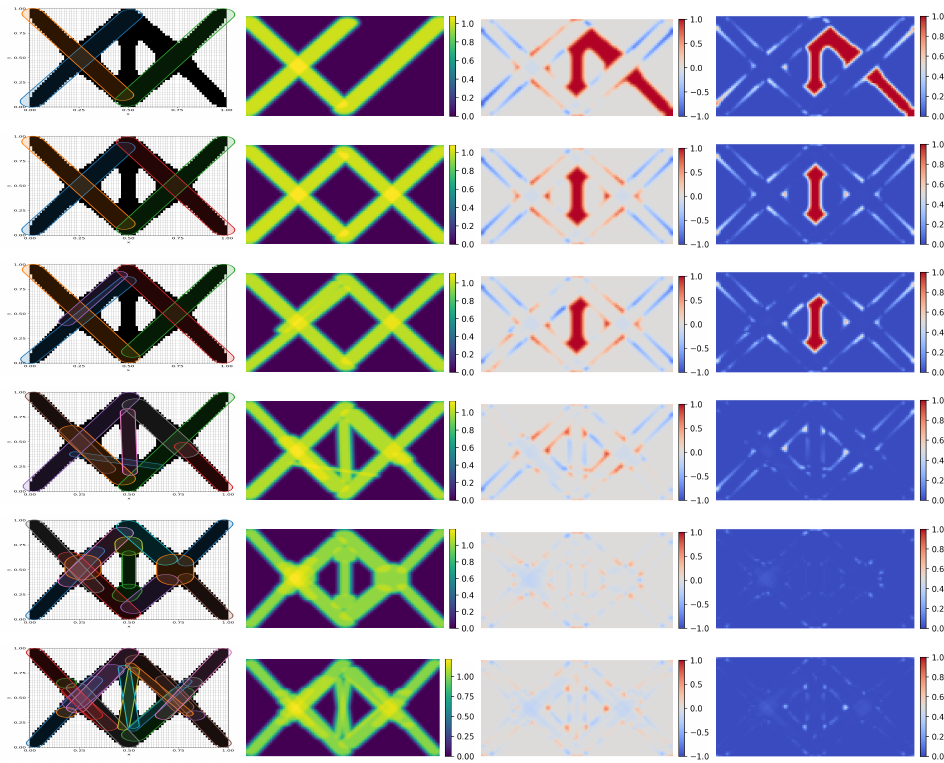}

  \caption{Reconstructions with different pill counts (rows: $3$, $4$, $5$, $8$, $13$, $18$).
  Columns show pill placement, aggregated density, residual and absolute residual.
  Low counts lead to underfitting, while high counts cause redundancy and blocking effects.}
  \label{fig:featurecount-grid}
\end{figure}

The \cref{fig:featurecount-grid} provides an overview of the resulting reconstructions, 
with each row corresponding to one pill count 
and the four columns showing pill placement, aggregated density, residual field, 
and absolute residual. 
The progression from under- to overparameterization is clearly visible. 
With only three or four pills, large parts of the target remain uncovered: 
the central vertical bar is missing entirely, 
and the residual maps show extended coherent errors. 
The corresponding absolute residual fields confirm that these are not localized mismatches 
but systematic failures to represent the geometry. 
Even with five pills, where coverage improves, 
the optimizer still fails to capture the central bar, 
illustrating the limitations of very low counts. 

A marked improvement occurs at eight pills. 
Coverage of the target becomes substantially more consistent, 
and the residuals are concentrated mainly at intersections and transition zones. 
The central bar is now recovered, but the representation of varying bar thickness remains coarse, 
as each pill carries only a single radius. 
With thirteen pills, nearly all gaps close, 
and the reconstructed density aligns closely with the target across the domain. 
Both residual and absolute residual maps show only minor deviations, 
reflecting the ability of additional pills to adapt to local thickness variations. 
At eighteen pills, the improvement is marginal: 
fine-scale details are refined, 
but overlapping and partially redundant pills introduce new artifacts. 
Residual maps reveal local over- and under-shoots 
caused by mutual interference between nearby pills. 
This marks the onset of overfitting, 
where additional degrees of freedom no longer provide substantial gains 
but instead destabilize the approximation. 

\begin{figure}[H]
  \centering
  \includegraphics[width=0.85\textwidth]{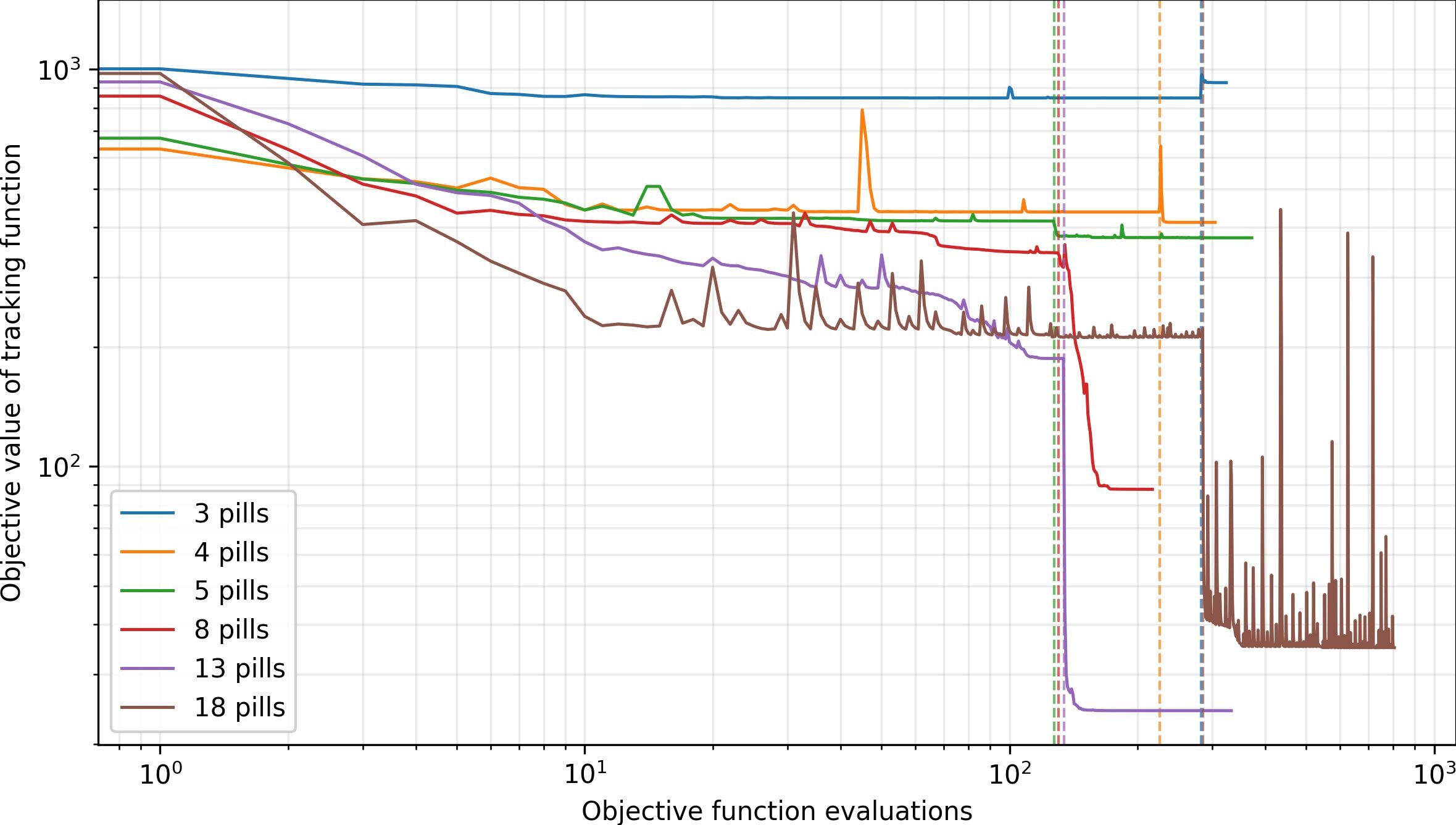}
  \caption{Objective value trajectories for different pill counts. 
  Logarithmic scaling emphasizes the stagnation of low counts (3–5), 
  the marked improvement at eight and thirteen, 
  and the limited additional gain at eighteen, where oscillations and instability appear.}
  \label{fig:featurecount-objval}
\end{figure}

The convergence histories in \cref{fig:featurecount-objval} provide the quantitative counterpart. 
All curves are plotted on a logarithmic scale, which emphasizes both early improvements 
and long-term convergence plateaus. 
Runs with three to five pills stagnate at high objective values, 
consistent with the poor reconstructions observed in the visualizations. 
Increasing to eight pills yields substantial improvement, 
but the objective remains elevated relative to higher counts. 
Thirteen pills reach considerably lower values, 
representing the best compromise between accuracy and stability. 
Eighteen pills achieve no further reduction and instead exhibit stronger oscillations, 
reflecting the interference and redundancy already visible in the residual maps. 
The log-scaling highlights that while the gap between very low and moderate counts is pronounced, 
the difference between thirteen and eighteen pills is comparatively minor. 

Overall, the experiments confirm that pill count governs the balance between
underfitting and overfitting. 
Very small counts cannot represent the target structure, 
while moderate redundancy (around twelve to thirteen pills for the five-bar benchmark) 
improves approximation quality and ensures stable convergence. 
Beyond this range, additional pills yield only diminishing returns 
and may introduce instability due to overlap and redundancy. 
Rather than treating pill count as fixed, the following sections explore strategies 
to dynamically adjust the number of pills during optimization. 
Iterative refinement (\cref{sec:refinement}) introduces missing pills to counteract undercoverage, 
while heuristic procedures (\cref{sec:heuristics}) focus on removing or merging redundant pills. 
Together, these approaches address the limitations observed here 
and provide mechanisms to balance model capacity in a controlled and adaptive manner. 

% ============================================================
% 6) Iterative Refinement
% ============================================================
\subsection{Iterative refinement}
\label{sec:refinement}

The preceding analysis of pill count (\cref{sec:pill-count}) has 
demonstrated the fundamental trade-off between underfitting and overfitting: 
too few pills cannot capture the structural complexity of the target, 
while excessive redundancy leads to instability and overlapping artifacts. 
A natural extension of this observation is adaptive refinement: 
instead of fixing the number of pills at the outset, 
new pills are added iteratively only where they are still required. 
This approach seeks to combine the efficiency of small initial pill sets 
with the accuracy of richer representations, 
while avoiding the disadvantages associated with redundancy. 

The refinement strategy builds directly on the staged optimization framework 
introduced in the single-pill case 
(cf.~\cref{sec:single_feature_optimization}). 
Optimization is still carried out in exploration, bridging and convergence stages, 
but is now supplemented by an outer loop in which new pills are inserted 
based on residual information. 
Specifically, after a baseline optimization with an initial pill count, 
the residual between the aggregated density and the target is evaluated. 
Regions where this residual exceeds a prescribed threshold are identified, 
and a new pill is seeded in one such region. 
This pill undergoes orientation and local convergence steps against the residual 
before being combined with the existing set and jointly optimized against the full target. 
If the objective decreases, the new pill is retained; 
otherwise, the process terminates. 

\begin{figure}[H]
  \centering
  \includegraphics[width=0.75\textwidth,height=0.375\textwidth]{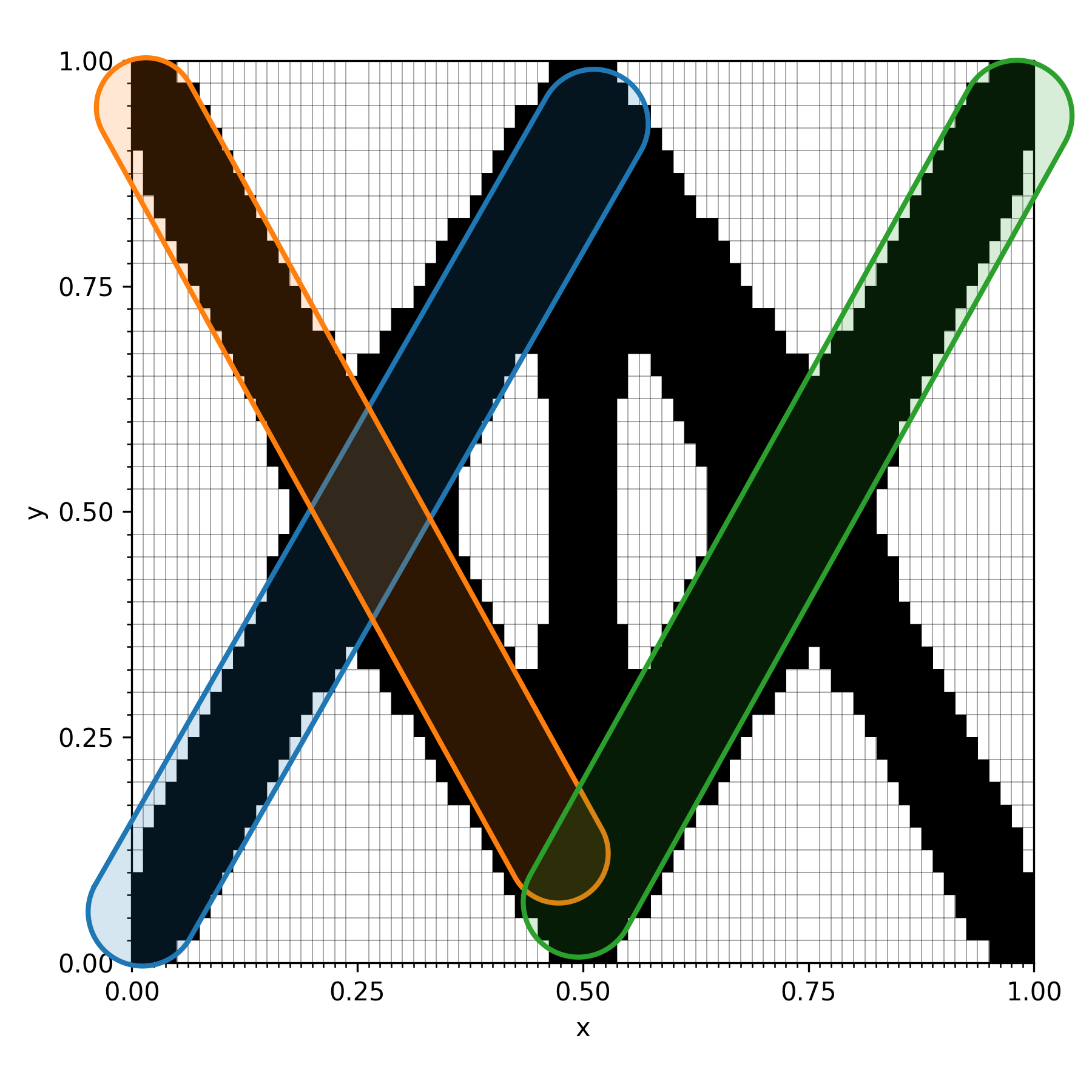}
  \caption{Pill placement after initial optimization with 
  three pills. Large parts of the target remain uncovered, 
  resulting in pronounced residuals.}
  \label{fig:iter-init-pill}
\end{figure}

\begin{figure}[H]
  \centering
  \includegraphics[width=0.75\textwidth]{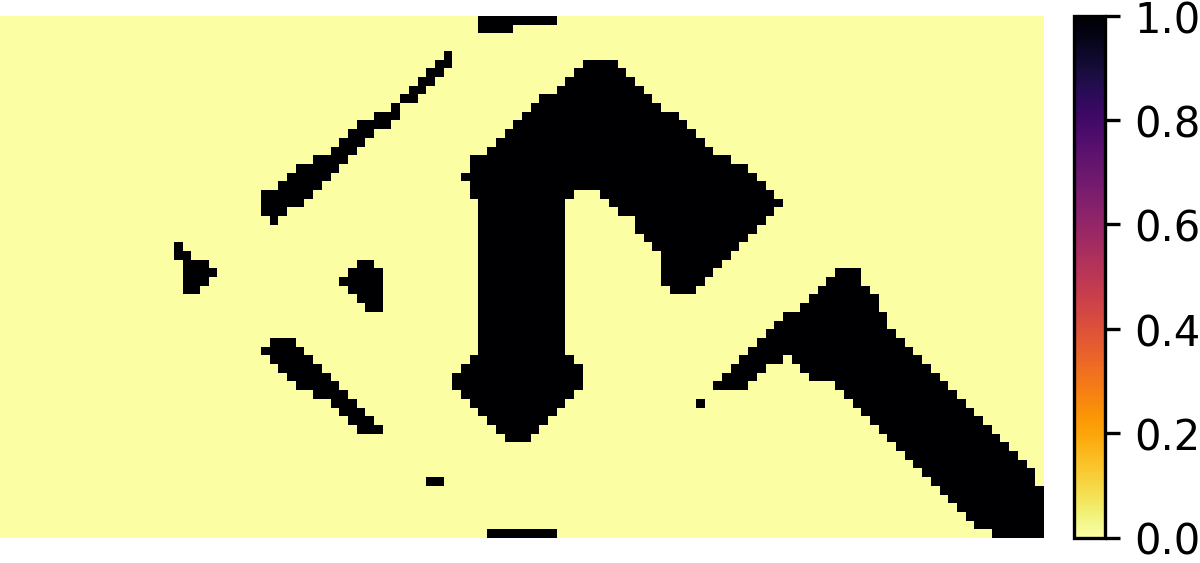}
  \caption{Residual field after three-pill optimization, 
  filtered with a threshold to suppress minor transition errors. 
  The remaining high residuals identify target regions that 
  are not captured by the current pills.}
  \label{fig:iter-init-residual}
\end{figure}

To illustrate this procedure, optimization is first performed with three pills. 
The resulting pill placement is shown in \cref{fig:iter-init-pill}. 
The target structure is only partially represented, leaving entire bars uncovered. 
The corresponding residual, filtered by a threshold to suppress minor deviations 
at pill boundaries, is displayed in \cref{fig:iter-init-residual}. 
Large coherent residuals highlight structural elements of the target that remain unrepresented, 
and these serve as insertion zones for new pills. 

\begin{figure}[H]
\centering
\includegraphics[width=0.75\textwidth,height=0.375\textwidth]{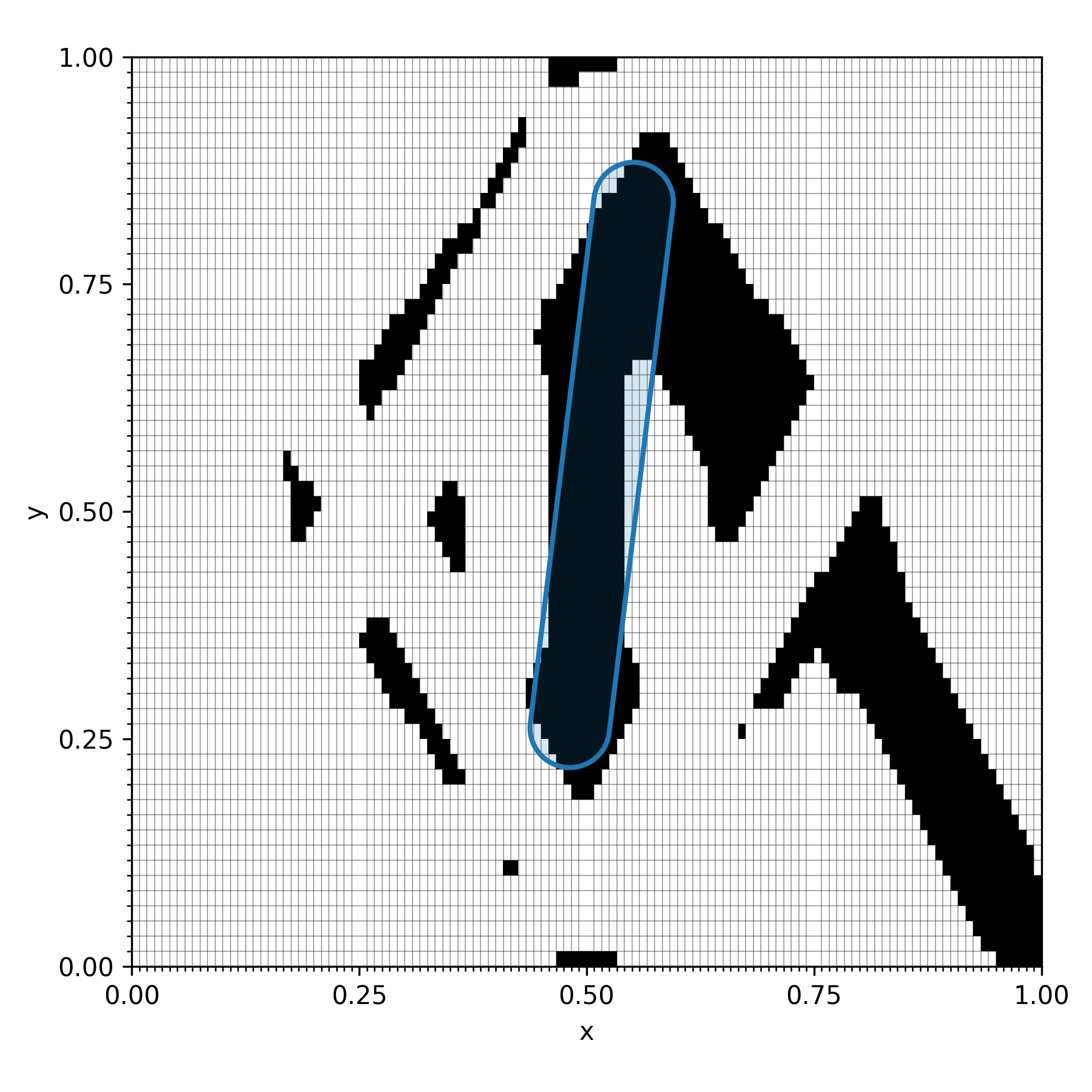}
\caption{Placement of the first added pill after local 
exploration and convergence against the residual. 
The new pill aligns with the vertical bar previously missing 
in the reconstruction.}
\label{fig:iter-first-pill}
\end{figure}

In the first refinement cycle, a new pill is inserted in the residual 
region corresponding to the central vertical bar. 
After exploration and local convergence against the residual, 
the pill aligns closely with the missing target component, 
as shown in \cref{fig:iter-first-pill}. 
This demonstrates that the residual-guided insertion mechanism 
effectively detects and fills structural gaps that cannot be represented 
by the initial configuration. 

\begin{figure}[H]
  \centering
  \includegraphics[width=0.75\textwidth,height=0.375\textwidth]{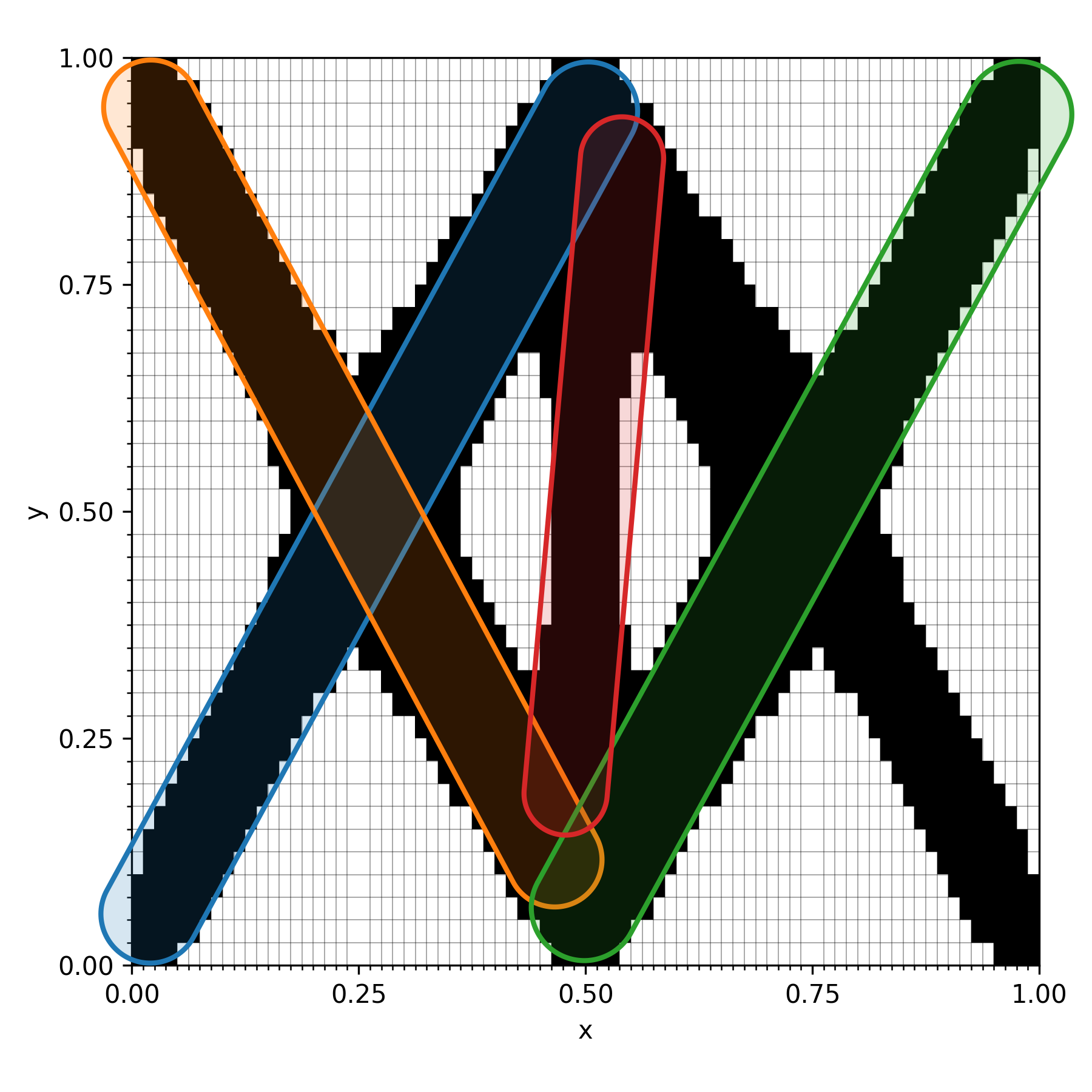}
  \caption{Pill placement after the first addition and subsequent combined 
  convergence. The red pill occupies the previously missing vertical bar, 
  while the pre-existing pills adjust their positions to improve overall coverage.}
  \label{fig:iter-first-converged}
\end{figure}

After this local optimization, the new pill is integrated with the existing set, 
and a combined convergence stage is performed against the full target. 
The resulting configuration is shown in \cref{fig:iter-first-converged}. 
The newly added red pill remains aligned with the vertical bar, 
while the previously existing pills adjust their positions. 
In particular, the blue pill shifts upward, 
whereas the orange and green pills move downward at their intersections 
with the central bar. 
These adjustments allow the original pills to relinquish parts of the target 
that are now adequately represented by the red pill, 
leading to a more balanced reconstruction. 
The refinement process therefore not only improves coverage by inserting missing components 
but also facilitates a reorganization of the full pill set, 
enhancing global consistency.

The refinement then continued with further insertions. 
After the second addition and subsequent convergence, 
a new segment emerged along the diagonal from the lower right support 
toward the mid–upper region. 
This pill successfully captured the corresponding target bar, 
which had previously been insufficiently represented. 
Its exploration aligned closely with the target, 
and the existing pills adapted slightly to redistribute coverage, 
thereby reducing the residual in this region.  

The third insertion improved the representation of the central cross–section. 
The added pill was placed in the mid–left region, 
filling density that had previously remained uncovered along the intersecting bar. 
This contribution reduced the residual in the crossing zone, 
where the existing pills alone could not fully capture the structure. 
Together, the second and third insertions closed two critical gaps of the earlier configuration: 
the diagonal trajectory on the right and the central cross–section. 

\begin{figure}[H]
  \centering
  \includegraphics[width=0.48\textwidth,height=0.24\textwidth]{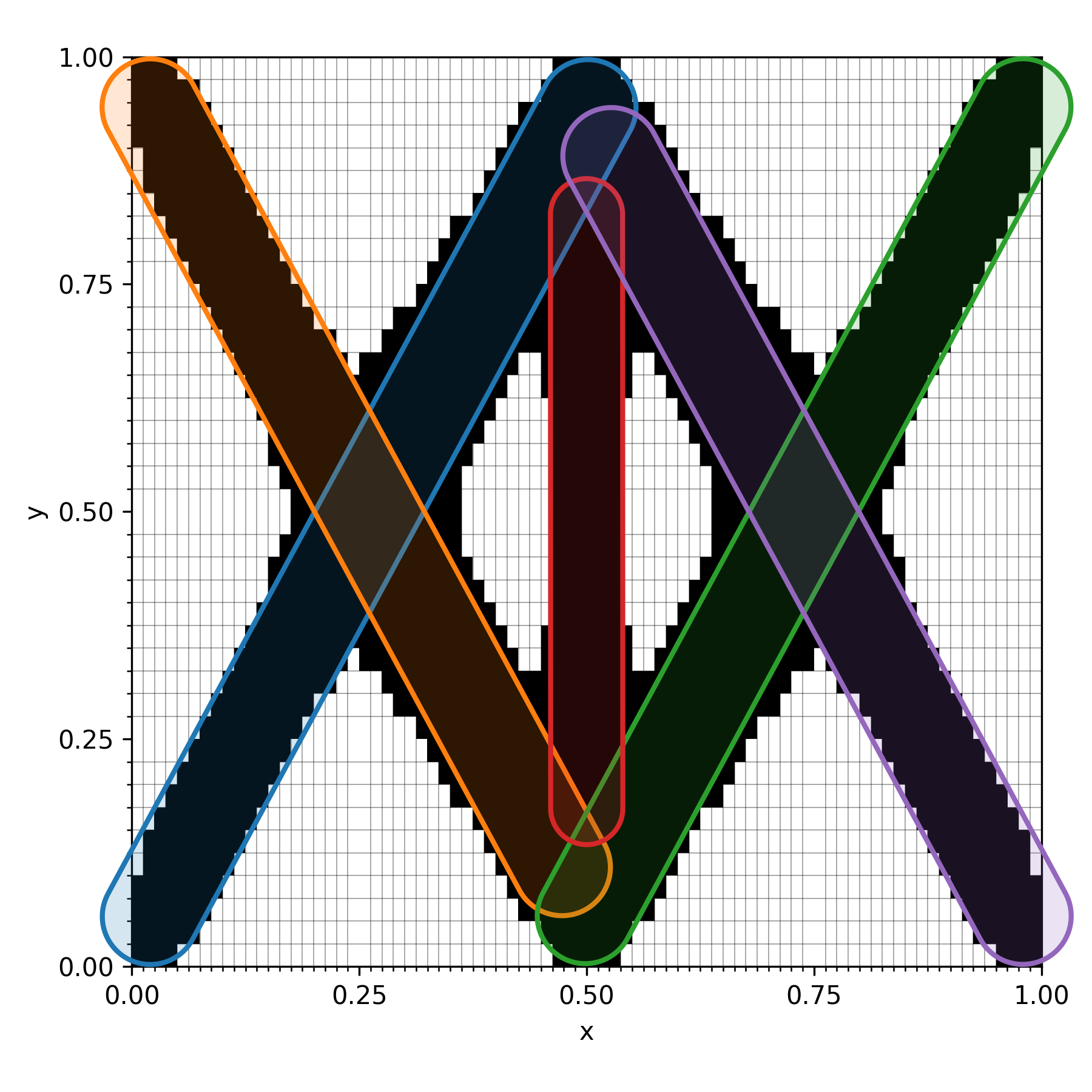}
  \includegraphics[width=0.48\textwidth,height=0.24\textwidth]{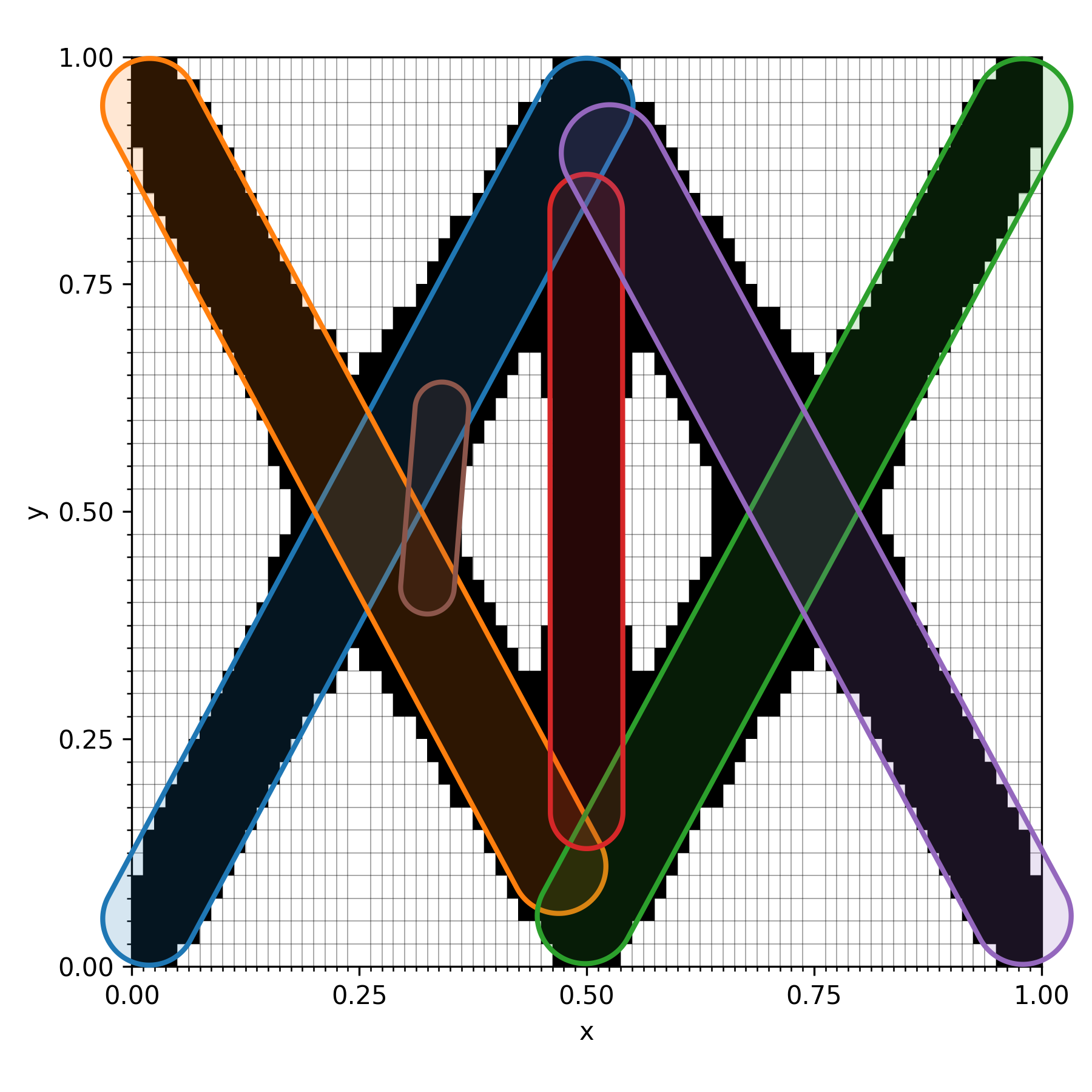}
  \caption{Pill placements after the second (left) and third (right) additions, 
  each followed by convergence. 
  The second pill captures the diagonal bar from the lower right, 
  while the third strengthens the central cross–section in the mid–left region.}
  \label{fig:iter-second-third}
\end{figure}

At this stage, a fourth insertion was attempted. 
The residual field had fragmented into several small and spatially disconnected patches, 
each corresponding to fine details not yet captured by the current pills. 
During exploratory optimization, the candidate pill attempted to span multiple regions simultaneously. 
After convergence, however, it reduced its radius and settled into a position 
that crossed several void areas without aligning cleanly with any residual patch. 
Its contribution was therefore ineffective: 
rather than closing a gap, the new segment introduced superfluous coverage 
and increased the objective value. 
Since no significant improvement could be achieved, 
the additive process terminated at this point, 
discarding the unsuccessful insertion and retaining the configuration after the third addition. 
The three accepted pills had already resolved the most critical deficiencies of the initial setup, 
while further insertions proved counterproductive.

\begin{figure}[H]
  \centering
  \includegraphics[width=0.48\textwidth,height=0.24\textwidth]{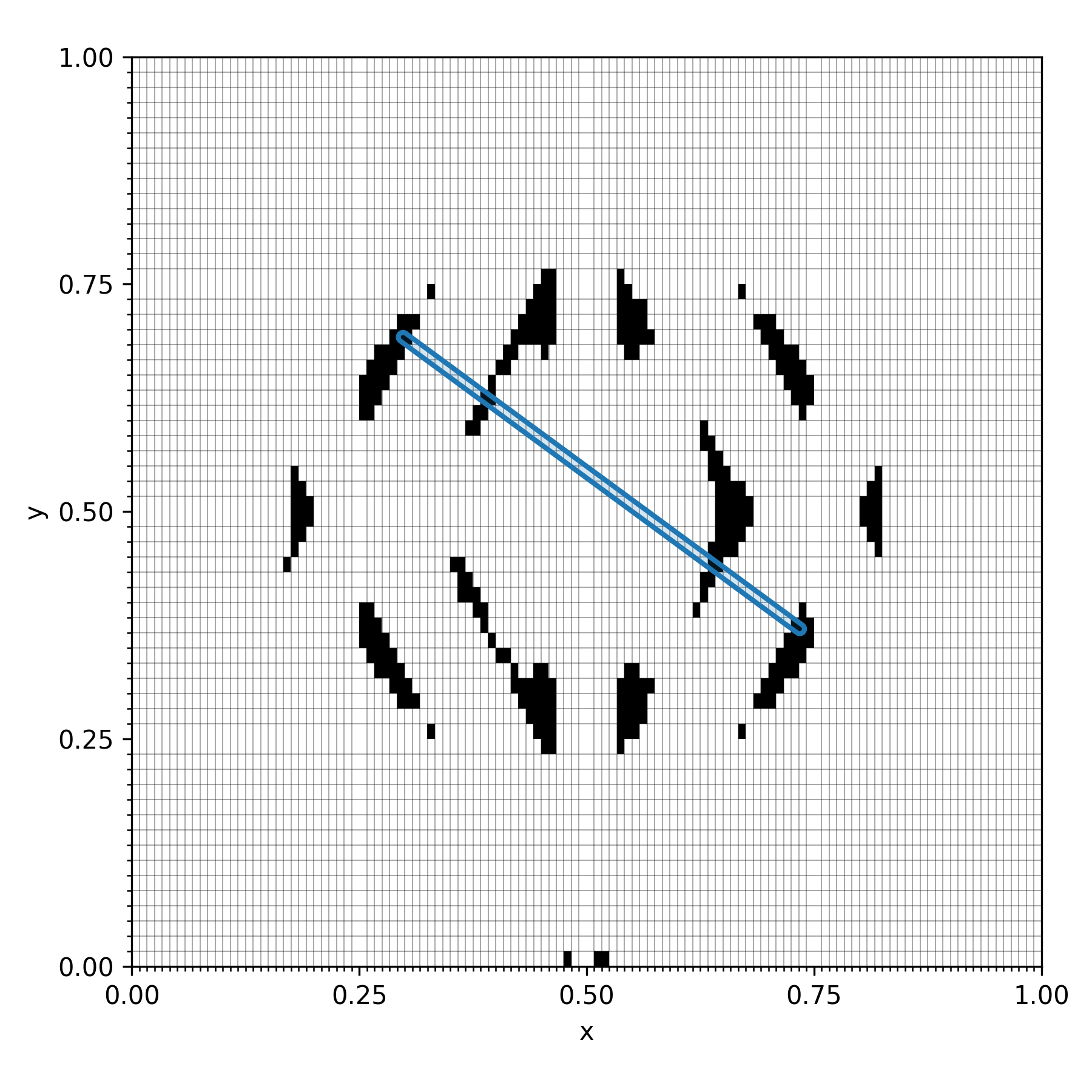}
  \includegraphics[width=0.48\textwidth,height=0.24\textwidth]{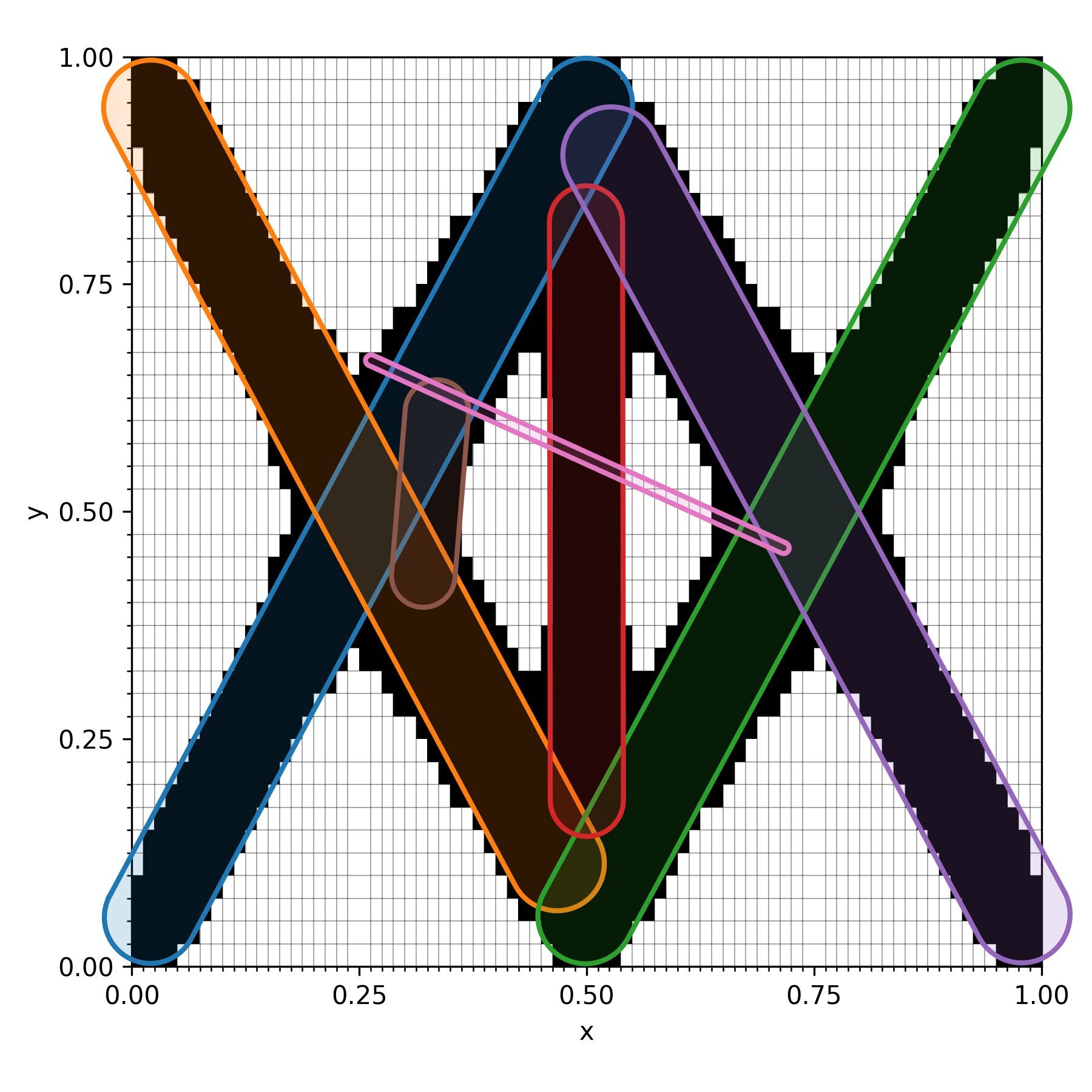}
  \caption{Attempted fourth insertion. 
  Left: residual field and candidate pill during optimization. 
  Right: converged placement of the added pill, which spans multiple void regions 
  without effectively covering the residual density. 
  As the objective value increased, the insertion was rejected and refinement terminated.}
  \label{fig:iter-final-failed}
\end{figure}

% ============================================================
% 7) Heuristics
% ============================================================
\subsection{Heuristics}
\label{sec:heuristics}

While the preceding analysis of pill count has shown that both
underfitting and overfitting may arise, 
the strategies of iterative refinement and heuristic reduction 
serve as optional mechanisms to address these issues. 
Iterative refinement mitigates undercoverage by adding pills in regions 
where residuals remain large, 
whereas heuristics counteract overcoverage by removing or consolidating 
redundant pills. 
The following results illustrate how such heuristics can be applied in practice. 

\begin{figure}[H]
  \centering
  \includegraphics[width=0.75\textwidth,height=0.375\textwidth]{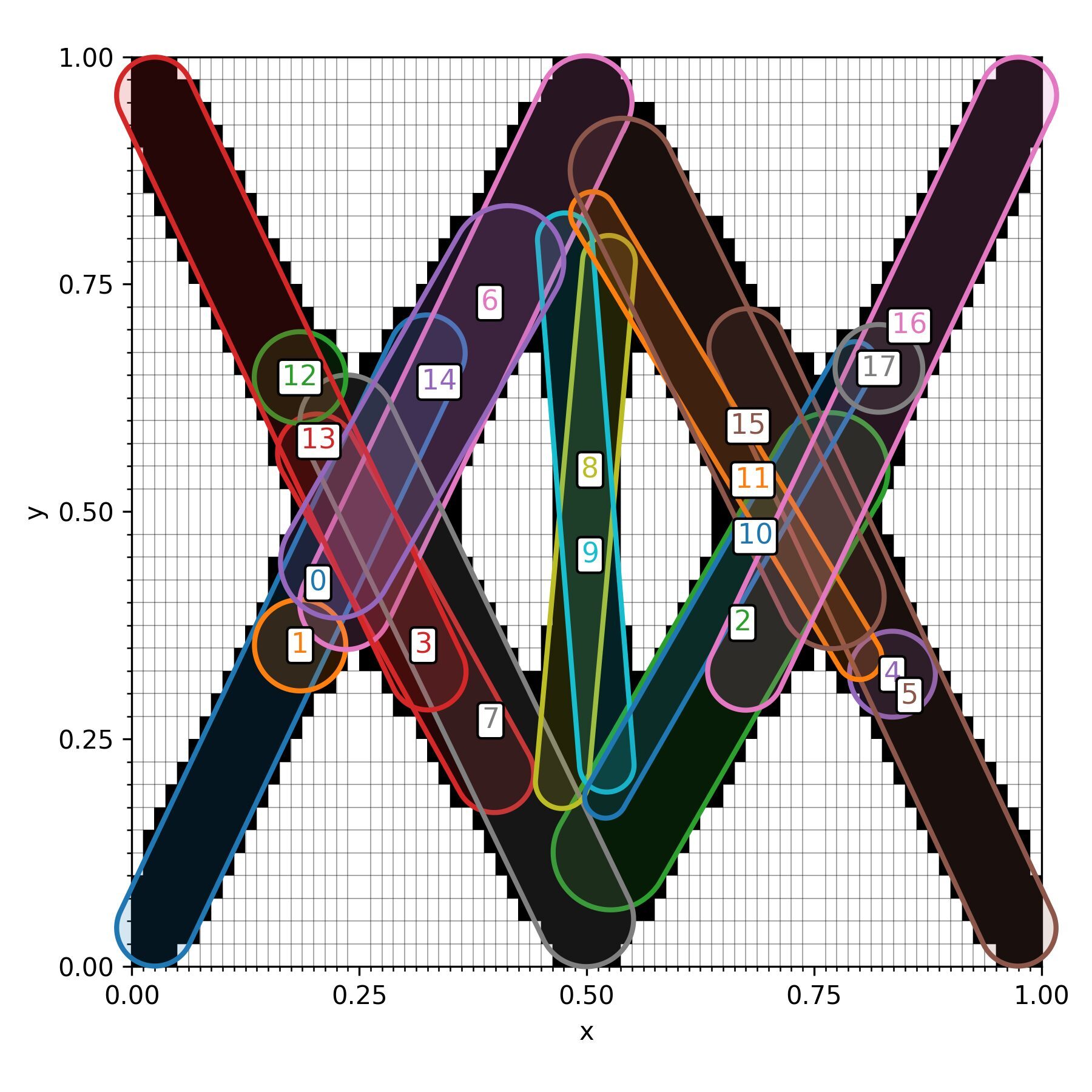}
  \caption{Initial pill placement with eighteen pills before applying
  heuristics. Pill IDs are displayed for reference.}
  \label{fig:heuristics-initial}
\end{figure}

As the starting point for illustration, 
the final configuration with eighteen pills is considered (\cref{fig:heuristics-initial}). 
This placement exhibits clear redundancy: several pills overlap 
or contribute only marginally to the representation. 
To make the effect of the heuristics more transparent, 
labels are included to identify individual pills. 

Three heuristic mechanisms are considered (see \cref{subsec:heuristics} for definitions). 
The first removes pills with very small area ratio (AR), which typically arise when 
a segment has shrunk in radius or length and contributes only marginally to coverage. 
The second discards pills with low unique area ratio (UR), indicating that their 
contribution is largely redundant with other overlapping pills. 
Care must be taken here, since overly strict thresholds may also remove all contributors 
from a structurally relevant region. 
The third heuristic consolidates groups of nearly parallel and adjacent pills by 
retaining one representative per group and reoptimizing the survivors in a subsequent run. 
This reduces redundancy while preserving essential coverage.

Table~\ref{tab:ar-ur-values} lists the AR and UR values of all eighteen pills
before applying any heuristic. The distribution of these values explains the
different outcomes obtained for varying thresholds.

\begin{table}[H]
  \centering
  \begin{tabular}{c|c|c}
    Pill ID & AR & UR \\
    \hline
    0 & 0.855 & 0.3760 \\
    1 & 0.110 & 4.8e-05 \\
    2 & 1.000 & 0.0707 \\
    3 & 0.551 & 5.5e-05 \\
    4 & 0.097 & 9.5e-06 \\
    5 & 0.859 & 0.325 \\
    6 & 0.946 & 0.0729 \\
    7 & 0.967 & 0.0485 \\
    8 & 0.490 & 4.4e-05 \\
    9 & 0.510 & 5.5e-04 \\
    10 & 0.368 & 1.1e-29 \\
    11 & 0.406 & 1.1e-10 \\
    12 & 0.110 & 5.6e-05 \\
    13 & 0.857 & 0.375 \\
    14 & 0.795 & 4.97e-05 \\
    15 & 0.948 & 0.0794 \\
    16 & 0.857 & 0.359 \\
    17 & 0.101 & 1.5e-08 \\
  \end{tabular}
  \caption{Area ratio (AR) and unique area ratio (UR) for all eighteen pills
  prior to applying heuristics.}
  \label{tab:ar-ur-values}
\end{table}

\begin{figure}[H]
    \centering
    \includegraphics[width=0.7\textwidth,height=0.35\textwidth]{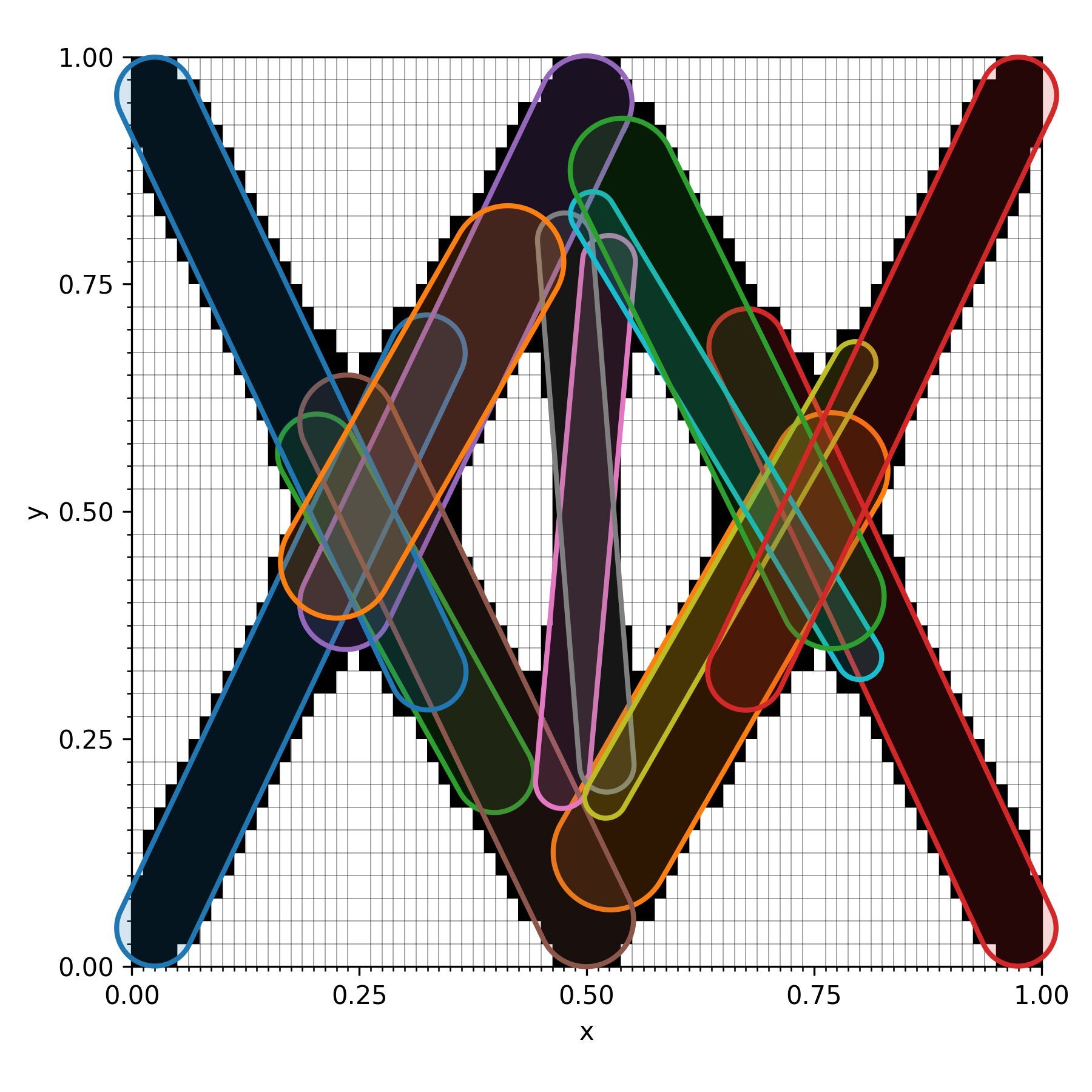}
    \caption{Pill placement after AR-based removal with a threshold of $0.15$. 
    Features with very small AR values (1, 4, 12, 17) are eliminated.}
    \label{fig:heuristics-AR}
\end{figure}

The \cref{fig:heuristics-AR} illustrates the effect of applying AR-based removal
with a threshold of $0.15$. 
This choice targets precisely those pills with very small AR values (IDs 1, 4, 12, 17), 
corresponding to segments of minimal radius or short length. 
Their elimination simplifies the representation without compromising coverage of the target. 

For UR-based removal, the situation is more delicate. 
Thresholds of $10^{-4}$ and $10^{-3}$ were tested. 
In both cases, several pills with extremely small UR values are removed, 
including IDs 1, 3, 4, 8, 10, 11, 12, 14 and 17. 
At the more aggressive threshold of $10^{-3}$, even pill~9 is eliminated. 
The consequences are most apparent in the central vertical bar: 
with the higher threshold, both overlapping candidates are discarded, 
leaving this part of the target uncovered. 
With the more conservative threshold, one pill remains in place and continues to capture the region. 
This demonstrates that UR-based removal can easily become too aggressive if applied without care, 
since overlapping pills are penalized simultaneously. 

\begin{figure}[H]
  \centering
  \includegraphics[width=0.48\textwidth,height=0.24\textwidth]{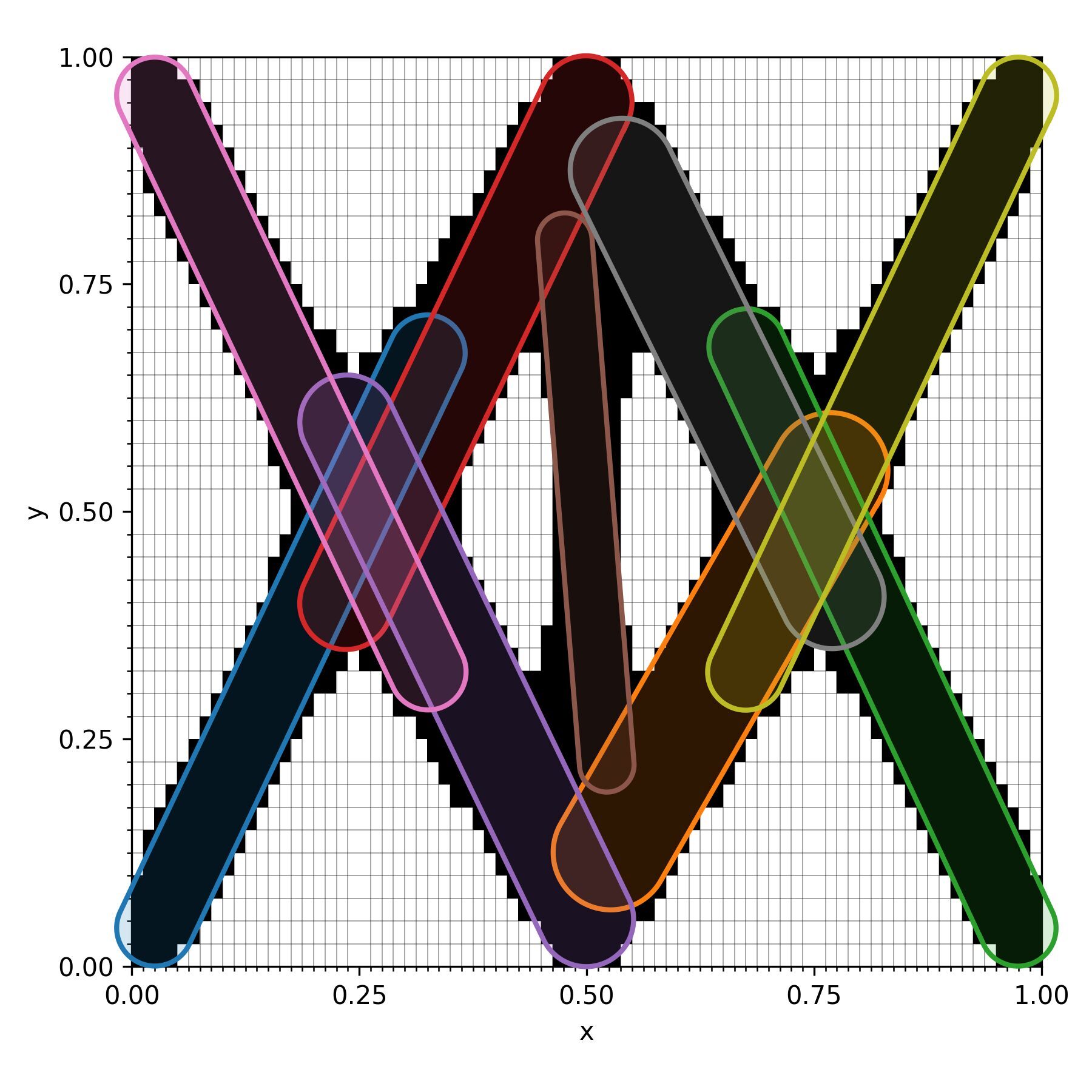}
  \includegraphics[width=0.48\textwidth,height=0.24\textwidth]{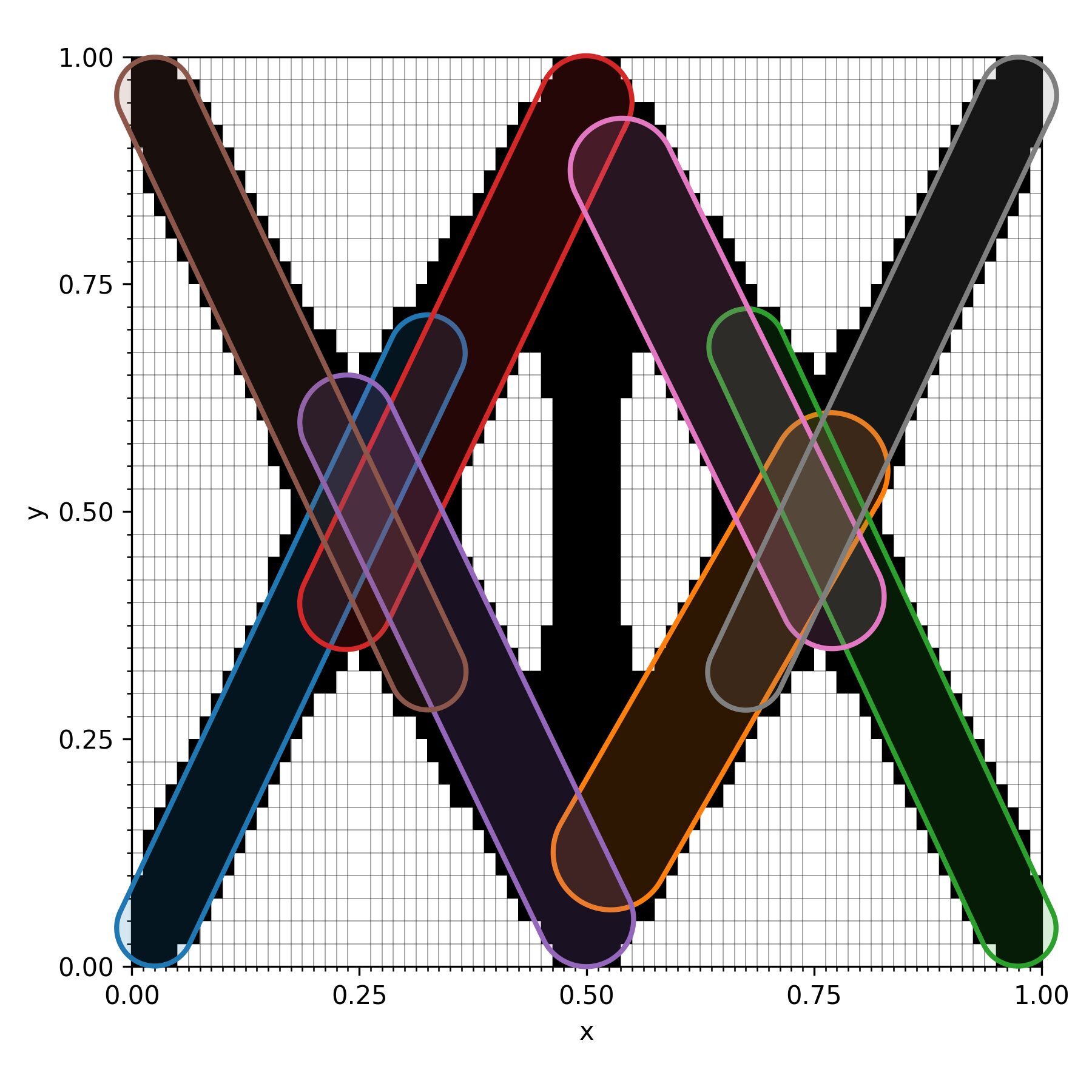}
  \caption{Feature placement after UR-based removal with thresholds of
  $10^{-4}$ (left) and $10^{-3}$ (right). 
  The aggressive threshold removes both candidates covering the central vertical bar, leaving it uncovered.}
  \label{fig:heuristics-UR-heuristics}
\end{figure}

Finally, the effect of pill consolidation is considered. 
Groups of pills that are nearly parallel and in close spatial proximity 
are identified using an angular tolerance of $\theta_{\mathrm{lim}}=10^\circ$ 
and a minimum separation distance of $d_{\min}=0.15$. 
For the configuration with eighteen pills, five such groups are obtained:
\[
\{0,6,14\}, \quad \{2,10,17\}, \quad \{3,7,13\}, \quad \{8,9\}, \quad \{11,15\}.
\]

Within each group, the members exhibit similar exploration patterns 
and are positioned so closely that their endpoints nearly touch. 
In such constellations, a single longer segment can often approximate 
the combined extent of the group. 
The consolidation heuristic therefore retains the longest element in each group 
and removes the companions; no explicit geometric fusion is performed. 
During the subsequent convergence run, 
the surviving representatives adapt their positions and radii, 
and thereby recover most of the coverage previously provided by the full group, 
up to differences caused by radius variation or slight angular offsets. 

\begin{figure}[H]
  \centering
  \includegraphics[width=0.75\textwidth,height=0.375\textwidth]{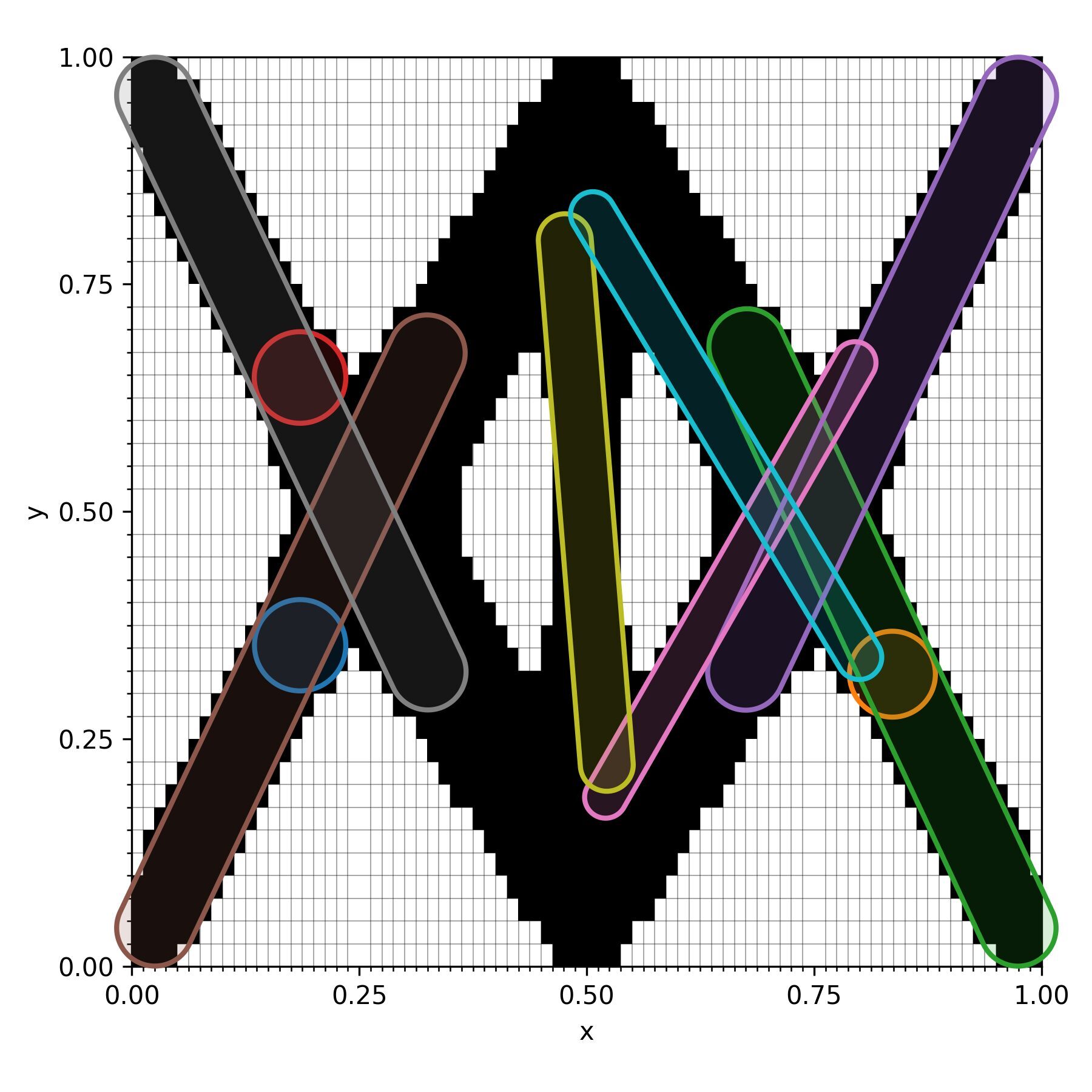}
  \caption{Feature placement after applying the combine heuristic with
  $\theta_{\mathrm{lim}}=10^\circ$ and $d_{\min}=0.15$. 
  In each group of nearly parallel and adjacent pills, the longest element is retained 
  while the others are removed. 
  The surviving representatives are reoptimized in a convergence run, 
  thereby reproducing the coverage of the full group with fewer pills.}
  \label{fig:heuristics-combine}
\end{figure}

The resulting pill placement is shown in \cref{fig:heuristics-combine}. 
For example, in the nearly vertical group $\{0,6,14\}$ only the longest element is retained, 
while for the pairs $\{8,9\}$ and $\{11,15\}$ one representative remains in each case. 
This reduces the number of pills while largely preserving target coverage; 
some local flexibility is traded for compactness. 

After the application of AR-, UR-, or combined heuristics, 
a convergence stage was performed against the full target 
to reoptimize the surviving representatives 
and restore consistency with the objective. 
For UR and its more aggressive variant UR\_agg, 
the converged feature allocations are presented together with their objective values. 
UR attained a final objective of $4.73 \times 10^{1}$, 
whereas UR\_agg converged to $2.55 \times 10^{2}$. 
For comparison, the pre-heuristic configuration achieved a lower objective of $3.50 \times 10^{1}$, 
but only at the expense of a substantially larger number of pill features. 
These results indicate that heuristic consolidation reduces feature count 
while maintaining competitive convergence behavior, albeit with trade-offs 
between objective quality and model compactness.

\begin{figure}[H]
  \centering
  \includegraphics[width=0.48\textwidth,height=0.24\textwidth]{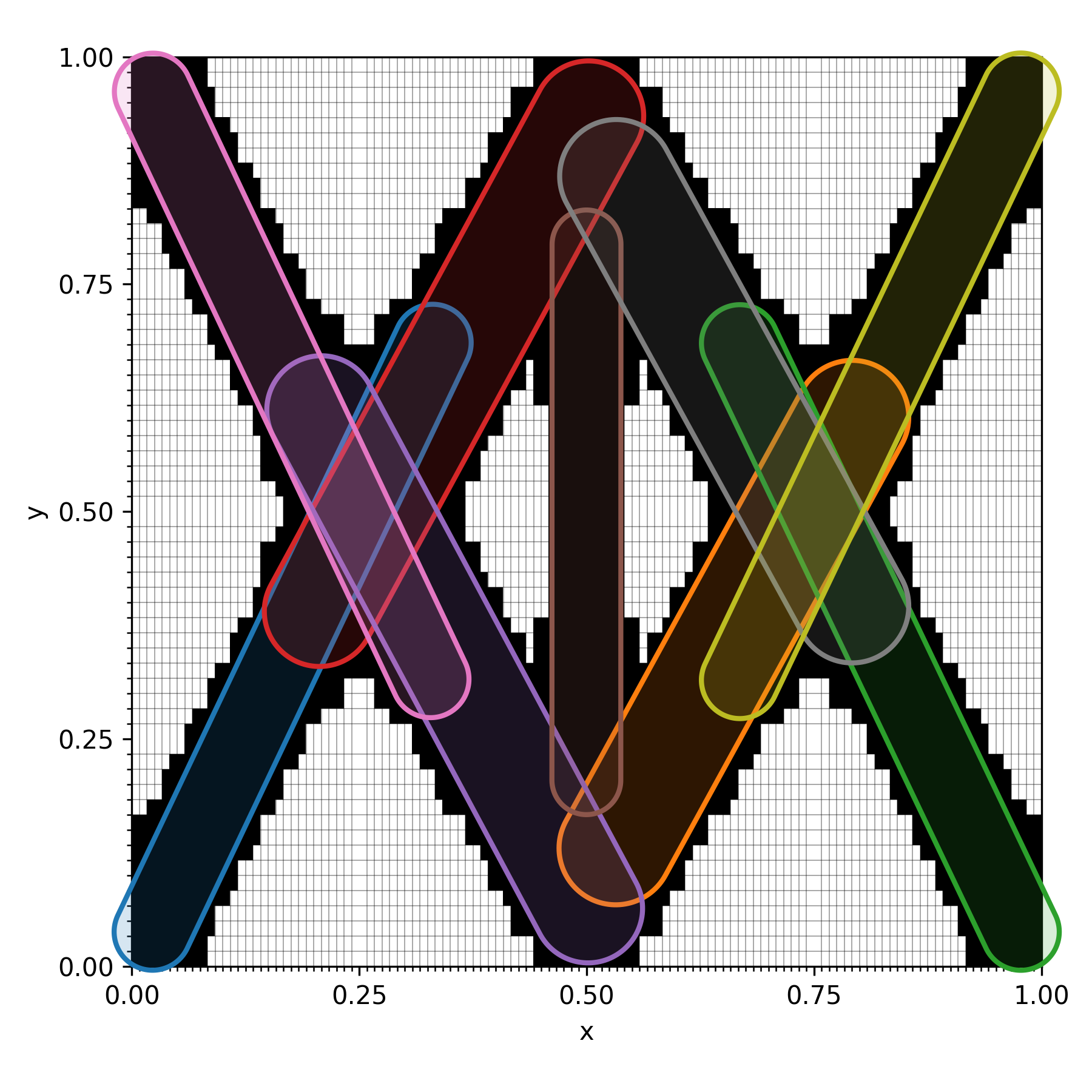}
  \includegraphics[width=0.48\textwidth,height=0.24\textwidth]{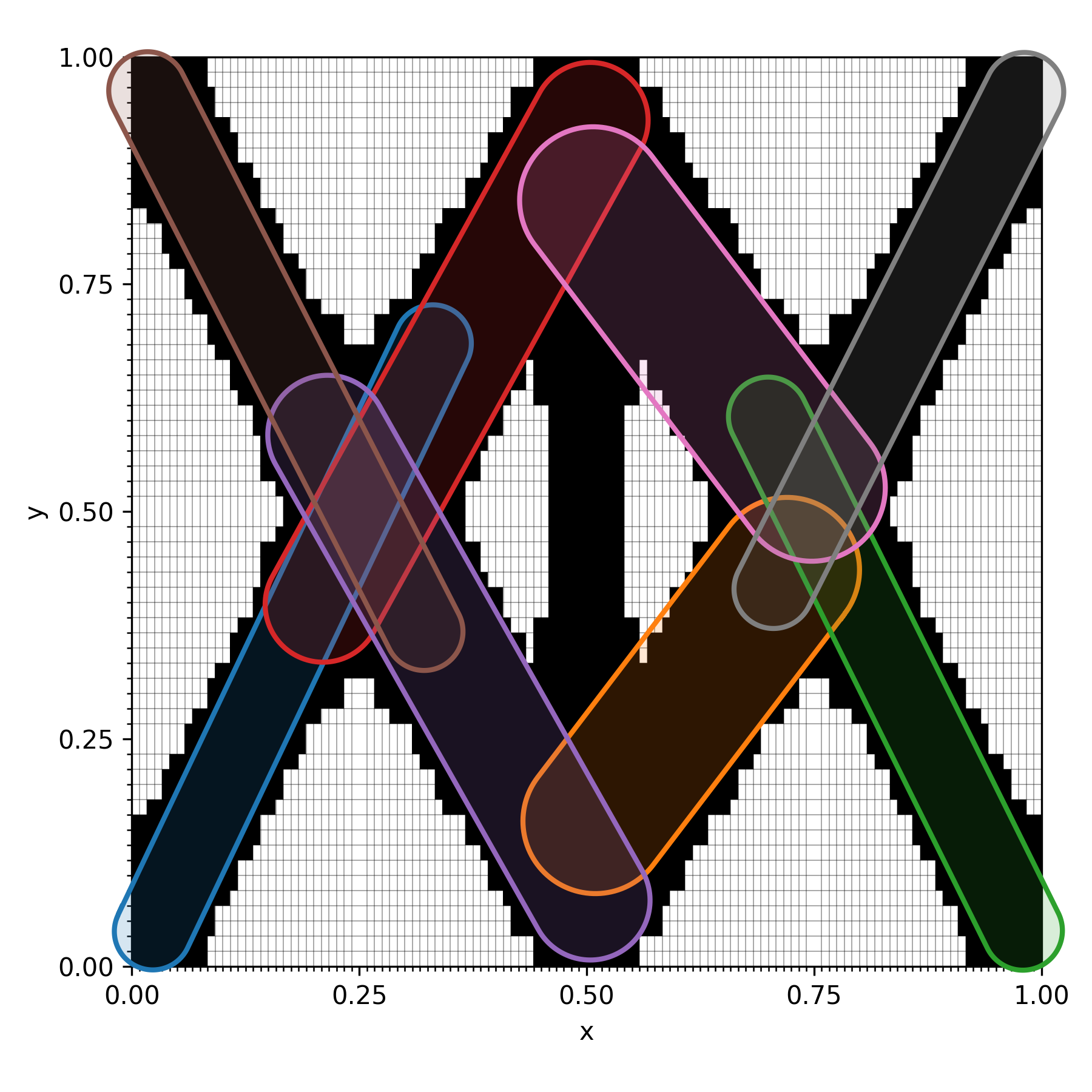}
  \caption{Final feature allocations after convergence for UR and UR\_agg heuristics. 
  The corresponding objective values are $4.73 \times 10^{1}$ (UR) and $2.55 \times 10^{2}$ (UR\_agg).}
  \label{fig:heuristics-UR-final}
\end{figure}

% ============================================================
% 8) Cantilever Test
% ============================================================

\subsection{Cantilever test}
\label{sec:cantilever}

\begin{figure}[H]
    \centering
    \includegraphics[width=0.6\textwidth]{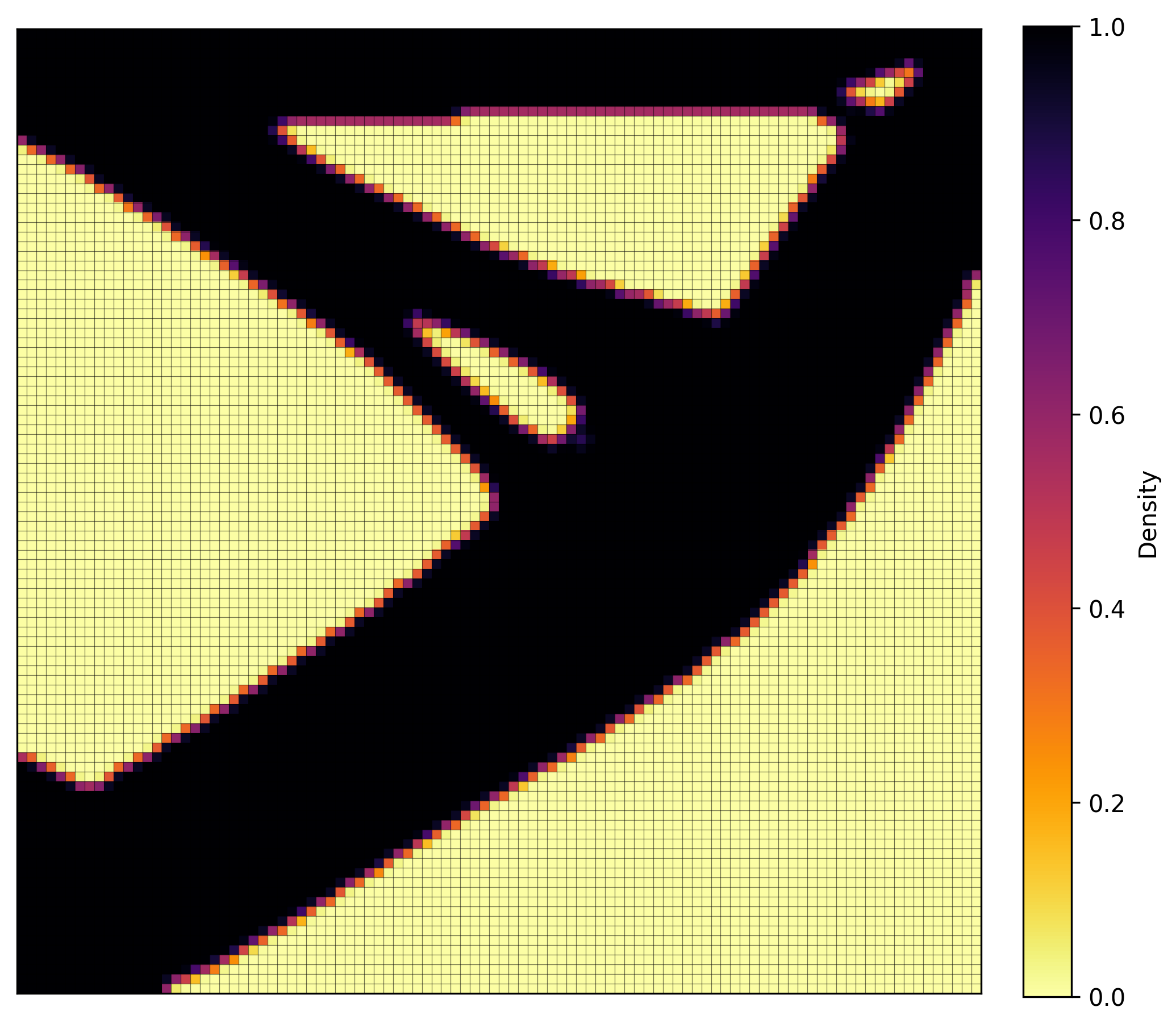}
    \caption{Target density for the cantilever beam, generated with SIMP.}
    \label{fig:cantilever_target}
\end{figure}

To demonstrate generalization beyond the five--bar benchmark, 
the reconstruction framework is applied to a classical cantilever beam. 
The domain is discretized with $100\times 100$ finite elements under plane--stress conditions; 
displacements are fixed along the left boundary, 
and a unit point load acts downward at the free right edge. 
The target density distribution is generated with the \textsc{CFS++} framework using SIMP, 
as shown in \cref{fig:cantilever_target}. 

Both aggregation rules established in the previous sections are tested: 
the $p$--norm scheme with exponent $p=9$ 
and the softmax scheme with parameter $\beta=18$. 
All other settings follow the defaults defined earlier: 
cross--initialization with eighteen pills of radius $r=0.05$, 
third--order quadrature and staged optimization with exploration, bridging and convergence phases. 
Heuristics are applied with thresholds $\mathrm{AR}_{\min}=0.15$ and $\mathrm{UR}_{\min}=10^{-3}$, 
while iterative refinement is governed by the residual threshold $\rho=0.2$ 
and improvement criteria $\Delta_{\mathrm{rel}}>10^{-3}$, $\Delta_{\mathrm{abs}}>0.0$. 

\begin{figure}[H]
  \centering
  \setlength{\tabcolsep}{2pt}
  \renewcommand{\arraystretch}{0.5}
  \begin{tabular}{cccc}
    \includegraphics[width=0.24\textwidth]{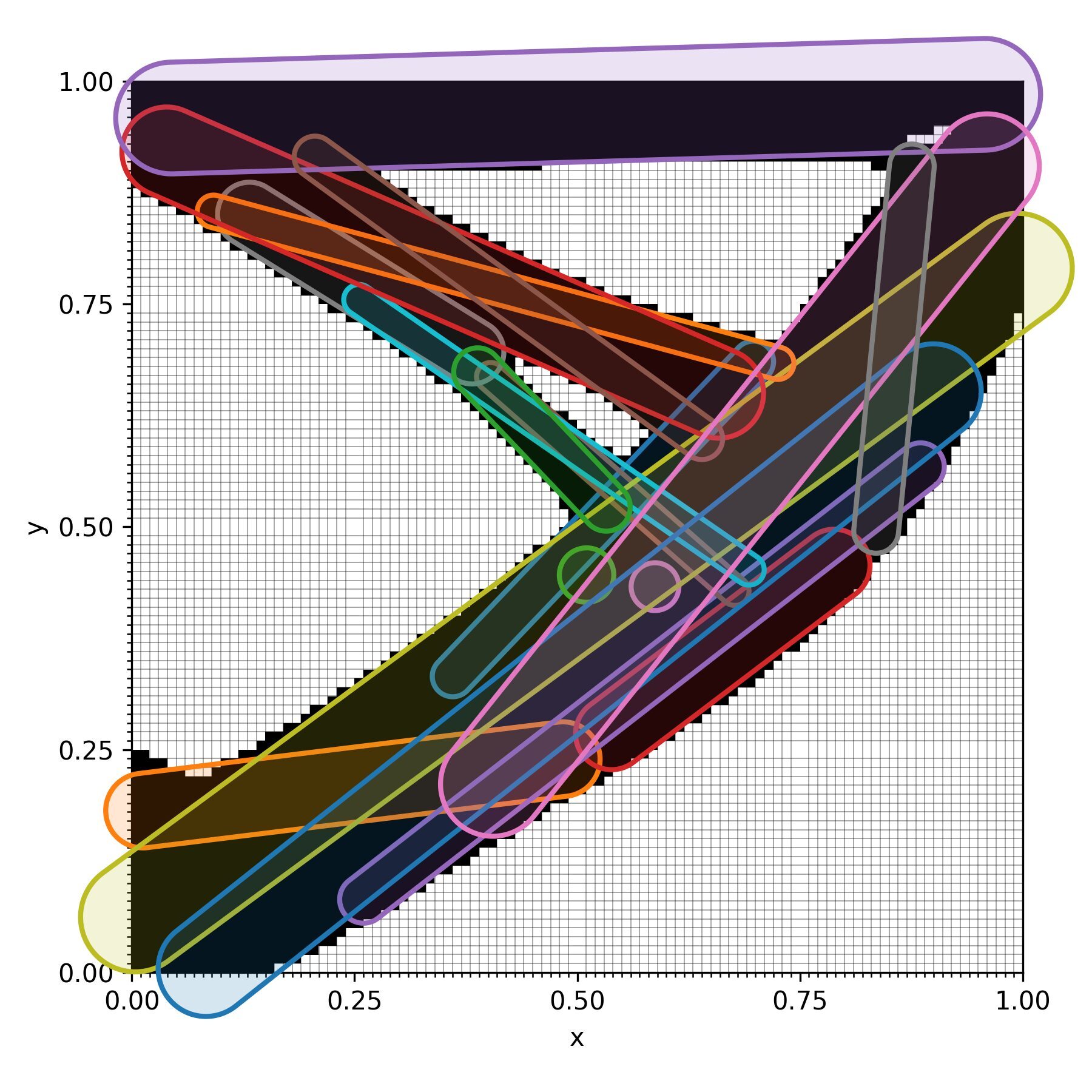} &
    \includegraphics[width=0.24\textwidth]{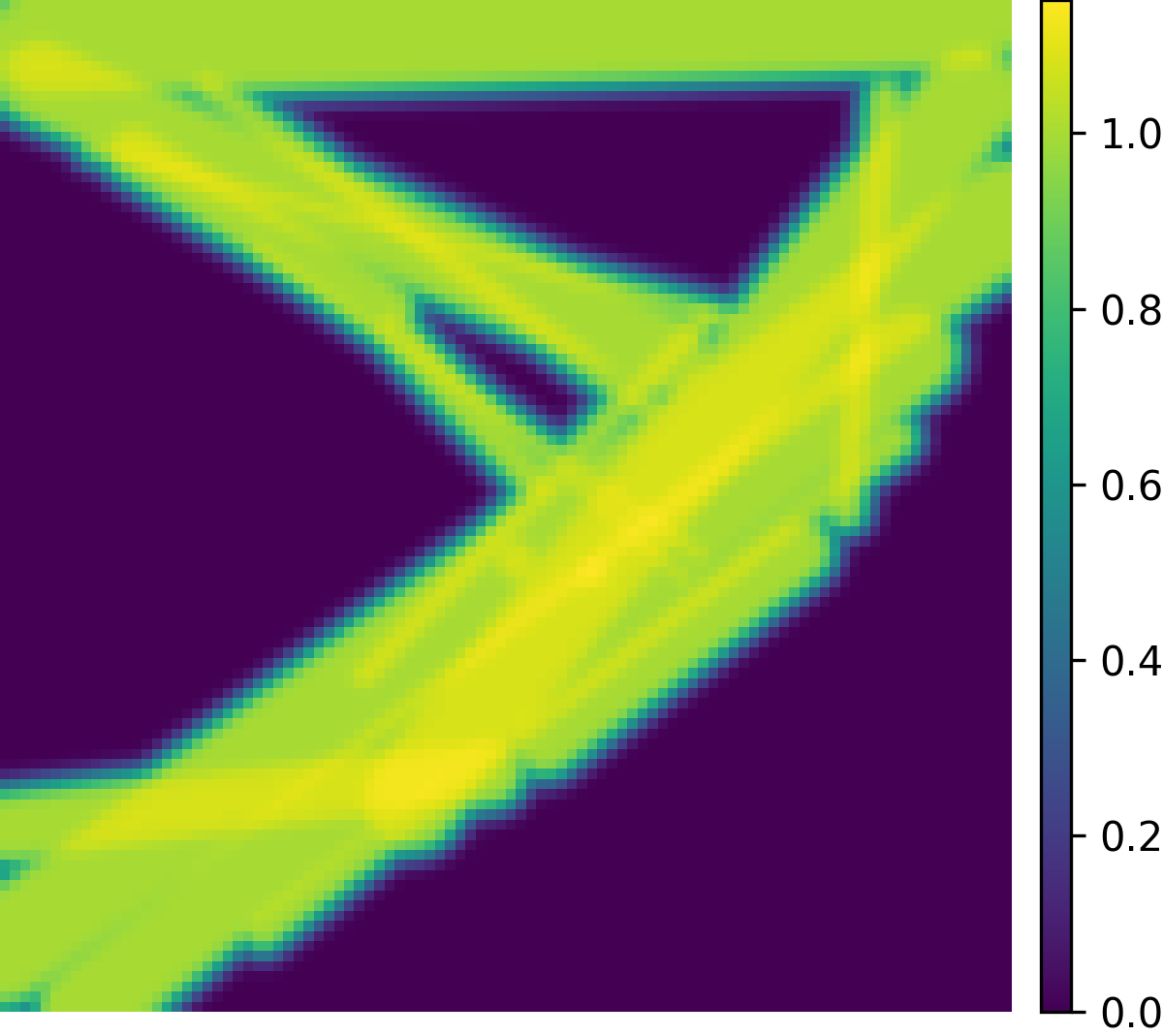} &
    \includegraphics[width=0.24\textwidth]{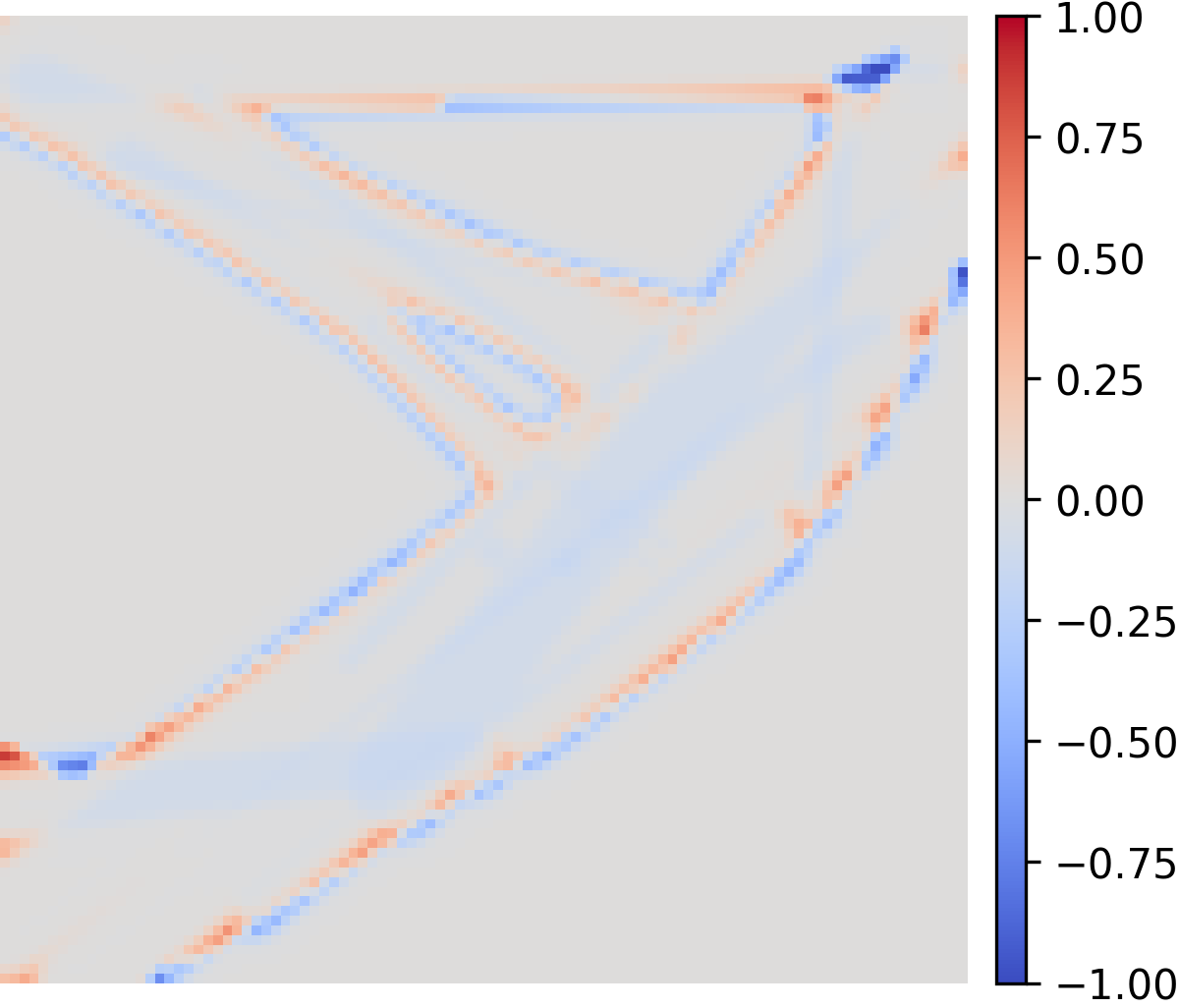} &
    \includegraphics[width=0.24\textwidth]{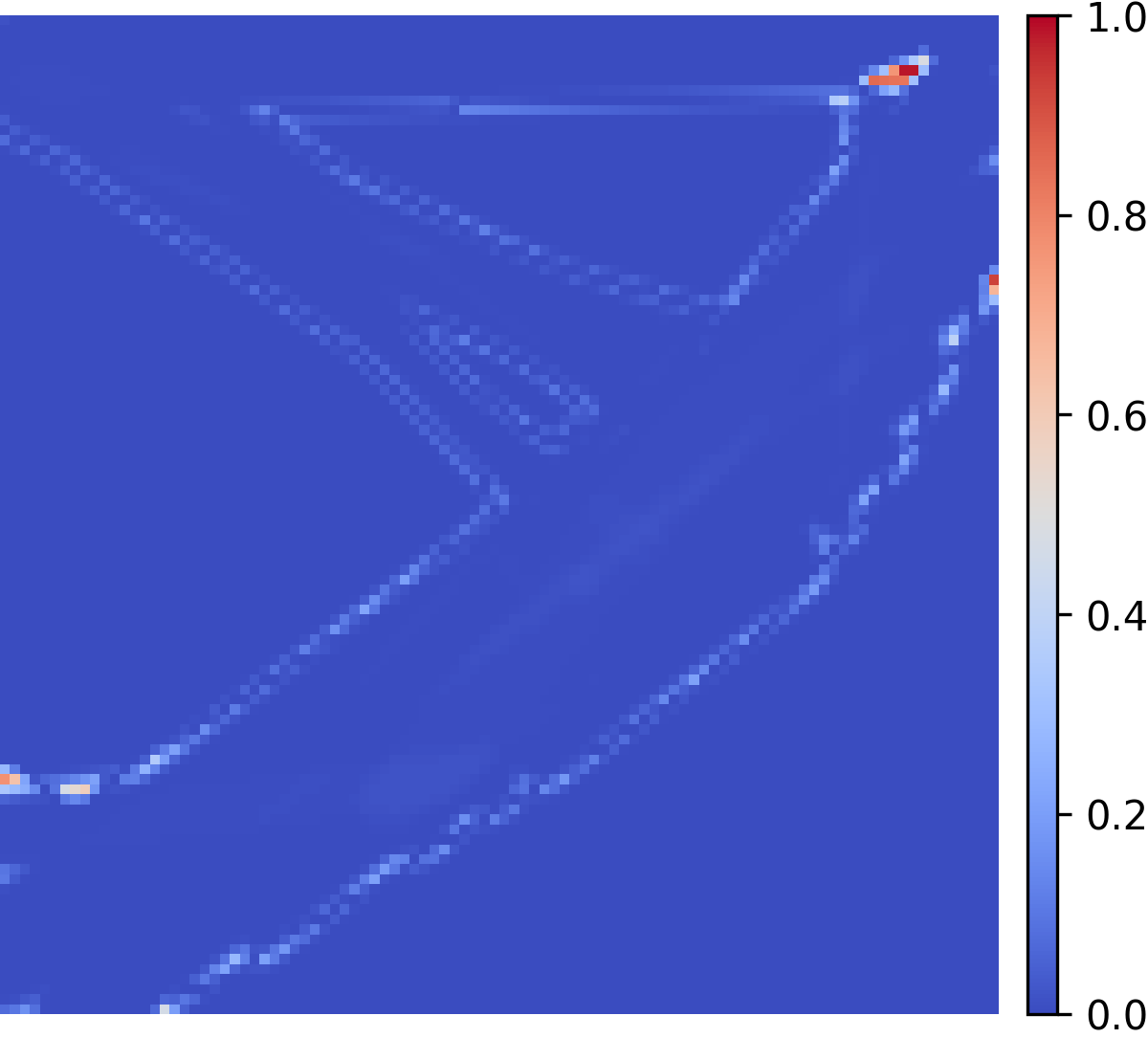} \\
    \includegraphics[width=0.24\textwidth]{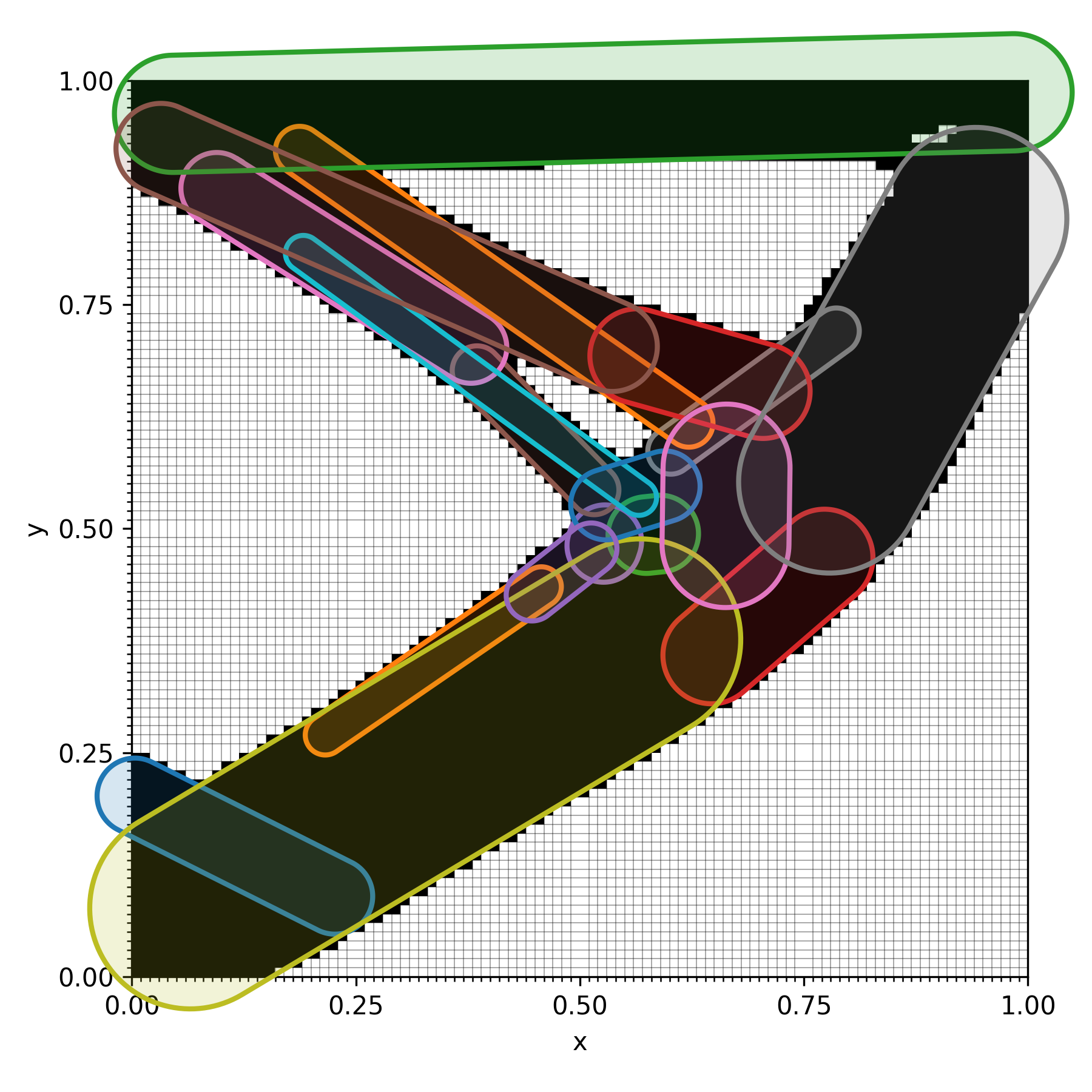} &
    \includegraphics[width=0.24\textwidth]{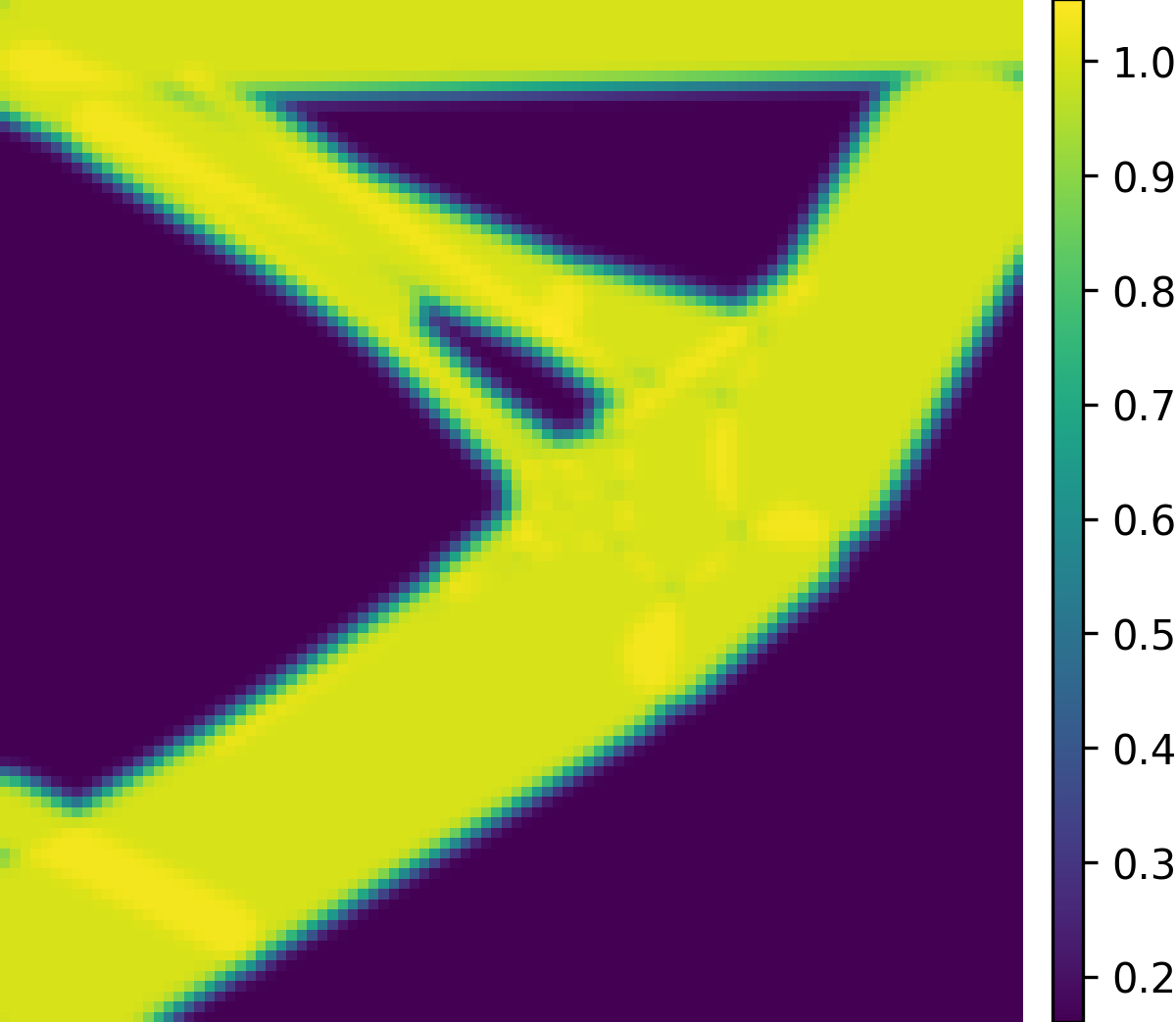} &
    \includegraphics[width=0.24\textwidth]{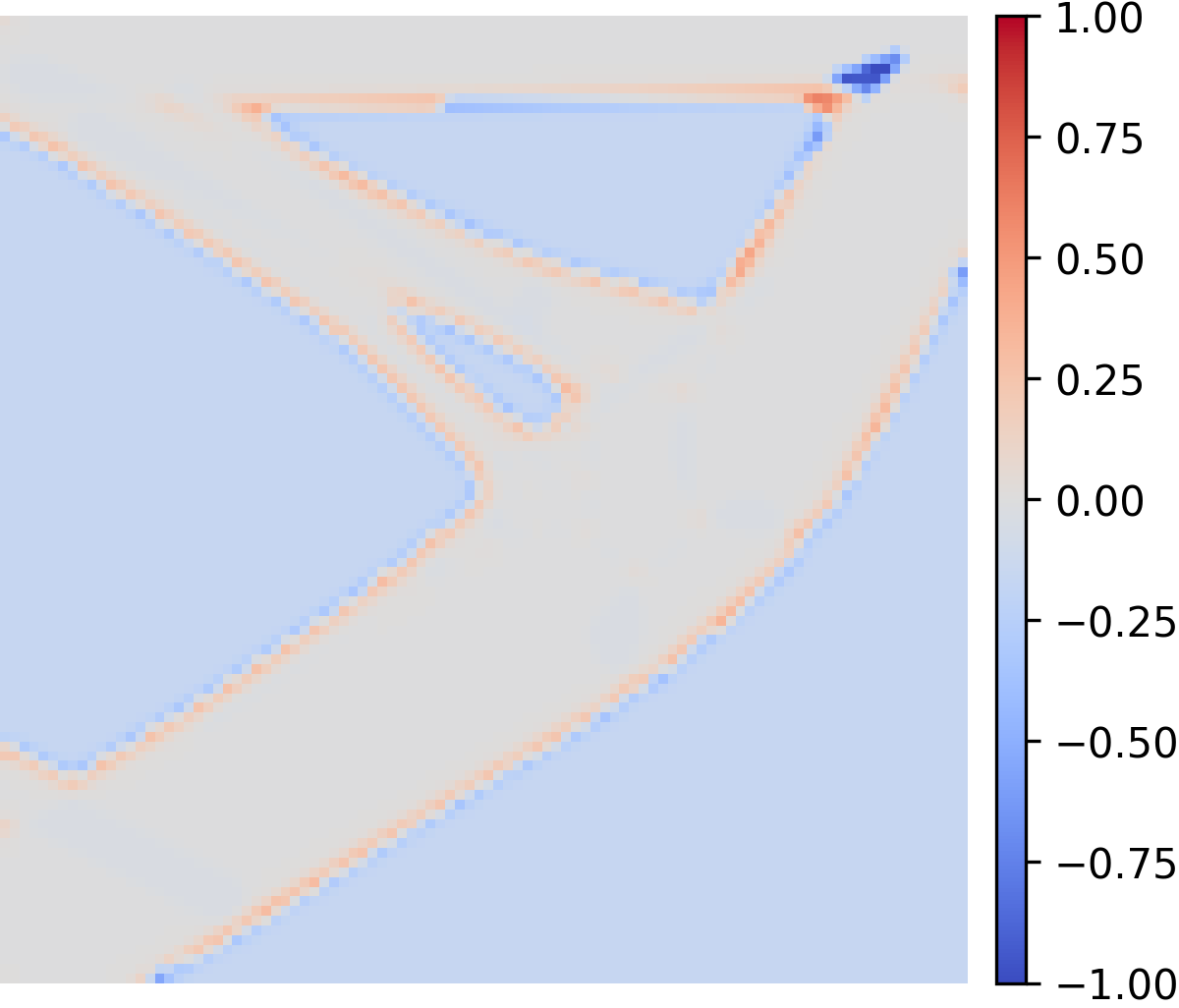} &
    \includegraphics[width=0.24\textwidth]{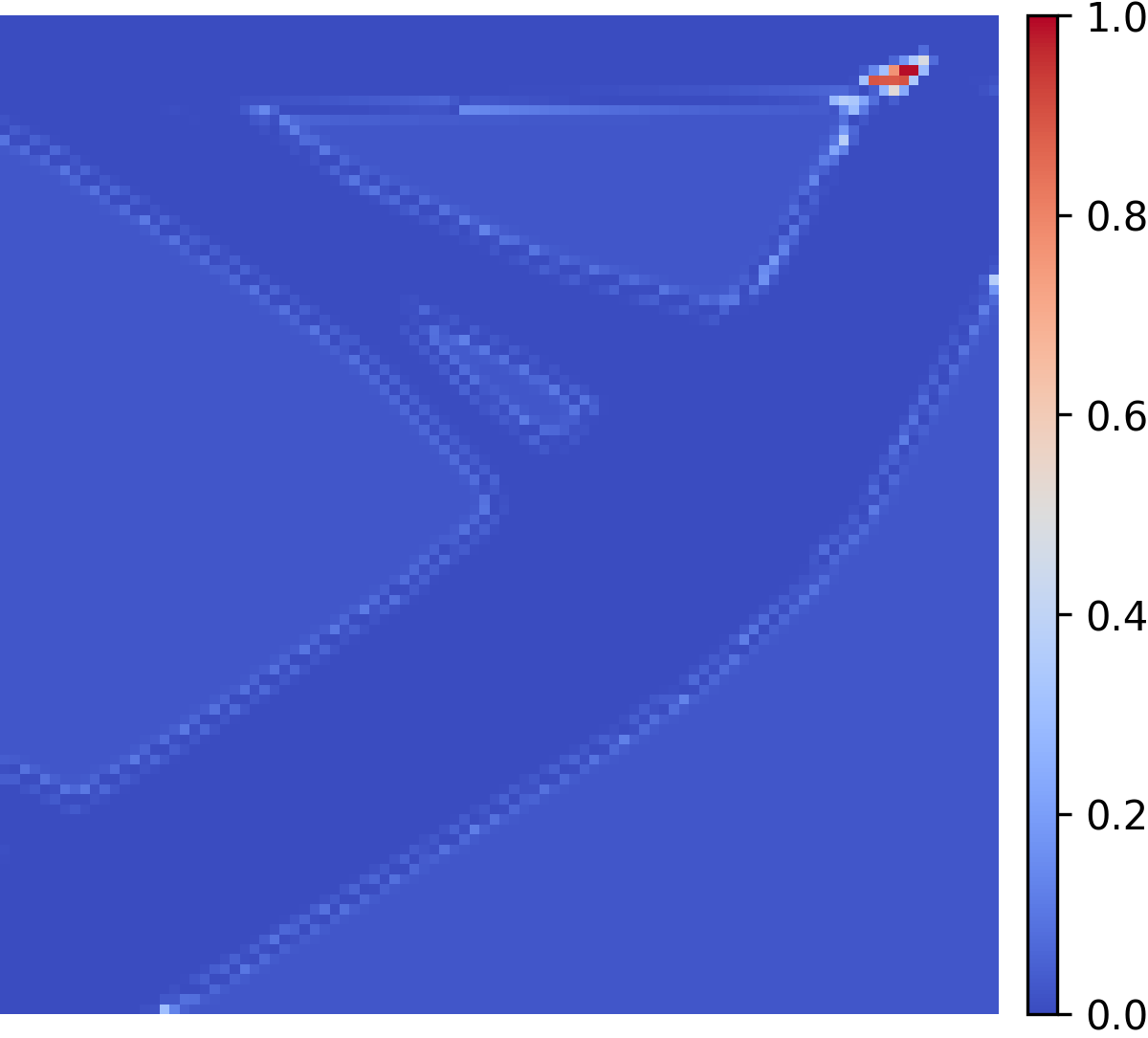} \\
  \end{tabular}
  \caption{Cantilever reconstructions with $p$--norm (top) and softmax
  (bottom) aggregation prior to heuristics. 
  Shown are pill placement, density field, residual and absolute residual.}
  \label{fig:canti-final}
\end{figure}

The cantilever reconstructions with $p$--norm and softmax aggregation 
are first considered in their final pre--heuristic form (\cref{fig:canti-final}). 
Both pipelines recover the main load--carrying bars and capture the overall target layout. 
In the $p$--norm case, the high exponent $p=9$ does not penalize overlap strongly, 
so that the large curved bar from the lower left support to the upper right 
is realized by several partially parallel and overlapping pills. 
This redundancy is functional: the varying width and curvature of the bar 
require multiple pills for accurate approximation. 
Coverage is complete apart from a small mismatch at the lower left corner, 
where one pill simultaneously attempts to represent a small density pocket 
and part of the adjacent curvature. 
This prevents correct alignment and leaves a visible discrepancy. 

In contrast, the softmax variant relies on fewer overlaps 
and exhibits greater variation in pill sizes. 
A few large pills dominate the structure, 
while smaller ones refine curved or otherwise uncovered regions. 
This yields a smoother and more uniform density distribution, 
with fewer irregularities in the aggregated field. 

\begin{figure}[H]
  \centering
  \setlength{\tabcolsep}{2pt}
  \renewcommand{\arraystretch}{0.5}
  \begin{tabular}{cccc}
    \includegraphics[width=0.24\textwidth]{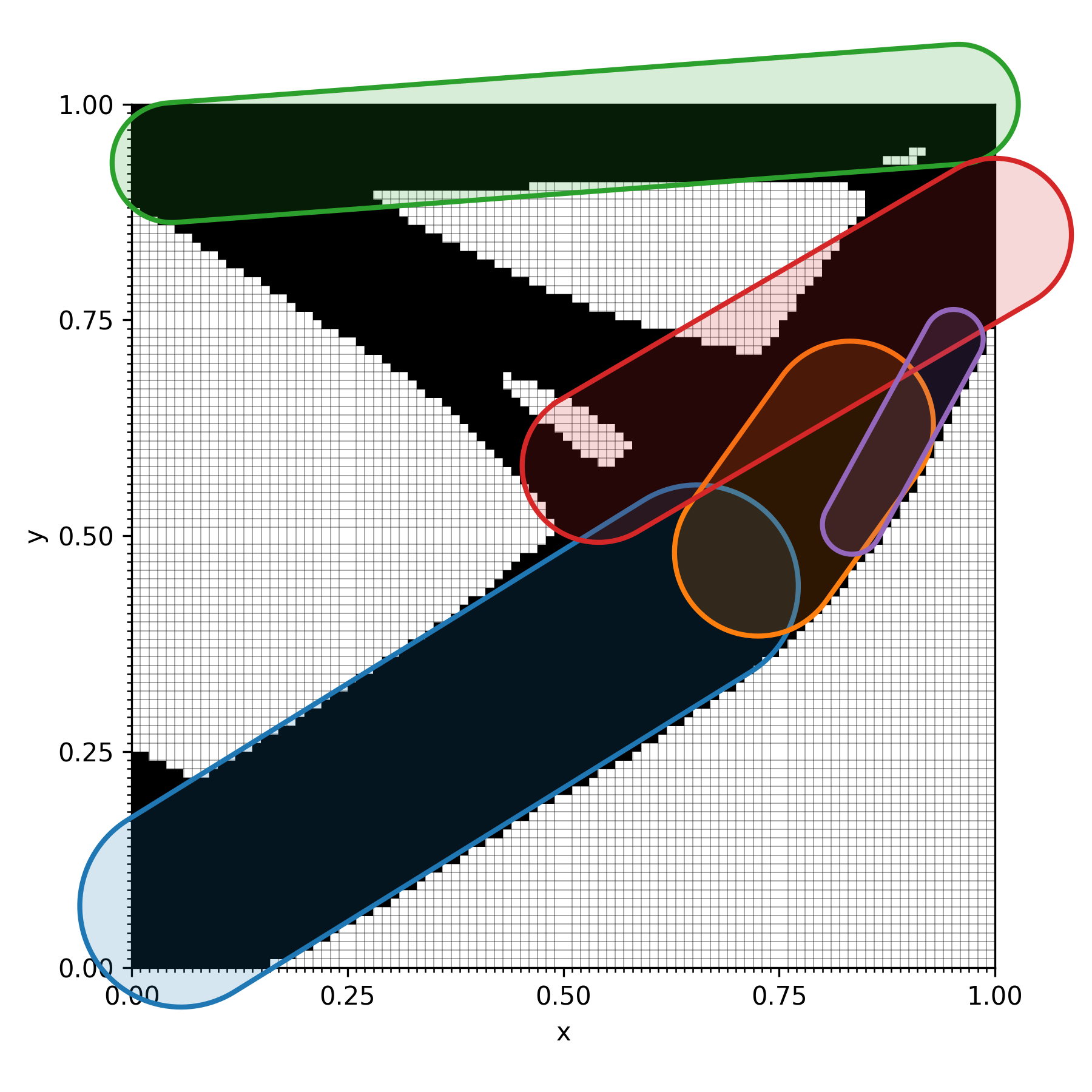} &
    \includegraphics[width=0.24\textwidth]{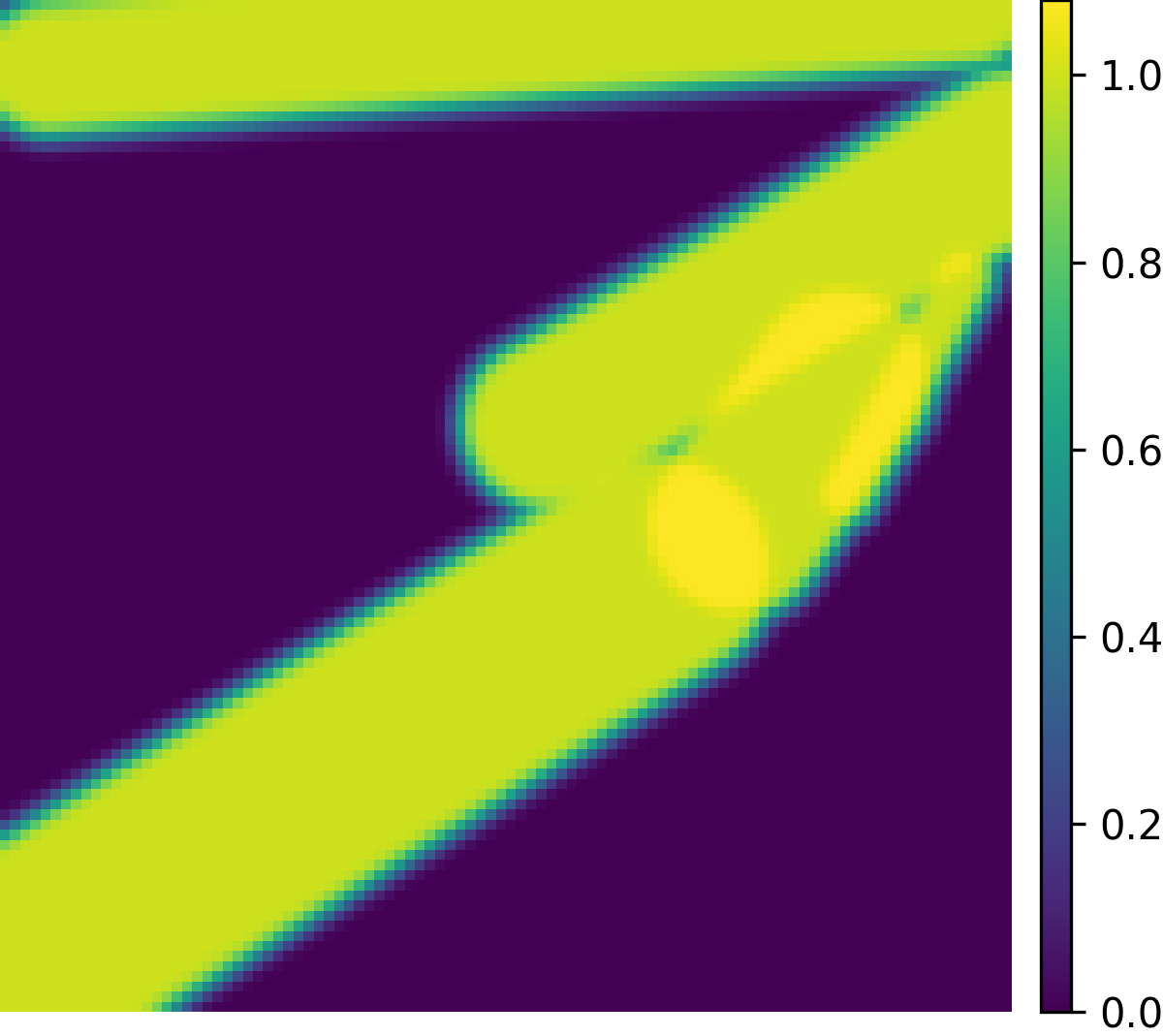} &
    \includegraphics[width=0.24\textwidth]{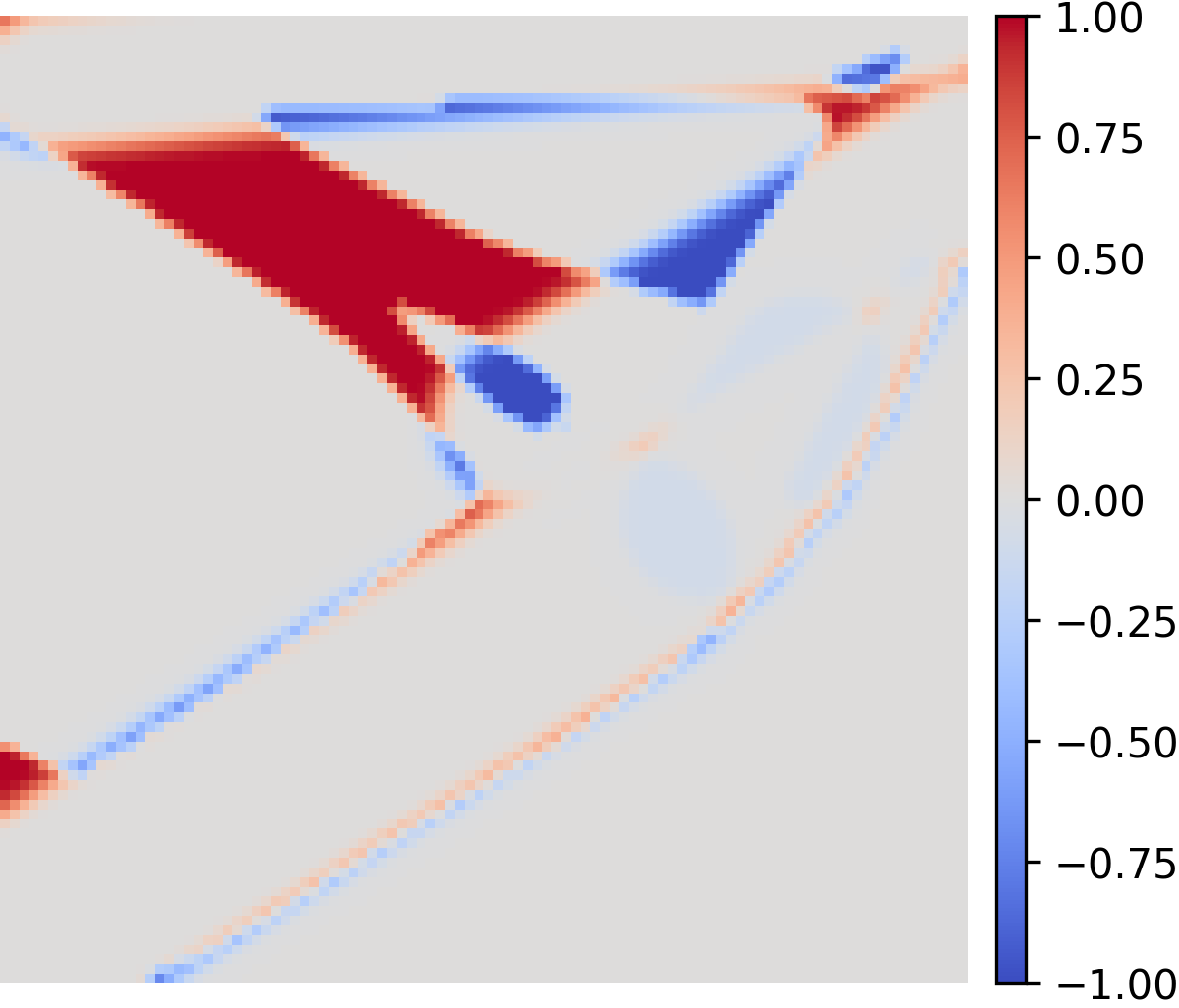} &
    \includegraphics[width=0.24\textwidth]{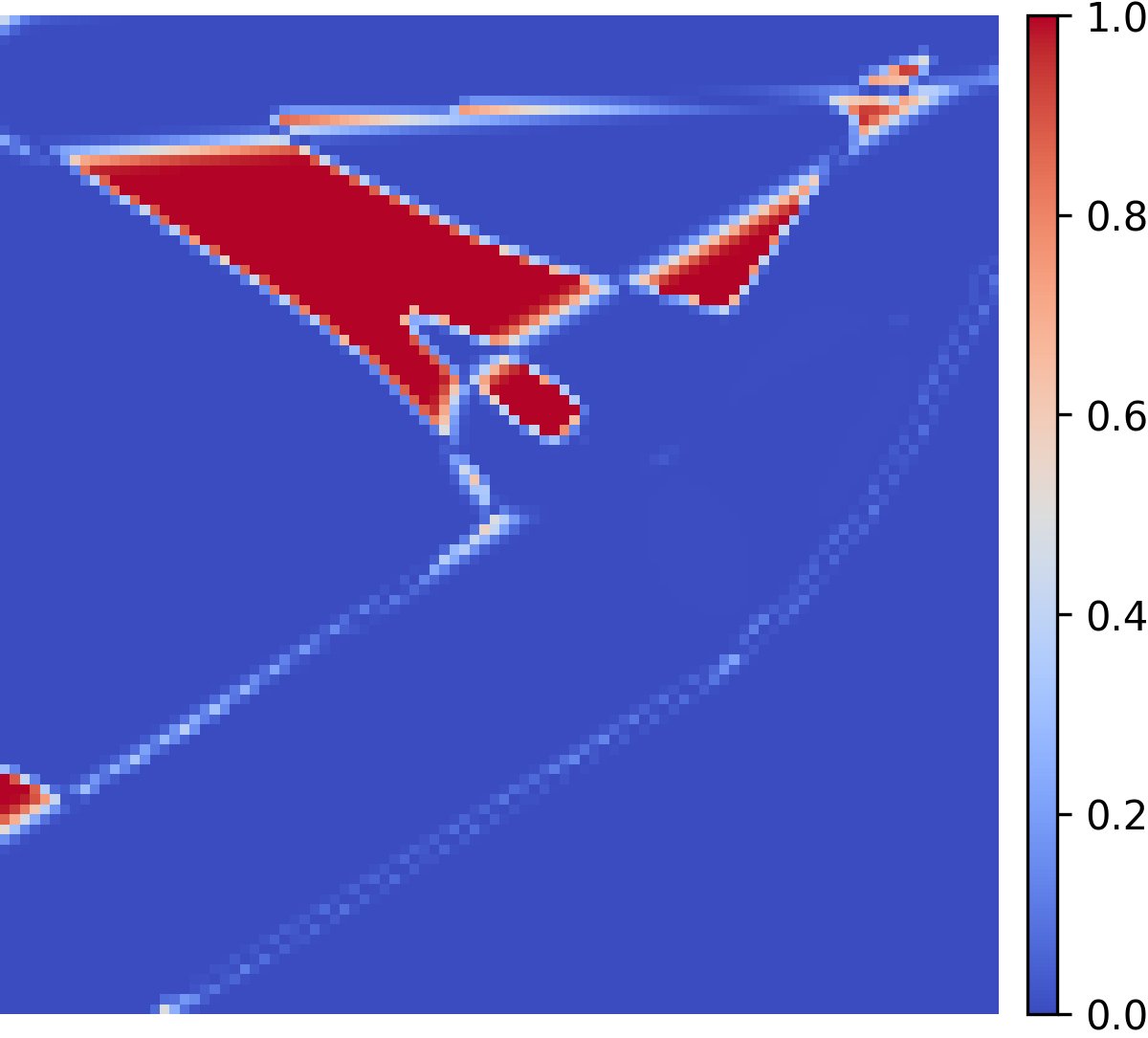} \\
    \includegraphics[width=0.24\textwidth]{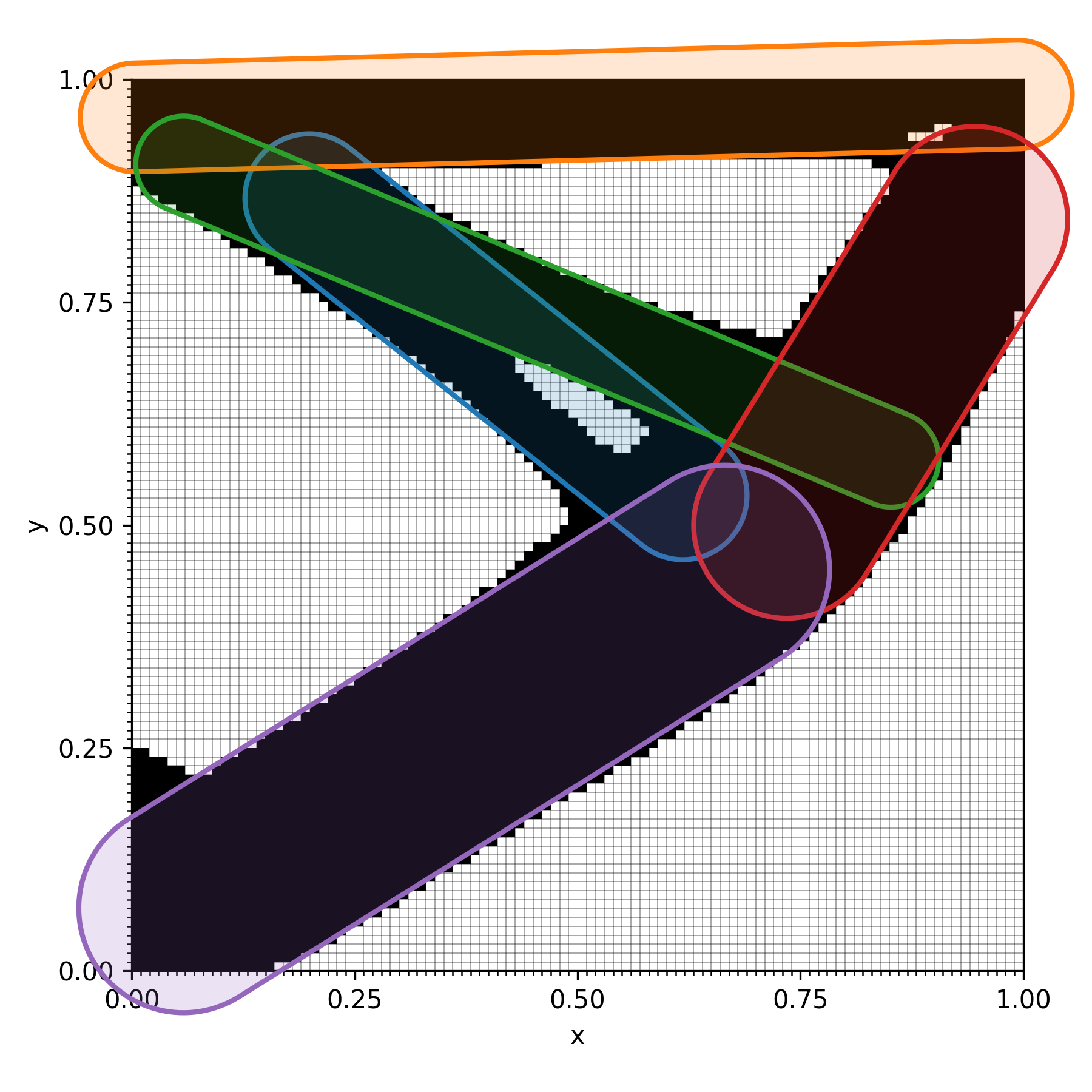} &
    \includegraphics[width=0.24\textwidth]{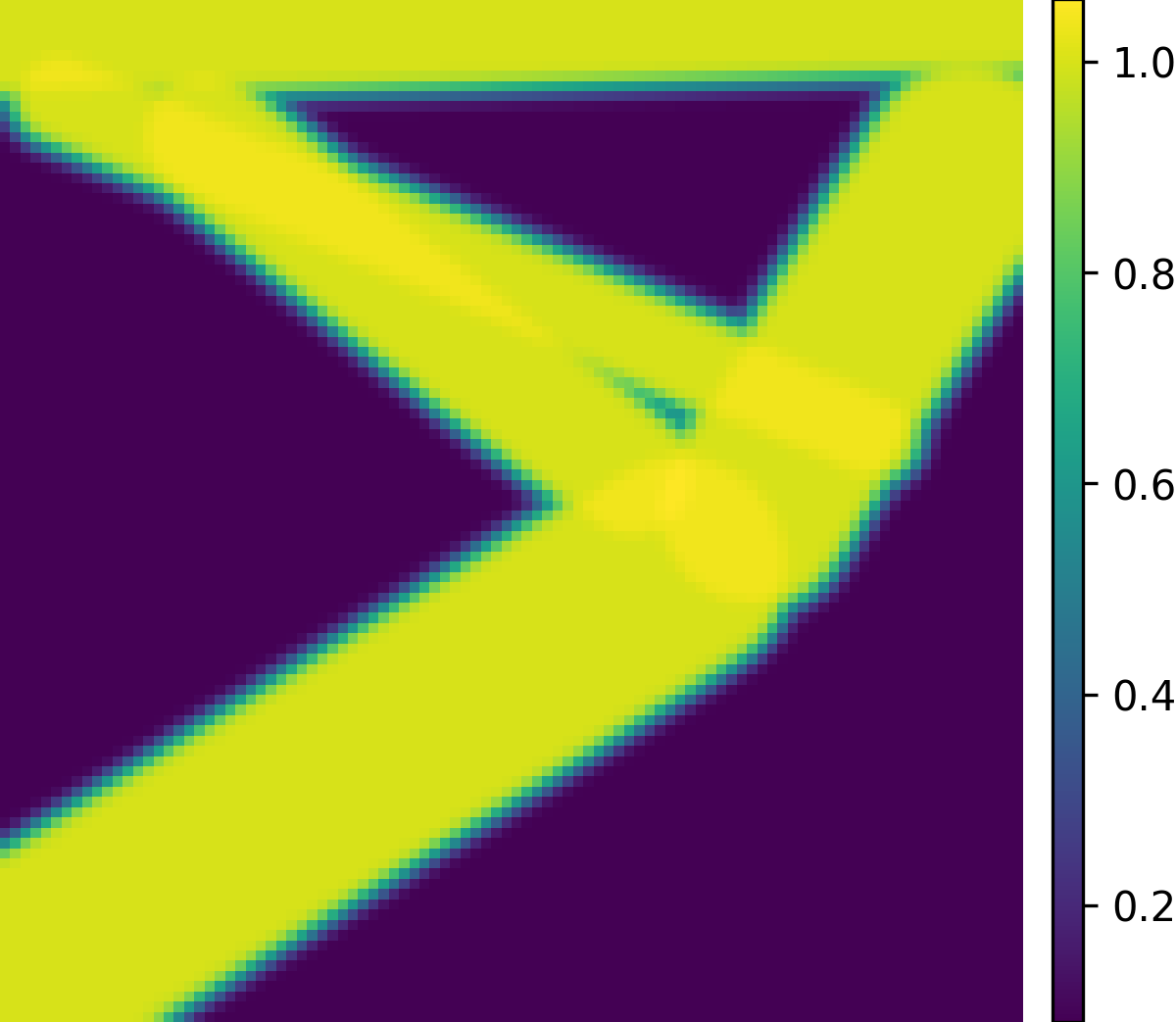} &
    \includegraphics[width=0.24\textwidth]{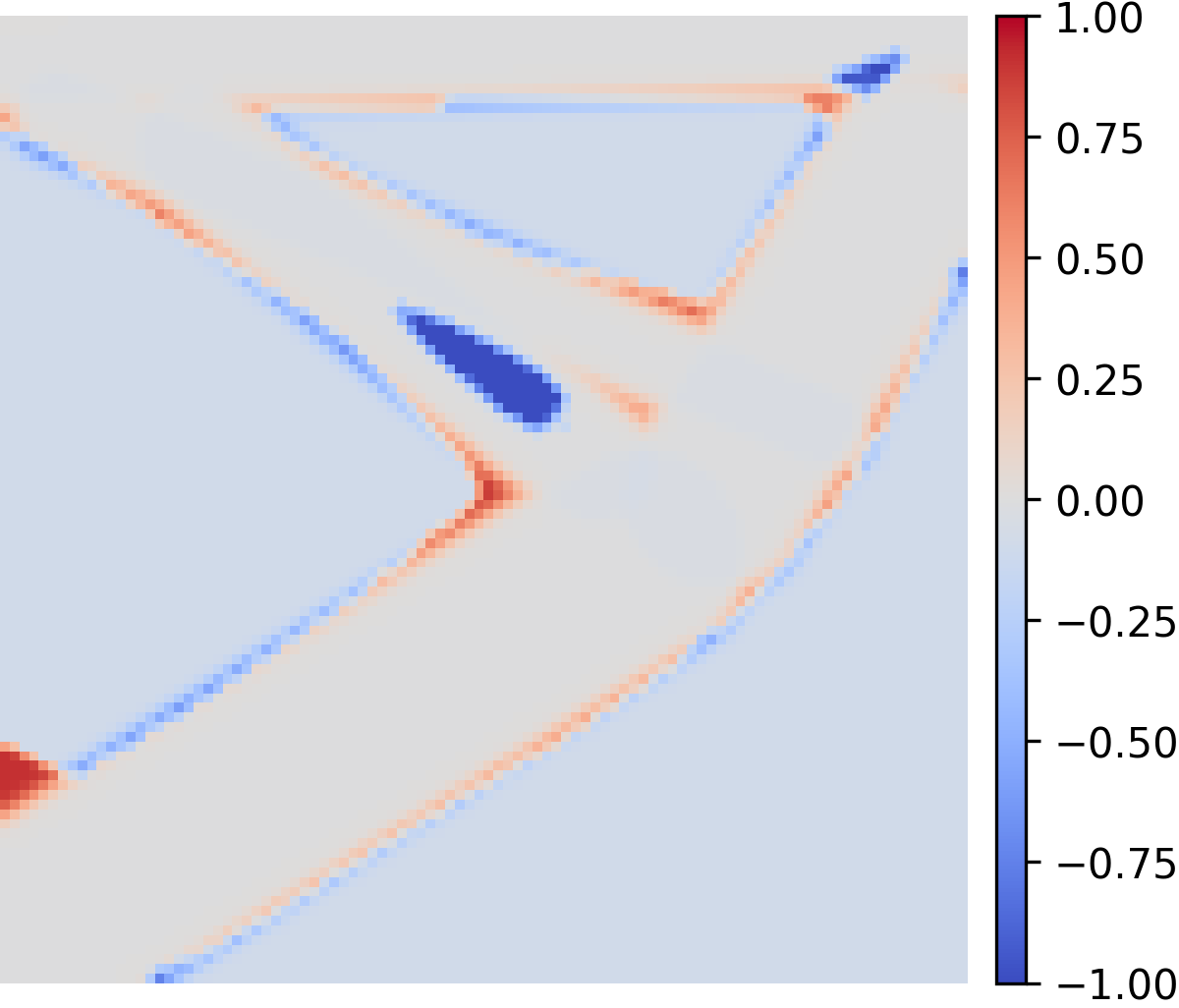} &
    \includegraphics[width=0.24\textwidth]{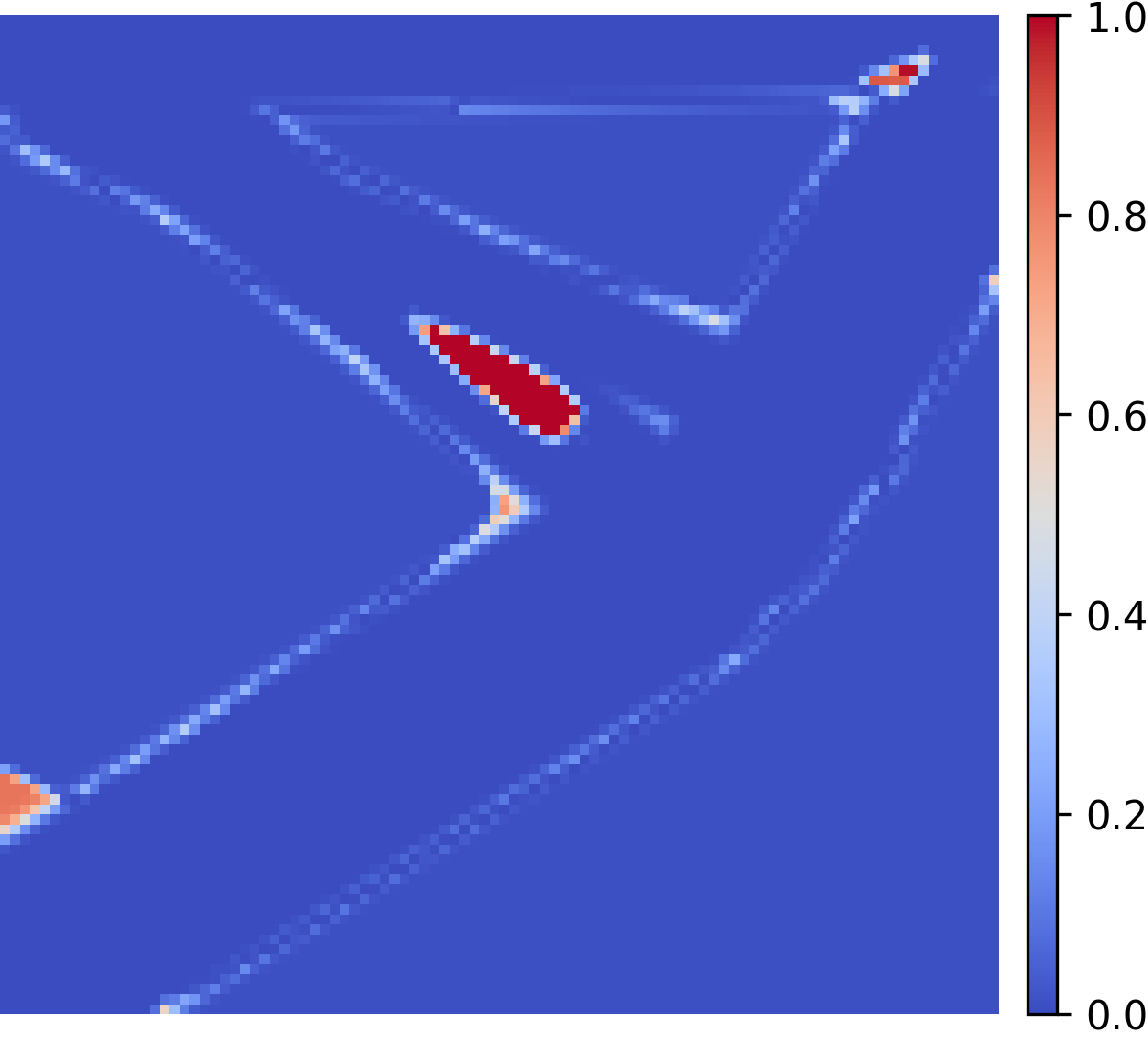} \\
  \end{tabular}
  \caption{Reconstructions after heuristic pruning with $p$--norm (top) and
  softmax (bottom). 
  The UR threshold removes auxiliary pills, 
  leading to coverage losses particularly in the central and upper regions.}
  \label{fig:canti-heuristics}
\end{figure}

When heuristics are applied, both reconstructions degrade, albeit in different ways 
(\cref{fig:canti-heuristics}). 
For the $p$--norm case, aggressive UR thresholds remove several pills, 
including the orange segment at the lower left. 
While this pruning simplifies the representation, 
it temporarily leaves parts of the target underrepresented. 
For softmax, the consequences are more pronounced: 
many of the smaller auxiliary pills around the central joint and in curved regions are discarded. 
Coverage that had previously been achieved by coordinated smaller pills 
is then replaced by a single oversized segment. 
This leads to gaps and distortions, particularly in the middle of the domain 
and near the upper right joint, where fidelity to the target geometry is visibly reduced. 
The issue arises from setting the UR threshold too high, 
which removes not only redundant but also structurally relevant pills. 

\begin{figure}[H]
  \centering
  \setlength{\tabcolsep}{2pt}
  \renewcommand{\arraystretch}{0.5}
  \begin{tabular}{cccc}
    \includegraphics[width=0.24\textwidth]{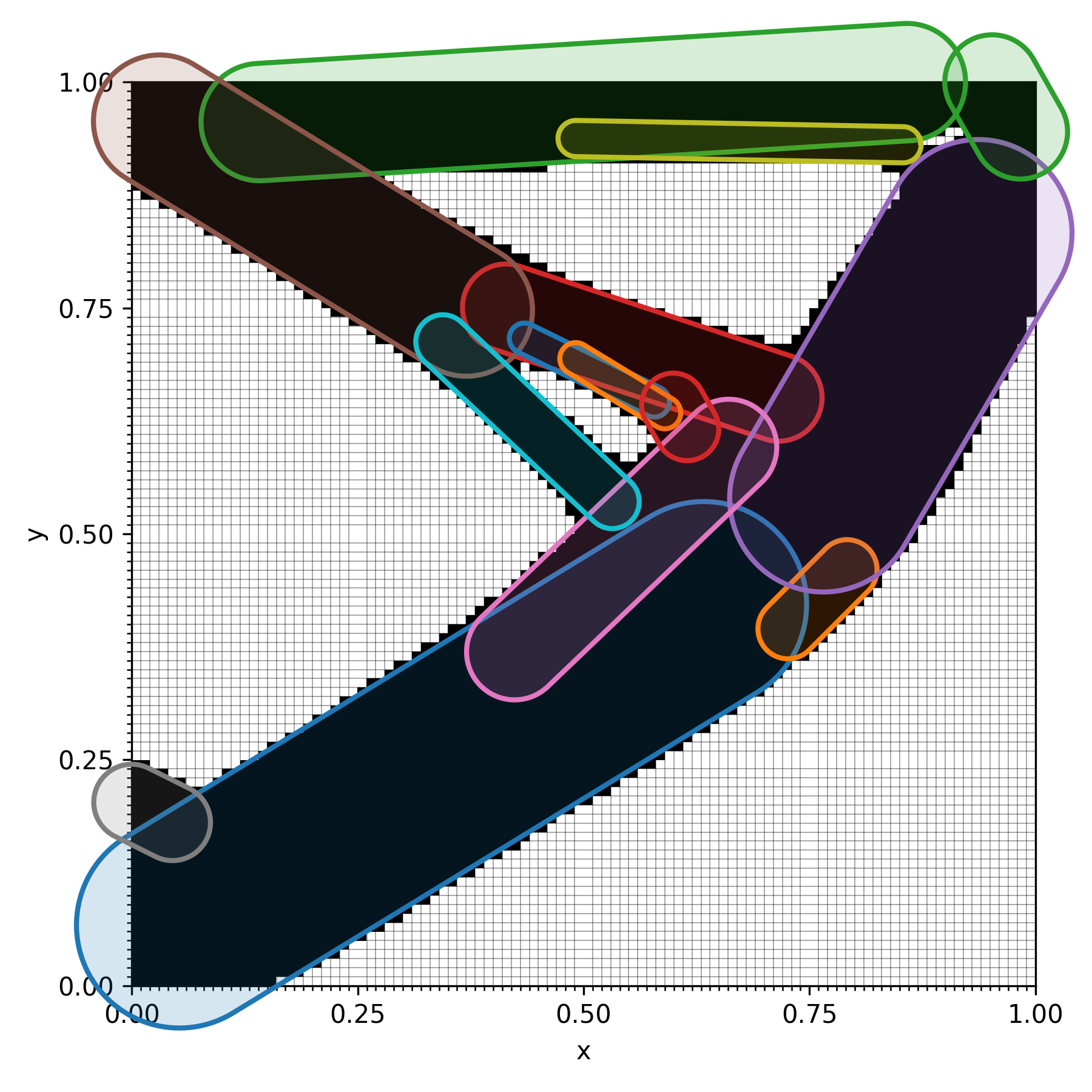} &
    \includegraphics[width=0.24\textwidth]{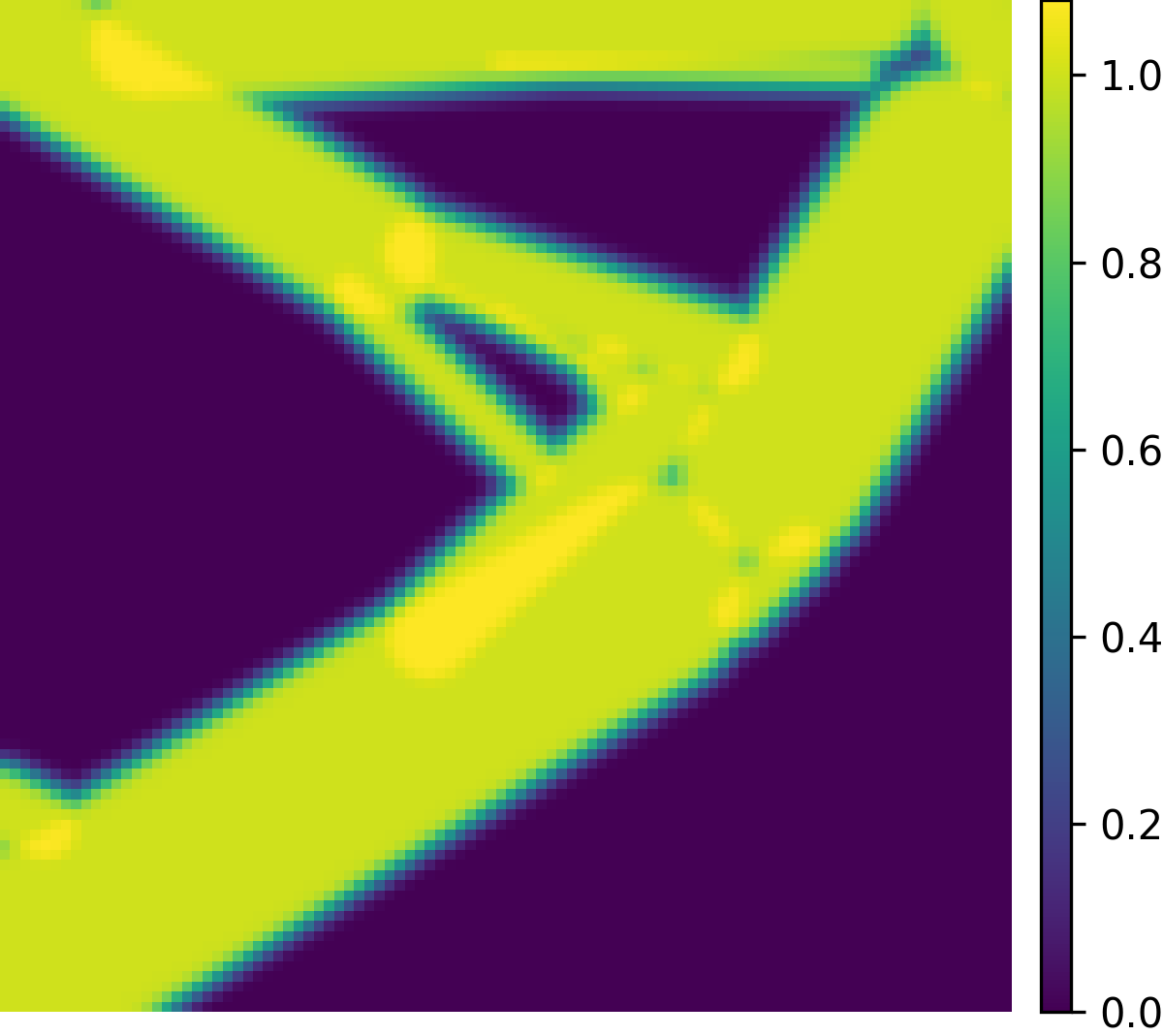} &
    \includegraphics[width=0.24\textwidth]{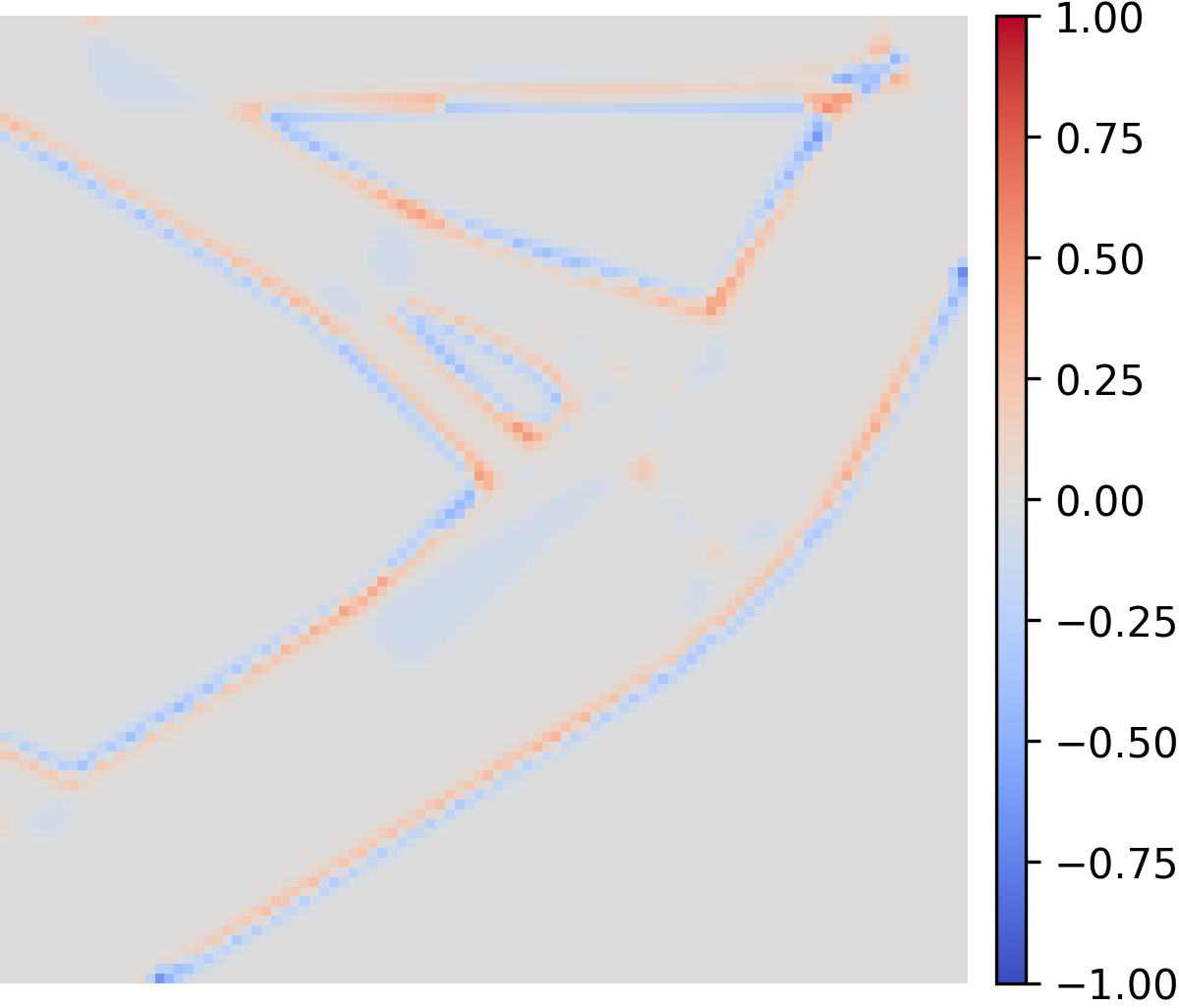} &
    \includegraphics[width=0.24\textwidth]{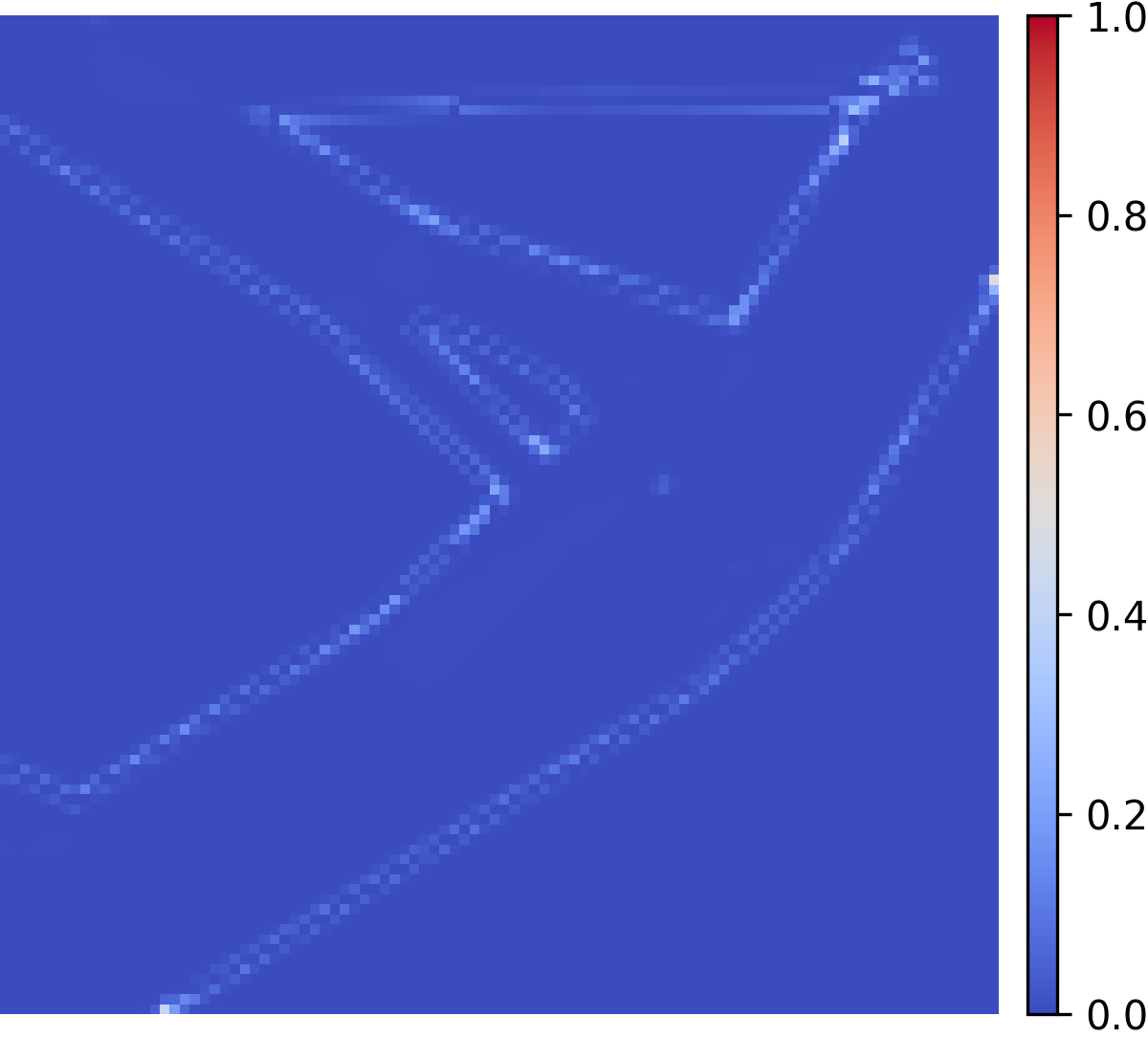} \\
    \includegraphics[width=0.24\textwidth]{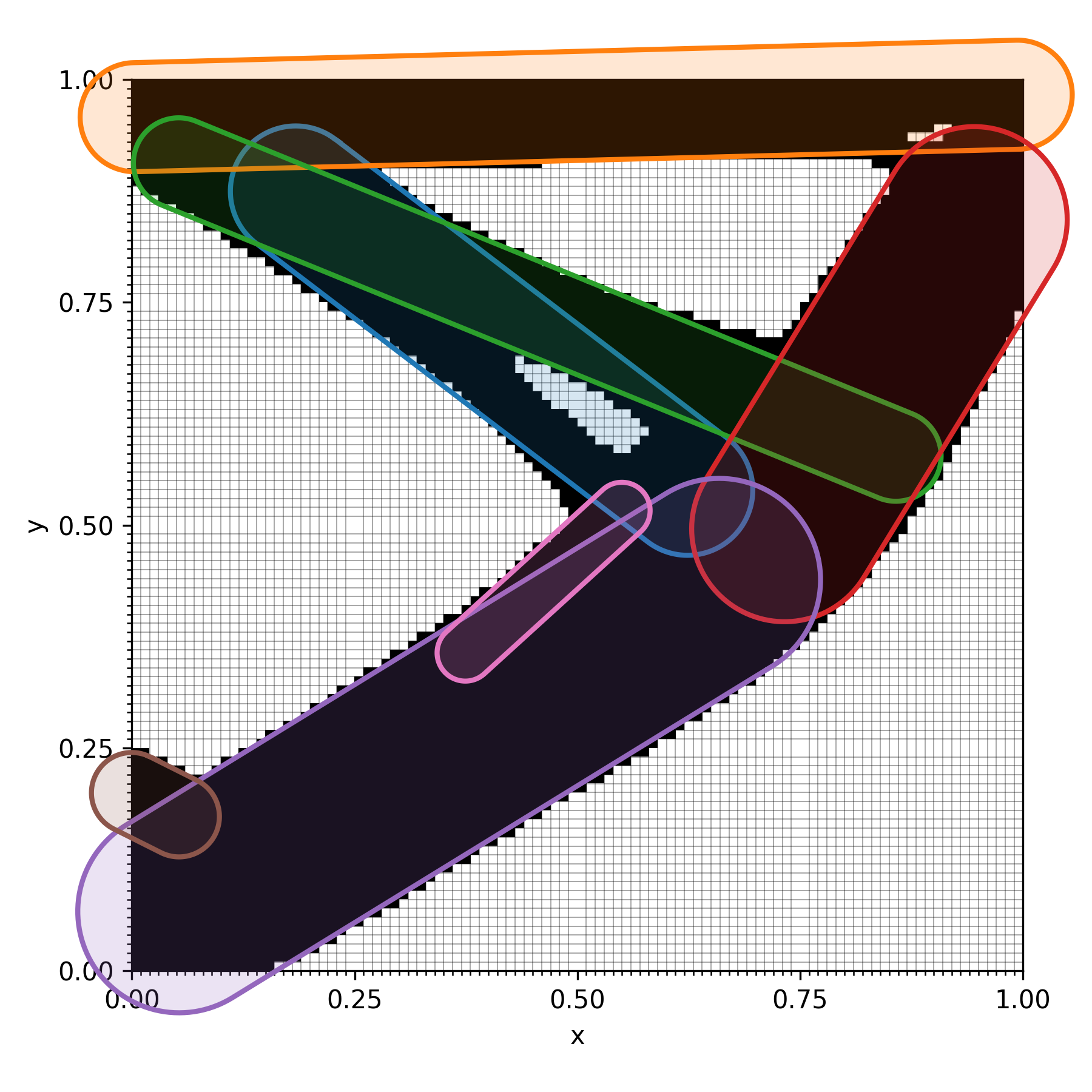} &
    \includegraphics[width=0.24\textwidth]{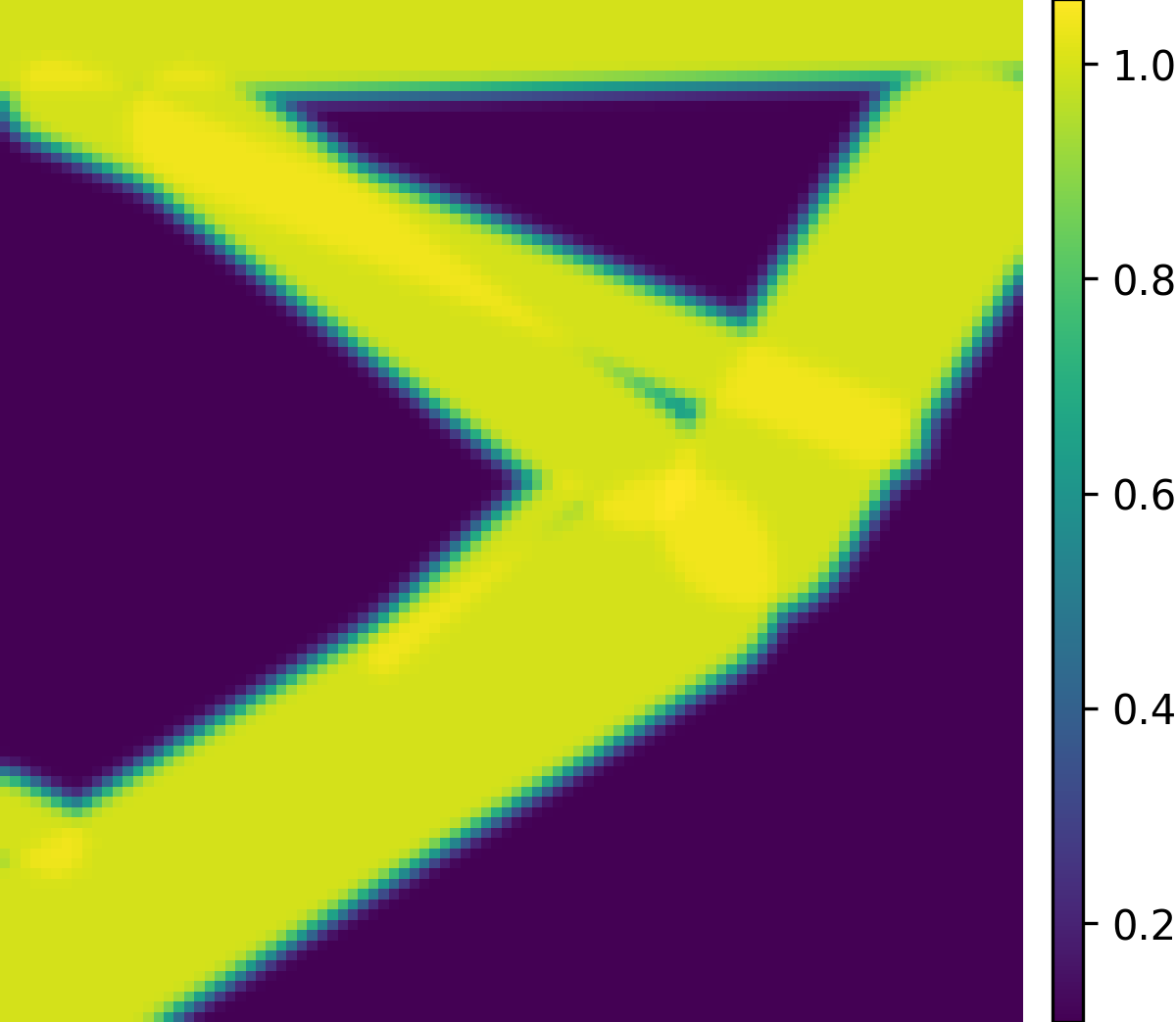} &
    \includegraphics[width=0.24\textwidth]{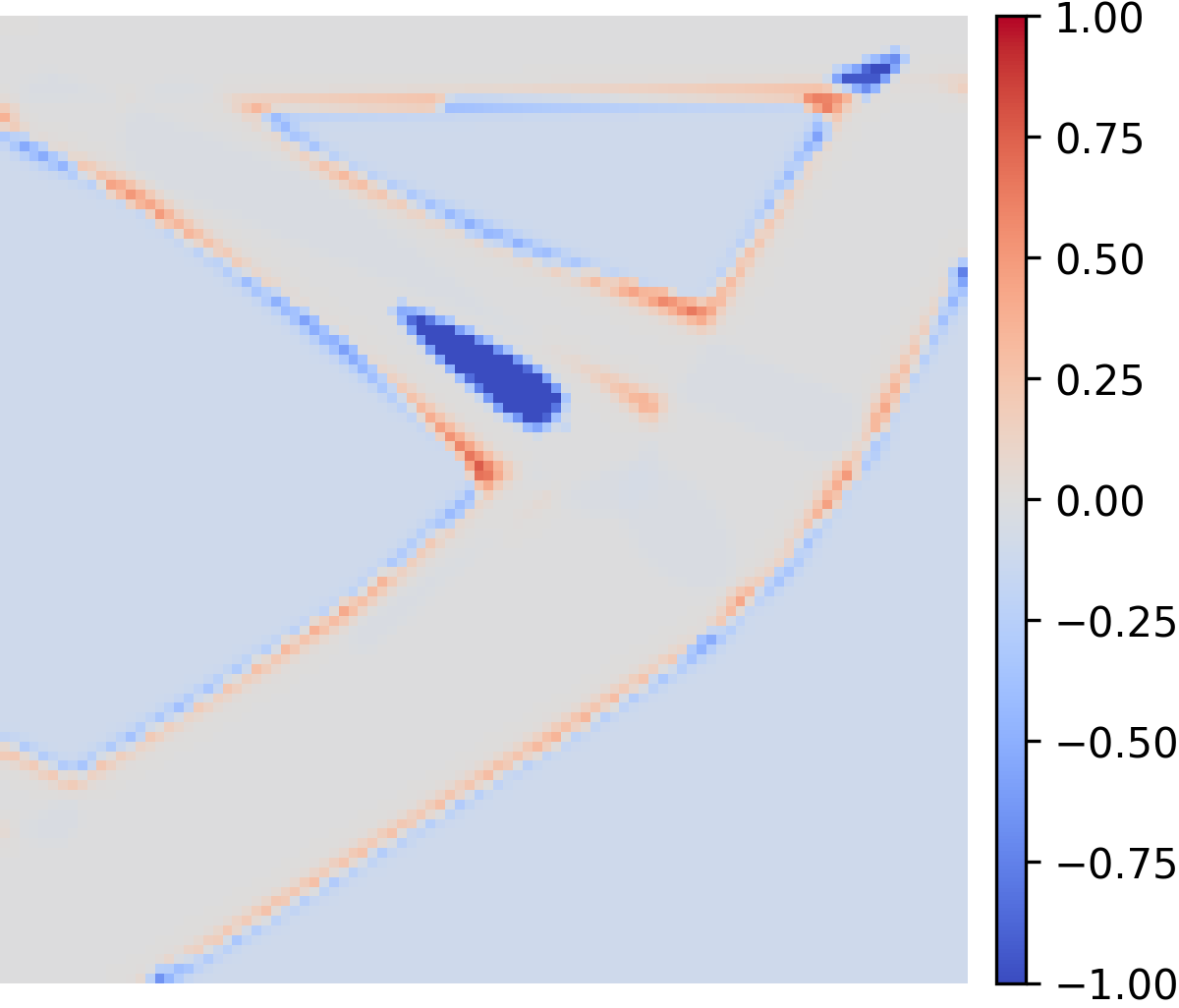} &
    \includegraphics[width=0.24\textwidth]{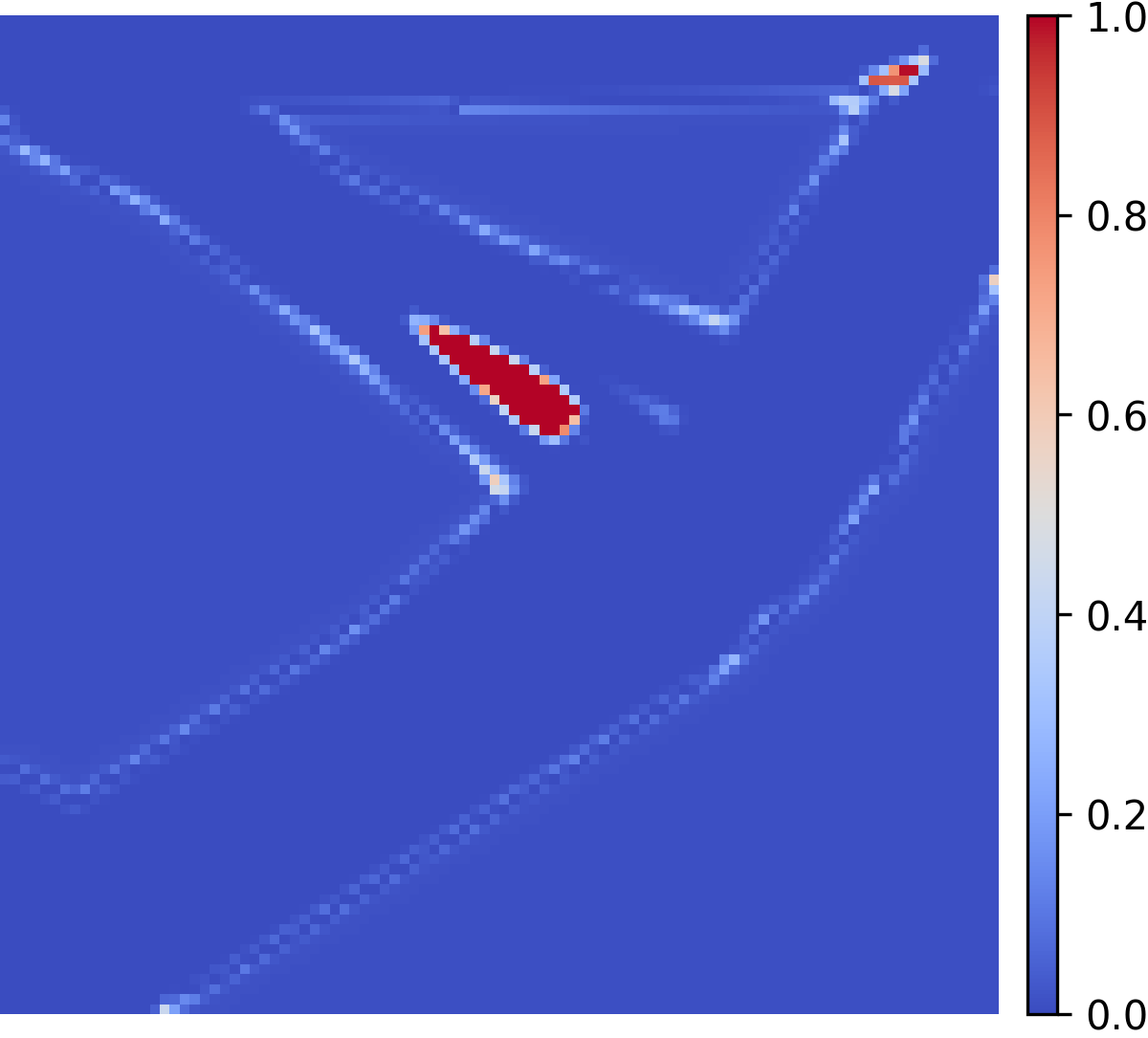} \\
  \end{tabular}
  \caption{Reconstructions after iterative refinement with $p$--norm (top) and
  softmax (bottom). 
  Refinement restores coverage in the $p$--norm case but is less effective for softmax.}
  \label{fig:canti-additive}
\end{figure}

Iterative refinement is able to compensate for these losses to varying degrees 
(\cref{fig:canti-additive}). 
In the $p$--norm pipeline, heuristic pruning removed thirteen of the initial eighteen pills, 
after which iterative refinement reintroduced nine new pills. 
Among these was a pill placed in the lower left corner, 
which focused exclusively on the small density pocket 
and aligned correctly with the target. 
Additional insertions concentrated around the central joint, 
where finer structures were needed to reduce residual errors, 
and two further capsules resolved the void pocket near the upper right boundary. 
This configuration is noteworthy in that it was the only case 
where the spurious filling of the upper void was successfully avoided. 
Compared to the initial state with many partially overlapping capsules, 
the final solution relies on fewer but larger pills for the main bars, 
complemented by smaller ones in critical curved regions. 
The outcome is a more compact yet higher-fidelity reconstruction, 
with redundant overlaps eliminated and finer details recovered through targeted additions. 

For the softmax pipeline, by contrast, refinement proves less effective. 
Once a large pill has stretched across a region to cover gaps left by pruning, 
the optimizer tends to preserve it rather than reintroduce finer structures. 
The reconstruction therefore remains less accurate than before pruning, 
despite the smoother appearance of the initial solution. 

A recurring limitation is visible in the upper right region. 
In most solutions, a large pill extends along the top boundary 
and simultaneously absorbs the small void pocket near the corner. 
Because such holes are only weakly penalized, 
they tend to be ``swallowed'' by nearby pills, 
and once this configuration is established it is rarely corrected during optimization.

The numerical results support these qualitative findings. 
Before heuristics and refinement, the objective values are 
$7.13\times 10^{1}$ for the $p$--norm and 
$1.71\times 10^{2}$ for softmax. 
After heuristics, they increase to 
$8.10\times 10^{2}$ and $2.02\times 10^{2}$, respectively, 
reflecting the loss of relevant coverage due to removal of auxiliary pills. 
Additive refinement restores accuracy, reducing the objectives to 
$4.33\times 10^{1}$ for the $p$--norm 
and $1.89\times 10^{2}$ for softmax. 

\begin{figure}[H]
    \centering
    \includegraphics[width=0.48\textwidth]{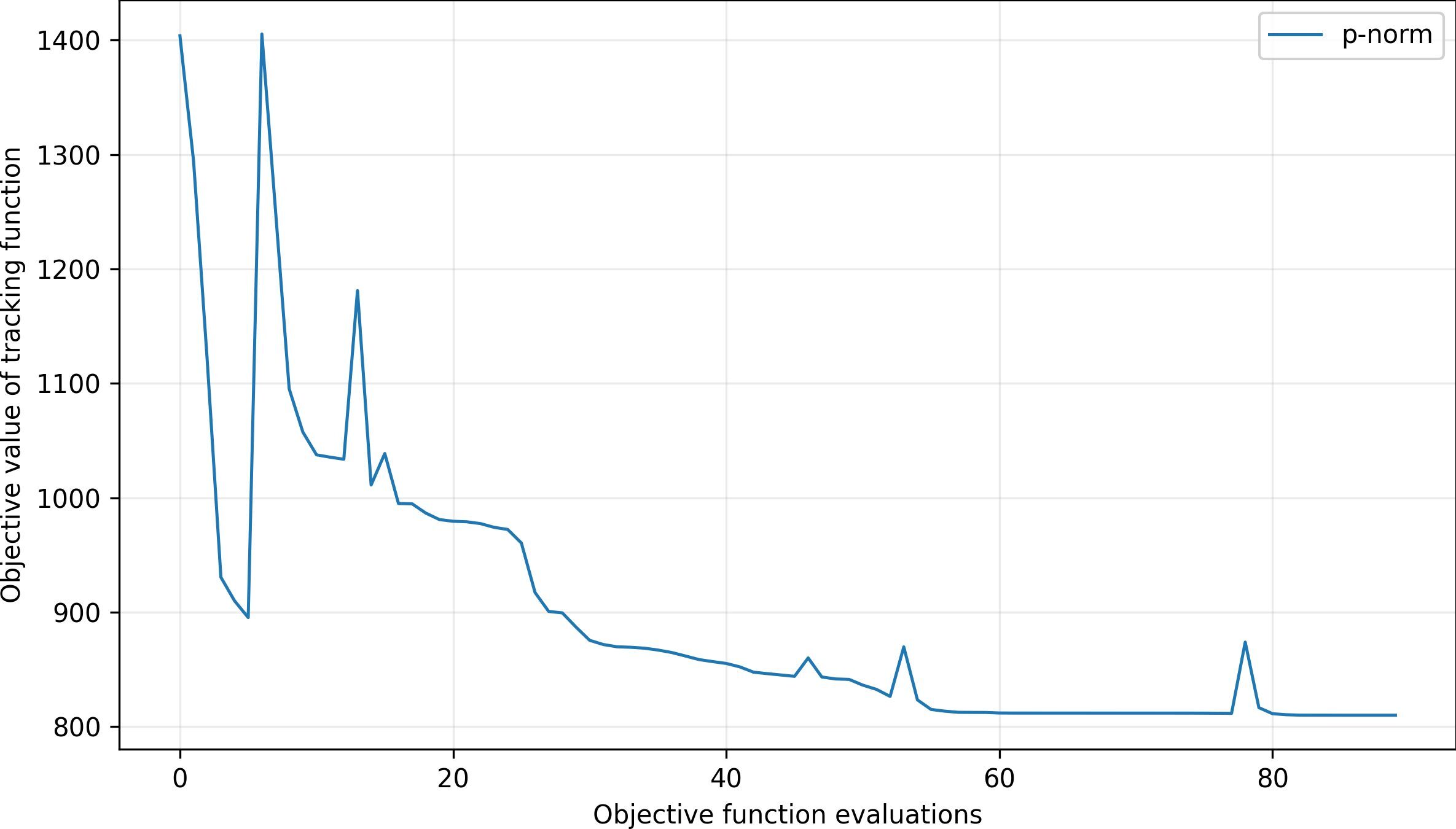}
    \includegraphics[width=0.48\textwidth]{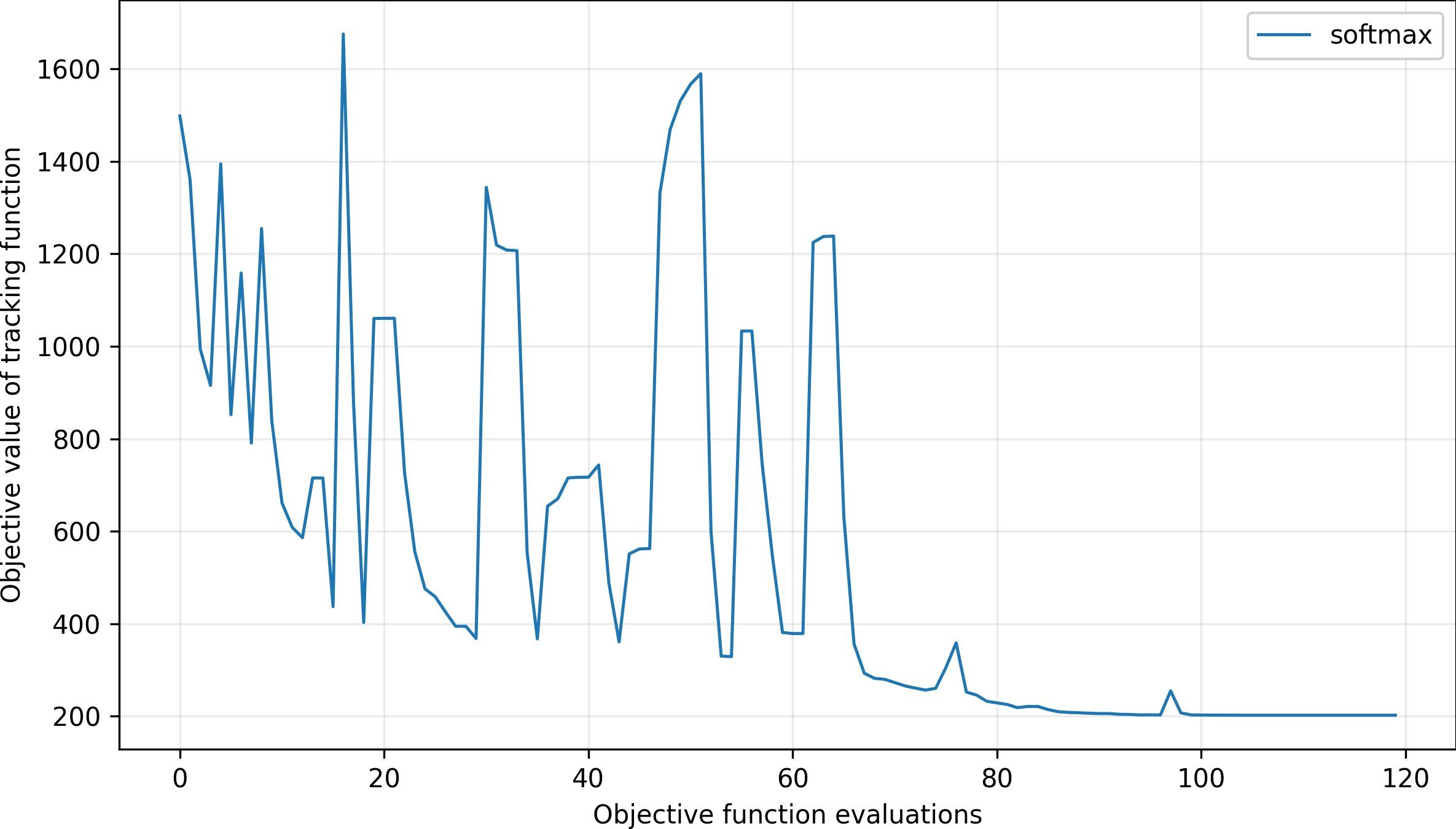} \\
    \includegraphics[width=0.48\textwidth]{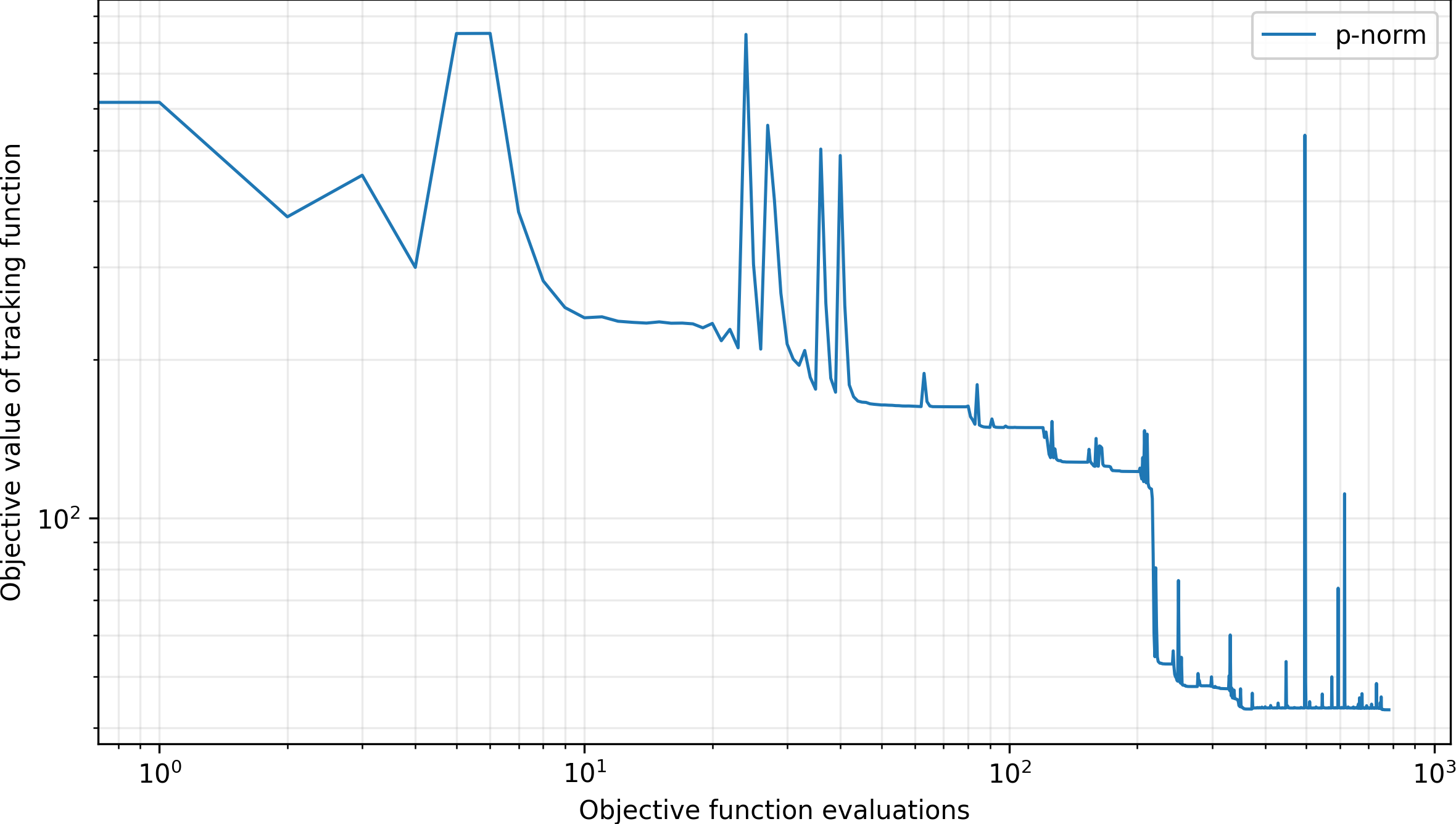}
    \includegraphics[width=0.48\textwidth]{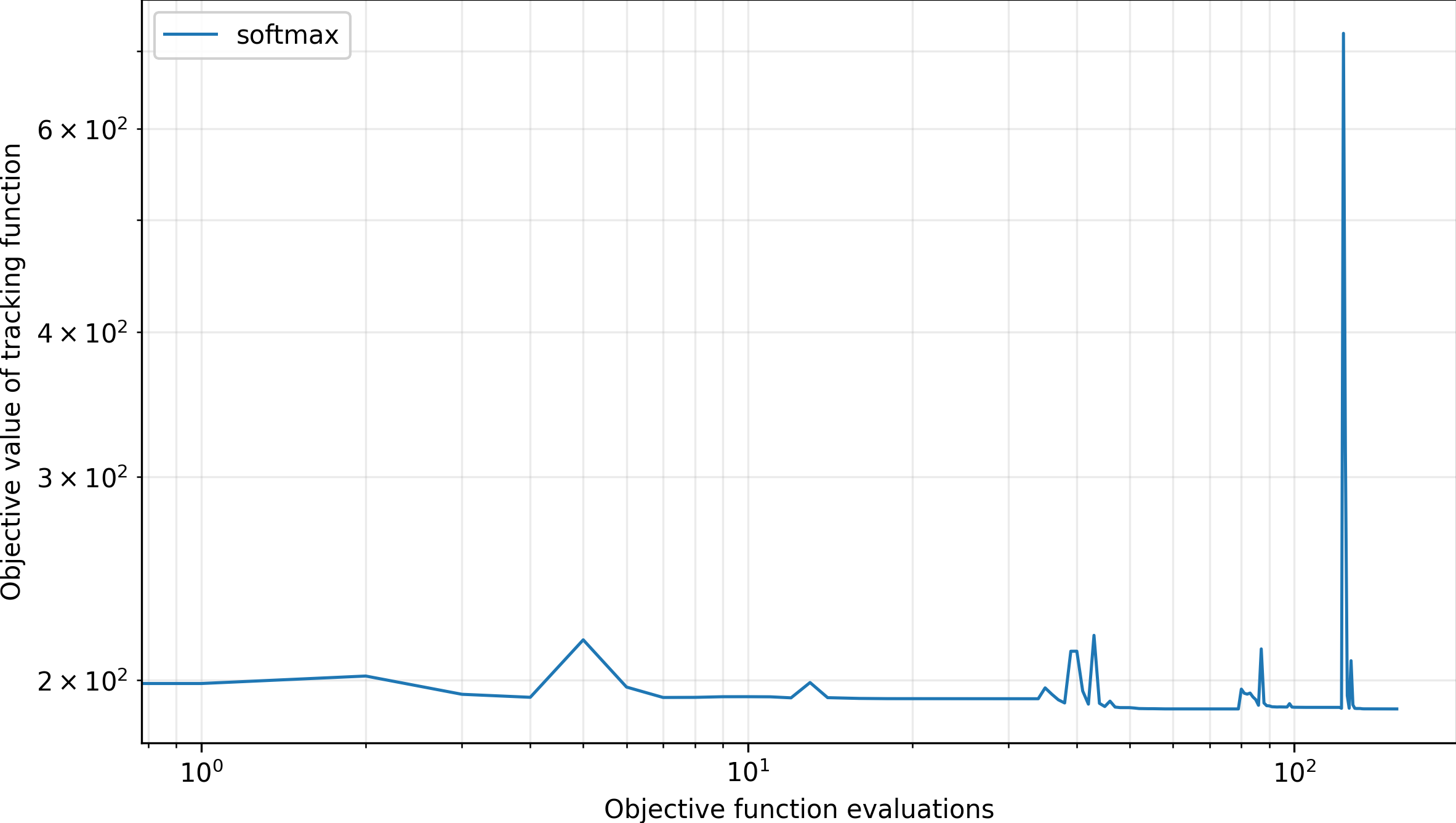} \\
    \includegraphics[width=0.7\textwidth]{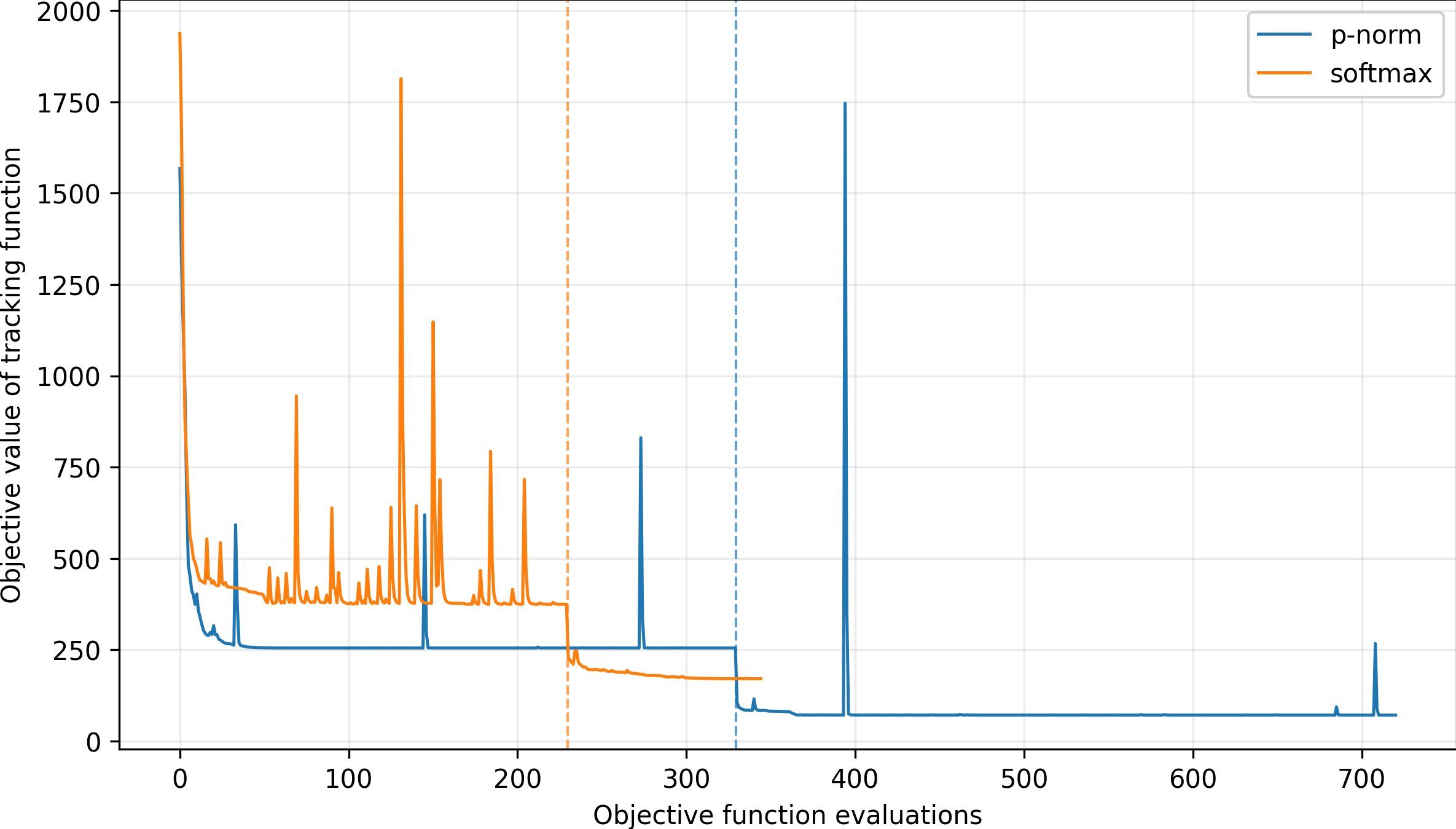}
    \caption{Objective trajectories for the cantilever test. Top: effect of 
    heuristics for $p$--norm (left) and softmax (right). Middle: effect of additive 
    refinement, with $p$--norm shown on logarithmic scale (left) and softmax on 
    linear scale (right). Bottom: combined comparison of both pipelines.}
    \label{fig:canti-obj-stages}
\end{figure}

The corresponding trajectories of the objective function provide further insight. 
\cref{fig:canti-obj-stages} compares the progression across the different stages. 
After the initial runs, both pipelines display rapid drops in the objective, 
but the curves typically plateau rather than terminating smoothly at tolerance criteria. 
This behavior reflects the limited ability of pills to contract or slide once 
a local minimum has been reached, a phenomenon already observed in the five--bar case. 

During the heuristic stages, the $p$--norm curve (top left) shows a pronounced increase 
caused by aggressive pruning of overlapping pills, 
followed by a relatively stable plateau. 
The softmax counterpart (top right) rises more moderately but remains at a higher level, 
consistent with the qualitative observation that the central joint was inadequately reconstructed after pruning. 
In both cases, oscillations indicate pills shrinking or reorienting erratically 
rather than converging smoothly. 

Additive refinement demonstrates its corrective effect most clearly in the $p$--norm case. 
The log--scaled trajectory (bottom left) highlights sharp objective decreases 
whenever a new pill is added and successfully aligned with the residual. 
The final objective is significantly lower than in the pre--pruning state, 
reflecting improved coverage with fewer but more effectively placed pills. 
The softmax trajectory (bottom right), plotted on a linear scale, 
exhibits a more gradual descent and does not achieve the same level of improvement. 
Once large pills have spanned residual holes, 
refinement cannot fully reintroduce the finer structures lost during pruning. 
This asymmetry between the two aggregation schemes explains the differences 
observed in the final reconstructions. 

In summary, the cantilever test provides quantitative confirmation of the effects 
observed qualitatively. 
Both pipelines were initialized with eighteen pills. 
In the $p$--norm case, heuristic pruning reduced this number substantially, 
and after iterative refinement the final configuration stabilized with fourteen pills. 
Despite the reduction, the achieved objective value 
($4.33\times 10^{1}$) is lower than in the original run, 
showing that a more compact representation can yield a superior reconstruction 
when redundant overlaps are removed and new pills are inserted in critical regions. 

The softmax pipeline follows a similar trend, though with less pronounced recovery. 
After pruning, the configuration contained only seven surviving pills, 
and the objective initially deteriorated. 
Iterative refinement improved the coverage, lowering the objective to 
$1.89\times 10^{2}$. 
While this remains slightly worse than the pre--pruning state, 
the solution is achieved with less than half of the original number of pills. 
Thus, both pipelines demonstrate that staged pruning and refinement 
can substantially reduce model complexity while retaining— 
and in the case of the $p$--norm even improving— 
the accuracy of the reconstruction. 

\chapter{Summary and Outlook\markboth{Summary and Outlook}{}}
\label{sec:outlook}

This final chapter reflects on the main contributions of the thesis and 
discusses possible avenues for further development. 
The first part outlines directions for extending the framework in terms of 
modeling flexibility and algorithmic robustness. 
The second part summarizes the findings of the numerical studies, 
highlighting the role of single- and multi-pill optimization, 
aggregation behavior, initialization strategies, discretization effects, 
second-order information and adaptive mechanisms.

\section*{Outlook}

The present study has established a pill-based framework for density representation 
and optimization in two dimensions. Building on this foundation, several directions 
for further development emerge, addressing both modeling flexibility and algorithmic robustness.  

A first avenue concerns the definition of geometric primitives. In the current formulation, 
each pill is parameterized by a single radius, enforcing uniform thickness along its medial axis. 
This restriction could be relaxed by assigning independent radii to the two endpoints, 
thereby enabling tapered geometries that approximate triangular or otherwise nonuniform structures 
with greater fidelity. Such an extension would increase adaptability of the representation 
while retaining analytic tractability of distances and sensitivities.  

A natural continuation is the extension to three dimensions. Conceptually, this requires 
replacing the two-dimensional distance evaluations with their spatial counterparts, 
while the overall pipeline of distance evaluation, pseudo-density mapping and staged optimization 
remains applicable. Such a generalization would allow volumetric pill assemblies and bring the 
framework closer to practical applications in structural design and additive manufacturing.  

The formulation of the exploration objective offers another potential refinement. 
The current reward-only objective favors overlap and may promote redundant coverage. 
Alternative formulations that penalize excessive overlap while preserving long-range 
directional cues could improve exploration efficiency and reduce stagnation in congested regions.  

Initialization strategies likewise admit further enhancement. The cross initialization adopted here 
ensures domain coverage but introduces pronounced symmetries, forcing pills to repel one another 
before adapting to the target. More flexible schemes, such as randomized perturbations, 
hierarchical seeding or radius scaling, could reduce redundancy and accelerate convergence 
while preserving coverage.

Further potential lies in refining transition functions. In particular, the construction 
of smooth asymmetric transitions based on incomplete beta functions provides an attractive option. 
The proposed Beta-CDF formulation produces compactly supported, strictly monotone $\Ccal^{2}$ profiles 
with adjustable exterior reach. It attains the exact midpoint value at $d=0$, 
while ensuring that both the value and its first two derivatives vanish at the outer boundary. 
Analytical expressions confirm strict monotonicity, convex–concave shape and endpoint regularity. 

Taken together, these outlook directions indicate how the framework may be advanced both 
in modeling capabilities—through tapered and three-dimensional primitives, refined transition 
functions and explicit void sensitivity—and in algorithmic robustness—through improved objectives 
and initialization. These developments form a clear roadmap toward a more versatile and effective 
pill-based topology optimization methodology.  
\begin{figure}[H] \centering \includegraphics[width=0.8\textwidth]{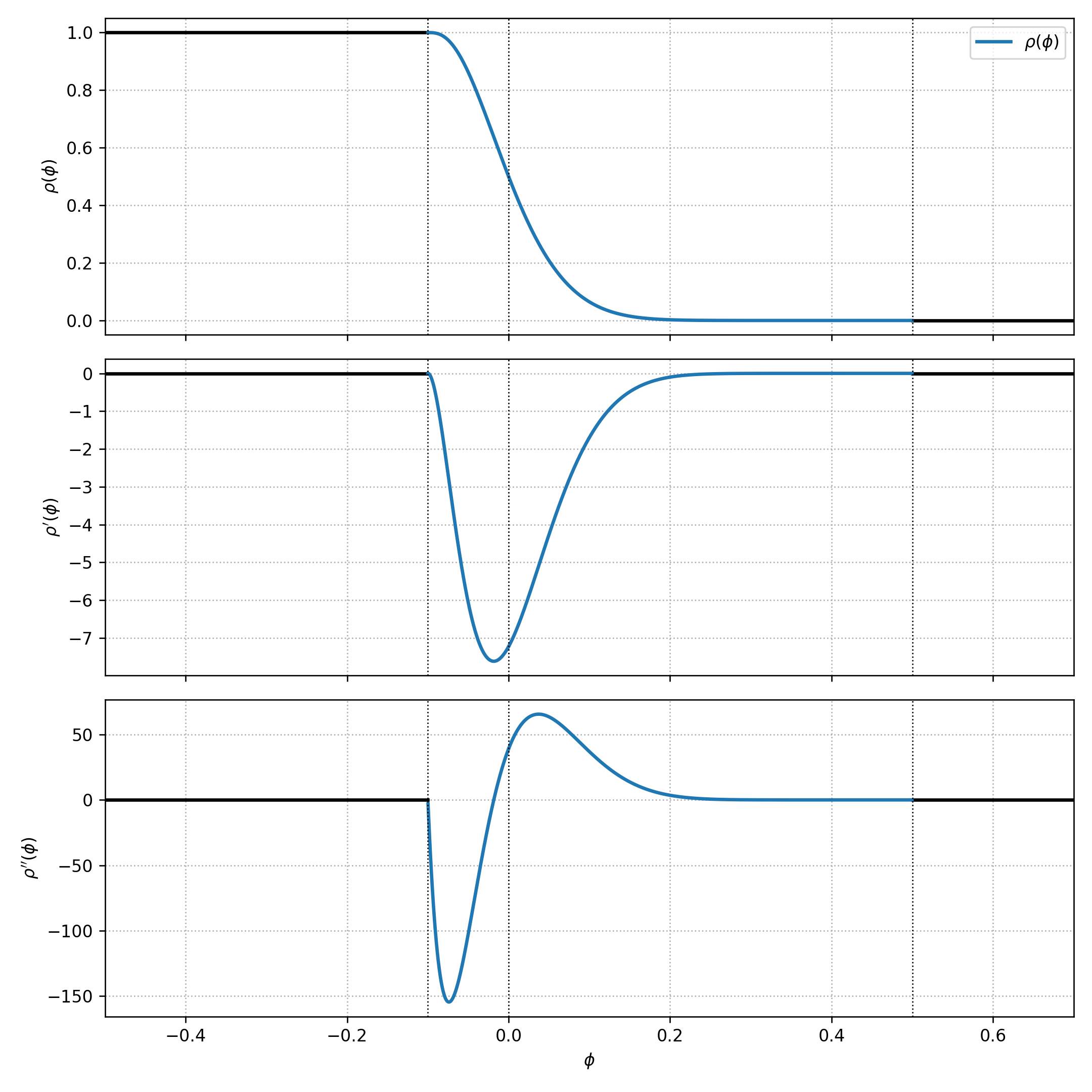}
   \caption{Beta-CDF transition with parameters $h=0.1$ and $\mathrm{ext}=0.4$. Top: density $\rho(d)$, middle: first derivative $\rho'(d)$, bottom: second derivative $\rho''(d)$. Both derivatives vanish at $d=-h$ and $d=L$, confirming compact support and $\Ccal^{2}$ regularity. The monotone negative slope of $\rho'(d)$ and the sign change of $\rho''(d)$ reflect the smooth convex-concave transition.}
   \label{fig:beta_mid_0104}
\end{figure}
\section*{Summary}

This thesis has developed and evaluated a pill-based reconstruction framework 
for density fields obtained from density-based topology optimization. 
The representation employs compact geometric primitives in the form of pill-shaped capsules, 
for which closed-form signed distances and analytic derivatives up to second order were derived. 
These ingredients enabled differentiable pseudo-densities through smooth transition functions 
and aggregation operators, mapped consistently onto fixed Cartesian grids.  

A central methodological contribution is the consistent derivation of analytic gradients and Hessians. 
This eliminated inaccuracies of semi-analytical sensitivities and permitted the use of exact 
second-order information in Newton-type solvers. Based on this analytic foundation, 
a staged optimization pipeline was constructed, comprising exploration, bridging and final convergence. 
Robustness was further enhanced by asymmetric transition extensions that activate far-field sensitivities, 
a reward-only exploration objective for reliable pose recovery and geometric constraints on segment length.  

Beyond the single-pill case, the thesis systematically investigated 
multi-pill optimization. Comparisons of aggregation families and parameter choices 
clarified their influence on coverage, stability and convergence. 
$p$--norm aggregation with $p=7$–$9$ and softmax aggregation with $\beta=10$–$14$ 
proved most effective, balancing compactness and accuracy. 
Initialization strategies were shown to strongly influence performance: 
cross initialization guarantees coverage but induces symmetry bias, 
while randomized seeding provides greater flexibility at the risk of inconsistency.  
Discretization quality was analyzed through mesh resolution and integration order, 
demonstrating that higher spatial resolution improves accuracy, 
while integration order beyond two yields little additional benefit.  
The role of second-order information was clarified by comparing exact Hessians 
with limited-memory BFGS updates: the exact Hessian consistently delivered 
higher fidelity reconstructions, particularly in cases with strong pill interactions.  

The number of pills was identified as a decisive factor in balancing underfitting and overfitting. 
Experiments with counts ranging from three to eighteen revealed that very low numbers 
leave major target structures uncovered, while excessive redundancy introduces instability. 
Around twelve to thirteen pills emerged as a robust compromise for the five-bar benchmark.  
Adaptive mechanisms were proposed to address these challenges:  
iterative refinement adds pills in residual regions to counteract undercoverage, 
while heuristics based on area and uniqueness ratios, as well as pill consolidation, 
reduce redundancy by pruning or merging overlapping elements.  

Finally, the framework was applied to a cantilever benchmark to demonstrate 
generalization beyond the five-bar case. 
Here, staged pruning and refinement showed that a more compact representation 
can in fact improve reconstruction quality when redundant overlaps are removed 
and new pills are introduced in structurally critical regions. 

Overall, the study demonstrated that pill-based reconstruction of voxelized SIMP densities 
is feasible, efficient and interpretable. 
The approach provides a more interpretable parametric representation of SIMP densities 
and thereby offers a potential pathway toward CAD-amenable formulations, 
although additional mechanisms such as additional heuristics would be required 
to ensure compatibility with practical design requirements.

% Bibliography
\printbibliography

@article{Norato2015GeometryProjection,
  author  = {Norato, Juan A. and Bell, Benjamin and others},
  title   = {A Geometry Projection Method for Continuum-Based Topology Optimization},
  journal = {Structural and Multidisciplinary Optimization},
  year    = {2015},

}

@article{WeinStingl2018ParamErsatz,
  author  = {Wein, Florian and Stingl, Michael},
  title   = {Parametric Ersatz Modeling for Feature-Based Topology Optimization},
  journal = {Journal },
  year    = {2018},

}

@article{Wein2020Review,
  author  = {Wein, Florian and Stingl, Michael and others},
  title   = {Feature-Based or Geometry-Aware Topology Optimization: A Review},
  journal = {Journal },
  year    = {2020},

}

@article{Greifenstein2023Spaghetti,
  author  = {Greifenstein, Niklas and others},
  title   = {Spaghetti/Bar-like Feature Decomposition for Topology Optimization},
  journal = {Journal},
  year    = {2023},

}

@article{Shannon2023PostProcessing,
  author  = {Shannon, John and others},
  title   = {Post-Processing of Density Fields via Bar Assemblies},
  journal = {Journal },
  year    = {2023},

}

@article{Hoang2017MMB,
  author  = {Hoang, Viet-Tien and others},
  title   = {Moving Morphable Bars for Structural Topology Optimization},
  journal = {Structural and Multidisciplinary Optimization},
  year    = {2017},
  volume  = {56},
  number  = {3},
  pages   = {761--773},
  doi     = {10.1007/s00158-017-1702-7},

}

@article{Zhang2016MMC,
  author  = {Zhang, Wenbin and others},
  title   = {Moving Morphable Components Based Structural Topology Optimization},
  journal = {Computer Methods in Applied Mechanics and Engineering},
  year    = {2016},
  volume  = {310},
  pages   = {393--414},
  doi     = {10.1016/j.cma.2016.07.015},

}

@article{Zhang2017CBS,
  author  = {Zhang, Wenbin and others},
  title   = {Consistent Bars: A Feature-Based Approach for Structural Topology Optimization},
  journal = {Structural and Multidisciplinary Optimization},
  year    = {2017},
  volume  = {55},
  number  = {3},
  pages   = {953--968},
  doi     = {10.1007/s00158-016-1551-3},

}

% Declaration of Authorship
\chapter*{Declaration of authorship\markboth{Declaration of authorship}{}}
\addcontentsline{toc}{chapter}{Declaration of authorship}

I confirm that I have written this thesis unaided and without using sources other than those listed and that this thesis has never been submitted to another examination authority and accepted as part of an examination achievement, neither in this form nor in a similar form. All content that was taken from a third party either verbatim or in substance has been acknowledged as such.

\vspace{.5cm}
\vspace*{1.5cm} \par

\begin{flushright}
    \line(1,0){200} \par
\end{flushright}
\docCity, \docDeadline 
\hfill\pbox{\textwidth}{
\docAuthor
}

\end{document}